\documentclass[headsepline]{scrbook}

\usepackage{amssymb,amsfonts,amsthm,amsmath} 
\usepackage{fullpage}
\usepackage{wasysym}
\usepackage{makeidx}
\usepackage{comment}
\usepackage{multicol}
\usepackage{scrlayer-scrpage}

\usepackage{background}
\SetBgContents{}

\makeindex

\usepackage{graphicx}
\usepackage[dvipsnames,svgnames,table]{xcolor}
\usepackage[colorlinks=true, pdfstartview=FitV,linkcolor=ForestGreen,citecolor=ForestGreen, urlcolor=blue]{hyperref}

\usepackage{esint} 
\usepackage{enumitem}

\usepackage[utf8]{inputenc}

\newcommand{\R}{\mathbb{R}}
\newcommand{\N}{\mathbb{N}}

\newcommand{\eps}{\varepsilon}
\newcommand{\etad}{{\eta_2}}

\newcommand{\cyril}[1]{\textcolor{red}{#1}}

\newcommand{\dd}{\, \mathrm{d} }
\newcommand{\dv}{\, \mathrm{d} v}
\newcommand{\ds}{\, \mathrm{d} s}
\newcommand{\dr}{\, \mathrm{d} r}
\newcommand{\dt}{\, \mathrm{d} t}
\newcommand{\dx}{\, \mathrm{d} x}

\newcommand{\dy}{\, \mathrm{d} y}
\newcommand{\dz}{\, \mathrm{d} z}
\newcommand{\dw}{\, \mathrm{d} w}

\newcommand{\phitrunc}{\varphi_{\mathrm{loc}}}
\newcommand{\ftrunc}{f_{\mathrm{loc}}}
\newcommand{\mtrunc}{\mathfrak{m}_{\mathrm{loc}}}
\newcommand{\Ttrunc}{\mathfrak{S}_{\mathrm{loc}}}
\newcommand{\Strunc}{S_{\mathrm{loc}}}
\newcommand{\bStrunc}{\bar{S}_{\mathrm{loc}}}
\newcommand{\Gammax}{\Gamma^{(x)}}
\newcommand{\Gammav}{\Gamma^{(v)}}

\newcommand{\Calp}{C_{\alpha,p}}

\newcommand{\Cee}{C_{\mathrm{DG}}^\pm}
\newcommand{\Ceem}{C_{\mathrm{DG}}^-}
\newcommand{\Cexp}{C_{\mathrm{exp}}}
\newcommand{\Cpdgp}{C_{\mathrm{pDG}^+}}
\newcommand{\Ckdgp}{C_{\mathrm{kDG}^+}}
\newcommand{\Cpdgm}{C_{\mathrm{pDG}^-}}
\newcommand{\Ckdgm}{C_{\mathrm{kDG}^-}}
\newcommand{\Cc}{C_{\mathrm{c}}}

\newcommand{\Csob}{C_{\mathrm{Sob}}}
\newcommand{\Cpw}{C_{\mathrm{PW}}}

\newcommand{\Civl}{C_{\mathrm{IVL}}}
\newcommand{\Cdg}{C_{\mathrm{DG}}}
\newcommand{\Clmp}{C_{\mathrm{LMP}}}
\newcommand{\bClmp}{\bar{C}_{\mathrm{LMP}}}
\newcommand{\eClmp}{{C}_{\mathrm{LMP},\eps}}
\newcommand{\pClmp}{{C}_{\mathrm{LMP},p}}
\newcommand{\Cdgun}{C_{\mathrm{DG1}}}

\newcommand{\Charnack}{C_{\mathrm{H}}}
\newcommand{\Cwhi}{C_{\mathrm{whi}}}
\newcommand{\Ceps}{C_{\mathrm{int}}}
\newcommand{\Cnu}{C_{\mathrm{ls}}}
\newcommand{\Cslsd}{C_{\mathrm{lsd}}}

\newcommand{\Cpar}{C_{\mathrm{par}}}
\newcommand{\Ckin}{C_{\mathrm{kin}}}

\newcommand{\domain}{\Omega}
\newcommand{\Ox}{\Omega_x}
\newcommand{\Ov}{\Omega_v}
\newcommand{\Breita}{B_{r_\eta}}
\newcommand{\Qmax}{Q_{\mathrm{max}}}
\newcommand{\Qwhi}{Q_{\mathrm{whi}}}
\newcommand{\Qharnack}{Q_{\mathrm{harn}}}
\newcommand{\Qpast}{Q_{\mathrm{past}}}
\newcommand{\Qfuture}{Q_{\mathrm{future}}}
\newcommand{\Qext}{Q_{\mathrm{ext}}}

\newcommand{\Qint}{Q_{\mathrm{int}}}

\newcommand{\Qexp}{Q_{\mathrm{exp}}}

\newcommand{\Qstack}{Q_{\mathrm{stack}}}

\newcommand{\Qgu}{Q_{\mathrm{pos}}}
\newcommand{\Ngu}{\mathcal{N}_{\mathrm{pos}}}
\newcommand{\bQgu}{Q^{\mathrm{crop}}_{\mathrm{pos}}}
\newcommand{\Qmid}{Q_{\mathrm{mid}}}
\newcommand{\Tmid}{T_{\mathrm{mid}}}
\newcommand{\dkin}{d_{\mathrm{kin}}}
\newcommand{\astkin}{*_{\mathrm{kin}}}
\newcommand{\weak}{\mathrm{weak}}

\newcommand{\un}{\mathbf{1}}

\newtheorem{thm}{Theorem}[section]
\newtheorem{lemma}[thm]{Lemma}
\newtheorem{prop}[thm]{Proposition}
\newtheorem{cor}[thm]{Corollary}
\newtheorem{defi}{Definition}
\theoremstyle{remark}
\newtheorem{remark}{Remark}
\newtheorem*{ellipticity}{Ellipticity}

\newtheorem{example}{Example}

\DeclareMathOperator{\PV}{PV}
\DeclareMathOperator{\osc}{osc}
\DeclareMathOperator{\dive}{div}

\DeclareMathOperator{\esssup}{ess-sup}
\DeclareMathOperator{\essinf}{ess-inf}
\DeclareMathOperator{\supp}{supp}

\begin{document}

\setlength{\headsep}{5mm}

\begin{titlepage}
   \centering
  \phantom{trick}
    \vspace{3cm}

    \huge \textbf{\textsf{ De Giorgi's regularity theory }}  \\
      \textbf{\textsf{for elliptic, parabolic and kinetic equations}}
    \vskip 5mm
    
    {\Large \textbf{\textsf{ Local  diffusions}}}
    
    \vskip 2cm
    {\Large Cyril Imbert\par}
    \vskip 1cm
    \date{\today}
    
    \vfill
\begin{figure}[h]
  \centering{\includegraphics[height=6cm]{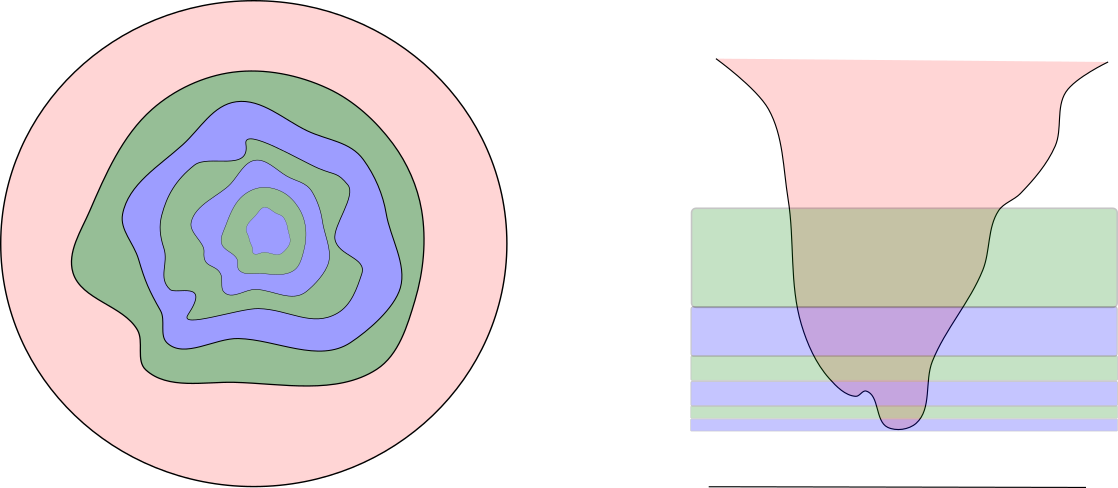}}
\end{figure}


    \vfill
\end{titlepage}


\setcounter{tocdepth}{1}
\tableofcontents

\frontmatter

\chapter{Foreword}

This book is based on lectures I delivered on De Giorgi's regularity theory. A preliminary version of the material in Chapters~\ref{c:elliptic} and \ref{c:parabolic} was first presented at \'Ecole normale sup\'erieure (Paris Sciences et Lettres) in spring 2025. I later expanded and delivered these lectures at the University of California, Berkeley, in fall 2025.
\smallskip

I first heard about De Giorgi's regularity theory by discussing with A.~F.~Vasseur during my  visits at the University of Texas at Austin in the early 2000s.
I very often referred to his two concise lecture notes about De Giorgi's methods \cite{MR3525875,MR2660718}, the second co-authored with L.~Caffarelli. I am very grateful to Alexis and Luis. 
\smallskip

In 2010, I began my research on kinetic equations together with Clément Mouhot while we were both at ENS. Clément was already a significant figure in the field and he proposed
the study of some toy nonlinear models. After attempting to apply the Krylov-Safonov approach, we eventually recognized that De Giorgi's methods could be strategically leveraged to our benefit for studying the regularity of kinetic Fokker-Planck equations. Our research led us to a paper by W. Wang and L. Zhang, published a couple of years earlier, building upon ideas from S.~N.~Kruzhkov. Over the ensuing years, Clément and I devoted ourselves to mastering various techniques and tools derived from the theory of elliptic and parabolic equations in divergence form. Clément also wrote recently  lecture notes with G.~Brigati \cite{brigati2025introductionquantitativegiorgimethods} about quantitative De Giorgi's methods.

Since 2008, Luis Silvestre and I have been collaborating on various challenges posed by elliptic and parabolic equations in non-divergence form. Around the year  2015, I have engaged in extensive discussions with Luis, François Golse, Clément Mouhot, and Alexis Vasseur regarding Harnack's inequality for kinetic Fokker-Planck equations. It is worth mentioning that Luis could have been a co-author of the paper we produced during that time \cite{zbMATH07050183}. After posting this first work, and following the publication of his seminal paper \cite{MR3551261} about  the Boltzmann equation, we embarked on the conditional regularity program for the space-in-homogeneous Boltzmann equation without cut-off.

I am deeply grateful to Clément and Luis for their years of steadfast and supportive collaboration.
\smallskip

I am also grateful for the many reactions, insights, and questions raised during lectures, both by students from ENS and UC Berkeley, and by colleagues from the two research programs hosted at the Simons-Laufer Mathematical Research Institute (SL Math) during the fall of 2025. Additionally, two colleagues from UC Berkeley attended the lectures and provided valuable feedback.
I would like to express my sincere gratitude to all of them for their attention and patience.
\smallskip

This book was made possible by the Chancellor's Professorship at the University of California, Berkeley. My time in this exceptional department was truly enjoyable. I also greatly benefited from the remarkable environment provided by SL Math (funded by NSF Grant No. DMS-1928930) during my stay in California, particularly the outstanding staff who took such excellent care of us.

\begin{flushright}
Berkeley, December 20, 2025
\end{flushright}

\chapter{From De Giorgi to Boltzmann}
\label{c:history}

This book presents a comprehensive regularity theory for solutions of elliptic, parabolic, and kinetic equations. The foundation of this theory was laid by E. De Giorgi's groundbreaking resolution of Hilbert's nineteenth problem in 1956. The innovative tools he developed to tackle this problem proved to be remarkably versatile. In 1957, just one year later, J. Nash independently developed analogous techniques for parabolic equations, concurrently with De Giorgi's research. By the year 2000, these techniques had been extended to address elliptic and parabolic equations featuring integral diffusion, such as the fractional Laplacian. More recently, the theory has evolved to encompass kinetic equations, accommodating both local and integral diffusion processes. This book aims to present these results in a unified and coherent manner, beginning with the classical elliptic framework and progressing through to the most recent advancements in kinetic equations.
\bigskip

\noindent \fbox{\parbox{.987\linewidth}{
\textbf{Disclaimer.}
The first version of these lecture notes does not contain results about equations with integral diffusions. In particular linear kinetic equations related to the Boltzmann equation are not addressed. Hopefully, they will be covered in a second version. 
}}

\section{Hilbert's 19\textsuperscript{th} problem}

The field of regularity theory for elliptic equations has seen remarkable advancements, notably with the resolution of Hilbert's 19\textsuperscript{th} problem. This problem was one of 23 posed by David Hilbert to the mathematical community during the International Congress of Mathematicians held in Paris.

\subsection{Statement}

In the English translation of Hilbert's problems published in 1902 in the bulletin of the American Mathematical Society \cite{MR1557926} (see also this \href{https://en.wikipedia.org/wiki/Hilbert's_problems}{Wikipedia's page}), the problem is stated in the following way,
\begin{quote}
\textit{  Are the solutions of regular problems in the calculus of variations always necessarily analytic?}
\end{quote}
It is motivated by the classical fact that  minimizers of $\int |\nabla_x u|^2 \dx$ are harmonic, and thus analytic.
On the one hand, D. Hilbert makes the notion of regular problem precise,
\begin{quote}
 \textit{ If one assumes that $L$ is uniformly convex and analytic, are the minimizers of a functional of the form $E(u)=\int_{\Omega} L (\nabla_x u(x)) \dx$ always necessarily analytic?}
\end{quote}
On the other hand, he neither defines the domain of integration $\Omega$ nor specifies the set of functions $u \colon \Omega \to \R$ over which the minimum is considered.
We also mentions that he considers more general functionals of the form $\int L(x,u(x),\nabla_x u(x)) \dx$ but we stick to the previous framework for simplicity. 
The uniform convexity assumption of the Lagrangian $L$ can be understood in the following sense,
\begin{equation}
\label{convexity}  \tag{Convexity}
  \text{There exist $\lambda, \Lambda >0$ such that for all $p \in \R^d$}, \quad \lambda |\xi|^2 \le D^2L (p) \xi \cdot \xi \le \Lambda |\xi|^2. 
\end{equation}
In this statement, $D^2L (p)$ denotes the Hessian matrix of the function $L$ at point $p$. It is $d \times d$, real and symmetric.

\subsection{The Euler-Lagrange equation}

Let $v$ denote a minimizer of the functional $E (v)$. If $\varphi$ is smooth and compactly supported in $\Omega$ and $\eps >0$, then $E (v + \eps \varphi)
\ge E (v)$. Using the definition of $E$, dividing by $\eps$ and letting it goes to $0$ yields,
\[ \int_\Omega \nabla L (\nabla_x v) \cdot \nabla_x \varphi \dx =0 \] 
where $\nabla L (p)$ denotes the gradient of $L$. Since $L$ is analytic and if $v$
is twice differentiable (in the classical Fréchet sense), then we can integrate by parts
and obtain that \[ -\int_\Omega \dive_x (\nabla L( \nabla_x v)) \varphi \dx =0.\] Since
$\varphi$ is an arbitrary smooth and compactly supported function on $\Omega$, this implies that,
\begin{equation}
\label{e:el}  \tag{Euler-Lagrange}
  - \dive_x(\nabla L (\nabla_x v)) =0.
\end{equation}
In the previous formula, $\dive_x$ denotes the divergence operator with respect to the $x$ variable. We draw the attention of the reader towards the fact that this equation is nonlinear. Studying such equations is a priori very challenging.

\subsection{Schauder's theory}

Schauder's theory from the years 1930 applies to elliptic and parabolic equations with H\"older continuous
coefficients. For instance, it applies to elliptic equations under divergence form, 
\begin{equation}
  \label{e:parabolic0}
  - \sum_{i,j=1}^d a^{ij}(x) \frac{\partial^2 u}{\partial x_i \partial x_j} + \sum_{i=1}^d b^i (x) \frac{\partial u}{\partial x_i}= S(x) \quad \text{ for } x \in \R^d 
\end{equation}
where coefficients $a^{ij}$ and $b^i$ together with the function $f$ are assumed to be H\"older continuous and bounded in $\R^d$. 
Under the following ellipticity condition on coefficients $a^{ij}$,
\begin{equation}
\label{ellipticity}  \tag{Ellipticity}
  \text{There exists $\lambda, \Lambda >0$ s.t. for all $x \in \R^d$}, \quad \forall \xi \in \R^d, \lambda |\xi|^2 \le A (x) \xi \cdot \xi \le \Lambda |\xi|^2
\end{equation}
J.~Schauder showed that it is possible to construct  solutions that are twice differentiable in the space variable $x$. In this first setting, the regularity
of the function and of its derivatives of order $1$ and $2$ is measured with a  modulus of continuity of H\"older type. 

\subsection{De Giorgi's theorem}

A first important idea of De Giorgi in \cite{DeG56} is to consider a spatial derivative $u = \frac{\partial v}{\partial x_i}$ for some $i \in \{1,\dots,d\}$ and to remark that it satisfies
$-\dive_x( D^2 L (\nabla_x v) \nabla_x u)=0$ because the differential operator $\frac{\partial}{\partial x_i}$ commutes with the divergence one. If we consider the matrix-valued function
\[A (x) = D^2 L (\nabla_x v (x)), \]
we remark that Assumption~\eqref{convexity} implies that the ellipticity condition~\eqref{ellipticity} is satisfied. Let us recall it,
\begin{equation}
  \label{ellipticity-bis}
  \tag{Ellipticity}
  \text{There exists $\lambda, \Lambda >0$ s.t. for all $x \in \R^d$}, \quad \forall \xi \in \R^d, \lambda |\xi|^2 \le A (x) \xi \cdot \xi \le \Lambda |\xi|^2
\end{equation}
But in stark contrast with Schauder setting, the map $x \mapsto A (x)$ is not known (yet) to be  H\"older continuous.
E.~De Giorgi proposes to forget about the nonlinear equation~\eqref{e:el} and to focus on the study of the linear equation
\begin{equation}
  \label{e:ell}
  - \dive_x ( A(x) \nabla_x u ) = 0, \quad x \in \Omega.
\end{equation}
He proposes that we only retain from the nonlinear problem that the Euler-Lagrange equation is elliptic, that is to say $A$ satisfies \eqref{ellipticity}.
In doing so, E.~De~Giorgi has to work with elliptic equations with \emph{rough coefficients}: this means that the function $A$
has no regularity assumption but \eqref{ellipticity}. 
\begin{thm}[E.~De~Giorgi -- \cite{DeG56}]
  Let $\Omega$ be an open set of $\R^d$. Assume that $x \mapsto A (x)$ satisfies \eqref{ellipticity} over $\Omega$.
  Then any (weak) \emph{solution} of \eqref{e:ell} is H\"older continuous in the interior of $\Omega$. 
\end{thm}
This statement is not complete because we did not make the notion of (weak) solution precise. Moreover, De~Giorgi's theorem has
a stronger conclusion. In particular, the H\"older exponent only depends on dimension and ellipticity constants $\lambda,\Lambda$.
Such constants will be called universal. 

\subsection{Resolution of Hilbert's 19\textsuperscript{th} problem}

Thanks to De~Giorgi's theorem, we now know that all derivatives $\frac{\partial v}{\partial_i}$ of the minimizer $u$ are H\"older continuous.
In particular, the map $A(x) = D^2 L (\nabla_x v)$ is H\"older continuous in the interior of the open set $\Omega$. It is then possible to
get a Schauder theory in $\Omega$ and conclude that $v$ is twice differentiable with H\"older continuous second order derivatives. Indeed, one can write
\eqref{e:el} under the following form,
\[ - \sum_{i,j=1}^d a^{ij}(x) \frac{\partial^2 v}{\partial x_i \partial x_j} (x) =0 \text{ in } \Omega. \]
Then Schauder theory can be localized to prove that $v$ is indeed $C^2$ in $\Omega$ and second derivatives of $v$ are H\"older continuous. 
But this implies now that first order derivatives  of $a^{ij} (x)$ are H\"older continuous. In particular, first derivatives $u=\frac{\partial v}{\partial x_k}$ of $v$ satisfy,
\[ - \sum_{i,j=1}^d a^{ij}(x) \frac{\partial^2 u}{\partial x_i \partial x_j}  =S \text{ in } \Omega \]
with the H\"older continuous source terms $S = - \sum_{i,j =1}^d \frac{\partial a^{ij}}{\partial x_k} \frac{\partial^2 v}{\partial x_i \partial x_j}$. 
In particular, functions $u$ are $C^2$ with H\"older continuous second derivatives. This means that $v$ is $C^3$. We can iterate this reasoning and finally
reach the conclusion that $v \in C^\infty$.

Then proving that $C^\infty$ solutions are analytic was known at that time. We do not discuss this point since it is remotely concerned
with the regularity theory we are interested in. 

\section{Parabolic equations with rough coefficients}

\subsection{Nash's contribution}

In \cite{nash}, J.~Nash proved  De~Giorgi's theorem for parabolic equations one year later. It is somewhat a generalization of De Giorgi's theorem since solutions of elliptic equations
can be seen as time-independent solutions of parabolic equations with time-independent coefficients. J.~Nash was not aware of De~Giorgi's result before writing his paper.
Let us give more details. He considered parabolic equations in divergence form,
\[ 
\frac{\partial u}{\partial t} - \sum_{i,j=1}^d \frac{\partial}{\partial x_i} \left( a^{ij}(x) \frac{\partial u}{\partial x_j} \right) = S, \quad t>0, x \in \R^d 
\]
or equivalently,
\begin{equation}
  \label{e:para} \tag{Parabolic}
  \frac{\partial u}{\partial t} - \dive_x (A(t,x) \nabla_x u) = S, \quad t>0, x \in \R^d .
\end{equation}
We can say that J.~Nash followed the same path as De~Giorgi. In particular, he derived a modulus of continuity for solutions that does not depend on the regularity
of the coefficients $a^{ij}$. But, even if the two articles share some similarities, J.~Nash's reasoning is quite different from De~Giorgi's one.
In particular, he took inspiration from statistical mechanics (see below) and considered the logarithm of positive solutions.
This idea will be further explored by J.~Moser shortly afterwards \cite{moser}.
Moreover, J.~Nash  worked with fundamental solutions, allowing him to transform integral estimates into
point-wise ones. 

\subsection{Fluid dynamics}

O.~A.~Ladyženskaja, V.~A.~Solonnikov and N.~N.~Ural'ceva wrote in 1968 a book entitled ``Linear and quasi-linear equations of parabolic type'' \cite{L1}. 
This book will be very influential in the second half of the twentieth century. It presents in particular the notion of parabolic De~Giorgi's classes and
expose the proof of De Giorgi \& Nash's theorem for their elements.

O.~A.~Ladyženskaja is also renown for a contribution to the study of the Navier-Stokes system for in-compressible fluids.
This system has a parabolic structure since it is written,
\begin{equation}
\label{ns}  \tag{Navier-Stokes}
  \begin{cases}
    \frac{\partial U}{\partial t} + (U \cdot \nabla_x) U = \Delta U + \nabla_x P, \quad t>0, x \in \R^3, \\
    \dive_x U = 0.
  \end{cases}
\end{equation}
This being said, it is quite different from the linear and scalar parabolic equations studied by J.~Nash.
First because it is a system (the function $U$ is valued in $\R^3$). Second because of the presence of the non-linear convection term. 
The understanding of the regularity of the solutions of this system is still largely open
and one of the Millennium prize problems is devoted to it. A remarkable contribution about this important mathematical question
was made in a series of papers by Scheffer, starting with \cite{schaffer565turbulence}, and later improved by L.~Caffarelli, R.~Kohn and L.~Nirenberg \cite{MR673830}. Their results quantify
how big (or small) is the set where the solution might not be smooth.\footnote{It gives an upper bound on its (parabolic) Hausdorff dimension.} An interesting remark
for us is that  A.~F.~Vasseur \cite{MR2374209} gave an alternative proof of Caffarelli-Kohn-Nirenberg's theorem
by following closely De~Giorgi's ideas. 

Later on, L.~Caffarelli and A.~F.~Vasseur applied to other models from fluid mechanics De~Giorgi's regularity methods. 
Among these models, the surface quasi-geostrophic equation attracted a lot of attention of the mathematical community. L.~Caffarelli and A.~F.~Vasseur
managed to prove in \cite{CV1} that solutions remain smooth for all times by applying ideas of De Giorgi's type. Let us mention that another team, made up of A.~Kiselev, F.~Nazarov and A.~Volberg \cite{MR2276260}, proved simultaneously the same result by a completely different method. 

\subsection{Parabolic equations with integral diffusion}

The diffusion term in the surface quasi-geostrophic equation is the square root of the Laplacian. It can be defined through Fourier analysis, but it can also be defined thanks to a singular integral. More precisely, for $s \in (0,1)$, the fractional Laplacian $(-\Delta)^s$ of a function $u\colon \R^d \to \R$ is defined by
\[ (-\Delta)^s u (x) =  c_{d,s} \int_{\R^d} (u(x)-u(y)) \frac{\dy}{|x-y|^{d+2s}} \]
for some positive constant $c_{d,s}$ depending dimension $d$ and $s \in (0,1)$. 

The study of general parabolic equations with integral diffusion were a very active field of research for years, in particular around 2010. 
We would like to mention two landmark contributions in this direction: the article by M.~Kassmann \cite{MR2448308}  and the one by L.~Caffarelli, Chan and A.~F.~Vasseur \cite{MR2784330}. We will present results for the following class of equations,
\begin{equation}
  \label{para-int} \tag{Parabolic with integral diffusion}
  \frac{\partial f}{\partial t} = \int_{\R^d} (f(t,w) - f(t,v)) K(v,w) \dw 
\end{equation}
for some positive kernels $K$ whose structure will be discussed in due time.
It will be assumed that it is comparable to the kernel of the fractional Laplacian $K_s (v,w) = c_{d,s} |v-w|^{-d-2s}$ in
a sense that will be made precise in the chapter devoted to parabolic (and kinetic) equations with integral diffusion.

\section{Equations from kinetic theory of gases}

In this section, is quickly exposed a naive point of view on the importance of the Boltzmann equation in mathematical physics.
It is a way to introduce this central nonlinear model with which the last chapters deal. It is also a way to share with the reader the enthusiasm for its study. 

\subsection{Irreversibility}

In the 19\textsuperscript{th} century, scientists were interested in thermodynamics. This branch of physics  ``[...] developed
out of a desire to increase the efficiency of early steam engines'' (\href{http://en.wikipedia.org/wiki/Thermodynamics}{Wikipedia}). 
Some phenomena were known to be irreversible. Such a principle was first stated by Sadie Carnot and further developed by Rudolf Clausius.
It is now known as the \emph{second law of thermodynamics}. 

Irreversibility cannot be (easily) understood from classical mechanics since Newton's law are reversible in time.
A classical example of this apparent paradox is given by a gas in a box with two compartments. If the gas is initially confined
in the left compartment and the wall between the two compartments is removed, the gas will quickly occupy the entire box.
It is not expected that he could be confined again in the left part of the box for later times. However, this is compatible with
Newton's laws. 

\subsection{Statistical mechanics}

Under the influence of the both the emergence of statistics in social sciences and the development of atomistic theories in physics,
J.~C.~Maxwell derived in 1867 an equation for the probability density function $f(x,v)$ of a dilute gas. It
accounts for the number of particles with velocities $v$ at position $x$. Such a function can evolve in time and J.~C.~Maxwell showed that
it satisfies an equation of the form,
\begin{equation}
  \label{e:boltzmann0}
  \tag{Boltzmann}
  \frac{\partial f}{\partial t} (t,x,v) + v \cdot \nabla_x f (t,x,v) = Q_B(f(t,x,\cdot),f(t,x,\cdot)) (v), \quad t>0, x, v \in \R^3.
\end{equation}
The left hand side of the previous equation encodes the fact that particles travel along straight lines at velocity $v$ when they
do not collide with other particles. The operator in the right hand side accounts for collisions between particles. 
During a collision at time $t$ and position $x$, the velocities of the two colliding particles are modified. This is
the reason why the collision operator acts only on the velocity variable. It is thus applied to $f(t,x,\cdot) \colon v \mapsto f(t,x,v)$.
The notation $Q_B(f,f)$ also draw one's attention towards the fact that the operator is quadratic. For some interaction potentials
between particles, this operator $Q_B$ has a diffusive effect. The reader can be surprised that the tag of the equation is
Boltzmann and not Maxwell. And that the subscript for the collision operator is $B$. 

In 1872, L.~Boltzmann studied the long time behavior of solutions of the equation derived by J.~C.~Maxwell. He
made the seminal observation that the quantity
\begin{equation}
  \tag{Entropy} - \int (\log f) f \dx \dv 
\end{equation}
increases with time. This fact is now known as Boltzmann's H-theorem and the equation originally derived by Maxwell nowadays bears the name of Boltzmann.
The H-theorem can be thought of a quantitative version
of the second law of thermodynamics that we mentioned above. 

\subsection{The Landau equation}

Boltzmann's collision operator $Q_B$ depends on the choice of the potential from which the inter-particle force derives.
These potentials should be less singular than the Coulombian one: indeed, in the latter case, $Q_B$ contains a singularity
that is too strong for the operator to be well defined. For this reason, Lev Landau proposed in \cite{landau1936kinetische} (see also \cite{landau1980statistical}) another collision
operator $Q_L (f,f)$ in order to take into account Coulomb interactions. It has the following form,
\[ Q_L (f,f) = \dive_v ( A_f \nabla_v f) + \dive_v (f b_f) \]
for some  matrix $A_f$ and some vector field $b_f$ depending on the solution $f$. 
Since $A_f$ is semi-definite, the operator $Q_L(f,f)$ has a diffusive structure.
It is reminiscent of the class of parabolic equations in divergence form considered by J.~Nash (see \eqref{e:parabolic0} above). 
This being said, the lower order term $b_f$ turns out to be very singular, and quadratic, just like for \eqref{ns}.

\section{The Kolmogorov equation \& kinetic Fokker-Planck equations}

There are many more kinetic equations beyond the Boltzmann and the Landau equations. Moreover, from the perspective
of the study of the regularizing effect of their collision operators, an equation introduced by A.~Kolmogorov in the early years 1930 played
an important role. 

\subsection{Kolmogorov's seminal observation}

In 1934, A.~Kolmogorov computed explicitly in \cite{MR1503147} solutions $f(t,x,v)$ of the following linear equation,
\begin{equation}
  \tag{Kolmogorov}
  \frac{\partial f}{\partial t} + v \cdot \nabla_x f = \Delta_v f 
\end{equation}
thanks to Fourier analysis. He showed that this equation has a smoothing effect: rough
initial data are immediately smoothed out for positive times. This important observation
will inspire L.~H\"ormander to develop his theory of hypoellipticity in his seminal paper \cite{hormander} from 1967.

\subsection{Kinetic equations with local diffusion}

In this monograph, we will present a regularity theory for kinetic equations of the form
\begin{equation}
  \label{e:kfp0} 
  \frac{\partial f}{\partial t} + v \cdot \nabla_x f = L_v f 
\end{equation}
associated with a linear and ``diffusive'' operator $L_v$ acts only in the velocity variable. 
The presentation will follow recent contributions on this topic that we quickly review in the next two paragraphs.

In view of the discussion above, it is relevant to  consider local diffusions  of the form,
\[ L_v f =  \dive_v (A \nabla_v f) \]
for some diffusion matrix satisfying an ellipticity condition.

\paragraph{Ultraparabolic equations}
The study of such equations was first addressed by considering a more general class of equations, called
ultraparabolic equations. The Italian school contributed in an essential way to the study of these equations,
which was launched by the paper by E.~Lanconelli and S.~Polidoro \cite{MR1289901}. The latter author played a key role in
the Italian community working on these questions. Among other things, he proved with A. Pascucci that weak solutions are locally bounded \cite{pp}.
The first  De Giorgi-type result was proved by W.~Wang and L.~Zhang \cite{wz09}. Obtaining such a result was a breakthrough.
The proof relied on not-so-classical ideas due to S.~N.~Kruzhkov about classical parabolic equations. A proof closer to De Giorgi's
paper was later devised by F.~Golse, C.~Mouhot, A.~F.~Vasseur and the author of this monograph \cite{zbMATH07050183}. 

\paragraph{Conditional regularity for the in-homogeneous Landau equation}
There were many further developments to the theory that we will review in a dedicated section of the chapter
related to the study of \eqref{e:kfp0}. We would like to mention that the conditional regularity program initiated
in \cite{zbMATH07050183} was continued in \cite{MR3778645} and completed in \cite{MR4072211}. This program consists in proving that, if physically relevant density functions $\rho (t,x)$, $E(t,x)$ and $H(t,x)$ are assumed to be bounded and $\rho$ is bounded from below, the the solution $f$ of the space-in-homogeneous Landau equation is smooth. These density functions are $\rho (t,x) = \int_{\R^d} f(t,x,v) \dv$, $E(t,x) = \int_{\R^d} f(t,x,v)|v|^2 \dv$ and $H(t,x) = \int_{\R^d} (\log f) f (t,x,v) \dv$. 

\subsection{Kinetic equations with integral diffusion}

A similar conditional regularity program was completed by L.~Silvestre and the author of this monograph for the space-in-homogeneous Boltzmann equation.  
In order to obtain the final $C^\infty$ estimate on solutions in \cite{zbMATH07537725}, the authors first needed to derive a theorem à la De Giorgi, by gaining
a control on the modulus of continuity of solutions only from ``ellipticity''. 
This was achieved in \cite{MR4049224} where a general class of kinetic equations with integral diffusion is introduced,
\[ L_v f = \int_{\R^d} (f(w) -f(v)) K (v,w) \dw \]
The kernel is assumed to be ``sufficiently'' elliptic in a sense that we will make precise.

We will see that the class of kernels that one has to work with in order to deal with the Boltzmann equation are much more difficult
to handle than the one corresponding to the fractional Laplacian. But it is too early to dive into this.

\mainmatter

\chapter{De Giorgi's approach to regularity}
\label{c:dg-big}

In this first chapter, we describe the overall picture that emerges when one follows the path initially taken by E.~De Giorgi,
revealing a broader perspective beyond  elliptic, parabolic and kinetic equations.
\bigskip 

\noindent
\fbox{\begin{minipage}{\linewidth}{%
      \textit{ This book is self-contained. Moreover, it is written in such a way  that each chapter can be red independently. In particular, we will always go through
        proofs again, in full details, even when they are very similar, if not copied/pasted from a previous chapter.
      This being said, there is a progression from Chapter~\ref{c:elliptic}
      to Chapter~\ref{c:kin}, in the sense that new theoretical challenges and technical difficulties are faced when passing from the elliptic to the parabolic,
      and from the parabolic to the kinetic framework.
      In the parabolic chapter, properties related to the time variable are to be taken into account.
      In the kinetic chapter, the main new phenomenon to be understood is the transfer of regularity from the velocity variable to the spatial one. 
      If the reader is not familiar with De Giorgi's techniques and they want to get the big picture, for instance to be able to adapt these techniques
      to their framework, it is probably useful to read the first chapters, or the sections of first chapters that are related to the estimate of interest (local maximum principle,
    H\"older estimate, Harnack's inequality etc.).}
    }\end{minipage}}

\section{H\"older regularity and  oscillations}

De Giorgi's regularity consists in local H\"older estimate of solutions of partial
differential equations, and more generally of functions satisfying appropriate families of
local (energy) estimates. 

\paragraph{Point-wise H\"older regularity.}
The H\"older regularity of a function $u$ on a set $\Omega$
is characterized in terms of the oscillations of the function around a point $x_0$. 
Given a measurable set $\mathcal{N}$ (for neighborhood),
the \emph{oscillation} of the function $u$ over  $\mathcal{N}$ is defined as,
\[ \osc_{\mathcal{N}} u = \esssup_{\mathcal{N}} u - \essinf_{\mathcal{N}} u.\]
\begin{figure}[h]
\begin{minipage}{.4\linewidth}
  \centering{\includegraphics[width=4cm]{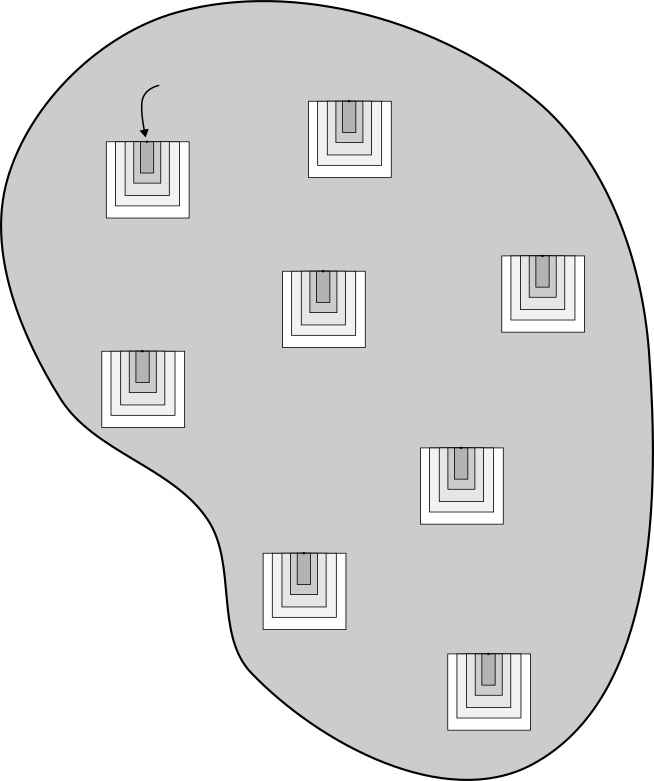}}
 \put(-83,120){$x_0$}
 \put(-120,130){\Huge $\Omega$}
\end{minipage}
\begin{minipage}{.6\linewidth}
  We aim at proving that a function $u \colon \Omega \to \R$ is H\"older continuous.
  This means that it is H\"older continuous around each point $x_0 \in \Omega$.
  \medskip
  
  Being H\"older continuous at $x_0$ for a function $u$ is equivalent to
  say that its \emph{oscillation} decays algebraically with the radius $r$ of 
  the \emph{neighborhood} $\mathcal{N}_r$.
  \medskip
  
  On this figure, neighborhoods around any point $x_0$ look like cylinders
  with the point $x_0$ at the top. This will be the case in the parabolic
  and kinetic frameworks. 
\end{minipage}
\caption{\textit{H\"older regularity at each point of $\Omega$.}}
   \label{fig:omega}
\end{figure}
We will see that H\"older regularity reduces to proving that around any point $x_0 \in \Omega$, we have
\[ \osc_{\mathcal{N}_r} u \le C r^\alpha \]
where $\mathcal{N}_r$ is a neighborhood of radius $r$ around $x_0$. 
We will see that it is enough to consider a sequence of shrinking cylinders:
\[ \forall k \ge 1, \qquad  \osc_{\mathcal{N}_{r_k}} u \le C {r_k}^\alpha  \qquad \text{ with } r_k \to 0 \text{ as } k \to \infty.\]

\paragraph{From micro to unit scale.}
A very important idea of De Giorgi's method is to reduce the proof of the
algebraic decay of the oscillation of a solution to an improvement of oscillation \emph{at unit scale}.
This improvement has to be independent of the solution, of any smoothness of coefficients, of the source term, \textit{etc}.
It shall only depend on dimension $d$ and constants characterizing the ellipticity
of the class of equations $\lambda,\Lambda$, or constants appearing in the local energy estimates. 
\begin{figure}[h]
  \centering{\includegraphics[height=5cm]{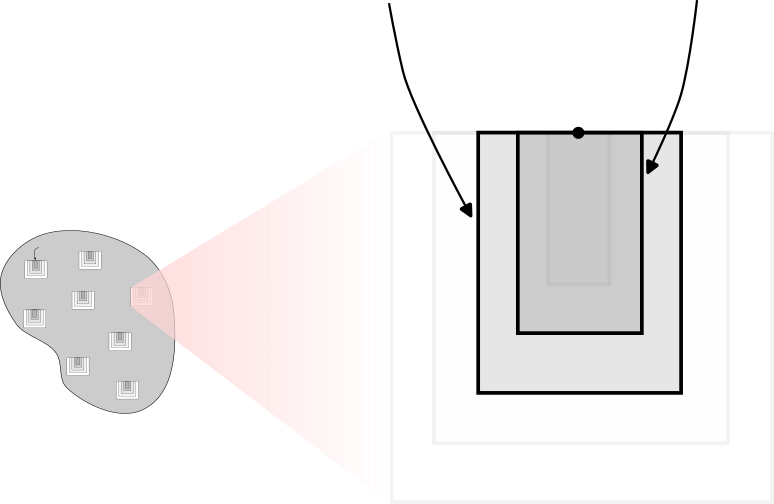}}
  \put(-19,150){$\mathcal{N}_{\frac12}$}
    \put(-115,150){$\mathcal{N}_1$}
\caption{\textit{Zooming in}: Two consecutive shrinking cylinders are scaled to $\mathcal{N}_{\frac12}$ and $\mathcal{N}_1$.} 
\label{fig:zoom}
\end{figure}

\noindent In order to prove the algebraic decay of the oscillation along a sequence, another key idea
is a consider two consecutive neighborhoods $\mathcal{N}_{k+1}$ and $\mathcal{N}_k$ and scale them so that $\mathcal{N}_k$ transforms into $\mathcal{N}_1$.
Let us assume for simplicity that then $\mathcal{N}_{k+1}$ scales into $\mathcal{N}_{\frac12}$. 
We aim at proving that 
\[ \osc_{\mathcal{N}_{\frac12}} u \le (1-\mu) \osc_{\mathcal{N}_1} u \]
for some $\mu = \mu(d,\lambda,\Lambda) \in (0,1)$. We say that the constant $\mu$ is \emph{universal}.
If we can prove such an improvement of oscillation at unit scale,
then, after scaling back, we get the algebraic decay of the oscillation of $u$ over $\mathcal{N}_{r_k}$ with $r_k = 2^{-k}$ and
$\alpha$ such that $2^{-\alpha} = (1-\mu)$. 

\paragraph{Shrinking neighborhoods: balls and cylinders.}
In order to make the previous reasoning applicable, we need neighborhoods around any point $x_0 \in \Omega$, at any scale $r >0$.
These neighborhoods shall encode two invariances of the class of equations under study:
translation (from the origin to $x_0$) and scale (by a factor $r$) invariances. 
\begin{figure}[h]
\centering{\includegraphics[height=3cm]{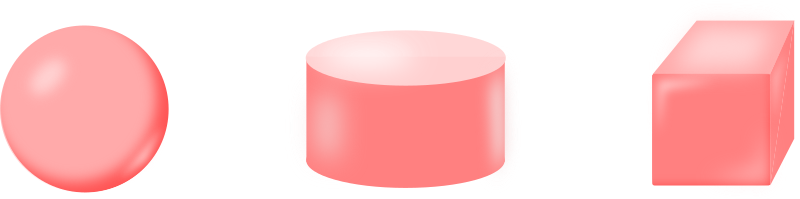}}
\caption{\textit{Neighborhoods}: balls (elliptic equations), straight cylinders (parabolic equations), slanted cylinders (kinetic equations)}
\label{fig:solides}
\end{figure}
Straight cylinders are invariant under the parabolic scaling $(t,x) \mapsto (r^2t,r x)$ while kinetic cylinders are invariant under the scaling $(t,x,v) \mapsto (r^2t,r^3x,rv)$. 
Balls and parabolic cylinders are invariant under translations $x \mapsto x_0+x$ and $(t,x) \mapsto (t_0+t,x_0+x)$, while kinetic cylinders are invariant under Galilean translation:
$(t,x,v) \mapsto (t_0+t,x_0+x + tv_0,v_0+v)$. This is the reason why they are slanted. 

\section{Local energy estimates and De Giorgi's classes}

\paragraph{From local energy estimates to improvement of oscillation.}
The local H\"older regularity (in fact the universal improvement of oscillation at unit scale) derives from
\emph{local energy estimates}. For this reason, De Giorgi's methods apply to equations from which we can derive
such local estimates. Such estimates are naturally associated with equations in divergence form.

\paragraph{De Giorgi's classes.}
In De Giorgi's original paper \cite{DeG56}, it is proven that not only solutions of elliptic equations are H\"older continuous,
but any function satisfying local energy estimates. This feature was later used in the context of calculus of variations to prove
regularity of quasi-minimizers of some functionals \cite{MR666107,MR778969}. 
These local energy estimates have to be satisfied after truncating the candidate $u$ (from below and from above) by an arbitrary constant $\kappa$:
one has to consider $(u-\kappa)_+$ and $(u-\kappa)_-$ with $a_\pm = \max(0, \pm a)$. 

\paragraph{Sub- and super-solutions.}
A sub-solution to the elliptic equation $-\dive_x (A \nabla_x u) = S$ is a function satisfying $-\dive_x (A \nabla_x u) \le S$ (in the sense of distributions). 
If a candidate $u$ is a sub-solution (resp. super-solution)
of the equation, then it satisfies the local energy inequalities after truncation from above $(u-\kappa)_+$ (resp. from below, $(u-\kappa)_-$). 
A general fact is that sub-solutions are contained in ``positive'' De Giorgi's classes ($\mathrm{DG}^+,\mathrm{pDG}^+, \mathrm{kDG}^+$) while
super-solutions  are in ``negative'' De Giorgi's classes ($\mathrm{DG}^-,\mathrm{pDG}^-, \mathrm{kDG}^-$).

\section{Local maximum principle and gain of integrability}

\paragraph{From measure to point-wise.}
We saw in the previous section that the relevant information to establish H\"older regularity is contained in two families of local energy estimates. 
This information is thus encoded in some estimates in Lebesgue spaces on the (truncated) function and some of its derivatives. Since the goal is to
estimate the oscillation of the candidate $u$, one has to derive point-wise information from local energy estimates: this is achieved through \emph{maximum principle}. 

\paragraph{Results in the book.} Propositions~\ref{p:lmp} (elliptic), \ref{p:lmp-parab} (parabolic), \ref{p:lmp-kinetic} (kinetic FP). 

\paragraph{Maximum principle from gain of integrability.}
The first step in establishing the improvement of oscillation of an element of the positive De Giorgi's class is to prove that it is locally bounded.
In order to do so, we consider the square of $(u-\kappa_k)_+$ on a shrinking neighborhood $\mathcal{N}_{r_k}$ with $\kappa_k$ increasing from $1$ to $2$ and
$r_k$ decreasing from $1$ to $1/2$.
\[ A_k := \int_{\mathcal{N}_{r_k}} (u-\kappa_k)_+^2.\]
De Giorgi's method consists in establishing the following nonlinear iterative estimate, \medskip

\noindent \fbox{\begin{minipage}{.99\linewidth}
 \centering{ \(\forall k \ge 1, \qquad A_{k+1} \le C^{k+1} A_k^\beta \qquad \) for two universal constants $C \ge 1$ and $\beta >1$. }
\end{minipage}}
\vspace{.25mm}

\noindent We recall that a constant is universal if it only depends on dimension and ellipticity constants. 
If the first term of the sequence $A_0$ is small, then the iterative estimate implies that $A_k \to 0$ as $k \to \infty$.
But $A_k \to \int_{\mathcal{N}_{\frac12}} (u -2)_+^2$ as $k \to \infty$ and $A_0 \le \int_{\mathcal{N}_1} u_+^2$.
We thus conclude that $u \le 2$ in $\mathcal{N}_{\frac12}$ as soon as the $L^2$-norm of $u_+$ in $\mathcal{N}_1$ is small enough.
\noindent \fbox{\begin{minipage}{.99\linewidth}
\begin{equation}
\label{e:lmp}
  \int_{\mathcal{N}_1} u_+^2 \le \delta_0 \quad \Rightarrow \quad u \le 2 \quad \text{ a.e. in } \mathcal{N}_{\frac12}.
\end{equation}
\end{minipage}
}

It can be surprising that we obtain a nonlinear estimate ($A_{k+1} \le C^{k+1} A_k^\beta$) for a linear equation.
This is made possible thanks to some local \emph{gain of integrability} and the truncation procedure: one way or the other, one has to prove that it is possible to control
\[ B_k := \left(\int_{\mathcal{N}_{r_k}} (u-\kappa_k)_+^p \right)^{\frac2{p}}\]
for some universal $p >2$. We can then use Bienaymé-Chebyshev's inequality and H\"older's inequality to make $A_k^\beta$ appear
for some universal $\beta >1$. More precisely, we prove that local energy estimates yield,
\[ B_{k+1} \le C^{k+1} A_k .\]

Then we  estimate the $L^2$-norm of the truncated function $(u-\kappa_{k+1})_+ = (u-\kappa_{k+1})_+ \un_{\{ u \ge \kappa_{k+1}\}}$
by the product of the $L^p$-norm of $(u-\kappa_{k+1})_+$  with the $L^{q}$-norm of $\un_{\{ u \ge \kappa_{k+1}\}}$ with $\frac1p + \frac1q =1$.
\begin{align*}
  A_{k+1} & = \int_{\mathcal{N}_{k+1}} (u-\kappa_{k+1})^2 \\
          & \le B_{k+1} \left( \int_{\mathcal{N}_{k+1}} \un_{\{ u \ge \kappa_{k+1}\}} \right)^{\frac2q} \\
      & \le C^{k+1} A_k  \left| \{ u \ge \kappa_{k+1}\} \cap \mathcal{N}_{k+1} \right|^{\frac2q}  .
\end{align*}
Here is the trick: write $u - \kappa_{k+1} = (u-\kappa_k) - (\kappa_k - \kappa_{k+1})$ and let $\delta_k$ denote $\kappa_k - \kappa_{k+1}>0$.
Then Bienaymé-Chebyshev's inequality implies that
\begin{align*}
  \left| \{ u \ge \kappa_{k+1}\} \cap \mathcal{N}_{k+1} \right| & \le \delta_k^{-2} \int_{\mathcal{N}_{k+1}} (u-\kappa_k)_+^2 \\
  & \le \tilde C^k A_k
\end{align*}
(we used that $\mathcal{N}_{k+1} \subset \mathcal{N}_k$). We thus get $A_{k+1} \le \bar C^{k+1} A_k^{1+ \frac2q}$ and we get the desired result with $\beta = 1 + \frac2q >1$ (universal). 

\section{Improvement of oscillation through expansion of positivity}

\paragraph{Expansion of positivity.}
The expansion of positivity states that any non-negative function $u$  from the (negative) De Giorgi's class  satisfies, \smallskip

\noindent \fbox{\begin{minipage}{.99\linewidth}
\[ |\{u \ge 1 \} \cap \mathcal{N}_{\mathrm{pos}} | \ge \frac12 |\mathcal{N}_{\mathrm{pos}}| \quad \Rightarrow \quad \{ u \ge \ell \quad \text{ a.e. in } \mathcal{N}_1\}\]
\end{minipage}} \vspace{0.1ex}

\noindent for some universal constant $\ell >0$. 
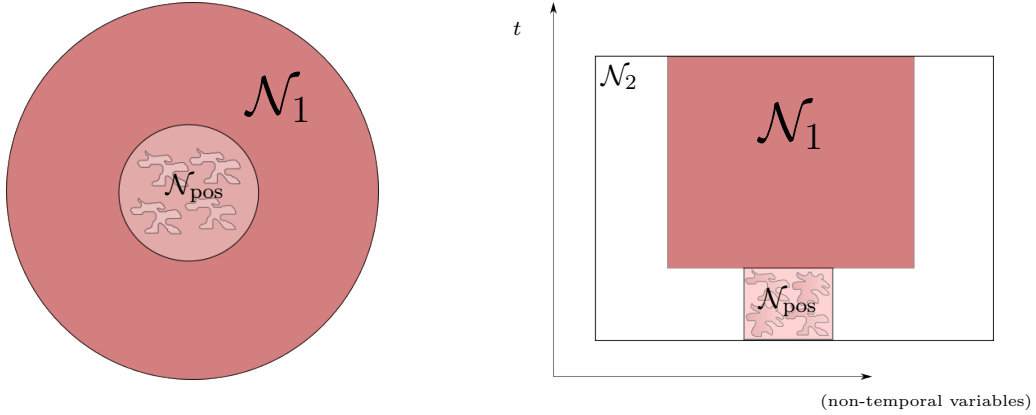
\begin{figure}[h]
  \centering{
    \includegraphics[height=5cm]{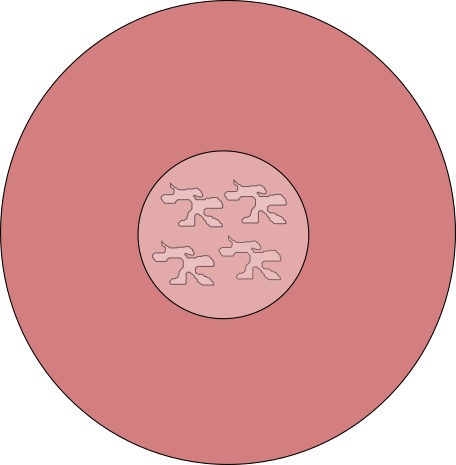} \hspace{2cm}
    \includegraphics[height=5cm]{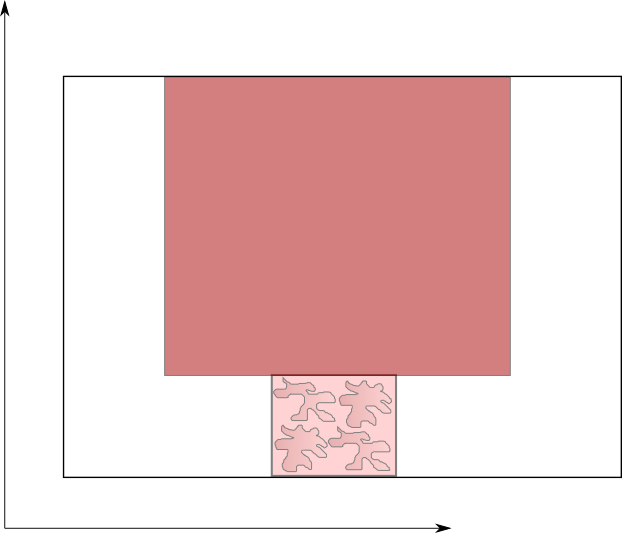}
  }
  \put(-280,100){\huge $\mathcal{N}_1$}
    \put(-310,70){$\mathcal{N}_{\mathrm{pos}}$}
  \put(-180,130){\scriptsize $t$}
 \put(-147,111){$\mathcal{N}_2$}
 \put(-88,90){\huge $\mathcal{N}_1$}
 \put(-88,27){$\Ngu$}
 \put(-65,-10){\tiny (non-temporal variables)}
 \caption{\textit{Expansion of positivity.}
 A lower bound on the super-level set on $\mathcal{N}_{\mathrm{pos}}$ implies a point-wise lower bound in $\mathcal{N}_1$.
   On the left, illustration of the expansion of positivity for elliptic equations.
   On the right, the parabolic and kinetic cases. For time-dependent equations,  positivity is expanded as time increases.}
\label{fig:expansion-positivity}
\end{figure}
We say that positivity is extended because if $u$ is positive in half of the small neighborhood $\mathcal{N}_{\mathrm{pos}}$ then it is
positive  \emph{(almost) everywhere} in the larger neighborhood $\mathcal{N}_1$. The price to pay is that the (universal) lower bound
deteriorates. 
\medskip

This property implies the decrease of oscillation of $u$. Indeed, the local maximum principle implies that the solution
is locally bounded in a set $\mathcal{N}_2$ containing $\mathcal{N}_1$. Using the linearity of the equation (encoded in the definition of De Giorgi's classes),
we reduce to the case where $0 \le u \le 2$ a.e. in $\mathcal{N}_2$.
In particular $\osc_{\mathcal{N}_2} u \le 2$. We distinguish two cases.
\begin{itemize}
\item If $|\{u \ge 1 \} \cap \mathcal{N}_{\mathrm{pos}} | \ge \frac12 |\mathcal{N}_{\mathrm{pos}}|$ then $\inf_{\mathcal{N}_1} u \ge \ell$. In particular, $\osc_{\mathcal{N}_1} u \le 2 -\ell$. 
\item If $|\{u \ge 1 \} \cap \mathcal{N}_{\mathrm{pos}} | < \frac12 |\mathcal{N}_{\mathrm{pos}}|$ then we can apply the previous case to $v = 2 -u$ and get
  $\inf_{\mathcal{N}_1} v \ge \ell$. In particular, $\osc_{\mathcal{N}_1} u \le 2 -\ell$ in this case too. 
\end{itemize}
We thus proved, \vspace{1ex}

\noindent \fbox{\begin{minipage}{.99\linewidth}
    \[ \osc_{\mathcal{N}_2} u \le 2 \quad \Rightarrow \osc_{\mathcal{N}_1} u \le 2 -\ell.\]
  \end{minipage}} \vspace{.01ex}

\noindent Making a long story short, by scaling the neighborhood by a factor $2$ (from $\mathcal{N}_1$ to $\mathcal{N}_{\frac12}$), we gain a universal factor $1-\ell/2 \in (0,1)$ on the oscillation
of the function. 

\paragraph{Results in the book.} Corollary~\ref{c:spreading} (elliptic), \ref{p:expansion-parab} (parabolic), \ref{p:expansion-kin} (kinetic FP).

\section{The intermediate value principle}

We finally explain how the expansion of positivity is obtained thanks to an intermediate value principle.

\paragraph{The local maximum principle, upside down.} The conclusion of the expansion of positivity is a point-wise lower bound on $u$.
Such a point-wise lower bound can be obtained by applying the local maximum principle to the function $1-u$. Indeed, we can apply
\eqref{e:lmp} to $v =1-u$ and get that $v \le 1/2$ as soon as $\int_{\mathcal{N}_1} v_+^2 \le \eps_1 |\mathcal{N}_1|$ for some universal constant $\eps_1$.
We then estimate this $L^2$-norm by the measure of the sub-level set of $u$,
\[ \int_{\mathcal{N}_1} v_+^2 \le |\{ v \ge 0 \} \cap \mathcal{N}_1 | = |\{ u \le  1 \} \cap \mathcal{N}_1 |. \]
In practice, we are going to apply this reasoning to $\theta^k u$ for some well chosen $\theta$ and an integer $k$. 
In conclusion, the local maximum principle (upside down) asserts that, \smallskip

\noindent \fbox{\begin{minipage}{.99\linewidth}
\begin{equation}
\label{e:lmp-usd}
  |\{ u \le  1 \} \cap \mathcal{N}_1 | \le \eps_1 |\mathcal{N}_1| \quad \Rightarrow \quad \bigg\{ u \ge \frac12 \text{ a.e. in } \mathcal{N}_{\frac12} \bigg\}.
\end{equation}
\end{minipage}}

\paragraph{The intermediate value principle for elliptic equations.} For elliptic equations, the intermediate value principle quantifies the fact
that a function with a square integrable gradient cannot jump for $0$ to $1/2$. More precisely, for $u \in H^1(B_1)$, 
\[  \Rightarrow |\textcolor{OliveGreen}{\{u \le 0 \}} \cap B_1| \times |\textcolor{OliveGreen}{\{ u \ge 1/2\}} \cap B_1| \le C \|\nabla_x u \|_{L^2 (B_1)}
  |\textcolor{BrickRed}{\{0 < u < 1/2 \}} \cap B_1|^{\frac12}  .\]
The important consequence of this estimate is that, if we have lower bounds $\textcolor{OliveGreen}{\delta_1,\delta_2}$
for the measures of sub- and super-level sets $\{ u \le 0\}$ and $\{ u \ge 1/2\}$,
then we have a lower bound $\textcolor{BrickRed}{\delta_{1,2}}$ on the intermediate value set $\{0 < u < 1/2 \}$. This is the way
the intermediate value principle is stated in general, as we shall see below  for kinetic Fokker-Planck equations. 

 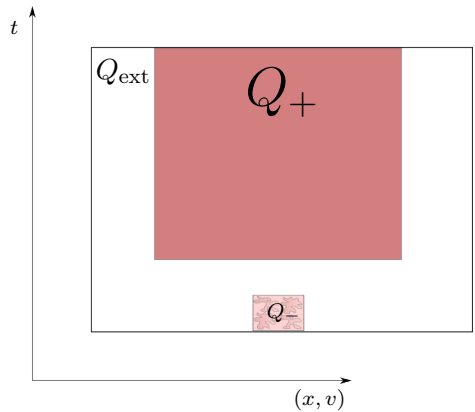
\begin{figure}[h]
 \centering{\includegraphics[height=5cm]{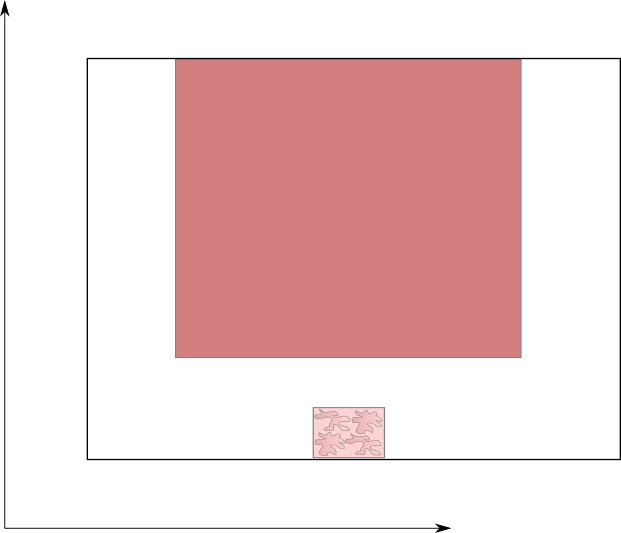}}
 \put(-173,132){\scriptsize $t$}
 \put(-141,115){$\Qext$}
 \put(-85,105){\huge $Q_+$}
 \put(-77,25){\scriptsize $Q_-$}
 \put(-67,-8){\scriptsize $(x,v)$}
 \caption{\textit{Geometric setting of the intermediate value principle.}}
   \label{fig:iv-principle-intro}
 \end{figure}

\paragraph{The intermediate value principle for kinetic Fokker-Planck equations.} For evolution equations, in particular for kinetic Fokker-Planck equations,
the geometric setting of the intermediate value principle is made of two cylinders, one sitting in the past, let us call it $Q_-$, and one sitting in the future, we call it $Q_+$.
It is necessary that they are separated by a time lap. Moreover, they are both contained in a large cylinder $\Qext$. The intermediate value principle asserts that, given the geometric setting (that is radii of cylinders $Q_\pm$ and $\Qext$ and time lap),  and for $\delta_1$ and $\delta_2$ given, there exists  constants $\theta \in (0,1)$ and $\delta_{1,2}$ such that
\vspace{1mm}

\noindent \fbox{\begin{minipage}{.99\linewidth}
\[
  \left.
  \begin{array}{r}
    | \{ f \ge 1 \} \cap Q_-| \ge \textcolor{OliveGreen}{\delta_1} |Q_-| \\
    |\{ f \le \theta \} \cap Q_+| \ge \textcolor{OliveGreen}{\delta_2} |Q_+|
  \end{array}
  \right\}
    \Rightarrow \quad |\{ \theta < f < 1 \} \cap \Qext | \ge \textcolor{BrickRed}{\delta_{1,2}} |\Qext|.
  \]
\end{minipage}
} \vspace{.1mm}

\noindent  In practice, all the constants $\theta,\delta_1,\delta_2, \delta_{1,2}$ are universal. 

\paragraph{How to apply it to prove the expansion of positivity.} We consider $f_k = \theta^{-k} f$. The assumption of the expansion of
positivity ensures that $|\{ f \ge 1 \} \cap Q_-| \ge \frac12|Q_-|$.
This implies that for all $k \ge 1$, we have $|\{ f_k \ge 1 \} \cap Q_-| \ge \frac12|Q_-|$.

If we can find $k \ge 1$, such that $|\{ f_k \le 1 \} \cap Q_-| \le \eps_1 |Q_-|$, then the local maximum principle upside down -- see \eqref{e:lmp-usd} --
implies that $f_k \ge \frac12$ in the ``interior'' of $Q_+$. 

But if $|\{ f_k \le 1 \} \cap Q_-| > \eps_1 |Q_-|$, or equivalently if $|\{ f_{k-1} \le \theta \} \cap Q_-| > \eps_1 |Q_-|$, since we already know that $|\{ f_{k-1} \ge 1 \} \cap Q_-| \ge \frac12|Q_-|$, then the intermediate value principle ensures that \( |\{ \theta < f_{k-1} <1  \} \cap \Qext| \ge \delta_{1,2} |\Qext|\)
for some universal $\delta_{1,2}$. In terms of intermediate value sets of $f$, this means that
\[ |\{ \theta^k < f < \theta^{k-1}  \} \cap \Qext| \ge \delta_{1,2} |\Qext|.\]

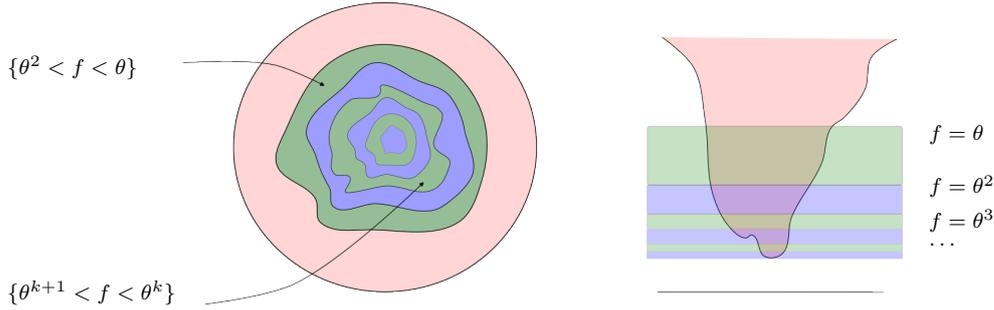
\begin{figure}[h]
\centering{\includegraphics[height=4cm]{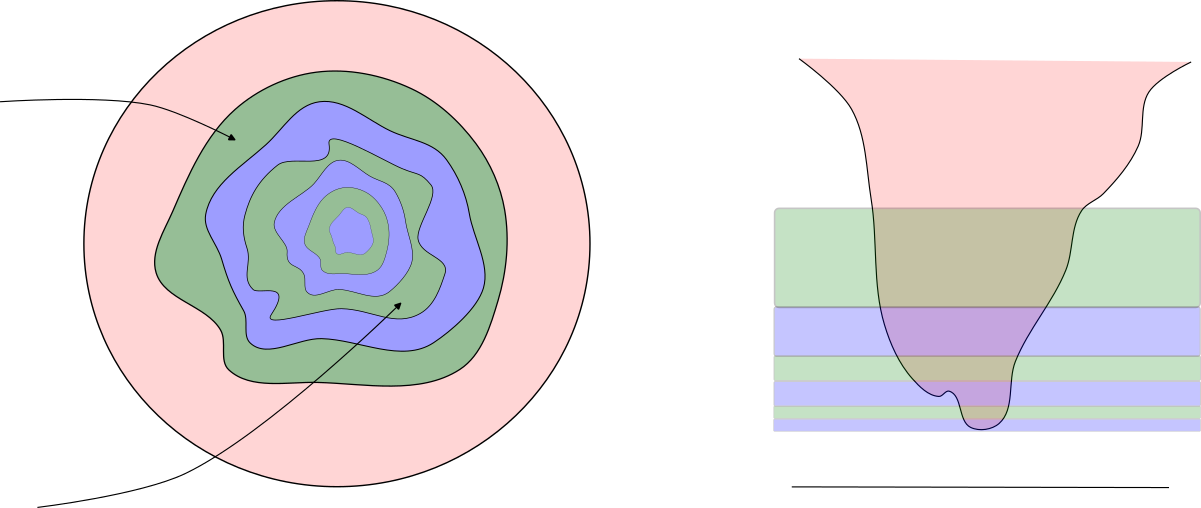}}
\put(-335,2){\scriptsize $\{ \theta^{k+1} < f < \theta^k\}$}
\put(-335,87){\scriptsize $\{ \theta^2 < f < \theta\}$}
\put(10,62){\scriptsize $f = \theta$}
\put(10,42){\scriptsize $f = \theta^2$}
\put(10,29){\scriptsize $f = \theta^3$}
\put(10,22){\scriptsize $\dots$}
\caption{\textit{Intermediate value sets}: On the left, is represented the neighborhood (ball or cylinder)
  where the function $u$ is studied. The green and blue rings correspond to
intermediate value sets $\{ \theta^{k+1} < f  < \theta^k \}$. 
}
\label{fig:ivs}
\end{figure}

But these intermediate value sets are distinct, and they occupy a universal proportion $\delta_{1,2}$ of the cylinder $\Qext$ (see Figure~\ref{fig:ivs}). 
This implies that there is only a finite number of them. In particular for $k$ large enough, we do have $|\{ f_k \le 1 \} \cap Q_-| \le \eps_1 |Q_-|$ and this produces
a lower bound on $f_k$ (and thus on $f$) in the interior of $Q_+$. 

\paragraph{Proof of the intermediate value principle.}
In the elliptic case, the intermediate value principle is a straightforward consequence of a Poincaré-Wirtinger's inequality.
In the parabolic case, the same inequality is used after freezing the time variable.
In the kinetic case, the proof requires to establish a Poincaré-Wirtinger's inequality for weak sub-solutions, involving the cylinders $Q_-$ and $Q_+$. 

\paragraph{Corresponding results in the book.} Lemma~\ref{l:ivl-elliptic} (elliptic), proof of Lemma~\ref{l:spreading} (parabolic), Proposition~\ref{p:iv-principle} (kinetic Fokker-Planck).

\section{Summary / conclusion}

De Giorgi's theorem is proved in two steps: \medskip

\noindent \fbox{\begin{minipage}{.99\linewidth}
\begin{itemize}
\item Square integrable solutions are  locally bounded (\textbf{Local maximum principle}). \vspace{-1ex}
\item Locally bounded solutions are  H\"older continuous (\textbf{improvement of oscillation}). 
\end{itemize}
\end{minipage}
}
\bigskip

\noindent Improvement of oscillation derives from two  principles: \medskip

\noindent \fbox{\begin{minipage}{.99\linewidth}
\[
  \left.\begin{array}{r} \text{Local maximum principle {\scriptsize (upside down)}} \\[1ex] \text{\textbf{Intermediate value principle}} \end{array} \right\}
  \Rightarrow \begin{array}{c} \text{Expansion} \\ \text{of positivity} \end{array} \Rightarrow \begin{array}{c}\text{Improvement} \\ \text{of oscillation} \end{array}
\]
\end{minipage}}

\section{Parameters}

There are a lot of parameters in this theory. We try to make notation as homogeneous as possible.
Let us present the main ones. 
\begin{itemize}
  \item $\lambda$ and $\Lambda$ denote the ellipticity constants. 
  \item $x_0$, $X_0$, $z_0$ denotes the ``center'' of a ``neighborhood'', that is to say of balls or cylinders. 
  \item $r$ and $R$ denote a small and a large radius for balls or cylinders (for instance in local energy estimates). 
  \item $\kappa$ denotes the truncation parameter.
  \item $\eps$ denotes a small parameter in approximation arguments.
  \item $p,q,p_\ast$ \textit{etc.} denotes Lebesgue exponents. 
  \item $\eps_0, \eps_{0,1},\bar \eps_0$ \textit{etc.} denotes the ``smallness'' of source terms.
  \item $\eps_1$ denotes the ``small'' proportion of a cylinder yielding a local point-wise (upper or lower) bound.
  \item $\eta_0,\eta_1,\eta_2$ and $\omega$ typically denote small radii in geometric settings.
  \item $\beta$ denotes a universal constant larger than $1$ in De Giorgi's iteration method for the local maximum principle.  
  \item $\omega_0$ denotes a parameter in the local maximum principle between a large cylinder $Q_R$ and a small cylinder $Q_r$.
\end{itemize}


\chapter{Elliptic equations}
\label{c:elliptic}

This chapter is devoted to the case of elliptic equations under divergence form.
It corresponds to De Giorgi's original framework. We aim at deriving a local H\"older estimate
for weak solutions.

We will also establish that the H\"older exponent, as well as the constant in the estimate,
are \emph{universal}: they depend on very few parameters, namely dimension and ellipticity constants. 

In the course of the reasoning, we will identify classes of functions that satisfy De Giorgi's theorem.
Their are coined as De Giorgi's classes. 

\section{Elliptic equations}

This section is devoted to the presentation of the class of elliptic equations
that we consider. We first introduce notation for Euclidean balls and differential operators (partial derivative, gradient, divergence). 

\subsection*{Balls, gradient and divergence}

\label{balls}
\begin{itemize}
\item
  Given $x_0 \in \R^d$ and $r>0$, \index{$B_r (x_0)$} denotes the open ball centered
  at $x_0$ of radius $r$. The closed ball is denoted by \index{$\bar B_r (x_0)$}.
  If $x_0=0$, we simply write $B_r$ and $\bar B_r$. 
  \item
    Given a function $u \colon \Omega \to \R$: $\partial_i u$ denotes the
    partial derivative of $u$ with respect to the real variable $x_i$, 
and $\nabla_x u$ denotes the gradient of  $\nabla_x u = (\partial_{1} u, \cdots, \partial_{d} u) \in \R^d$. 
\item
Given a vector field $F \colon \Omega \to \R^d$, 
$\dive_x  F$ denotes its divergence: $\dive_x F = \sum_{i=1}^d \partial_i F_i$.
\end{itemize}

\subsection*{Ellipticity}

We consider the following class of elliptic equations: \label{elliptic}
\[ - \dive_x  ( A \nabla_x u ) = S, \qquad x \in \Omega \]
where $\Omega$ denotes an open set of $\R^d$.
The function $u \colon \Omega \to \R$ is called the \emph{solution} of the equation while $S \colon \Omega \to \R$ is called the \emph{source term}.
The function $S$ is given and we aim at studying the function $u$. The function $A$ is also defined in $\Omega$ but it is matrix-valued. More precisely,
it takes values in  the set $\mathcal{S}_d (\R)$ of real symmetric $d \times d$ matrices. \label{symmetric}
It satisfies the following \emph{ellipticity condition},
\begin{ellipticity} There exists $\lambda, \Lambda >0$ such that
  for a.e. $x \in \Omega$, 
  \[ \forall \xi \in \R^d,  \quad \lambda |\xi|^2 \le A (x) \xi \cdot \xi \le \Lambda |\xi|^2.\]
\end{ellipticity} 
\noindent It is convenient to simply write $A \in \mathcal{E}(\lambda, \Lambda)$ for the set of $A$'s satisfying the previous condition.
We remark that the ellipticity condition is equivalent to impose that for a.e. $x \in \Omega$, the eigenvalues of the real symmetric matrix $A(x)$ lie in the interval $[\lambda,\Lambda]$. 

When studying these equations, we will see that the lower bound on eigenvalues will allow us to control the gradient of the solution in the set  $L^2$ of square integrable functions.
The upper bound will be used in order to get such a control, but localized on balls (See the Caccioppoli estimate contained in Proposition~\ref{p:DGclass}). 

De Giorgi's theorem asserts that under this mere assumption on the coefficients (and suitable source terms), one can control the modulus of continuity of solutions. 

\subsection*{Scaling and translation}

Let $u$ be a solution of $- \dive_x ( A \nabla_x u ) = S$ in
a ball $B_R (x_0)$. Then the function $v (y) = u ( r y)$ is
a solution of $-\nabla_y \cdot ( \bar A \nabla_y v) = \bar S$ in the ball $B_{rR}(rx_0)$
with $\bar A(y) = A (ry)$ and $\bar S (y) = r^2 S (ry)$. Remark in particular that if $A$
is elliptic ($A \in \mathcal{E}(\lambda,\Lambda)$) then so is $\bar A$:
we have $\bar A \in \mathcal{E}(\lambda,\Lambda)$. In other words, the class
of elliptic equations that we work with is invariant under scaling: ellipticity
constants are conserved (and source terms are scaled).

The class of elliptic equations is also translation invariant.
Indeed, the function $w(y) = u (y_0+y)$ satisfies $-\nabla_y \cdot ( \tilde A \nabla_y v) = \tilde S$ in the ball $B_{rR}(rx_0)$
with $\tilde A(y) = A (y_0+y)$ and $\tilde S (y) =  S (y_0+y)$. In particular,
$\tilde A \in \mathcal{E}(\lambda,\Lambda)$ if $A \in \mathcal{E}(\lambda,\Lambda)$. 

\subsection*{H\"older continuity and oscillation}

A function $u$ is $\alpha$-H\"older continuous in a  set $F \subset \R^d$ if for any $x,y \in F$,
\( |u(x)-u(y)| \le C |x-y|^\alpha.\) The set of $\alpha$-H\"older continuous functions on $F$ are
denoted by $C^\alpha (F)$. The semi-norm $[\cdot]_{C^\alpha (F)}$ is defined by
\[ [u]_{C^\alpha (F)} = \sup_{\stackrel{x,y \in F}{ x \neq y}} \frac{|u(x)-u(y)|}{|x-y|^\alpha}.\]
This space is equipped with the norm $\|u\|_{C^\alpha (F)} = \|u\|_{C(F)} + [u]_{C^\alpha (F)}$
where $\|u\|_{C(F)}$ denotes $\sup_F |u|$.

The H\"older continuity of an essentially bounded function at a point $x_0$ can be established by studying its oscillation around $x_0$.
\begin{defi}[Oscillation] \label{d:osc}
  Let $\Omega$ be an open set and $u \in L^\infty (\Omega)$ be real valued and an open set $\omega \subset \Omega$. The \emph{oscillation} of $u$ in $\omega$
  is defined by
  \[ \osc_\omega u = \esssup_\omega u - \essinf_\omega u .\]
\end{defi}
\begin{prop}[Characterization of H\"older continuity] \label{p:holder-osc}
  Let $B$ be an open ball and $u \in L^\infty (B)$. Assume that for all $x_0 \in B$ and all $r>0$,
  \[ \osc_{B_r (x_0) \cap B} u \le C r^\alpha.\]
  Then $u$ is $\alpha$-H\"older continuous in $B$ and
  \( [u]_{C^\alpha (B)} \le  C .\)
\end{prop}
\begin{proof}
The proof proceeds in three steps. 
  
\paragraph{Step 1.}  Assume that $u \colon B \to \R$ is continuous and consider $x,y \in B$.
  We now consider $r = |x-y|$ and remark that $y \in \bar B_r(x)$. In particular,
  \[ u(y) - u(x) \le \osc_{\bar B_r(x) \cap B} u = \osc_{B_r(x) \cap B} u \le C r^\alpha = C |x-y|^\alpha.\]
  We conclude that   \( [u]_{C^\alpha (B)} \le  C \) for $u$ continuous.

  If now $u$ is merely essentially bounded in $B$, we argue by approximation and consider
  a mollifier $\rho \colon \R^d \to \R$ with $\rho \ge 0$, $\rho \in C^\infty (\R^d)$, compactly supported in $B_1$ and $\int_{\R^d} \rho (x) \dx =1$.
  We then define for any $\eps >0$ and $x \in \R^d$,
  \[ u^\eps (x) = \int_{B} u(y) \rho_\eps (x-y) \dy . \]

\paragraph{Step 2.}
  Let $B=B_R (x_R)$ for some $x_R \in \R^d$ and $R>0$. We claim that for all $x_0 \in B^\eps = B_{R-\eps}(x_R)$ and all $r>0$,
  \[ \underset{B_r (x_0) \cap B^\eps}{\osc} u^\eps \le C r^{\alpha}.\]
  Remark that for $x \in B^\eps$, we have
  \[ u^\eps (x) = \int_{B}  u(y) \rho_\eps (x-y) \dy = \int_{\R^d}  u (x-\eps z) \rho (z) \dz \]
  since $x - \eps z \in B$. 
  In particular,
   \begin{align*}
    \underset{B_r (x_0) \cap B^\eps}{\osc} u^\eps & = \sup_{x \in B_r (x_0) \cap B^\eps} \int_{\R^d} u(x-\eps z) \rho (z) \dz - \inf_{x \in B_r (x_0) \cap B^\eps} \int_{\R^d} u(x-\eps z) \rho (z) \dz \\
                                         & \le \int_{\R^d} \left( \osc_{B_r (x_0-\eps z) \cap B} u \right) \rho (z) \dz 
   \end{align*}
   and the claim follows.

   \paragraph{Step~3.}
   Let $\eps_0>0$ and $\eps \in (0,\eps_0)$. We conclude from Step~1 that for all $x,y \in B^{\eps_0}$,
   \[ |u^\eps (x) - u^\eps (y) | \le C |x-y|^\alpha.\]
   By dominated convergence, we have that $u^\eps \to u$ a.e. in $B^{\eps_0}$, and we conclude that
   \[ |u(x) - u(y) | \le C |x-y|^\alpha.\]
   Since $\eps_0>0$ is arbitrarily small, we conclude that $[u]_{C^\alpha (B)} \le C$. 
\end{proof}

\subsection*{Energy and integration by parts}

The elliptic equations of the form $-\dive_x (A \nabla_x u) = S$ are
called to be in divergence form, or conservative form. The reason is that
a natural ``energy'' is associated to them,
\[ E(u) :=\int_{\Omega} A \nabla_x u \cdot \nabla_x u .\]
The ellipticity assumption can be interpreted as a condition under which
the energy behaves like the $L^2$ norm of the gradient.

\paragraph{Local energy.}
Let us understand why this quantity is naturally associated with the
elliptic equations that we presented above. 
In order to do so, we consider local energies by considering a ball $B \subset \Omega$,
\[ E_B (u) := \int_{B} A \nabla_x u \cdot \nabla_x u. \]
If a function $u \colon B \to \R$ is such that this local energy is minimal when $u$ is perturbed by a $C^\infty_c$ function
$\varphi$ supported in $B$, then 
\begin{equation} \label{weak}
  \tag{weak} \int_B A \nabla_x u \cdot \nabla_x \varphi =0.
\end{equation}
Indeed, the real function
\[ E_B (u+ t \varphi) = \int_B A (\nabla_x u + t \nabla_x \varphi) \cdot (\nabla_x u + t \nabla_x \varphi) \]
has to achieve a minimum at $t=0$, in particular its derivative vanishes
at $t=0$. This implies \eqref{weak}.  Integrating by parts (if possible),
we recover the elliptic equation in $B$.

\paragraph{Towards a weak formulation.} Recall that we do not want to make
any smoothness assumption on $A$, so that (among other things) we can scale
solutions and stay in the class of solutions of elliptic equations. We
do not want to differentiate $A \nabla_x u$, at any cost. Integrating by
parts, we can think of the elliptic equation under the form~\eqref{weak}. Such a form
makes sense as soon as $\nabla_x \varphi$ and $\nabla_x u$ are
square integrable function (because eigenvalues of $A$ are bounded from
above). 

\section{The  space $H^1(\Omega)$ and weak solutions}

In this section, we recall the definition of the Sobolev space $H^1(\Omega)$ that is needed in order to define weak solutions.
We will also state and prove a result about the composition of $H^1$ functions with Lipschitz ones. This result is classical
but a little bit less than the other ones. We thus state it and prove it. It is useful for the reader that is not completely
at ease with the functional setting. We then recall some classical functional inequalities and give precise references for proofs. 

\subsection{Two Hilbert spaces}

The $H^1$ space is made of square integrable functions whose gradient is also square integrable. 
\begin{defi}[The space $H^1 (\Omega)$]
  The space $H^1(\Omega)$ is the vector subspace of $L^2 (\Omega)$ made of functions $u$ admitting first order derivatives (in the sense of distributions) lying in $L^2(\Omega)$.
  The (weak) gradient of an element $u$ of $H^1(\Omega)$ lies in $L^2 (\Omega)^d$. 
\end{defi}
We want to pass to the limit in test functions $\varphi$. This is the reason
why we introduce the following subset of $H^1(\Omega)$. It is a convenient way to impose $u=0$ at the boundary
in a weak sense.
\begin{defi}[The space $H^1_0(\Omega)$]
We denote by  $H^1_0(\Omega)$  the closure of $C^\infty_c (\Omega)$ with respect to the topology induced by the
  $H^1$ norm (associated with  $(\cdot,\cdot)$). 
\end{defi}

We recall that smooth functions are dense in $H^1(\Omega)$. See for instance \cite[Corollary~9.8]{MR2759829}.
\begin{prop}[Density of smooth functions in $H^1(\Omega)$] \label{p:density}
  Let $\Omega$ be open with $C^1$ boundary and $u \in H^1 (\Omega)$.  Then there exists a sequence $U_n \in C_c^\infty(\R^d)$ such that
  $ u_n = U_n \vert_\Omega \to u$ in $H^1(\Omega)$. 
\end{prop}
This proposition immediately implies the following technical lemma.
\begin{lemma} \label{l:localization}
  Let $\Omega$ be open with $C^1$ boundary and $u \in H^1(\Omega)$ and $\varphi \in C^\infty_c (\Omega)$. Then
  $\varphi u \in H^1_0 (\Omega)$. 
\end{lemma}
\begin{proof}
  Consider the sequence $u_n$ from Proposition~\ref{p:density} and consider $v_n = u_n \varphi$. We easy verify that  $v_n \to u \varphi$ in $H^1(\Omega)$.
  Since $v_n \in C^\infty_c (\Omega)$, we conclude that $u \varphi \in H^1_0(\Omega)$.
\end{proof}

\subsection{Weak solutions}

We recall that $\mathcal{E}(\lambda, \Lambda)$ denotes the set of all elliptic ``matrices'', see \eqref{weak} on page~\pageref{weak}. 
\begin{defi}[Weak solutions] \label{d:weak-sol}
  Let $\Omega$ be an open set of $\R^d$ and $A \in \mathcal{E}(\lambda, \Lambda)$ and $S \in L^1 (\Omega)$. A function $u \colon \Omega \to \R$ is
  a \emph{weak solution} of $-\dive_x (A \nabla_x u ) = S$ in $\Omega$ if $u \in H^1(\Omega)$ and if for all $\varphi \in C^\infty_c (\Omega)$,
  \[ \int_\Omega A \nabla_x u \cdot \nabla_x \varphi = \int_\Omega S \varphi. \]
\end{defi}
By using the definition of $H^1_0 (\Omega)$, we see that we can extend the class of functions $v$ used in the weak formulation.
Let us put into a lemma whose proof is left to the reader.
\begin{lemma}[Extending the set of test functions] \label{l:test}
  If $u$ is a weak solution of $-\dive_x (A \nabla_x u ) = S$ in $\Omega$ with $S \in L^2 (\Omega)$, then for all $v \in H^1_0 (\Omega)$, we have
  \[ \int_\Omega A \nabla_x u \cdot \nabla_x v = \int_\Omega S v. \]
\end{lemma}

\subsection{Composition in $H^1 (\Omega)$}

The following result is classical but we state and prove it for the reader's convenience. It is sometimes called Stampacchia's theorem. 
\begin{prop}[Composition]\label{p:composition}
  Let $\Omega$ be open. We assume that either $\Omega$ is bounded or $T(0)=0$. Consider $u \in H^1(\Omega)$ and $T \colon \R \to \R$.
  \begin{enumerate}
  \item \label{one} If $T$ is $C^1$ and $T'$ is bounded in $\R$ (that is to say globally Lipschitz continuous),
    then $T(u) \in H^1 (\Omega)$ and $\nabla_x (T(u)) = T'(u) \nabla_x u$ in $L^2 (\Omega)$.
  \item \label{two} If $T (r) = r_+ = \max(r,0)$, then $T(u) \in H^1 (\Omega)$ and $\nabla_x u_+ = \un_{\{ u \ge 0 \}} \nabla_x u = \un_{\{ u > 0 \}} \nabla_x u$. 
  \end{enumerate}
\end{prop}
\begin{proof}
We prove each item successively. 
  
\textsc{Proof of \ref{one}.}  We first consider $T \in C^1$ with $T'$ bounded in $\R$ by a constant $L>0$.
  Then we know that
  \[ |T(u)| \le |T(0)| + L |u|.\]
  In particular $T(u) \in L^2 (\Omega)$ (either because $\Omega$ is bounded or because $T(0)=0$). Moreover, $T'(u) \nabla_x u$ is square integrable because
  $|T'(u) \nabla_x u | \le L |\nabla_x u|$.

  We prove next that $T'(u) \nabla_x u$ is the weak gradient of $T(u)$.
  By Proposition~\ref{p:density}, we know that there exists $u^n \vert_\Omega \in C^\infty_c (\R^d)$ such that $u^n \to u$ in $H^1 (\Omega)$.
  For clarity, we simply write $u^n$ for $u^n \vert_\Omega$. 
  In particular $u^n \to u$ in $L^2 (\Omega)$ and, up to a sub-sequence,
  we have that $u^n \to u$ almost everywhere in $\Omega$.

  Since
  \[ |T(u^n)-T(u)| \le L |u^n - u|, \]
  we first get that $T(u^n) \to T(u)$ in $L^2 (\Omega)$.

  Second we claim that $\nabla_x T(u^n) \to T'(u) \nabla_x u$ in $L^2 (\Omega)$.
  Indeed $\nabla_x T(u^n) = T'(u^n) \nabla_x u^n$ and we can write,
  \begin{align*}
    |\nabla_x T(u^n) - T'(u) \nabla_x u | &\le |T'(u^n)-T'(u)| |\nabla_x u | + |T'(u^n)||\nabla_x u^n - \nabla_x u| \\
  &\le |T'(u^n)-T'(u)| |\nabla_x u | + L |\nabla_x u^n - \nabla_x u|.
  \end{align*}
  The second term converges to $0$ in $L^2 (\Omega)$ and one can apply dominated convergence to prove that the first one also goes to $0$ in $L^2 (\Omega)$ since
  \[ |T'(u^n)-T'(u)| |\nabla_x u | \le 2 L |\nabla_x u|.\]
  
Thanks to the uniqueness of distribution limits, we conclude that $\nabla_x (T(u)) = T'(u) \nabla_x u$ in $L^2 (\Omega)$. 
  \medskip
  
\textsc{Proof of \ref{two}.} We argue by approximation. More precisely, we consider two
bump functions in $\R$: $\bar \rho$ supported in $[0,1]$ and $\underline{\rho}$ supported in $[-1,0]$.
We then consider for $\eps \in (0,1)$ the function $\bar \theta^\eps$ such that $(\bar \theta^\eps)'' (t)= \bar \rho_\eps (t)= \frac1\eps \bar \rho \left(\frac{t}\eps\right)$ and
$\bar \theta^\eps (-\eps) = (\bar \theta^\eps)'(-\eps) =0$. Similarly, we consider $\underline \theta^\eps$ such that $(\underline \theta^\eps)''(t) = \underline \rho_\eps (t) =
\frac1\eps \underline \rho \left( \frac{t}{\eps} \right)$ and $\underline  \theta^\eps (0) = (\underline \theta^\eps)'(0) =0$. Easy computations show that
\[ \un_{[0,\infty)} \le (\bar \theta^\eps)' \le \un_{[-\eps,\infty)} \le \un_{[\eps,+\infty)} \le (\underline \theta^\eps)' \le \un_{(0,+\infty)}.\]
In particular,
\begin{align*}
  (\bar \theta^\eps)' \to \un_{[0,+\infty)} \quad &\text{ and } \quad (\underline \theta^\eps)' \to \un_{(0,+\infty)}   \quad \text{(pointwise)}\\
\forall r, \quad  r_+ \le  (\bar \theta^\eps) (r) \le (r+\eps)_+   \quad &\text{ and } \quad (r-\eps)_+ \le (\underline \theta^\eps) (r) \le r_+
\end{align*}
and both $\bar \theta^\eps(u)$ and $\underline \theta^\eps(u)$ converge to $u_+$ in $L^2 (\Omega)$ (by dominated convergence). 

By \ref{one}, we know that $\bar \theta^\eps(u)$ and $\underline \theta^\eps(u)$ are in $H^1(\Omega)$ and
\[ \nabla_x (\bar \theta^\eps(u)) = (\bar \theta^\eps)'(u) \nabla_x u \quad
  \text{ and } \quad\nabla_x (\underline \theta^\eps (u)) = (\underline \theta^\eps)'(u) \nabla_x u.\]
Since $(\bar \theta^\eps)'$ and $(\underline \theta^\eps)'$ are both bounded by $1$, we can apply dominated convergence
and conclude that 
\[ \nabla_x (\bar \theta^\eps (u)) \to \un_{[0,+\infty)} (u) \nabla_x u \quad
  \text{ and } \quad\nabla_x (\underline \theta^\eps (u)) \to \un_{(0,+\infty)} (u) \nabla_x u  \quad \text{ in } L^2 (\Omega).\]
Then uniqueness of distribution limits implies that both limits coincide with $\nabla_x u_+$. 
\end{proof}

\subsection{Functional inequalities}

In this section, we state without proofs two functional inequalities.
The first one corresponds to Sobolev's embedding \cite[Corollary~9.14]{MR2759829}.
\begin{prop}[Sobolev's inequality] \label{p:sobolev}
  There exists a positive constant $\Csob$, depending on dimension $d$, such that for all $r>0$ and $u \in H^1(B_r)$,
  \[ \| u \|_{L^{p^\ast} (B_r)}^2 \le \Csob \left( \|\nabla_x u\|^2_{L^2(B_r)} + r^{-2} \|u\|^2_{L^2(B_r)} \right) \]
  with $\frac1{p^\ast} = \frac12 - \frac1d$ if $d \ge 3$ and any $p^\ast$ if $d=1,2$.
  The constant $\Csob$ also depends on $p^\ast$ if $d=1,2$. 
\end{prop}

We state the next functional inequality on the unit ball \cite[Theorem~3.2]{MR2021262}. We will use it in Chapter~\ref{c:kin}
with some $q \in (1,2]$. 
\begin{prop}[Poincaré-Wirtinger's inequality]\label{p:poincare}
Let $q \in [1,2]$. 
  There exists a constant $\Cpw$, only depending on $d$ and $q$, such that for all  $u \in L^1 (B_1)$ with $\nabla_x u \in L^1 (B_1)$.
  Then,
  \[ \int_{B_1} \left| u - \fint_{B_1} u \right|^q \dx \le \Cpw \int_{B_1} |\nabla_x u |^q \dx \]
    where $\fint_{B_1} u = \frac1{|B_1|} \int_{B_1} u \dx$. 
\end{prop}
\begin{remark}
  In the previous proposition, the functions $u$ and their weak derivatives $\nabla_x u$ are only in $L_q$ on the unit ball:
  the appropriate setting for such a statement is the Sobolev space $W^{1,q}(B_1)$, see for instance \cite{MR2759829}. 
\end{remark}

We use the Poincaré-Wirtinger's inequality to get a result
that relates the size of the sets where an $H^1$ function is respectively above $\frac12$ and below $0$.
The regularity of the weak derivative implies that the set of intermediate values (i.e. between $0$ and $\frac12$) cannot be too small.
\begin{lemma}[Intermediate value] \label{l:ivl-elliptic}
  There exists $\Civl>0$ only depending on dimension such that for all $u \in H^1 (B_1)$,
  \[ |\{ u \le 1/2 \} \cap B_1 | \cdot |\{u \ge 1\}\cap B_1| \le \Civl \|\nabla_x u\|_{L^2 (B_1)} |\{1/2 < u < 1\} \cap B_1|^{\frac12}.\]
\end{lemma}
\begin{proof}
  Let $v = T(u)$ with $T(r) = \max(\min (r,1),1/2)$. We claim that $v \in H^1(B_1)$ and
  \[ \nabla_x v = \un_{\{1/2 < u < 1\}} \nabla_x u \text{ in } L^2 (B_1).\]
  In order to justify the claim, we remark that $\min (u,1) = - \max (-u ,-1) = 1 - (1 -u)_+$.
  Since we know that $1-u \in H^1(B_1)$, we have that $w= \min (u,1) \in H_1 (B_1)$ by Proposition~\ref{p:composition}
  and $\nabla_x w = \un_{\{ u < 1\}} \nabla_x u$. Now we can apply the proposition again and conclude that $v$ is indeed in $H^1(B_1)$
  with a gradient supported in the intermediate value set.

  We next apply the Poincaré-Wirtinger's inequality (see Proposition~\ref{p:poincare}) to the function $v$. Letting $\bar v$ denote $\fint_{B_1} v$, we write
  \begin{align*}
    \int_{B_1} |v -\bar v | \dx & \le \Cpw \int_{B_1} |\nabla_x v | \dx \\
                                & = \Cpw \int_{B_1} |\nabla_x u | \un_{\{0 < u < 1/2\}} \dx \\
    & \le \Cpw \|\nabla_x u \|_{L^2 (B_1)} \left| \{ 1/2 < u < 1\} \cap B_1\right|^{\frac12}. 
  \end{align*}
  We now get a lower bound on the left hand side of the first inequality.
  \begin{align*}
    \int_{B_1} |v -\bar v| \dx &\ge \int_{ \{ v = 1/2 \} \cap B_1} |1/2 - \bar v| \dx \\
                               & = (\bar v - 1/2) |\{v =0\} \cap B_1 |  \\
    & = \frac1{|B_1|} \left( \int_{B_1} (v-1/2) \dx \right)|\{v =1/2\} \cap B_1 | \\
    & \ge  \frac1{2|B_1|}  |\{v = 1\} \cap B_1| |\{v =1/2\} \cap B_1 |.
  \end{align*}
  We get the announced inequality with $\Civl = \Cpw 2 |B_1|$. 
\end{proof}
\begin{remark}
  Exponents of set measures are not optimal. We lost information when we used the Poincaré-Wirtinger's inequality, that is sub-optimal too.
  Exponents  can be improved in order to get the functional inequality originally proved by E.~De Giorgi in \cite{DeG56} and nowadays known as De Giorgi's isoperimetric inequality. The interested reader is also referred to  \cite{MR3525875}.
\end{remark}

\section{Elliptic De Giorgi's classes  \& the local maximum principle}
s
In this section, we derive a local maximum principle for a class of functions satisfying some local energy estimates. 

\subsection{Elliptic De Giorgi's classes}

Our next goal is to derive a family of inequalities satisfied by a truncated weak solution.
\begin{defi}[Elliptic De Giorgi's classes]\label{d:DGclass}
  Let $B$ be an open ball of $\R^d$ and $S \in L^2 (B)$.   A function $u \colon B \to \R$ belongs to the \emph{elliptic De Giorgi's class} $\mathrm{DG}^\pm (B,S)$ if $u \in H^1 (B)$
  and there exists a  constant $\Cee \ge 1$ such that 
  for all ball $B_R (x_0) \subset B$, all $\kappa \in \R$ and all $r \in (0,R)$,
  \begin{equation}\label{e:EE}
    \int_{B_r (x_0)} \left| \nabla_x (u-\kappa)_\pm \right|^2 \dx \le \frac{\Cee}{(R -r)^2} \int_{B_R (x_0)} \left| (u-\kappa)_\pm \right|^2 \dx + \Cee \int_{B_R(x_0)} |S| (u-\kappa)_\pm \dx
  \end{equation}
  where $(u-\kappa)_+ = \max (u-\kappa,0)$ and $(u-\kappa)_- = \max (-(u-\kappa),0)$.

  The class $\mathrm{DG} (B,S)$ is the intersection of $\mathrm{DG}^+ (B,S)$ and $\mathrm{DG}^- (B,S)$.
\end{defi}
\begin{remark}[Universal constants]
  Before introducing De Giorgi's classes, a constant was called \emph{universal} if it only depends on dimension $d$ and ellipticity parameters $\lambda,\Lambda$.
  Now we extend this definition to include constants that only depend on dimension and on the constants $\Cee$ used to define De Giorgi's classes.  
\end{remark}
The classes are invariant under translation and scaling.
\begin{lemma}[Invariance of De Giorgi's classes] \label{l:invariance-dg}
  If $B$ be an open ball of $\R^d$ and $u \in \mathrm{DG}^\pm(B)$ and $B_r (x_0) \subset  B$. Then the function $v = \lambda u (\frac{x-x_0}r)$
  lies in $\mathrm{DG}^\pm(B_1, \mathfrak{S})$  with $\mathfrak{S} (x)= \frac{\lambda}{r^2} S (\frac{x-x_0}r )$. 
\end{lemma}

We now check that the De Giorgi's class $\mathrm{DG}$ contains all weak solutions of the elliptic equations that are considered in this chapter. 
\begin{prop}[Weak solutions and DG classes]\label{p:DGclass}
  Let $u$ be a weak solution of an elliptic equation $-\dive_x ( A \nabla_x u) = S$ in an open ball $B$ with $S \in L^2 (B)$.
Then $u \in \mathrm{DG}(B,S)$.
\end{prop}
\begin{remark}[Caccioppoli's estimate]
  Such an estimate is some times called a Caccioppoli's estimate (for the functions $(u-\kappa)_\pm$). It can be described as a reverse Poincaré's inequality, with
  a key difference: the $L^2$-norm of the gradient is taken with respect to a ball $B_r (x_0)$ that is strictly contained in
  the ball $B_\rho (x_0)$ where the $L^2$-norm of the function is computed. 
\end{remark}
\begin{remark}[Local energy estimates]
The inequalities from the proposition are some times referred to as local energy estimates. 
\end{remark}
Before turning to the proof of this proposition, we state and prove an elementary technical lemma. 
\begin{lemma}[Truncation function] \label{l:trunc}
  For any positive numbers $r,R$ with $r< R$ and $x_0 \in \R^d$, there exists
  a $C^\infty$ function $\rho\colon \R^d \to [0,1]$ with
  compact support in $B_R(x_0)$, equal to $1$ in $B_r(x_0)$ and such that
  \[ \forall x \in \R^d, \qquad \left|\nabla_x \rho(x) \right| \le \frac{2}{R-r}.\]
\end{lemma}
\begin{proof}
Consider  a $C^\infty$ function $\theta \colon \R \to [0,+\infty)$ supported in $[-1,1]$, such that $\theta (0) =1$ and $|\theta'(r)| \le 2$ for all $r \in \R$. Then consider $\rho (x) = \theta \left(\frac{|x|-r}R \right)$. Such a function is smooth since $\nabla_x \rho$ and higher derivatives is supported in $B_R \setminus B_r$ where $|x|$ is smooth. 
\end{proof}
We can now prove the local energy estimates. 
\begin{proof}[Proof of Proposition~\ref{p:DGclass}]
We want to use $v = (u-\kappa)_+ \rho^2$ as a test function in Definition~\ref{d:weak-sol} of a weak solution for $u$. 
We know that $u-\kappa \in H^1(B)$ and by Proposition~\ref{p:composition}, we know that $(u-\kappa)_+ \in H^1 (B)$.
Then we can localize thanks to $\rho$ by using Lemma~\ref{l:localization} and get that $(u-\kappa)_+ \rho^2 \in H^1_0 (B)$. 

We now use Lemma~\ref{l:test} and the fact that $\nabla_x ((u-\kappa)_+ \rho^2) = \rho^2 \nabla_x (u-\kappa)_+ +2 (u-\kappa)_+ \rho \nabla_x \rho$
in order to write,
\[
  \int_B \bigg[ A \nabla_x u \cdot \nabla_x (u-\kappa)_+ \bigg] \rho^2 = \int_B S (u-\kappa)_+ \rho^2 - 2 \int_B \bigg[ A \nabla_x u \cdot \nabla_x \rho \bigg] (u-\kappa)_+ \rho.
\]
We now use the fact that $\nabla_x (u-\kappa)_+ = \un_{\{u > \kappa\}} \nabla_x u$ and $A$ can be written $(\sqrt{A})^2$ with $\sqrt{A}$ symmetric and definite in order
to write the previous equality as follows,
\begin{align*}
  \int_B \bigg[ A \nabla_x (u-\kappa)_+ &\cdot \nabla_x (u-\kappa)_+  \bigg] \rho^2 \\
  &= \int_B S (u-\kappa)_+ \rho^2 - 2 \int_B \bigg[ \sqrt{A} \nabla_x (u-\kappa_+) \cdot \sqrt{A} \nabla_x \rho \bigg] (u-\kappa)_+ \rho \\
  & \le \int_{B_R (x_0)} |S| (u-\kappa)_+ + \frac12    \int_B \bigg[ A \nabla_x (u-\kappa)_+ \cdot \nabla_x (u-\kappa)_+ \bigg] \rho^2 \\
  &  + 2 \int_B (A \nabla_x \rho \cdot \nabla_x \rho ) (u-\kappa)_+^2.
\end{align*}
This implies,
\[
\frac12  \int_B \bigg[ A \nabla_x (u-\kappa)_+ \cdot \nabla_x (u-\kappa)_+  \bigg] \rho^2 
   \le \int_{B_R (x_0)} |S| (u-\kappa)_+  + 2 \int_B (A \nabla_x \rho \cdot \nabla_x \rho ) (u-\kappa)_+^2.
\]
We now use ellipticity of $A$ and the fact that $\rho \equiv 1$ in the smaller ball $B_\rho (x_0)$ in order to finally get,
\[
\frac\lambda 2  \int_{B_\rho (x_0)} |\nabla_x (u-\kappa)_+|^2  
   \le \int_{B_R (x_0)} |S| (u-\kappa)_+  + \frac{8 \Lambda}{(R-r)^2} \int_{B_R (x_0)}  (u-\kappa)_+^2.
\]
We thus proved the local energy estimate with $\Cee = \max (2/\lambda, 16 \Lambda/\lambda)$. 
\end{proof}

\subsection{Local maximum principle}

We saw above that our first task is to control the $L^\infty$-norm of weak solutions satisfying
an elliptic equation in $B_1$ in the interior ball $B_{1/2}$. As a matter of fact, we can prove
it in any interior ball. It will be convenient for us to get it in $B_{3/4}$. And we will derive
a point-wise upper bound rather than a two-sided bound. Here is a precise statement. We recall that
a constant is \emph{universal} if it only depends on the constants appearing in the definition of
the De Giorgi's class. 
\begin{prop}[Local maximum principle] \label{p:lmp}
  There exists a universal constant  $\Clmp >0$ such that for any $u \in \mathrm{DG}^+ (B_1,S)$,  
 \[ \| u_+ \|_{L^\infty (B_{3/4})} \le \Clmp \left( \|u_+\|_{L^2 (B_1)} + \|S\|_{L^\infty (B_1)} \right) .\]
\end{prop}
In order to prove this lemma, we will need the following technical result about sequences of real numbers.
\begin{lemma} \label{l:induc}
  Let $(A_k)_k$ be a sequence of positive real numbers such that there exists $\beta >1$ and $C >1$ such that
  \[ \forall k \ge 0, \quad  A_{k+1} \le C^{k+1} A_k^\beta .\]
  If $A_0 < C^{-\frac\beta{(\beta-1)^2}}$, then  $A_k \to 0$ as $k \to +\infty$. 
\end{lemma}
\begin{proof}
  The proof is elementary, we just iterate the estimate on $A_k$ in order to get
  \[ \forall k \ge 1, \quad  A_k \le C^{p_k} A_0^{\beta^k}\]
  with $p_k = \sum_{i=0}^{k} (k-i) \beta^i$. This estimate can be proved by induction.
  Indeed, $p_1 = 1$ and $A_1 \le C A_0^\beta$. Now if $A_k \le C^{p_k} A_0^{\beta^k}$ for some $k \ge 1$, then
  \[ A_{k+1} \le C^{k+1} A_k^{\beta} \le C^{k+1} C^{p_k \beta} A_0^{\beta^{k+1}} \]
  and we do have
  \[p_k \beta + (k+1) = \sum_{i=0}^k (k-i) \beta^{i+1} + (k+1)=  p_{k+1}.\]
  
  It is now possible to compute explicitly $p_k$ by remarking that
  \[ p_k = \frac{\partial P}{\partial \alpha} (1,\beta) \quad \text{ with } \quad P (\alpha,\beta ) = \sum_{i=1}^{k} \alpha^{k-i} \beta^i
    = \frac{\beta^{k+1}-\alpha^{k+1}}{\beta -\alpha}.\]
  We compute and estimate the partial derivative of $P$ with respect to $\alpha$ with $0 < \alpha < \beta$.
  \[ \frac{\partial P}{\partial \alpha} (\alpha, \beta ) = \frac{\beta^{k+1}-\alpha^{k+1}}{(\beta -\alpha)^2} - (k+1) \frac{\alpha^k}{\beta - \alpha}
    \le \frac{\beta^{k+1}}{(\beta -\alpha)^2}.\]
  In particular,
  \[ \forall k \ge 1, \quad  p_k \le \frac{\beta^{k+1}}{(\beta -1)^2}.\]
  This implies that $A_k$ satisfies
  \[ \forall k \ge 1, \quad  A_k \le \left(C^{\frac{\beta}{(\beta-1)^2}} \right)^{\beta^k} A_0^{\beta^k} = \left(C^{\frac{\beta}{(\beta-1)^2}} A_0 \right)^{\beta^k} . \qedhere\]
\end{proof}
We are ready to prove the local maximum principle.
\begin{proof}[Proof of Proposition~\ref{p:lmp}]
  We first prove that there exists some universal constant $\delta_0 \in (0,1)$ such that, 
  if $\|S\|_{L^\infty(B_1)} \le 1$ and if  $\|u_+ \|_{L^2 (B_1)} \le \delta_0$,  then $u \le 2$ a.e. in $B_{3/4}$. 

\paragraph{Iterative truncation.}
  De Giorgi's original idea for getting an upper bound on the weak solution under study is to truncate it by an increasing sequence $\kappa_k$
  and integrate it on shrinking balls $B^k = B_{r_k}$. Precisely, we consider
  \[ A_k = \int_{B^k} (u-\kappa_k)_+^2 \dx. \]
  with
  \[ \forall k \ge 0, \quad \kappa_k = 2 - 2^{-k}, \quad r_k = \frac34 + \frac{2^{-k}}{4} .\]
  In order to obtain an upper bound on $u$ in $B_{3/4}$, we have to find two universal constants $\beta>1$ and $C>0$ such that, for all $k \ge 1$, we have
  $A_{k+1} \le C^k A_k^\beta$. Indeed, in this case, Lemma~\ref{l:induc} implies that $A_k \to 0$ as soon as $A_0 < C^{-\frac\beta{\beta-1}}$.
  Since
  \[ A_0 = \int_{B_1} (u-1)_+^2 \dx \le \|u_+ \|_{L^2 (B_1)}^2 \le \delta_0^2, \]
  we see that we can choose $\delta_0 = (1/2) C^{-\frac{\beta}{2(\beta-1)}}$. Such a constant is universal since so are $C$ and $\beta$.
  Since the limit of $A_k$ as $k \to +\infty$ is $\|(u-2)_+\|^2_{L^2 (B_{3/4})}$, the fact that $A_k \to 0$ yields $u \le 2$ almost everywhere in $B_{3/4}$.

  \paragraph{Local energy estimates.}
  We use the definition of elliptic De Giorgi's classes to write the corresponding inequality for $u$ with $x_0=0$, $R= r_k$ and $r = r_{k+1}$.
  In particular, $R-r = r_k - r_{k+1}= 2^{-k-3}$, and recalling that $\|S\|_{L^\infty (B_1)} \le \delta_0$, we obtain,
  \[ \int_{B^{k+1}} |\nabla_x (u - \kappa_{k+1})_+|^2 \dx \le \Cee 4^{k+3} \int_{B^k} (u-\kappa_{k+1})_+^2 \dx + \Cee  \int_{B^k}  (u-\kappa_{k+1})_+ \dx .\]

\begin{itemize}
\item \textsc{(Gain of integrability)}
  We then use Sobolev's inequality in $B^{k+1}$, see Proposition~\ref{p:sobolev}, and get
  \begin{equation*}
\begin{aligned}    \Csob^{-1} \|(u - \kappa_{k+1})_+\|_{L^{p_*}(B^{k+1})}^2 \le & \left(\Cee 4^{k+3} +(3/4)^{-2} \right) \int_{B^k} (u-\kappa_{k+1})_+^2 \dx \\
    & + \Cee  \int_{B^k}  (u-\kappa_{k+1})_+ \dx .
\end{aligned}
\end{equation*}
We also used both that $B^{k+1} \subset B^k$ and $r_{k+1} \ge 3/4$. We now use Cauchy-Schwarz inequality and $\kappa_{k+1} \ge \kappa_k$,
  \begin{equation}\label{e:dg1}
    \|(u - \kappa_{k+1})_+\|_{L^{p_*}(B^{k+1})}^2 \le  \Cdgun \left( 4^{k+3}  A_k +   A_k^{\frac12} |\{ u \ge \kappa_{k+1} \} \cap B^k |^{\frac12} \right)
\end{equation}
with a universal constant $\Cdgun = \Csob (\Cee +1)+1$. It is convenient to add $1$ in order to ensure that $\Cdgun \ge 1$.
\item \textsc{(Nonlinearization procedure)}
  We now use that $\{ u \ge \kappa_{k+1} \} = \{ (u-\kappa_k)_+ \ge 2^{-k-1} \}$ and
  Bienaymé-Chebyshev's inequality in order to estimate the norm of the indicator function,
  \begin{align}
    \nonumber
    |\{ u \ge \kappa_{k+1} \} \cap B^{k} | & =  |\{ (u-\kappa_k)_+ \ge 2^{-k-1}\} \cap B^{k} | \\
    & \nonumber \le 4^{k+1} \| (u-\kappa_k)_+\|_{L^2 (B^{k})}^2 \\
    & \le 4^{k+1} A_k \label{e:dg3}
  \end{align}
  We then can combine \eqref{e:dg1} and \eqref{e:dg3} and get,
  \begin{align}
\nonumber       \|(u - \kappa_{k+1})_+\|_{L^{p_*}(B^{k+1})}^2 &\le  \Cdgun    \left( 4^{k+3}  A_k +   2^{k+1} A_k \right) \\
\label{e:dg4}        &\le 2^{2k+7} \Cdgun A_k .
  \end{align}
\item \textsc{(Nonlinear iteration)}
  We now estimate $A_{k+1}$ from above by using H\"older's inequality with $q \in (1,2)$ such that $\frac12 = \frac1{p_*} + \frac1q$,
  \begin{align*}
    A_{k+1} &\le \|(u - \kappa_{k+1})_+\|_{L^{p_*}(B^{k+1})}^2 \left\| \un_{\{ u \ge \kappa_{k+1}\}}\right\|_{L^q (B^{k+1})}^2 \\
            & \le  \|(u - \kappa_{k+1})_+\|_{L^{p_*}(B^{k+1})}^2 \left|\{ u \ge \kappa_{k+1}\} \cap B^k\right|^{\frac2q} \\
            & \le 2^{2k+7} 4^{\frac{2(k+1)}q}\Cdgun  A_k^{1+\frac2q}  
  \end{align*}
\end{itemize}
This implies in particular that $A_{k+1} \le C^{k+1} A_k^\beta$ with the universal exponent $\beta = 1 +\frac2q >1$
and the universal constant $C \ge 1$ only depending on $q$ and $\Cdgun$. In particular, Lemma~\ref{l:induc}
implies that $A_k$ converges to $0$ as soon as $A_0 < C^{-\frac{\beta}{(\beta-1)^2}}$. Since $A_0 \le \delta_0^2$ (see the beginning of the proof)
  we  pick $\delta_0 \in (0,1)$ such that
\[ \delta_0^2 = \frac12  C^{-\frac{\beta}{(\beta-1)^2}}.\]
Such an $\delta_0$ is indeed universal. 

\paragraph{The general case.} We now  remark that if we do not assume anymore that $\|S\|_{L^\infty(B_1)} \le 1$ and $\|u_+\|_{L^2 (B_1)} \le \delta_0$,
either $u \le 0$ a.e. in $B_1$ or $\|u_+ \|_{L^2 (B_1)} > 0$. In the latter case, we consider
\[ \tilde u  = \frac{u}{ \delta_0^{-1} \|u_+ \|_{L^2 (B_1)}+\|S\|_{L^\infty(B_1)} }.\]
This function $\tilde u$ is a weak solution of $-\dive_x (A \nabla_x \tilde u ) = \tilde S$ with
\[ \tilde S = \frac{S}{ \delta_0^{-1} \|u_+ \|_{L^2 (B_1)}+\|S\|_{L^\infty(B_1)} } \le 1. \]
Since $\|\tilde u_+ \|_{L^2(B_1)} \le \delta_0$, we conclude that \( \|\tilde u_+ \|_{L^\infty (B_{\frac34})} \le 2 \), that is to say
 \[ \| u_+ \|_{L^\infty (B_{3/4})} \le 2 \delta_0^{-1} \|u_+\|_{L^2 (B_1)} + 2 \|S\|_{L^\infty (B_1)} . \qedhere \]
\end{proof}

\subsection{More about the local maximum principle}

This subsection can be skipped unless the reader is interested in the derivation of Harnack's inequality. 
In order to establish it, we first need to adapt  the previous proof to get a local maximum principle between two balls of arbitrary radii $r$ and $R$.
We recall that a constant is \emph{universal} if it only depends on the constants appearing in the definition of
the De Giorgi's class. 
\begin{prop}[Local maximum principle] \label{p:lmp-plus}
  There exist two universal constants  $\bClmp >0$ and $\omega_0 >0$ such that for any $u \in \mathrm{DG}^+ (B_R,S)$ and $r \in (0,R)$,  
  \[ \| u_+ \|_{L^\infty (B_r)} \le \bClmp \left( \left(1 + \frac1{r^2} + \frac{1}{(R-r)^2} \right)^{\omega_0}  \|u_+\|_{L^2 (B_R)} + \|S\|_{L^\infty (B_R)} \right) .\]
\end{prop}
\begin{remark}
The constant $\omega_0= \frac{\beta}{2 (\beta-1)^2}$ for $\beta = 1 + \frac2q = \frac32 + \frac1d$ if $d \ge 3$. 
\end{remark}
\begin{proof}
  We reduce the proof to the case where  there exists some constant $\delta_0 \in (0,1)$ (depending on $d,\lambda,\Lambda, r,R$) such that, 
  if $\|S\|_{L^\infty(B_R)} \le 1$ and if  $\|u_+ \|_{L^2 (B_R)} \le \delta_0 $,  then $u \le 2$ a.e. in $B_r$.
  We define $A_k$ as before but with the shrinking radii defined as
  \[ r_k = r + (R-r)2^{-k}.\]
  In particular, $r_k -r_{k+1} = (R-r)2^{-k-1}$ and $r_k \ge r$. In particular, \eqref{e:dg1} is replaced with,
  \begin{equation}\label{e:dg1-plus}
    \|(u - \kappa_{k+1})_+\|_{L^{p_*}(B^{k+1})}^2 \le  \Cdgun \left( 4^{k+1} ((R-r)^{-2} + r^{-2})  A_k +   A_k^{\frac12} |\{ u \ge \kappa_{k+1} \} \cap B^k |^{\frac12} \right).
\end{equation}
Next, \eqref{e:dg3} is unchanged but \eqref{e:dg4} is replaced
  \begin{align}
\nonumber       \|(u - \kappa_{k+1})_+\|_{L^{p_*}(B^{k+1})}^2 &\le  \Cdgun    \left( 4^{k+1} ((R-r)^{-2} + r^{-2})  A_k +   2^{k+1} A_k \right) \\
\label{e:dg4-plus}        &\le 2^{2k+6} (r^{-2}+(R-r)^{-2} +1) \Cdgun A_k .
  \end{align}
  This implies that $A_{k+1} \le C^{k+1} A_k^\beta$ with $C = \bar C (1 + r^{-2}+(R-r)^{-2})$ with $\bar C$ universal.
  Then we can conclude if $\delta_0 = \frac1{\sqrt{2}} C^{-\frac{\beta}{2(1-\beta)^2}} = \tilde C (1 + r^{-2}+(R-r)^{-2})^{-\omega_0}$ with $\tilde C$ universal. 
  \end{proof}
  We first state a straightforward consequence of the local maximum principle.
  \begin{cor}[Upside down maximum principle]\label{c:lmp-usd}
There exist two universal constants  $\bar \eps_0,\eps_1 \in (0,1)$  such that for any non-negative $u \in \mathrm{DG}^- (B_1,S)$ with $\|S\|_{L^\infty(B_1)} \le \bar \eps_0$, we have
  \[ |\{ u \ge 1 \} \cap B_1| \ge (1-\eps_1)|B_1| \quad \Rightarrow \quad \{ u \ge 1/2  \text{ in } B_{1/2} \} .\]
\end{cor}
\begin{proof}
  We  apply the local maximum principle from Proposition~\ref{p:lmp-plus} to $v = 1-u$ with $r=1/2$ and $R=1$ and $z_0 = 0$.
  Remarking that $v \le 1$ a.e. in $B_1$ (because $u$ is non-negative), we can write,
  \begin{align*}
    \|v_+ \|_{L^\infty (B_{1/2})} &\le \bClmp \left((1+4+4)^{\omega_0} \| v_+ \|_{L^2 (B_1)} + \|S\|_{L^\infty (B_1)}  \right) \\
                                  & \le \bClmp \left(9^{\omega_0} |\{ v \ge 0\} \cap B_1|^{\frac12}  + \|S\|_{L^\infty (B_1)}  \right) \\
    & \le \bClmp \left(9^{\omega_0} \eps_1^{\frac12}  + \bar \eps_0  \right). 
  \end{align*}
  We now pick $\bar \eps_0 = \frac{\bClmp}{4}$ and $\eps_1 = \left(\frac{\bClmp}{9^{\omega_0} 4}\right)^2$ and conclude that $v \le \frac12$ in $B_{1/2}$.
  This means $u \ge \frac12$ in $B_{1/2}$. 
\end{proof}
We can next show that the  $L^2$-norm in the right hand side of the local maximum principle can be replaced with any $L^\eps$-``norm'' for any $\eps \in (0,2)$.  
\begin{cor}[Local maximum principle - again] \label{c:lmp}
 Given a universal constant $\eps \in (0,2)$, there exists a constant  $\eClmp >0$, only depending on $d,\lambda,\Lambda$ and $\eps$, such that for any $u \in \mathrm{DG}^+ (B_1,S)$,  
  \[ \| u_+ \|_{L^\infty (B_{1/2})} \le  \eClmp \left( \|u_+\|_{L^\eps (B_1)} + \|S\|_{L^\infty (B_1)} \right) \]
  where $\|u_+\|_{L^\eps (B_1)} = \|u_+^\eps\|_{L^1 (B_1)}^{\frac1\eps}$. 
\end{cor}
\begin{proof}
  This corollary is a consequence of the interpolation of $L^2$ between $L^\eps$ and $L^\infty$. If $\eps < 1$, we interpolate $L^{2/\eps}$ between $L^1$ and $L^\infty$.
We start by applying Proposition~\ref{p:lmp-plus} for $r,R \in (0,1)$, 
\begin{align*}
  \| u_+ \|_{L^\infty (B_r)}
  & \le \Clmp \left( \left(1 + \frac1{r^2} + \frac{1}{(R-r)^2} \right)^{\omega_0}  \|u_+\|_{L^2 (B_R)} + \|S\|_{L^\infty (B_R)} \right) \\
  & \le \Clmp \left( \left(1 + \frac1{r^2} + \frac{1}{(R-r)^2} \right)^{\omega_0}  \|u_+\|_{L^\eps (B_R)}^{\eps/2} \|u_+\|_{L^\infty (B_R)}^{1-\eps/2} + \|S\|_{L^\infty (B_R)} \right) \\
    & \le \frac12 \|u_+\|_{L^\infty (B_R)} + K_\eps  \left(1 + \frac1{r^2} + \frac{1}{(R-r)^2} \right)^{\omega_\eps}  
\end{align*}
with $K_\eps = 2^{\frac{2-\eps}{\eps}} \Clmp^{\frac2{\eps}}  \|u_+\|_{L^\eps (B_1)} + \Clmp  \|S\|_{L^\infty (B_1)}$ and $\omega_\eps = \frac{2\omega_0}{\eps}$. We now consider $r_0 = \frac12$ and $r_{n+1}= r_n + \delta (n+1)^{-2}$ with
$\delta = \frac12 \left(\sum_{k=1}^\infty k^{-2}\right)^{-1} =  \frac{3}{\pi^2}$. In particular, $\frac12 \le r_n \le 1$ for all $n \ge 0$.
Letting $N_n$ denote $  \| u_+ \|_{L^\infty (B_{r_n})}$, we thus have,
\[ N_n \le \frac12 N_{n+1} + K_\eps (1+4+ \delta^{-2} (n+1)^4)^{\omega_\eps} \le \frac12 N_{n+1} + K_\eps (3/\delta^2)^{\omega_\eps} (n+1)^{4 \omega_\eps}.\]
By induction, we thus get for all $n \ge 1$,
\[ N_0 \le \frac1{2^{n}} N_{n} +  K_\eps (3/\delta^2)^{\omega_\eps} \left( \sum_{k=1}^n \frac{k^{4\omega_\eps}}{2^{k-1}} \right).\]
Letting $n \to \infty$, we conclude that
\begin{align*}
\| u_+ \|_{L^\infty (B_{1/2})} & = N_0 \\
   &\le K_\eps (3/\delta^2)^{\omega_\eps} \left( \sum_{k=1}^\infty \frac{k^{4\omega_\eps}}{2^{k-1}} \right) \\
  & =  \eClmp \left(  \|u_+\|_{L^\eps (B_1)} + \Clmp  \|S\|_{L^\infty (B_1)} \right)
\end{align*}
with $\eClmp = (3/\delta^2)^{\frac{2\omega_0}{\eps}} \left( \sum_{k=1}^\infty \frac{k^{\frac{8\omega_0}{\eps}}}{2^{k-1}} \right)  \left( 2^{\frac{2-\eps}{\eps}} \Clmp^{\frac2{\eps}} + \Clmp \right)$ and  $\delta =  \frac{3}{\pi^2}$.
Since $\omega_0$ and $\Clmp$ are universal, the constant $\eClmp$ only depends on $d,\lambda,\Lambda$ and $\eps$. 
\end{proof}

\section{Improvement of oscillation \& De Giorgi's theorem}

In this section, we now state and prove De Giorgi's theorem. 
\begin{thm}[De Giorgi]\label{t:dg-ell}
  Let   $A \in \mathcal{E}(\lambda,\Lambda)$ with $\Omega = B_1$ and $\lambda, \Lambda >0$.
  There exist two universal constants $\alpha \in (0,1]$ and $\Cdg >0$ such that 
  any weak solution $u \in \mathrm{DG} (B_1,S)$ with $S \in L^\infty (B_1)$
  is $\alpha$-H\"older continuous and
 \[ \| u \|_{C^\alpha (B_{1/2})} \le \Cdg \left( \|u\|_{L^2 (B_1)} + \|S\|_{L^\infty (B_1)} \right) .\]
\end{thm}
\begin{remark}[Universal constants]
We recall that a constant is \emph{universal} if it only depends on the constant appearing in the definition of
the De Giorgi's classes $\mathrm{DG}^\pm$. 
\end{remark}

We already proved that solutions, and more generally functions in the De Giorgi's class $\mathrm{DG}^+$
are locally bounded: this is a consequence of the local maximum principle. It is sometimes called
the De Giorgi's first lemma
. We now will establish that the oscillation
improves while zooming in: this will be achieved by establishing the infimum lift.

\subsection{Infimum lift}

We now know from the local maximum principle (Proposition~\ref{p:lmp}) that elements of $\mathrm{DG}^+$
(and in particular weak solutions of the class of elliptic equations treated in this chapter)
are essentially bounded from above in the interior of the domain.
With such an information in hand, we can now study how the oscillation of functions in $\mathrm{DG}$ in a ball $B_r (x_0)$
behaves with the radius $r$. We aim at proving that it decays as $r^\alpha$ for some universal exponent $\alpha \in (0,1]$.
We  indeed saw earlier (Proposition~\ref{p:holder-osc}) that it is equivalent to being $\alpha$-H\"older continuous. 

To get such a decay, we aim at proving that the oscillation of a weak solution, and more generally of elements of the De Giorgi's class,
improves by a universal factor when zooming in by another universal factor.
In order to establish such a result, we first prove that we can lift the essential infimum of an non-negative element of $\mathrm{DG}^-$
above some universal constant $\theta$ in $B_{1/2}$ if the measure of its $1$-super-level set in $B_1$ is universally bounded from below. 
\begin{prop}[Lifting the infimum] \label{p:lower}
  Let $\iota \in (0,1)$ be a universal constant. There exist two other universal constants $\theta \in (0,1)$ and $\eps_0 \in (0,1)$ such that, 
  if $u \in \mathrm{DG}^- (B_2,S)$ and  $\|S\|_{L^\infty(B_2)} \le \eps_0$ and $u \ge 0$ a.e. in $B_2$,
  then
  \[ |\{ u \ge 1 \} \cap B_1 | \ge (1-\iota) |B_1| \quad \Rightarrow \quad \left\{ u \ge \theta \text{ a.e. in } B_{\frac12} \right\}.\]   
\end{prop}
\begin{remark}[About the parameter $\iota$]
  In order to prove De Giorgi's theorem, we only need to consider $\iota=1/2$. But when proving
  Harnack's inequality, it will be convenient to consider a larger $\iota$: indeed we will consider  $\iota = 1-4^{-d}$. 
\end{remark}
\begin{proof}
    We consider the sequence of scaled functions $u_k = 2^k u \in \mathrm{DG}^- (B_2,S_k)$  with source terms $S_k = 2^k S$ (see Lemma~\ref{l:invariance-dg}).
  We consider $\bar \eps_0,\eps_1 \in (0,1)$ from the upside down maximum principle (Corollary~\ref{c:lmp-usd}). 

  We use the definition of $\mathrm{DG}^-(B_2,S)$ in order to write,
  \begin{align*}
    \| \nabla_x (u_k-1)_- \|_{L^2(B_1)}^2 &\le \Ceem \left( \|(u_k-1)_- \|^2_{L^2 (B_2)} + \| S_k\|_{L^2 (B_2)}^2 \right) \\
                                        & \le \Ceem \left( |B_2| + 4^k  \eps_0^2 |B_2| \right) \\
    & \le 2 \Ceem |B_2|
  \end{align*}
  as long as $2^k \eps_0 \le 1$ (we will see that this condition is reached when we will choose $\eps_0$). We notice that we also used that $(u_k-1)_- \le 1$.

Because we have for all $k \ge 0$,
  \[ |\{ u_k \ge 1 \} \cap B_1| \ge |\{ u \ge 1\} \cap B_1 | \ge (1-\iota)|B_1|,\]
  we can apply the intermediate value lemma (Lemma~\ref{l:ivl-elliptic}) and deduce that
  \[ \frac{(1-\iota)^2|B_1|^2}{2 \Ceem |B_2|} |\{ u_k \le 1/2\} \cap B_1|^2 \le  |\{ 1/2 < u_k < 1 \} \cap B_1|\]
  In particular, for $k \ge 0$ such that
  \[ |\{ u_k \le 1/2\} \cap B_1| \ge \eps_1 |B_1|,\]
  we have
  \[ \alpha |B_1|   \le  |\{ 1/2 < u_k < 1 \} \cap B_1| \quad \text{ with } \quad \alpha :=\frac{(1-\iota)^2|B_1|^3}{2 \Ceem |B_2|} \eps_1^2.\]
  We now pick the largest integer $N \ge 1$ such that $N \alpha \le 1$ and consider
  \[ \eps_0 = 2^{-N} .\]

  Let
  \[ E = \bigg\{ k \in \{1,\dots, N+1\} : |\{ u_k \le 1/2 \} \cap B_1| \ge \eps_1 |B_1| \bigg\}. \]
  We just proved that for all $k \in E$,
  \[ \alpha |B_1|   \le  |\{ 1/2 < u_k < 1 \} \cap B_1| = |\{ 2^{-k-1} < u < 2^{-k} \} \cap B_1|.\]
  In particular,
  \[ (\# E) \alpha |B_1| \le \sum_{k \in E} |\{ 2^{-k-1} < u < 2^{-k} \} \cap B_1| \le |B_1|?\]
  We conclude that $\# E \le N$. In particular, there exists $k_0 \in \{1,\dots,N+1\} \setminus E$. For this integer, we have,
  \( |\{ u_{k_0+1} \le 1 \} \cap B_1| < \eps_1 |B_1| \) or equivalently,
  \[ | \{ u_{k_0+1} > 1 \} \cap B_1| > (1-\eps_1)|B_1|.\]
  Now the upside down maximum principle from Corollary~\ref{c:lmp-usd} implies that $u_{k_0+1} \ge 1/2$ a.e. in $B_{1/2}$ if $2^{k_0+1} \eps_0 \le \bar \eps_0$.
  We thus choose $\eps_0 = 2^{-N-1} \bar \eps_0$. 
  We get $u \ge 2^{-k_0-2}$ a.e. in $B_{1/2}$. We reached the desired conclusion with $\theta = 2^{-k_0-2}$. 
\end{proof}

\subsection{Improvement of oscillation}

An immediate consequence of this proposition is the fact that the oscillation of elements of the De Giorgi's class improves with a universal factor $(1-\mu)$
when zooming from $B_2$ to $B_{\frac12}$. 
\begin{prop}[Improvement of oscillation] \label{p:improve-osc}
  Let   $A \in \mathcal{E}(\lambda,\Lambda)$ with $\Omega = B_1$ and $\lambda, \Lambda >0$. Let $\eps_0 \in (0,1)$
  be given by Proposition~\ref{p:lower} for $\iota = 1/2$. 
  There exist a universal constant $\mu \in (4^{-2},1)$ and   such that, if $u \in \mathrm{DG} (B_2,S)$
  with $S \in L^\infty (B_2)$ such that $\| S\|_{L^\infty(B_2)} \le \eps_0$, and  $u \in L^\infty (B_2)$, then
  \[ \osc_{B_2} u \le 2 \quad \Rightarrow \quad \osc_{B_{\frac12}} u \le 2 \mu .\]
\end{prop}
\begin{remark}[Why do we care about a lower bound on $\mu$?]
  The fact that we pick $\mu > 4^{-2}$ is irrelevant for this proof. It is just a convenient condition
  for the proof of De Giorgi's theorem. 
\end{remark}
\begin{proof}
We let $M$  and $m$ denote the essential supremum and essential infimum of $u$ on $B_{2}$. 
In particular, $\osc_{B_{2}} u = M -m$. 

We  reduce to the case where $-1 \le u \le 1$  by considering
\( \tilde u = u-\frac{M+m}2.\)
The function $\tilde u$ takes values in $[-1,1]$ and lies in $\mathrm{DG} (B_2,\tilde S)$ with $\|\tilde{S}\|_{L^\infty (B_2)} \le \eps_0$. 

We now distinguish two cases.
\begin{itemize}
\item If $|\{ \tilde u \le 0 \} \cap B_1 | \ge \frac12 |B_1|$, then Proposition~\ref{p:lower} implies that
  $\tilde u \le 1 -\theta$ a.e. in $B_{\frac12}$. But since $\tilde u \ge -1$, we conclude that
  $\osc_{B_{\frac12}} \tilde u \le 2 -\theta$. 
\item If $|\{ \tilde u \le 0\} \cap B_1 | < \frac12|B_1|$, then $|\{ -\tilde u \ge 0\} \cap B_1 | < \frac12|B_1|$
  and the function $v = -\tilde u$ is smaller than $1$ a.e. in $B_2$ and lies in $\mathrm{DG} (B_2,-\tilde S)$
  with $\| - \tilde S\|_{L^\infty(B_2)} \le \eps_0$ and satisfies $|\{ v \le 0 \} \cap B_1| \ge |\{ v < 0 \} \cap B_1| > \frac12|B_1|$.
  We conclude that $v \le 1$ a.e. in $B_{\frac12}$, that is to say $\tilde u \ge -(1-\theta)$ a.e. in $B_{\frac12}$. Since $\tilde u \le 1$ a.e. in $B_2$,
  we conclude that   $\osc_{B_{\frac12}} \tilde u \le 2 -\theta$ in this case too. 
\end{itemize}
We thus proved that in both cases, $\osc_{B_{\frac12}} \tilde u \le 2 -\theta$.
We reached the desired conclusion with $\mu = \max (4^{-1}, (1-\theta)/2)$.
\end{proof}

\subsection{Proof of De Giorgi's theorem}

We are now ready to prove De Giorgi's theorem.
\begin{proof}[Proof of Theorem~\ref{t:dg-ell}]
The proof proceeds in several steps. 
  
\paragraph{Reduction.}
  The local maximum principle (Proposition~\ref{p:lmp}) ensures that $u$ is essentially bounded in $B_{\frac34}$,
  \[ \|u\|_{L^\infty (B_{\frac34})} \le \Clmp \left( \|u\|_{L^2 (B_1)} + \| S\|_{L^\infty (B_1)} \right).\]
  The essential upper bound on $u$ is obtained from Proposition~\ref{p:lmp} applied to $u$ and while the essential lower bound comes from
  its application to $-u$.

  We are thus left with proving that
  \[ [u]_{C^\alpha (B_{\frac12})} \le \bar{\Cdg} \left( \|u\|_{L^\infty (B_{\frac34})} + \| S\|_{L^\infty (B_1)} \right).\]
  If $\|u\|_{L^\infty (B_{\frac34})}=0$, we are done. If not, by considering
  \[ \tilde u = \frac{u}{\|u\|_{L^\infty (B_{\frac34})} + \eps_0^{-1} \| S\|_{L^\infty (B_1)}},\]
  we get that $ \|\tilde u \|_{L^\infty (B_{\frac34})} \le 1$ and $\|\tilde S \|_{L^\infty (B_1)} \le \eps_0$
  and we want to prove that
  \begin{equation}
    \label{e:goal}
    [\tilde u]_{C^\alpha (B_{\frac12})} \le \Cdg 
  \end{equation}
  for some universal constants $\alpha \in (0,1)$ and $\Cdg>0$. It is enough to study the oscillation
  of $\tilde u$ around points $x_0 \in B_{1/2}$. We thus consider such a point $x_0 \in B_{\frac12}$.
  We know that $\|\tilde u\|_{L^\infty (B_{1/4} (x_0))} \le 1$.

  \paragraph{Infinite iteration.}
  We now want to scale $\tilde u$ from $B_{1/4} (x_0)$ to $B_2$ in order to apply the result about the improvement of  oscillation, recall Proposition~\ref{p:improve-osc}. We thus consider for $x \in B_2$,
  \[ \bar u (x) = \tilde u ( x_0 + \frac18 x) .\]
  We have $\bar u \in \mathrm{DG} (B_2,\bar S)$  with $\bar A (x) = A (x_0 + \frac18 x)$ and
  $\bar S (x) = \left( \frac18 \right)^2 \tilde S (x_0 + \frac18 x)$. In particular, $\|\bar u\|_{L^\infty(B_2)} \le 1$ and $\|\bar S \|_{L^2 (B_2)} \le \eps_0$.
  Then Proposition~\ref{p:improve-osc} implies that
  \[ \osc_{B_{\frac12}} \bar u \le  2 \mu  \]
  for $\mu \in (4^{-1},1)$ universal. 
  
 Now we consider $\bar u_1 = \mu^{-1} \bar u (x/4)$. In particular, $\osc_{B_2} \bar u_1 \le 2$ and it satisfies an elliptic equation with the source term $S_1 = (4^2\mu)^{-1} \bar S (x/4)$.  
 Since $4^2 \mu \ge 1$ and $\|\bar S\|_{L^\infty(B_2)} \le \eps_0$, we also have $\| S_1 \|_{L^\infty (B_2)} \le \eps_0$. We thus can apply Proposition~\ref{p:improve-osc} and conclude
 that $\osc_{B_{\frac12}} \bar u_1 \le 2 \mu$. We iterate this procedure by consider $\bar u_{k+1} (x) = \mu^{-1} \bar u_k (x/4) = \mu^{-k-1} \bar u (4^{-k-1} x)$. These re-scaled functions $\bar u_k$ satisfy
 $\osc_{B_2} \bar u_k \le 2$ and they satisfy an elliptic equation with a source term $S_k$ such that $\|S_k\|_{L^\infty (B_2)} \le \eps_0$. We thus conclude that for all $k \ge 1$, 
 \( \osc_{B_2} \bar u_k \le 2,\) which translates into
 \[ \osc_{B_{r_k}} \bar u \le 2 \mu^k \]
 with $r_k = (1/2)4^{-k}$. We now consider $\alpha \in (0,1)$ such that $\mu^k = \left(4^{-k}\right)^\alpha$, that is to say $\alpha = \ln (1/\mu) / \ln 4 > 0$.
 We conclude that
 \[ \osc_{B_{r_k}} \bar u \le 2 (2 r_k)^\alpha =2^{1+\alpha} r_k ^\alpha.\]

 \paragraph{Conclusion.}
  We are almost done. We need to check that we control the oscillation of $\bar u$ over balls of arbitrary radius $r >0$.
  In order to do so, we first deal with $r \in (0,1/4]$ by considering $k \ge 1$ such that $r_k \le r \le r_{k-1}$. In this case, we write
    \[ \frac{\osc_{B_r} \bar u}{r^\alpha} \le \frac{\osc_{B_{r_{k-1}}} \bar u}{{r_{k-1}}^\alpha} \times \frac{r_{k-1}^\alpha}{r^\alpha}
      \le 2^{1+\alpha} \frac{r_k^\alpha}{r^\alpha} 4^{\alpha} \le 2^{1+3 \alpha}.\]
  We thus proved that for any $x_0 \in B_{\frac12}$ and any $r \in (0,\frac14]$, we have
  \[ \osc_{B_r (x_0)} \tilde u \le \left(2^{1+3\alpha}\right) r^\alpha.\]
  For a radius $r \ge 1/4$, we simply write
  \[ \osc_{B_r (x_0) \cap B_{1/2}} \tilde u \le \osc_{B_{1/2}} \tilde u \le 2 \le 2 (4^\alpha) r^\alpha = 2^{1+2 \alpha} r^\alpha.\]
  We now conclude from Proposition~\ref{p:holder-osc} that \eqref{e:goal} holds true with $\Cdg = 2^{1+3\alpha}$. 
\end{proof}

\section{(Weak) Harnack's inequality}

In this section, we will see that non-negative weak solutions of elliptic equations in divergence form are such that their supremum over a unit ball (after scaling) is controlled \emph{from above} by their infimum over the same ball, up to some universal constant. Such an inequality was first considered by C.~Harnack in 1887 in the study of convergence of sequences of harmonic functions. In view of the discussions of the previous sections, it is natural to expect that the result holds as well for non-negative elements of the De Giorgi's class $\mathrm{DG}$. 
\begin{thm}[Harnack's inequality]\label{t:harnack}
  There exist a universal constant $\Charnack$ such that for any $u \in \mathrm{DG} (B_1,S)$ 
  with $S \in L^\infty (B_1)$ and $u \ge 0$, we have
  \[ \sup_{B_{1/2}} u \le \Charnack \left( \inf_{B_{1/2}} u + \|S\|_{L^\infty (B_1)} \right) .\]
\end{thm}

Since $u$ is non-negative, its supremum coincides with its $L^\infty$-norm.
By the local maximum principle, we know that we can control it by its $L^2$-norm.
But the $L^2$-norm can be interpolated between an $L^\eps$-``norm'' for some small $\eps$ 
and the $L^\infty$-norm.\footnote{If $\eps \in (0,1)$, the space $L^{\frac2\eps}$ is
  interpolated between $L^1$ and $L^\infty$, see below.}

For this reason, the proof of Harnack's inequality reduces to the control of the mass of $f^\eps$ in $B_{1/2}$.
Such a result is known as a \emph{weak Harnack's inequality}. It is not weaker than Harnack's inequality,
it is in fact more general since it applies to any  element of the De Giorgi's class $\mathrm{DG}^-$.
\begin{remark}
We recall that this latter class contains all super-solutions of the elliptic equations with work with.
\end{remark}
\begin{thm}[Weak Harnack's inequality]\label{t:whi}
   There exist two universal constants $\Cwhi >0$ and $\eps>0$ such that for  any $u \in \mathrm{DG}^- (B_2,S)$
  with $S \in L^\infty (B_2)$ and $u \ge 0$, we have
  \[ \left( \int_{B_{\frac12}} u^\eps (x) \dx \right)^\eps \le \Cwhi \left( \inf_{B_\frac12} u + \|S\|_{L^\infty (B_2)} \right) .\]
\end{thm}
\begin{remark}[Universal constants]
We recall again for the reader's convenience that a constant is \emph{universal} if it only depends on the constant appearing in the definition of
the De Giorgi's classes $\mathrm{DG}^\pm$. 
\end{remark}

The proof of this theorem relies on the covering argument that is presented in the next subsection. 

\subsection{Ink spots}

\begin{lemma}[Ink spots] \label{l:is}
  Let $E \subset F \subset B_{\frac12}$ be measurable sets of $\R^d$. Assume there is a constant $\iota >0$ such that
  \begin{itemize}
  \item $|E| < (1-\iota) |B_{\frac12}|$,
  \item  any open ball $B \subset B_{\frac12}$ satisfying $|E \cap B| > (1-\iota)|B|$ is contained in $F$.
  \end{itemize}
Then $|E| \le (1-c \iota)|F|$ for some constant $c$ only depending on the dimension $d$. 
\end{lemma}
\begin{proof}
  For a ball $B = B_r(x)$ and $\kappa >0$, we write $\kappa B$ for $B_{\kappa r} (x)$.  
  
  By Lebesgue's differentiation theorem \cite[Theorem II.4.5]{boyer2012mathematical} applied to the integrable function
  \[ \un_E (x) = \begin{cases} 1 & \text{ if } x \in E, \\ 0 & \text{ if not,} \end{cases} \]
  we know that for a.e. $x \in E$, there exists an open ball $B^x$ such that $|E \cap B^x | \ge (1-\iota) |B^x|$.
  Let us now choose a maximal open ball $B^x \subset B_{\frac12}$ containing $x$ and satisfying $|E \cap B^x | \ge (1-\iota) |B^x|$. It is of the form $\bar B^x = B_{\bar r} (\bar x)$.
  By assumption, we know that $B_{\bar r} (\bar x) \neq B_{\frac12}$ and $B_{\bar r}(\bar x) \subset F$.

  We now claim that $|E \cap \bar B^x | = (1-\iota) |\bar B^x|$. Otherwise, there would be a ball $\tilde B^x$ and a $\delta>0$
  such that $\bar B^x \subset \tilde B^x \subset (1+\delta) \tilde B^x$ with  $\tilde B^x \subset B_{\frac12}$ and $|E \cap \tilde B^x| > (1-\iota)|\tilde B^x|$,
  contradicting the maximality of $\bar B^x$.

  The set $E$ is covered by the closed balls $\bar B^x$. By Vitali's lemma \cite[Theorem~1.24]{MR3409135}, there exists a countable sub-collection of non-overlapping closed balls $\bar B^j = \bar B_{r_j} (x_j)$, $j \ge 1$,
  such that $E \subset \cup_{j=1}^\infty 5 \bar B^j$. Since $B^j \subset F$ and $|B^j \cap E| \ge (1-\iota)|B^j|$, this implies that $|B^j \cap (F \setminus E) | \ge \iota |B^j|$.
  \[
    |F \setminus E|  \ge \sum_{j=1}^\infty |B^j \cap (F \setminus E)| 
    \ge \sum_{j=1}^\infty \iota |B^j| 
    = 5^{-d} \sum_{j=1}^\infty \iota |5 B^j| 
    \ge 5^{-d} \iota |E|.
  \]
  We conclude that $|F| \ge (1+5^{-d} \iota ) |E|$, from which we get $|E| \le (1-c \iota) |F|$ with $c=5^{-d}$ since $c \iota < 1$.  
\end{proof}

\subsection{Proof of the (weak) Harnack's inequality}

If we consider a non-negative $u \in \mathrm{DG}^-$ and we apply Proposition~\ref{p:lower} to the function $1-u \le 1$, we readily get
the following result.
\begin{cor}[Generating a lower bound] \label{c:lower-bound}
Let $\iota \in (0,1)$ be universal.  There exist two universal constants $\theta \in (0,1)$ and $\eps_0 \in (0,1)$ such that, 
  if $u \in \mathrm{DG}^-(B_2,S)$ with  $\|S\|_{L^\infty(B_2)} \le \eps_0$ and $u \ge 0$ a.e. in $B_2$,
  then
  \[ |\{ u \ge 1 \} \cap B_1 | \ge (1-\iota) |B_1| \quad \Rightarrow \quad \left\{ u \ge \theta \text{ a.e. in } B_{\frac12} \right\}.\]   
\end{cor}

The fact that $u \ge \theta$ a.e. in $B_{\frac12}$ implies that $\theta^{-1} u \ge 1$ in a proportion $(1-\iota)$ of a larger ball. 
We thus can iterate this estimate by re-scaling the function $u$ at each iteration, up to getting a lower bound on $u$ in $B_1$. 
\begin{cor}[Expansion of positivity] \label{c:spreading}
   There exist  universal constants $\iota \in (0,1)$, $M >1$ and $\eps_{0,s} \in (0,1)$ such that, 
   if $u \in \mathrm{DG}^-(B_4, S)$ with  $\|S\|_{L^\infty(B_4)} \le \eps_{0,s}$ and $u \ge 0$ a.e. in $B_4$,
   then
   \[ |\{ u \ge M \} \cap B_1 | \ge (1-\iota) |B_1| \quad \Rightarrow \quad \left\{ u \ge 1 \text{ a.e. in } B_1  \right\}\]
   or equivalently,
   \[ \inf_{B_1} u \le 1 \quad \Rightarrow \quad |\{ u \ge M \} \cap B_1 | < (1-\iota) |B_1| .\]
\end{cor}
\begin{proof}
Let $\iota = 1 - 4^{-d}$ and $\theta$ and $\eps_0$ be given by Corollary~\ref{c:lower-bound}. Let $M= \theta^{-2}> 1$ and $\eps_{0,s} =  \eps_0/4 \in (0,1)$.
 We can apply Corollary~\ref{c:lower-bound} to $\frac{u}M$ since the corresponding source term $\frac{S}M$ is essentially bounded in $B_4$ (and thus in $B_2$) by $\eps_0$.
 We get that  $u \ge \theta M = \theta^{-1}$ a.e. in $B_{\frac12}$. We precisely chose $\iota$ such that $(1-\iota) |B_2|=|B_{\frac12}|$. In particular, the function $u$ satisfies,
 \[ |\{u \ge \theta^{-1}\} \cap B_2 | \ge (1-\iota)|B_2|.\]
The function $u_1 (x) := \theta u(2x)$ satisfies,
\[ -\dive_x (A_1 \nabla_x u_1) = S_1 \text{ in } B_2 \]
with $S_1 (x) = 4\theta S (2 x)$.
In particular, $\| S_1 \|_{L^\infty (B_2)} \le \theta \eps_0 \le \eps_0$ since $\theta \in (0,1)$.
We can apply  Corollary~\ref{c:lower-bound} again, but to $u_1$ this time, because we have,
\begin{align*}
  |\{ u_1 \ge 1 \} \cap B_1 | &= \int_{B_1} \un_{[1,+\infty)} (\theta u (2x)) \dx = \int_{B_2} \un_{[\theta^{-1}, +\infty)} (u(y)) 2^{-d} \dy \\
  &= 2^{-d} |\{ u \ge \theta^{-1}\} \cap B_2| \ge (1-\iota)2^{-d} |B_2| = (1-\iota)|B_1|.
\end{align*}
The conclusion of the corollary is that  $u_1 \ge \theta$ a.e. in $B_{\frac12}$, that is to say  $u \ge 1$ in $B_1$. 
 \end{proof}    

\begin{proof}[Proof of Theorem~\ref{t:whi} (Weak Harnack's inequality)]
The proof proceeds in several steps. We first reduce the proof to a universal estimate on the super-level sets $\{ u > t \}$ of the function $u$. 
We then prove the result for a universal value $t=M$. We finally get the result for all $t = M^k$ by a covering argument (ink spots). 

\paragraph{Reduction.}
  In order to prove the result, we are now used to reduce to prove that, if
  \begin{equation}
    \label{e:assum}
    \inf_{B_{\frac12}} u \le 1 \quad \text{ and } \quad \|S\|_{L^\infty (B_1)} \le  \eps_0 ,
  \end{equation}
  with   $\eps_0$ given by Corollary~\ref{c:lower-bound}. 
  then  \( \int_{B_{1/2}} u^\eps (x) \dx \le \Ceps \)  for some universal constants $\eps$ and $\Ceps$.

  In order to estimate this integral, it is sufficient to prove
  that there exists  universal constants $\Cnu$ and $\nu >0$ such that
  \begin{equation} \label{e:sls-cont}
    \forall t >1, \quad |\{ u > t \} \cap B_{\frac12} | \le \Cnu t^{-\nu}.
  \end{equation}
  Indeed, starting from the layer cake formula \cite[Theorem~1.13]{MR1817225},  we write,
  \begin{align*}
    \int_{B_{\frac12}} u^\eps (x) \dx &= \eps \int_0^\infty t^{\eps -1} |\{ u > t \} \cap B_{\frac12}| \dt \\
                                      & \le \eps \int_0^1 t^{\eps-1} |B_{\frac12}| \dt + \eps \int_1^\infty t^{\eps -1} |\{ u > t \} \cap B_{\frac12}| \dt \\
                                      & \le |B_{\frac12}| + \eps \Cnu \int_1^\infty t^{\eps -\nu -1} \dt \qquad \text{ choose } \eps = \nu /2  \\
    & \le |B_{\frac12}| +  \Cnu .
  \end{align*}
  We can further reduce the proof to the case $t=M^k$ for some universal constant $M>1$ and integers $k \ge 1$, 
  \begin{equation} \label{e:sls-disc}
    \forall k \ge 1, \quad |\{ u > M^k \} \cap B_{\frac12} | \le \Cslsd (1-\delta)^k
  \end{equation}
  with $\Cslsd>1$ and $\delta \in (0,1)$ universal too. 
  Indeed, if  \eqref{e:sls-disc} holds, then for $t >M$, we pick $k$ such that $M^k < t \le M^{k+1}$ and $\nu>0$ such that $(1-\delta) = M^{-\nu}$ and we write
  \[
    |\{ u > t \} \cap B_{\frac12} | \le |\{ u > M^k \} \cap B_{\frac12} | 
                                     \le \frac{\Cslsd}{1-\delta} (1-\delta)^{k+1} 
       \le \frac{\Cslsd}{1-\delta} \left( M^{k+1} \right)^{-\nu} 
     \le \frac{\Cslsd}{1-\delta} t^{-\nu}.
  \]
  
\paragraph{Ink spots.} In order to establish \eqref{e:sls-disc}, we show that we can apply the ink spot lemma~\ref{l:is}  with 
$E = \{ u \ge M^{k+1}\} \cap B_{\frac12}$ and $F = \{ u \ge M^k \} \cap B_{\frac12}$ and $\iota \in (0,1)$ given by Corollary~\ref{c:spreading}.
Let us verify the assumption of the lemma. 

These sets are measurable and  $E \subset F \subset B_{\frac12}$. 

By applying Corollary~\ref{c:spreading} to \( v(x) =  u \left(\frac{x}2 \right) \)
that lies in $\mathrm{DG}^- (B_2,S_v)$ 
whose source term $S_v$ is essentially bounded in $B_2$ by $\eps_{0,s}$, we get that
\[ |\{ u \ge M \} \cap B_{\frac12} \}| < (1-\iota) |B_{\frac12}|.\]
This implies that $E$ satisfies the first assumption of Lemma~\ref{l:is} since $\{ u \ge M^{k+1} \} \subset \{ u \ge M\}$.

We next check that $E$ and $F$ also satisfy the second assumption of Lemma~\ref{l:is}.
In order to do so, we consider an open ball $B_r (x_0) \subset B_{\frac12}$ such that $|E \cap B_r(x_0)| \ge (1-\iota)|B_r (x_0)|$.
This means
\begin{equation}\label{e:mk}
  |\{ u \ge M^{k+1} \} \cap  B_r (x_0)| \ge (1-\iota) |B_r (x_0)|.
\end{equation}
We aim at proving that this implies $ u \ge M^k$  a.e. in $ B_r (x_0)$.
For $x \in B_4$, we consider $\tilde u (x) = \frac{u (x_0 + rx)}{M^k}$. It satisfies $-\dive_x (\tilde A \nabla_x \tilde u) = \tilde S$ in $B_4$
with $\tilde S (x) = M^{-k}r^2 S (x_0 +rx)$. Since $M \ge 1$, $r \le 1$ and $\|S\|_{L^\infty(B_4)} \le \eps_{0,s}$, we conclude that
$\|\tilde S\|_{L^\infty(B_4)} \le \eps_{0,s}$. Then \eqref{e:mk} translates into
\[ |\{ \tilde u \ge M \} \cap B_1 | \ge (1-\iota) |B_1|.\]
Applying Corollary~\ref{c:spreading}, we conclude that $\tilde u \ge 1$ a.e. in $B_1$, that is to say, $u \ge M^k$ a.e. in $B_r(x_0)$, as desired. 

\paragraph{Conclusion.} Applying Lemma~\ref{l:is}, we conclude that
\[ |\{ u \ge M^{k+1}\} \cap B_{\frac12}| \le (1-c \iota) |\{ u \ge M^{k}\} \cap B_{\frac12}|.\]
This inequality implies \eqref{e:sls-disc} with $\delta = c \iota$ and $\Cslsd = |B_1|$. 
\end{proof}

\begin{proof}[Proof of Theorem~\ref{t:harnack} (Harnack's inequality)]
We simply combine the weak Harnack's inequality (Theorem~\ref{t:whi}) with the improved local maximum principle (Corollary~\ref{c:lmp}). 
\end{proof}

\section{Bibliographical notes}
\label{s:biblio-elliptic}

This chapter follows closely De Giorgi's original proof \cite{DeG56}. In particular, the class of functions satisfying local energy estimates are called $\mathcal{B}(E,\gamma)$ and they correspond to the elliptic De Giorgi's class $\mathrm{DG}^\pm$ (see Definition~\ref{d:DGclass}).
Lower order terms are later considered by O.~Ladyzenskaya and N.~Ural$'$tseva in their book \cite{L1}. 

The first difference between De Giorgi's original proof and the one presented in subsequent works (including the proof contained in this book) lies in the extra condition made on the functions E.~De Giorgi works with. He assumes that the functions $w(x)$ are absolutely continuous on ``almost all segments contained in in $E$ and parallel to the coordinate
axes''. Since $H^1$ functions in a ball of $\R^d$ are such that (for instance) $x_1 \mapsto \partial_{x_1} u$ is square integrable for almost every $(x_2,\dots,x_d)$ (by Fubini's theorem), they are absolutely continuous on line segments.

Another difference lies in the way that the \emph{intermediate value principle} is obtained. Let us make this vague statement more explicit.  One way or the other, the proof boils down to controlling from below the measure of the set of intermediate values of an $H^1$ function  by the natural super-level and lower-level sets.
In  \cite[Lemma~II]{DeG56}, a functional inequality is derived for elements of $\mathcal{B}(E,\gamma)$. This inequality is very close to what is nowadays known as \emph{De Giorgi's isoperimetric
inequality}, see for instance \cite[Lemma~10]{MR3525875}. It is more precise than a Poincaré-Wirtinger's inequality. A modern way to explain this difference is to compare a Poincaré inequality  with the Sobolev embedding (in $W^{1,1}$ or $W^{1,2}$). 

A Poincaré-Wirtinger's inequality (Proposition~\ref{p:poincare} in this chapter) is cooked up or used in most of (if not all) works in this trend of research:  it appears in Nash's original paper \cite[top of p.~936]{nash} under the form now known as ``Nash's inequality''.
Amusingly enough, J.~Nash mentioned that E.~M.~Stein gave him the proof. It also appears in J.~Moser's article \cite[Lemma~1]{MR159138} about elliptic equations as a true Poincaré-Wirtinger inequality: a mean is retrieved to the function before considering its $L^2$-norm.
He also used a weighted Poincaré inequality  in his work about parabolic equations \cite[Lemma~3]{moser}. We would like also to mention Kruzhkov's work \cite[Theorem~1.1]{MR171086} (announced in \cite{MR151703}) and G.~Lieberman's classical book \cite[Proposition~6.14]{MR1465184}, both on parabolic equations.

J.~Moser \cite{moser} and later N.~L.~Trudinger  \cite[Theorem 1.2]{MR226168}   established a weak Harnack's inequality for parabolic equations. The proof by J.~Moser uses a  iterative procedure that departs from De Giorgi's original one. It is now referred to as Moser's iteration. Let us briefly describe it. 
On the one hand, it was known that if $u$ is a solution of a parabolic equation, then $\varphi (u)$ is a sub-solution if $\varphi$ is convex. On the other hand, the local maximum principle allows the control of the $L^p$-norm of the solution by its $L^2$-norm for some $p>2$. This is what we called the \emph{gain of integrability} of sub-solutions. J.~Moser observed that this gain of integrability can be applied iteratively, by considering the convex function $\varphi (r) = r^{p/2}$. This leads to the local maximum principle. He also observed that the convex change of variables $\varphi (r) = 1/r$ allows one to control the infimum of a positive solution from below by its $L^2$-norm. Using again such change of variables, proving Harnack's inequality boils down to be able to relate the $L^\eps$-``norm'' of $f$ with the $L^\eps$-``norm'' of $1/f$ for an arbitrarily small $\eps >0$. For the insecure reader, we make precise that the correct statement is to relate  $L^1$-norms of $f^\eps$  $1/f^\eps$. In order to relate them, he considers the logarithm of the solution (like J.~Nash did in his original contribution) and observes that the equation that it satisfies contains a quadratic term. This quadratic term allows him to control the propagation of the level sets of the logarithm. 

G.~Lieberman \cite{MR1465184} notes that while N.~Trudinger was not the first to prove the weak Harnack's inequality, but he was the first to recognize its significance, despite it being a straightforward consequence of previously known results. G.~Lieberman also mentions that DiBenedetto and Trudinger \cite{MR778976} demonstrated that non-negative functions in the elliptic De Giorgi's class, which correspond to super-solutions of elliptic equations, satisfy a weak Harnack's inequality. Finally, G. L. Wang \cite{MR1032780} (see also \cite{MR1246215}) proved a weak Harnack's inequality for functions in the corresponding parabolic De Giorgi's classes.

Proofs of  weak Harnack's inequalities use a covering argument. N.~L.~Trudinger's proof \cite[Theorem 1.2]{MR226168} already relies on a measure lemma from N.~V.~Krylov and M.~V.~Safonov \cite{MR563790}. 
When proving it for parabolic equations, E.~DiBenedetto and N.~L.~Trudinger \cite{MR778976} uses the same procedure. G. L. Wang also mentions that he could use such a lemma in \cite{MR1032780}.
N.~V.~Krylov and M.~V.~Safonov  emphasize that their lemma is related to some works by E.~Landis \cite{MR1487894}.  Landis calls such results from measure theory ``(aptly in his opinion)'' \cite{MR563790} \emph{crawling of ink spots}.

\chapter{Parabolic equations}
\label{c:parabolic}

In this chapter, the regularity of solutions of parabolic equations under divergence
form is studied. It corresponds to Nash's original framework. As we will see, the
techniques that were introduced in the previous chapter for elliptic equations
naturally extends to parabolic ones. 

\section{Ellipticity, cylinders and H\"older continuity}

Let $I$ be a bounded interval of $\R$ of the form $(a,b]$ with $a,b \in \R$.
Let $\Omega$  be an open set of $\R^d$.
Let $\lambda, \Lambda$ be two positive constants with $\lambda \le \Lambda$.
We consider
\begin{equation}
\label{e:ellipticity}
  \mathcal{E} (\lambda, \Lambda) = \{ A \in L^\infty ( I \times \Omega, \mathbb{S}_d (\R)):
\text{ a.e. in } I \times \Omega, \quad \forall \xi \in \R^d, \lambda |\xi|^2 \le A \xi \cdot \xi \le \Lambda |\xi|^2 \}.
\end{equation}
To each $A \in \mathcal{E} (\lambda, \Lambda)$, we associate the following equation,
\begin{equation}
  \label{e:parabolic}
  \partial_t f  = \dive_x ( A \nabla_x f)  +S 
\end{equation}
posed in $I \times \Omega$ with $S \in L^1 (I \times \Omega)$. 

When $A$ is the identity matrix, equation~\eqref{e:parabolic} is simply the \emph{heat equation},
\begin{equation}
  \label{e:heat}
  \partial_t f  = \Delta_x f +S .
\end{equation}

\subsection{Invariances and cylinders}

\paragraph{Parabolic scaling.}
Let $R>0$. For $X = (t,x) \in \R\times \R^d \times \R^d$, we define the scaling operator $S_R$ by
\[ S_R (X) =  (R^2t,R x).  \]
If $f$ is a solution of the parabolic equation \eqref{e:parabolic} for some $A \in \mathcal{E} (\lambda,\Lambda)$,
then the function $f_R (X) = f (S_R (X))$ satisfies \eqref{e:parabolic} with $A$ is replaced with $A_R (X) = A (S_R(X))$.
Notice that $A_R \in \mathcal{E}(\lambda,\Lambda)$. 

\paragraph{Translation invariance.}
Given $X_0 = (t_0,x_0) \in \R \times \R^d$ and a solution $u$ of \eqref{e:parabolic} with $A\in \mathcal{E}(\lambda,\Lambda)$, the function
$v (X) = u(X_0 +X) = u(t_0+t,x_0+x)$ is a solution of \eqref{e:parabolic} with $A$ replaced with
$A_0 (X) = A(X_0 +X)$. Notice that $A_0 \in \mathcal{E}(\lambda,\Lambda)$. 

\paragraph{Parabolic cylinders.}
\label{pcyl}
For $R>0$ and $X_0 = (t_0,x_0) \in \R \times \R^d$, the parabolic cylinder $Q_R (X_0)$ is defined by
\[ Q_R (X_0) = (t_0-R^2,t_0] \times B_R (x_0). \]
This family of cylinders encodes the parabolic scaling and the translation invariance of the class
of parabolic equations of the form \eqref{e:parabolic}. Indeed, we can write $Q_R (X_0) = X_0 + S_R (Q_1)$
with $Q_1 = (-1,0] \times B_1$ (unit cylinder).

\begin{remark}
We also draw the attention of the reader towards the
fact that these cylinders are neither open nor closed and that the point $X_0$ lies at the top of it.
One could justify this choice by arguing for instance that  the information used to study parabolic equations should come from the past. 
\end{remark}

\subsection{Parabolic H\"older regularity}

In this section, we give a sufficient condition for a function $u (t,x)$ to be
H\"older continuous. It is expressed in terms of the oscillation of $u$ in parabolic cylinders. 
\begin{defi}[Parabolic H\"older spaces and semi-norms]
  Let $\alpha \in (0,1]$ and $Q$ a parabolic cylinder. A function $u \colon Q \to \R$ is in the
  space $\Cpar^\alpha (Q)$ if there exists a constant $C_0>0$ such that,
for all  $X=(t,x),Y=(s,y) \in Q$,
\[ |u(X)-u(Y)| \le C_0 \left( |t-s|^{1/2} + |x-y| \right)^\alpha .\]
The smallest $C_0$ such that the previous inequality holds true is denoted by $[u]_{\Cpar^\alpha (Q)}$.
The space $\Cpar^\alpha (Q)$ is equipped with the norm $\|u\|_{\Cpar^\alpha (Q)} = \sup_Q |u| + [u]_{\Cpar^\alpha (Q)}$.
\end{defi}
\begin{remark}[Classical H\"older regularity]
A  function $u \in \Cpar^\alpha (Q)$ is $\alpha/2$-H\"older
  continuous in $t$ while it is $\alpha$-H\"older continuous in $x$. This is a general ``fact''
  about the regularity of solutions of parabolic equations: the solution is twice more regular in $x$ than in $t$. 
\end{remark}
\begin{remark}[Parabolic distance]
  For $X=(t,x),Y=(s,y) \in \R \times \R^d$, the quantity $|t-s|^{1/2} + |x-y|$ defines
  a distance, sometimes referred to as the parabolic distance. 
\end{remark}

We first recall that for a function $u \colon A \to \R$ essentially bounded on a Borel set $A$,
its oscillation on $A$ is defined as
\[ \osc_A u = \esssup_A u - \essinf_A u \]
where $\esssup_A u$ and $\essinf_A u$ are the essential supremum and infimum of $u$ on $A$. 
\begin{prop}[Parabolic H\"older regularity via oscillations] \label{p:holder-para}
  Consider a parabolic cylinder $Q$ and a function $u \in L^\infty(Q)$.
  Assume that there exist $\alpha \in (0,1]$ and $C_0>0$ such that for all $X \in Q$ and all $r>0$,
  we have $\osc_{Q_r (X) \cap Q} u \le C_0 r^\alpha$. Then $u$ is H\"older continuous in $Q$. More precisely, for all
\end{prop}
\begin{proof}
The proof proceeds in three steps. 
  
  \paragraph{Continuous functions.} We first prove that if $u$ is continuous in $Q$ and satisfies the assumption for some constant $C_0$, then
  it satisfies the conclusion with the same constant $C_0$. In order to do so, we consider $X=(t,x)$ and $Y = (s,y)$ and
  we define $r = |t-s|^{1/2} + |x-y|$.
  We observe next that $Y$ lies in the closure of $Q_r (X)$. Indeed, $|t-s| \le r^2$ and $|x-y| \le r$.
  We use next that $u$ is continuous in $Q$ in order to write,
  \[ u (X) - u (Y) \le \underset{Q_r (X) \cap Q}{\esssup} u - \underset{Q_r (X) \cap Q}{\essinf} u \le C_0 r^\alpha.\]
  We can now exchange the role of $X$ and $Y$ and conclude that
  \[ |u(X)-u(Y)| \le C_0 \left( |t-s|^{1/2} + |x-y| \right)^\alpha .\]

  \paragraph{Regularization.} We now consider a merely essentially bounded function $u$. Given $\eps \in (0,1)$, we consider
  a smooth non-negative function $\theta \colon \R \to \R$ supported in $[-1,0]$ and such that $\int_{\R} \theta (t) \dt =1$.
  The location of the support of $\theta$ is important to ensure that the mean in time is computed from past times.
  We also consider a smooth non-negative function $\rho \colon \R^d \to \R$ supported in $\bar B_1$ and such that $\int_{\R^d} \rho (x) \dx =1$.
  We then re-scale these functions with a parameter $\eps \in (0,1)$ and consider $\rho_\eps (x) = \eps^{-d} \rho(\eps^{-1}x)$ and $\theta_\eps (t) = \eps^{-2} \theta (\eps^{-2}t)$ and
  \[ u^\eps (t,x) = \iint_{Q} u(s,y) \rho_\eps (x-y) \theta_\eps (t-s) \ds \dy.\]
  Now consider $X_0$ and $R>0$ such that $Q = Q_R (X_0)$. We now consider $X_1 =(t_1,x_1) \in Q_{R-\eps} (X_0)$.
  We remark that for $Z=(r,z) \in Q_\eps$, we have $X_1-\eps Z \in Q_R(X_0)$. In particular, we can write for $X_1 \in Q$, 
  \begin{align*}
    \underset{X \in Q_r (X_1) \cap Q_{R-\eps} (X_0)}{\osc} u^\eps = &  \sup_{Q_r (X_1) \cap Q_{R-\eps} (X_0)} \iint_{Q} u(s,y) \rho_\eps (x-y) \theta_\eps (t-s) \ds \dy  \\
    - & \inf_{X \in Q_r (X_1) \cap Q_{R-\eps} (X_0)} \iint_{Q} u(s,y) \rho_\eps (x-y) \theta_\eps (t-s) \ds \dy \\
    \le &   \iint_{Q} \sup_{Q_r (X_1-\eps Z) \cap Q_R (X_0)} u(t-\eps r,x-\eps z) \rho (z) \theta (r) \ds \dy  \\
    - &  \iint_{Q} \inf_{X \in Q_r (X_1) \cap Q_R (X_0)} u(t-\eps r,x-\eps z) \rho (z) \theta (r) \ds \dy \\
    \le &   \iint_{Q} \underset{Q_r (X_1-\eps Z) \cap Q_R (X_0)}{\osc} u(t-\eps r,x-\eps z) \rho (z) \theta (r) \ds \dy  \\
    \le &   C_0 r^\alpha.
  \end{align*}

  \paragraph{Conclusion.} Now we conclude as in the elliptic case (Proposition~\ref{p:holder-osc}).
     Let $\eps_0>0$ and $\eps \in (0,\eps_0)$. We conclude from Step~1 that for all $X,Y \in Q_{R-\eps_0} (X_0)$,
   \[ |u^\eps (X) - u^\eps (Y) | \le C_0 \left( |t-s|^{\frac12} + |x-y|\right)^\alpha.\]
   By dominated convergence, we see that  $u^\eps \to u$ a.e. in $Q_{R-\eps_0} (X_0)$. We conclude that for all $X,Y \in Q_{R-\eps_0} (X_0)$,
   \[ |u(X) - u(Y) | \le C_0 \left( |t-s|^{\frac12} + |x-y|\right)^\alpha.\]
   Since $\eps_0>0$ is arbitrarily small, we conclude that $[u]_{C^\alpha (Q)} \le C_0$. 
\end{proof}

\section{Weak  solutions and De Giorgi \& Nash's theorem}

\subsection{The space $H^{-1}(\Omega)$ and $L^2(I, H^1(\Omega))$.}

\begin{defi}
The space $H^{-1} (\Omega)$ is the topological dual space of $H^1_0 (\Omega)$. 
\end{defi}
It is a Banach space when equipped with the norm
\[\| f \|_{H^{-1} (\Omega)} = \sup_{\stackrel{v \in H^1_0 (\Omega)}{\|v\|_{H^1 (\Omega)} \le 1}}|\langle f, v \rangle|.\]
\begin{example}
Weak derivatives (in the sense of distributions) of square integrable functions are important examples of elements of $H^{-1} (\Omega)$. More generally,
 given functions $h_0, h_1,\dots, h_d \in L^2(\Omega)$, we define the linear application $F$ on $H^1_0 (\Omega)$ as
 \[ \forall v \in H^1_0(\Omega), \qquad \langle F, v\rangle = \int_{\Omega} h_0 v \dx - \sum_{i=1}^d \int_{\Omega} h_i \partial_i v .\]
  Then $F \in H^{-1} (\Omega)$. 
\end{example}
\begin{remark}
The converse is true, see \cite[p.~299]{MR2597943}.
\end{remark}

\subsection{Weak solutions}

Our notion of weak solutions is guided by the following formal computation. Consider a smooth test function $\varphi(x)$ supported $B$ and multiply the equation
by $u \varphi^2$. After integrating in $B$ and writing $A = (\sqrt{A})^2$ with $\sqrt{A}$ real symmetric and semi-definite, one gets,
\begin{equation}
  \label{e:energy}
  \frac12 \frac{\dd }{\dt} \int_B u^2 \varphi^2 \dx + \int_B (A \nabla_x u \cdot \nabla_x u) \varphi^2 = 2 \int_B (\sqrt{A} \nabla_x u) \varphi \cdot u \sqrt{A} \nabla_x \varphi
  + \int_B S u \varphi.
\end{equation}
In particular, using Cauchy-Schwarz inequality twice (for each term of the right hand side) and the ellipticity of $A$, we get
\[ \frac12 \frac{\dd }{\dt} \int_B u^2 \varphi^2 \dx + \frac\lambda2 \int_B |\nabla_x u|^2 \varphi^2 \le 2\Lambda \int_B  u^2  |\nabla_x \varphi|^2 \dx + \frac12 \int_B S^2 \dx + \frac12 \int_B u^2 \varphi^2 \dx.\]
This formal computation suggests that only assuming that $u$ and $\nabla_x u$ are square integrable, one can deduce that $u \in L^\infty (I, L^2 (\Omega))$ and
(thanks to the equation) that $\partial_t u \in L^2 (I,H^{-1}(\Omega))$. A classical result from functional analysis, known as Aubin-Lions's lemma,
then implies time continuity with values in $L^2 (\Omega)$, see for instance \cite[Theorem II.5.16]{boyer2012mathematical}. 
\begin{defi}[Weak solutions]\label{def:weak}
  Let $I$ be an interval of the form $(a,b]$ and $\Omega$ be an open set of $\R^d$ and $S \in L^2(\Omega)$.
  A function $u \colon I \times \Omega \to \R$ is a \emph{weak solution} of a parabolic equation of the form
  \( \partial_t u - \dive_x (A \nabla_x u ) = S \) in $I \times \Omega$ if,
  \[ u \in C (I, L^2 (\Omega)), \quad \nabla_x u \in L^2 (I \times \Omega), \quad \partial_t u \in L^2 (I,H^{-1}(\Omega)),\]  and for all $\varphi \in L^2 (I, H^1_0 (\Omega))$,
  \[ \int_I \langle \partial_t u , \varphi \rangle_{H^{-1}(\Omega),H^1_0 (\Omega)} \dt = -\iint_{I \times \Omega} A \nabla_x u \cdot \nabla_x \varphi \dt \dx + \iint_{I \times \Omega} S \varphi \dt \dx. \]
\end{defi}
\begin{remark}[Weaker solutions]
  We can use a weaker notion of solutions by only assuming that $u \in L^2 (I \times \Omega)$, that the weak gradient $\nabla_x u \in L^2(I \times \Omega)$ and
  that the equation is satisfied in the weak sense by only testing with  smooth and compactly supported functions (\textit{i.e.} in the sense of distributions). 
  It requires some extra work to prove that such weaker solutions  are in fact weak solutions in the sense of Definition~\ref{def:weak}.
\end{remark}
\begin{remark}[Weak sub- and super-solutions] \label{r:sub-super}
  We can also consider functions that are only sub-solutions or super-solutions of a parabolic equation. This would be subsets of those weaker solutions
  that are discussed in the previous remark.
\end{remark}

\subsection{Local energy estimates}

Now that we made precise the notion of solutions that are going to work with, we derive the local estimates that will be enough for us to establish their continuity.
\begin{prop}[Local energy estimates for weak solutions] \label{p:lee}
  Let $u$ be a weak solution of $\partial_t u - \dive_x (A \nabla_x u) = S$ in $I \times B$ with $I=(a,b]$ and $B$ an open ball.
  Consider $B_r (x_0 ) \subset B_R (x_0) \subset B$ and a truncation function $\rho$ as in Lemma~\ref{l:trunc}. 
  If $S \in L^\infty (I \times B)$ and $A \in \mathcal{E} (\lambda,\Lambda)$, then for all $t_1,t_2 \in I$, $r,R>0$ and $X_0 \in I \times B$ such that $r<R$ and $Q_R (X_0) \subset I \times B$, and all $\kappa \in \R$, 
\begin{align*}
   \int_{B_r(x_0)} (u-\kappa)_+^2 (t_2 ,x)  \dx &  + \lambda  \int_{t_1}^{t_2} \int_{B_r (x_0)}  | \nabla_x (u-\kappa)_+ |^2  \dt \dx    
  \le   \int_{B_r(x_0)} (u-\kappa)_+^2 (t_1,x)  \dx \\ & + \frac{16\Lambda}{(R-r)^2} \int_{t_1}^{t_2} \int_{B_R(x_0)}   (u-\kappa)_+^2 \dt \dx 
        + \int_{t_1}^{t_2} \int_{B_R(x_0)} S (u-\kappa)_+ \dt \dx. 
\end{align*}
\end{prop}
In order to prove this proposition, we need the following technical lemma that allows us to use $(u-\kappa)_+ \varphi$ as a test function in the definition of weak solutions, for any
compactly supported and smooth function $\varphi$. 
\begin{lemma} \label{l:tech-pdg}
  Let  $u \in L^2(I \times \Omega)$ be such that $\nabla_x u \in L^2 (I \times \Omega)$ and $\partial_t u \in L^2 (I, H^{-1}(\Omega))$. Then for all $\varphi \in C^\infty_c (I \times \Omega)$ and $\kappa \in \R$, 
  we have $(u-\kappa)_+ \varphi \in L^2 (I , H^1_0 (\Omega))$ and
  \[ \int_I \langle \partial_t u , (u-\kappa)_+ \varphi  \rangle_{H^{-1}, H^1_0} \dt = - \frac12 \iint_{I \times \Omega} (u-\kappa)_+^2 \partial_t \varphi \]
  where $\langle \partial_t u , (u-\kappa)_+ \varphi  \rangle_{H^{-1}, H^1_0}$ denotes
  $\langle \partial_t u (t,\cdot), (u-\kappa)_+ (t,\cdot) \varphi (t,\cdot) \rangle_{H^{-1}(\Omega), H^1_0(\Omega)}$. 
\end{lemma}
\begin{proof}
  We consider an even non-negative function $\rho \in C^\infty_c (\R \times \R^d)$ supported in $(-1,1) \times B_1$ and $\int_{\R \times \R^d} \rho (t,x) \dt \dx = 1$ and define for all $\eps \in (0,1)$ the function  $\rho_\eps (t,x) = \eps^{-2-d} \rho (\eps^{-2} t, \eps^{-1} x)$. For any  $\nu \in (0,1)$, we also consider a function $P_\nu (r)= (r-\kappa)_+ \ast \theta_\nu$.
  Then consider the smooth and compactly supported function $\psi_{\eps,\nu} = \bigg(  P_\nu (u_\eps ) \varphi \bigg) \star \rho_\eps$ with $u_\eps = u \ast \rho_\eps$.
  \begin{align*}
    \iint_{I \times B} u \partial_t \psi_{\eps,\nu}
    &= \iint_{I \times B} u \bigg[ \partial_t ( P_\nu (u_\eps) \varphi) \ast \rho_\eps \bigg] =  \iint_{I \times B} u_\eps \partial_t ( P_\nu (u_\eps) \varphi) \\
     &= - \iint_{I \times B} (\partial_t u_\eps)  P_\nu (u_\eps) \varphi       = - \frac12 \iint_{I \times B} \left( \partial_t Q_\nu (u_\eps) \right) \varphi  
  \end{align*}
  with $Q_\nu(r) = \frac12(r-\kappa)_+^2 \ast \theta_\nu$. 
  After integrating by parts again, we thus get
  \begin{equation}
    \label{e:ee-p-inter}
    \int_{I \times B} \langle \partial_t u,  \psi_{\eps,\nu} \rangle_{H^{-1}, H^1_0} \dt =  \frac12 \iint_{I \times B} Q_\nu (u_\eps) \partial_t \varphi.   
  \end{equation}
  
  In order to conclude, we now pass to the limit as $\eps \to 0$ and $\nu \to 0$. We use the three following facts:  
\begin{align*}
&    u_\eps \to u \text{ in } L^2 (I,H^1(\Omega)), \\ 
&  P_\nu (u_\eps) \to P_\nu (u) \text{ in } L^2(I,H^1(\Omega)), \\
&  Q_\nu (u_\eps) \to Q_\nu (u) \text{  in } L^1 (I \times \Omega).
\end{align*} 
The fact that $u_\eps \to u$ in $L^2 (I \times \Omega)$ is a general fact for mollifiers \cite[Proposition II.2.25]{boyer2012mathematical}. Then $\nabla_x u_\eps = (\nabla_x u) \ast \rho_\eps$
in the sense of distributions, from which we get $\nabla_x u_\eps \to \nabla_x u$ in $L^2 (I \times \Omega)$. This leads to the first fact.
The second fact uses that $P_\nu$ is $1$-Lipschitz. Indeed,
\[ |P_\nu (u_\eps) - P_\nu (u)| \le |u_\eps -u| \]
and this implies the convergence in $L^2 (I \times \Omega)$. Then $\nabla_x P_\nu (u_\eps) = P'_\nu (u_\eps) \nabla_x u_\eps$. In particular,
\begin{align*}
  |\nabla_x P_\nu (u_\eps) - \nabla_x P_\nu (u)| & \le |P'_\nu (u_\eps) - P'_\nu (u)| |\nabla_x u| + |P'_\nu (u_\eps)| |\nabla_x u_\eps - \nabla_x u| \\
  & \le |P'_\nu (u_\eps) - P'_\nu (u)| |\nabla_x u| +  |\nabla_x u_\eps - \nabla_x u|. 
\end{align*}
We can apply dominated convergence to prove the $L^2$-convergence of the first term and we already proved the $L^2$-convergence of the second one.
We thus established the second fact. Finally,
\begin{align*}
 & \int_{I \times \Omega} | Q_\nu (u_\eps) - Q_\nu (u)|
  \le \frac12 \int_{I \times \Omega} \int_\R  |(u_\eps - \nu s -\kappa)_+^2 -(u-\nu s-\kappa)_+^2| \theta (s) \ds , \\
 & \le \frac12 \int_{I \times \Omega} \int_\R  \bigg(|(u_\eps - \nu s -\kappa)_+ -(u-\nu s-\kappa)_+||(u_\eps - \nu s -\kappa)_+ +(u - \nu s-\kappa)_+|  \bigg) \theta (s) \ds \\
 & \le \frac12 \int_{I \times \Omega} |u_\eps -u| \int_\R  \bigg(|(u_\eps - \nu s -\kappa)_+ +(u - \nu s-\kappa)_+|  \bigg) \theta (s) \ds \\
  & \le \frac12 \|u_\eps -u\|_{L^2 (I\times \Omega} \bigg( \|u_\eps\|_{L^2 (I\times \Omega)} +\|u\|_{L^2 (I\times \Omega)} + 2 |\kappa| + 2 \nu \int_{\R} |s| \theta (s) \ds \bigg) 
\end{align*}
Since $u_\eps \to u$ in $L^2 (I \times \Omega)$, $\|u_\eps\|_{L^2 (I\times \Omega)}$ is bounded independently of $\eps \in (0,\eps_0)$. This completes
the justification of the third fact. 

Thanks to  those three facts, we can pass to the limit in $\eps \to 0$ in \eqref{e:ee-p-inter} and get
\[    \int_{I \times B} \langle \partial_t u,  P_\nu(u) \varphi \rangle_{H^{-1}, H^1_0} \dt =  \frac12 \iint_{I \times B} Q_\nu (u) \partial_t \varphi.   \]
Arguing as above, we see that $P_\nu (u) \varphi \to P(u) \varphi$ in $L^2 (I , H^1(\Omega))$ and $Q_\nu (u) \to (u-\kappa)^2_+$ in $L^1 (I \times \Omega)$
and we conclude. 
\end{proof}
\begin{proof}[Proof of Proposition~\ref{p:lee}]
  We proceed in several steps. We first use $(u-\kappa)_+ \varphi$ as a test function in the definition of weak solutions (Step~1)
  before considering a specific test function $\varphi$ (Step~2). 

  \paragraph{Step~1.}
  Let $\varphi \in C^\infty_c (I \times B)$. Using Lemma~\ref{l:tech-pdg}, we know that $(u-\kappa)_+ \varphi \in L^2 (I, H^1_0 (B))$ and we have
  \begin{align*}
    - \frac12 \iint_{I \times B} (u -\kappa)_+^2 \partial_t \varphi 
    = &  \int_I \langle \partial_t u (t), (u-\kappa)_+ \varphi (t) \rangle_{H^{-1},H^1_0} \dt  \\
    = &  - \iint_{I \times B} A \nabla_x u \cdot \nabla_x (u-\kappa)_+ \varphi  \\
    & - \iint_{I \times B} \bigg( A \nabla_x u \cdot  \nabla_x  \varphi \bigg) (u-\kappa)_+  + \iint_{I \times B} S (u-\kappa)_+ \varphi \\
        = &  - \iint_{I \times B} A \nabla_x (u-\kappa_+) \cdot \nabla_x (u-\kappa)_+ \varphi  \\
    & - \iint_{I \times B} \bigg( A \nabla_x (u-\kappa)_+ \cdot  \nabla_x  \varphi \bigg)  (u-\kappa)_+   + \iint_{I \times B} S (u-\kappa)_+ \varphi .
  \end{align*}
  Like in the elliptic case, we used the fact that $\nabla_x (u-\kappa)_+ = \un_{\{ u \ge \kappa\}} \nabla_x u$ and $\un_{\{ u \ge \kappa\}} = \un^2_{\{ u \ge \kappa\}}$. 
  After remarking that it is positive, we  move the first term  of the last equality in the right hand side to the left hand side, 
  \begin{align*}
    - \frac12 \iint_{I \times B} (u -\kappa)_+^2 \partial_t \varphi
    & + \iint_{I \times B} A \nabla_x (u-\kappa_+) \cdot \nabla_x (u-\kappa)_+ \varphi  \\
    & = - \iint_{I \times B} \bigg( A \nabla_x (u-\kappa)_+ \cdot   \nabla_x  \varphi \bigg)  (u-\kappa)_+   + \iint_{I \times B} S (u-\kappa)_+ \varphi.
  \end{align*}

  \paragraph{Step~2.}
  Let $\rho$ be the truncation function from Lemma~\ref{l:trunc} that is supported in $B_R (x_0)$ and equal $1$ in $B_r(x_0)$.
  Given $ t_1,t_2 \in I$ with $t_1 < t_2$,   we consider a 1D mollifier $\theta \colon \R \to \R$ (smooth and of unit mass) that is supported in $[-1,0]$, and
  the smooth function $\Theta_\eps \colon (0,+\infty) \to \R$ such that $\Theta_\eps (0) = 0$ and for all $t \in \R$, we have $\Theta_\eps'(t) = \theta_\eps (t-t_1) - \theta_\eps (t-t_2)$.
  We choose $\eps>0$ small enough so that $t_1-\eps \in I$. 
  
  We now use Step~1 with $\varphi (t,x) = \Theta_\eps (t) \rho^2 (x)$. Since $\partial_t \varphi (t,x) = (-\theta_\eps (t-t_2) +\theta_\eps (t-t_1)) \rho^2 (x)$, we can
  rearrange terms and get,
  \begin{equation}\label{e:lee-prelim}
    \begin{aligned}
    \frac12 & \iint_{I \times B} (u-\kappa)_+^2 (t,x) \rho^2 (x)  \theta_\eps (t-t_2)
     + \iint_{I \times B} |\sqrt{A} \nabla_x (u-\kappa_+) |^2  \rho^2 \Theta_\eps   \\
    = &     \frac12 \iint_{I \times B} (u-\kappa)_+^2 (t,x) \rho^2 (x) \theta_\eps (t-t_1) + \iint_{I \times B} S (u-\kappa)_+ \rho^2 \Theta_\eps \\
            & - 2 \iint_{I \times B} \bigg( \sqrt{A} \nabla_x (u-\kappa)_+ \cdot   \sqrt{A} \nabla_x  \rho \bigg) \rho (u-\kappa)_+ \Theta_\eps.
  \end{aligned}
\end{equation}
We next let $\eps$ go to $0$. Since $u \in C(I,L^2 (B))$, we claim that $t \mapsto \int_B (u-\kappa)_+^2 \rho^2 (x) \dx$ is continuous.
  Indeed, letting $M = \|\rho\|_{L^\infty(B)}$ and $m(\cdot)$ be the modulus of continuity on $I$, we write
  \begin{align*}
    \bigg| \int_B (u(t)-\kappa)_+^2 \rho^2 (x) \dx - &\int_B (u(s)-\kappa)_+^2 \rho^2 (x) \dx \bigg| \\
    & \le M^2 \int_B \left| (u(t)-\kappa)_+^2 - (u(s)-\kappa)_+^2 \right|  \dx \\
    & \le M^2 \int_B \left| (u(t)-\kappa)_+ - (u(s)-\kappa)_+ \right| \left| (u(t)-\kappa)_+ - (u(s)-\kappa)_+ \right|  \dx \\
    \intertext{We now use that $r \mapsto (r_\kappa)_+$ is $1$-Lipschitz,}
                                                     & \le M^2 \int_B \left| u(t)- u(s) \right| ( |u(t)| + |u(s)| + 2 |\kappa| )  \dx \\
                                                     & \le M^2 \left( \int_B (u(t)- u(s))^2 \dx \right)^{1/2} \left( 2\int_B ( |u(t)|^2 + |u(s)|^2 + 4 \kappa^2 )  \dx \right)^{1/2} \\
    & \le M^2 m (|s-t|) \left( 4 \sup_{\tau \in I} \|u(\tau)\|_{L^2(B)}^2 + 8 \kappa^2 |B| \right)^{1/2} .
  \end{align*}
  The time continuity then implies that for $i=1,2$,
  \[ \iint_{I \times B} (u-\kappa)_+^2 (t,x) \rho^2 (x)  \theta_\eps (t-t_i) \to \int_{B} (u(t_i) -\kappa)_+^2 (t,x) \rho^2 (x) \dx \qquad \text{ as } \eps \to 0.\]
  As far as terms containing $\Theta_\eps$ are concerned, they are much simpler since we can apply dominated convergence theorem and replace $\Theta_\eps$ with $\un_{[t_1,t_2]}$.
  Passing to the limit $\eps \to 0$ in \eqref{e:lee-prelim} leads to,
    \begin{align*}
    \frac12 & \int_{B} (u-\kappa)_+^2 (t_2,x) \rho^2 (x) \dx     + \int_{t_1}^{t_2} \int_{B} |\sqrt{A} \nabla_x (u-\kappa_+) |^2  \rho^2   \\
    = &     \frac12 \int_{B} (u-\kappa)_+^2 (t_1,x) \rho^2 (x) \dx + \int_{t_1}^{t_2} \int_{ B} S (u-\kappa)_+ \rho^2  \\
       -     & 2 \int_{t_1}^{t_2} \int_{B} \bigg( \sqrt{A} \nabla_x (u-\kappa)_+ \cdot   \sqrt{A} \nabla_x  \rho \bigg) \rho (u-\kappa)_+ .
  \end{align*}
  We now treat the error term on the last line like we did in the elliptic case, see also the formal computation preceding the definition of weak solutions,
    \begin{align*}
    \frac12 & \int_{B} (u-\kappa)_+^2 (t_2,x) \rho^2 (x) \dx     + \frac12 \int_{t_1}^{t_2} \int_{B} |\sqrt{A} \nabla_x (u-\kappa_+) |^2  \rho^2   \\
    \le &     \frac12 \int_{B} (u-\kappa)_+^2 (t_1,x) \rho^2 (x) \dx  + \int_{t_1}^{t_2} \int_{B} S (u-\kappa)_+ \rho^2  
          +     2 \int_{t_1}^{t_2} \int_B |\sqrt{A} \nabla_x  \rho |^2 (u-\kappa)_+^2 .
  \end{align*}
We finally use ellipticity constants of $A$ and properties of $\rho$ to get,
     \begin{align*}
      & \int_{B} (u-\kappa)_+^2 (t_2,x) \rho^2 (x) \dx     + \lambda \int_{t_1}^{t_2} \int_{B} | \nabla_x (u-\kappa_+) |^2  \rho^2   \\
     \le &    \int_{B} (u-\kappa)_+^2 (t_1,x) \rho^2 (x) \dx  + 2\int_{t_1}^{t_2} \int_{B} S (u-\kappa)_+ \rho^2  
           +     \frac{16\Lambda}{(R-r)^2}  \int_{t_1}^{t_2} \int_{B_R (x_0)}  (u-\kappa)_+^2 .
     \end{align*}
    We can now let $\rho$ converge to $\un_{B_r (x_0)}$ in order to get the result by dominated convergence.
\end{proof}

\subsection{De Giorgi \& Nash's theorem}

The goal of this chapter is to establish the parabolic counterpart of De Giorgi's theorem for elliptic equations by following De Giorgi's seminal path.
In the following statement, a universal constant refers to a constant that only depends on dimension $d$ and ellipticity constants $\lambda,\Lambda$ from \eqref{e:ellipticity}. 
\begin{thm}[De Giorgi \& Nash]\label{t:dg-parab}
  Let   $A \in \mathcal{E}(\lambda,\Lambda)$ with $\Omega = B_1$ and $\lambda, \Lambda >0$.
  There exist two universal constants $\alpha \in (0,1]$ and $\Cdg >0$ such that for any
  weak solution $u$ of $\partial_t u = \dive_x (A \nabla_x u) + S$ in $I \times B$ with $S \in L^\infty (I \times B)$, we have
  \[ \|u\|_{\Cpar^\alpha (Q_{\frac12})} \le \Cdg \left( \|u\|_{L^2 (Q_1)} + \|S\|_{L^\infty (Q_1)} \right).\]
\end{thm}
\begin{remark}[Local maximum principle and expansion of positivity]
  Like in the elliptic case, the proof of the H\"older regularity of solutions consists in establishing first a local maximum principle
  for elements of $\mathrm{pDG}^+$ and, second, in  proving that the 
  oscillation of solutions improves by a universal factor when zooming in with a universal factor thanks to an expansion of positivity
  property for elements of $\mathrm{pDG}^-$. 
\end{remark}
\begin{remark}[Parabolic De Giorgi's classes]
We state the theorem for weak solutions but we will prove it for functions that are in the parabolic De Giorgi's class. More precisely,
the local maximum principle  holds  for functions in the parabolic De Giorgi's class $\mathrm{pDG}^+$ while the the infimum lifting property
 holds  for functions in the parabolic De Giorgi's class $\mathrm{pDG}^-$.
\end{remark}

\section{The parabolic De Giorgi's class pDG\textsuperscript{+} \& the maximum principle}

We now define the parabolic De Giorgi's class that will ensure that the local maximum principle holds true. 
\begin{defi}[The parabolic De Giorgi's class $\mathrm{pDG}^+$]
  Let $Q$ be an open set of $\R \times \R^d$ of the form $I \times B$ where $I = (a,b]$ and $B$ is an open ball.
  A function $u \colon I \times B$ lies in the \emph{parabolic De Giorgi's class} $\mathrm{pDG}^+(I \times B,S)$
  if $u \in  L^\infty (I, L^2 (B))$ and $\nabla_x u \in L^2 (I \times B)$ and for all $X_0 = (t_0,x_0) \in I \times B$ and all $r,R>0$ such that
  $r< R$ and $Q_R (X_0) \subset I \times B$, and all $\kappa \in \R$,
  \begin{multline} \label{e:pdgc}
    \underset{t \in (t_0-r^2,t_0]}{\sup}  \int_{B_r(x_0)} (u-\kappa)_+^2 \dx + \int_{Q_r (X_0)} |\nabla_x (u-\kappa)_+|^2 \dt \dx \\
    \le \Cpdgp \left( \frac{1}{(R-r)^2} + \frac1{r^2} \right) \int_{Q_R (X_0)} (u-\kappa)_+^2 \dt \dx + \Cpdgp \int_{Q_R (X_0)} |S| (u-\kappa)_+ \dt \dx
  \end{multline}
  for some universal constant $\Cpdgp \ge 1$.
\end{defi}
\begin{remark}[Universal constants]
Like in the previous chapter, a constant is \emph{universal} if it only depends on dimension and the constants $\Cpdgp, \Cpdgm$ appearing in the definition of
the De Giorgi's classes $\mathrm{pDG}^\pm$ (see below for $\mathrm{pDG}^-$). 
\end{remark}
We will use repeatedly the following elementary observation that we record in a lemma. 
\begin{lemma}[Invariance by scaling and translation of the De Giorgi's class]\label{l:invariance-pdgp}
  If $u \in \mathrm{pDG}^+(I \times B, S)$ and $Q_r (X_0) \subset I \times B$, then the function $v = \lambda u (\frac{t-t_0}{r^2}, \frac{x-x_0}r)$
  lies in $\mathrm{pDG}^+(Q_1, \mathfrak{S})$
  with $\mathfrak{S} (t,x)= \frac{\lambda}{r^2} S (\frac{t-t_0}{r^2} , \frac{x-x_0}r )$. 
\end{lemma}

  \subsection{Weak solutions are in pDG\textsuperscript{+}}

\begin{prop}[Weak solutions are in the parabolic DG class] \label{p:weak-pdg}
  Let $u$ be a weak solution of $\partial_t u = \dive_x ( A \nabla_x u ) + S$ in $I \times B$ with $S \in L^2 (I \times B)$. Then
  $u \in \mathrm{pDG}^+ (I \times B,S)$.
\end{prop}
\begin{remark}
 The parabolic De Giorgi's class $\mathrm{pDG}^+$ also contains all sub-solutions mentioned in Remark~\ref{r:sub-super}.
\end{remark}
\begin{proof}
The local energy estimates satisfied by the weak solution $u$ given by Proposition~\ref{p:lee} can be expressed with 
  \[\left\{
  \begin{aligned}
    E (t) &= \int_{B_r(x_0)} (u(t)-\kappa)_+^2 \dx, \\
    D(t) & = \int_{B_r (x_0)} |\nabla_x (u-\kappa)_+|^2 \dx, \\
    \Sigma (t) &= \int_{B_R (x_0)} \left ( |S|(u-\kappa)_+ + \frac{16\Lambda}{(R-r)^2} (u-\kappa)_+^2 \right)\dx. 
  \end{aligned} \right.
\]
We obtain that for all $t_1,t_2 \in (t_0-R^2,t_0]$ with $t_1 < t_2$, we have
\[ E (t_2) + \lambda \int_{t_1}^{t_2} D (t) \dt \le E (t_1) + \int_{t_1}^{t_2} \Sigma (t) \dt. \]

\paragraph{Control of the dissipation.}
Taking $t_1 \in  [t_0-R^2,t_0-r^2]$ and $t_2 = t_0$ and discarding some part of the time interval in the left hand side, we get
\[ \lambda \int_{t_0-r^2}^{t_0} D (t) \dt \le E (t_1) + \int_{t_0-R^2}^{t_0} \Sigma (t) \dt. \]
We can now take the mean with respect to $t_1 \in (t_0-R^2,t_0-r^2]$ and get 
\[ \lambda \int_{t_0-r^2}^{t_0} D (t) \dt \le \frac{1}{R^2-r^2} \int_{t_0-R^2}^{t_0} E (t) \dt + \int_{t_0-R^2}^{t_0} \Sigma (t) \dt. \]
From the definition of $E,D$ and $\Sigma$, we get 
\[
  \lambda \int_{Q_r (x_0)} |\nabla_x (u-\kappa)_+|^2 \le
  \int_{Q_R (x_0)} \left ( |S|(u-\kappa)_+ + \left( \frac1{R^2-r^2} + \frac{16\Lambda}{(R-r)^2} \right) (u-\kappa)_+^2 \right).
\]

\paragraph{Control of the energy.}
Now taking $t_2 \in (t_0-r^2,t_0]$ and $t_1 \in (t_0-R^2,t_0 -r^2]$, we discard the dissipation and take a mean in $t_1$ and get, 
\[
\sup_{t_2 \in (t_0-r^2,t_0]}\int_{B_r(x_0)} (u(t_2)-\kappa)_+^2 \dx \le
  \int_{Q_R (x_0)} \left ( |S|(u-\kappa)_+ + \left( \frac1{R^2-r^2} + \frac{16\Lambda}{(R-r)^2} \right) (u-\kappa)_+^2 \right).
\]
Combining the two estimates yield the announced result with $\Cee =  (4 \Lambda +1)(1 + \frac1\lambda)$ by remarking that $R^2 -r^2 \ge (R-r)^2$.  
\end{proof}

\subsection{Gain of Integrability}

\begin{prop}[Gain of integrability in $\mathrm{pDG}^+$]\label{p:gain-parab}
   Let $u \in \mathrm{pDG}^+ (I \times B,S)$. Then for all cylinder $Q_R(X_0) \subset I \times B$ and $r \in (0,R)$ and $\kappa \in \R$,
   \[ \| (u-\kappa)_+ \|_{L^{p_c} (Q_r(X_0))}^2 \le \Cc \left( \frac{1}{(R-r)^2} + \frac1{r^2} \right) \| (u-\kappa)_+\|^2_{L^2 (Q_R(X_0))} + \Cc \int_{Q_R (X_0)} |S| (u-\kappa)_+ \dt \dx \]
 for some universal exponent $p_c =2 + \frac4d > 2$ for $d \ge 3$ and any $p_c >2$ for $d=1,2$. In the latter case, the constant $\Cc$ also depends on $p_*$.  
 \end{prop}
 The proposition derives from Sobolev's inequality and the following elementary lemma. 
 \begin{lemma}[Interpolation in Lebesgue spaces] \label{l:interpolation}
   Let $p_i,q_i \in [1,\infty]$ for $i=0,1$, we consider $p,q$ such that 
   \[ \frac1p = \frac\theta{p_0} + \frac{1-\theta}{p_1} \quad \text{ and } \quad \frac1q = \frac\theta{q_0} + \frac{1-\theta}{q_1} \]
   for some $\theta \in (0,1)$. Then
   \[ \| u \|_{L^p (I, L^q (B))} \le \| u \|^\theta_{L^{p_0}(I,L^{q_0} (B))}  \| u \|^{1-\theta}_{L^{p_1}(I,L^{q_1} (B))} .\]
 \end{lemma}
 This lemma is a consequence of H\"older's inequality. 
 \begin{proof}[Proof of Proposition~\ref{p:gain-parab}]
   Let $I^0= (t_0-r^2,t_0]$ and $B^0 = B_r(x_0)$.

   \paragraph{Integrability.} The fact that $u \in \mathrm{pDG}^+ (I \times B),S$ together with Sobolev's inequality implies that $u \in L^\infty (I^0,L^2 (B^0)) \cap L^2(I^0,L^{p_*} (B^0))$.
   Then the previous lemma implies that $u \in L^{p_c} (Q_r (X_0))$ with
   \[ \frac1{p_c} =   \frac{\theta}{\infty} + \frac{1-\theta}2  = \frac{\theta}2 + \frac{1-\theta}{p_*} \]
   where we recall that $p_*$ comes from Sobolev's inequality and satisfies $\frac1{p_*} = \frac12 -\frac1d$ for $d \ge 3$ and $p_* >2$ arbitrary for $d=1,2$.

   For $d \ge 3$, this implies $p_c = \frac2{1-\theta}$ and $\frac{1-\theta}2  = \frac{\theta}2 + \frac{1-\theta}{p_*}$. The second condition is equivalent to,
   \[ \frac12 -\frac1{p_*}= \theta \left( 1 - \frac1{p_*} \right) .\]
   In view of the definition of $p_*$ this is equivalent to $\frac1d = \theta (\frac12 + \frac1d)$. Hence $\theta = \frac{2}{d+2}$ and
   $p_c = \frac2{1-\frac2{d+2}} = 2 + \frac{4}d$.

   In the case where $d=1,2$, we can pick any $p_* >2$. Writing $\frac1{p_*} = \frac12 - \eps$, we get $\theta = \frac{2\eps}{1+2\eps}$ and $p_c = 2 + \frac1{\eps}$. 

 \paragraph{Estimate.} We now write the estimate that we obtain after using the interpolation lemma.  Recall that $I^0= (t_0-r^2,t_0]$ and $B^0 = B_r(x_0)$.
 Use the elementary inequality $a^\theta b^{1-\theta} \le \theta a+(1-\theta)b$ for all $a,b>0$ and Sobolev's inequality from Proposition~\ref{p:sobolev} in order to write,
 \begin{align*}
   &\| (u-\kappa)_+ \|_{L^{p_c} (Q_r(X_0))}^2 \\
   \le &  \| (u-\kappa)_+ \|_{L^{\infty}(I^0, L^2 (B^0))}^{2 \theta} \| (u-\kappa)_+ \|_{L^2 (I^0, L^{p_*} (B^0))}^{2 (1-\theta) }   \\
    \le &   \| (u-\kappa)_+ \|_{L^{\infty}(I^0, L^2 (B^0))}^2 +  \| (u-\kappa)_+ \|_{L^2 (I^0, L^{p_*} (B^0))}^2   \\
    \le &  \| (u-\kappa)_+ \|_{L^{\infty}(I^0, L^2 (B^0))}^2 +  \Csob \| \nabla_x  (u-\kappa)_+ \|_{L^2 (Q_r (X_0))}^2  + \Csob r^{-2} \| (u-\kappa)_+\|^2_{L^2 (Q_r (X_0))}  \\
    \le& \Cc \left( \frac{1}{(R-r)^2} + \frac{1}{r^2} \right) \| (u-\kappa)_+\|^2_{L^2 (Q_R(X_0))} + \Cc \int_{Q_R (X_0)} |S| (u-\kappa)_+ \dt \dx 
 \end{align*}
 with the universal constant $\Cc = \max (1,\Csob) \Cee$.
\end{proof}

\subsection{Local maximum principle}

\begin{prop}[Local maximum principle] \label{p:lmp-parab}
Let $Q \subset \R \times \R^d$ be a parabolic cylinder and $S \in L^\infty( Q)$. There exists a universal constant $\Clmp \ge 1$ such that for all $u \in \mathrm{pDG}^+(Q,S)$ and $Q_R(X_0) \subset Q$ and $r< R$, 
   \[ \| u_+ \|_{L^\infty (Q_r(X_0))} \le \Clmp \left( \left(1 + \frac1{r^2} + \frac{1}{(R-r)^2} \right)^{\omega_0} \|u_+ \|_{L^2 (Q_R(X_0))} + \| S\|_{L^\infty (Q_R(X_0))} \right) \]
  with $\omega_0=\frac{(d+4)(d+2)}8$.
 \end{prop}
 \begin{proof}
We proceed in several steps. 

\paragraph{Reduction.}
  We are are by now used to start by reducing to the case where the right hand side is universal. We can also assume $X_0=0$. 
  More precisely, we reduce to the case where $\|u_+ \|_{L^2 (Q_R)} \le \delta_0$ and $\| S\|_{L^\infty (Q_R)} \le 1$
  by considering
  \[ \tilde{u}(X) = \frac{u(X_0+X)}{\delta_0^{-1} \|u_+ \|_{L^2 (Q_R(X_0))} + \| S\|_{L^\infty (Q_R(X_0))}},\]
  and we aim at proving that $\tilde{u} \le 2$ a.e. in $Q_r$. The reduction is possible because we have  $\tilde{u} \in \mathrm{pDG}^+ (Q,S)$ with
  $S$ replaced with $\tilde S (X) = S (X_0+X) / (\delta_0^{-1} \|u_+ \|_{L^2 (Q_R(X_0))} + \| S\|_{L^\infty (Q_R(X_0))})$. 

\paragraph{Iterative truncation.}
We now assume that $\|u_+ \|_{L^2 (Q_R)} \le \delta_0$ and $\| S\|_{L^\infty (Q_R)} \le 1$ and we aim at finding a constant  $\delta_0>0$,
only depending on $d,\lambda,\Lambda, r,R$, such that $u \le 2$ a.e. in $Q_r$. In order to do so, we consider shrinking cylinders $Q^k = Q_{r_k}$ with $r_k = r + (R-r) 2^{-k}$ and 
\[ A_k = \int_{Q^k} (u -\kappa_k)^2_+ \dt \dx  \]
with $\kappa_k = 2 -2^{-k}$. We still have
\[ A_0 = \int_{Q_R} (u-1)_+^2 \dt \dx \le \|u_+\|_{L^2 (Q_R)}^2 \le \delta_0^2 \]
and we prove that there exist $C>1$ and $\beta >1$ such that for all $k \ge 0$, $A_{k+1} \le C^{k+1} A_k^\beta$ with $C$ depending on $d,\lambda,\Lambda, r,R$.
Then the iteration lemma~\ref{l:induc} from the previous chapter implies that $A_k \to 0$ as soon as
$\delta_0^2 = \frac12 C^{-\frac{\beta}{(\beta-1)^2}}$. We thus are left with finding two  constants $C>1$ and $\beta >1$
such that the iterative inequality $A_{k+1} \le C^{k+1} A_k^\beta$ holds true. The constant $\beta$ is universal while $C$ also depends on $r$ and $R$, and more precisely,
on $(r^{-2} + (R-r)^{-2})$. 

  \paragraph{Local energy estimates.}
   We use the local energy estimates defining the parabolic De Giorgi's class $\mathrm{pDG}^+ (B_1,S)$
   with $R=r_k$ and $r= r_{k+1}$ and $\kappa = \kappa_k$ and $X_0 = (t_0,x_0) = 0$. In particular, we recall that $r_k-r_{k+1} = (R-r) 2^{-k-1}$.
   Recalling that $S$ is bounded by $1$, and write $Q^{k+1} = I^{k+1} \times B^{k+1}$ for an interval and a ball, we can write,
\begin{multline*}
   \underset{t \in I^{k+1}}{\sup} \int_{B^{k+1}}  (u-\kappa_{k+1})_+^2 \dx    +  \int_{Q^{k+1}} |\nabla_x (u-\kappa_{k+1})_+|^2 \dt \dx \\
    \le \Cpdgp \left( 4^{k+1} (R-r)^{-2}  + r^{-2}\right) \int_{Q^k} (u-\kappa_k)_+^2 \dt \dx + \Cpdgp \int_{Q^k}  (u-\kappa)_+ \dt \dx.
  \end{multline*}
  \begin{itemize}
    \item
      \textsc{(Gain of integrability)} 
      We use Sobolev's inequality in $B^{k+1}$ for all time $t \in I^{k+1}$ in order to get,
\begin{align*}
     \|(u-\kappa_{k+1})_+ &\|^2_{L^\infty(I^{k+1}, L^2 (B^{k+1}))} +  \Csob^{-1}  \| (u-\kappa_{k+1})_+ \|^2_{L^2 (I^{k+1}, L^{p_*} (B^{k+1}))}  \\
     \le & \left( \Cpdgp 4^{k+1} (R-r)^{-2} + (\Cpdgp+1)r^{-2} \right) \|(u-\kappa_{k+1})_+\|^2_{L^2 (Q^k)}   \\
     + &\Cpdgp \int_{Q^k}  (u-\kappa_{k+1})_+ \dt \dx.
  \end{align*}
  Then use Proposition~\ref{p:gain-parab} to estimate the left hand side and use Cauchy-Schwarz's inequality
  to estimate the second term in the right hand side and get,
\[
   \|(u-\kappa_{k+1})_+\|^2_{L^{p_c} (Q^{k+1})}  
  \le \Cdgun \left(  4^k ((R-r)^{-2} + r^{-2})  A_k  +  A_k^{\frac12} |\{ u \ge \kappa_{k+1} \} \cap Q^k |^{\frac12} \right)
  \]
for some universal constant $\Cdgun \ge 1$. 
\item \textsc{(Nonlinearization procedure)}
  Applying  Bienaymé-Chebyshev's inequality, we obtain like in the elliptic case -- see Eq.~\eqref{e:dg3} -- that
  \( |\{ u \ge \kappa_{k+1} \} \cap Q^k | \le 4^{k+1}  A_k\)
  so that the previous estimate can be continued as,
  \[   \|(u-\kappa_{k+1})_+\|_{L^{p_c} (Q^{k+1})}^2    \le  2^{2k+1} \Cdgun  ((R-r)^{-2} + r^{-2}+1)  A_k .\]
  \item \textsc{(Nonlinear iteration)}
Like in the elliptic case and after replacing balls with cylinders, we write,
    \begin{align*}
    A_{k+1}  & \le  \|(u - \kappa_{k+1})_+\|_{L^{p_c}(Q^{k+1})}^2 \left|\{ u \ge \kappa_{k+1}\} \cap Q^k\right|^{\frac2q} \\
            & \le 2^{2k+1} 4^{\frac{2(k+1)}q}\Cdgun ((R-r)^{-2} + r^{-2}+1) A_k^{1+\frac2q}  
  \end{align*}
with $q \in (1,2)$ such that $\frac12 = \frac1{p_c} + \frac1q$. Since $p_c= \frac{2(d+2)}d$, we have $q = d+2$ and $\beta = 1 + \frac2q = 1 + \frac2{d+2}$. 
  \end{itemize}
We conclude that we have $A_{k+1} \le C^{k+1} A_k^\beta$ for $C = \bar C (1+ (R-r)^{-2} + r^{-2})$ and $\bar C>1$ and $\beta >1$, as desired. 

Since
\[\delta_0^2 = \frac12 C^{-\frac{\beta}{(\beta-1)^2}}= \frac12 \bar C^{-\frac{\beta}{(\beta-1)^2}} (1+ (R-r)^{-2} + r^{-2})^{-\frac{\beta}{(\beta-1)^2}},\]
the factor $\omega_0$ corresponds to $\frac{\beta}{2(\beta-1)^2} = \frac{d+4}{}$. 
\end{proof}
We next apply the local maximum principle to get some information about lower bounds.
\begin{cor}[Upside down local maximum principle] \label{c:lower-gen}
  There exist two universal constants $\eps_0, \eps_1 \in (0,1/2)$ such that for all $\rho \in (1,2]$, 
  all $-u \in \mathrm{pDG}^+ (Q_\rho,S)$ with $\|S \|_{L^\infty (Q_\rho)} \le \eps_0$ and $u \ge 0$ a.e. in $Q_\rho$, we have 
  \[ |\{ u \ge 1 \} \cap Q_\rho | \ge (1-\eps_1) |Q_\rho| \quad \Rightarrow \quad \left\{ u \ge \frac12 \text{ a.e. in } Q_1 \right\}.\]
\end{cor}
\begin{remark}
It is convenient to pick $\eps_1 < 1/2$ when we will use it in the expansion of positivity. 
\end{remark}
\begin{proof}
  We apply the local maximum principle to $v = 1-u \in \mathrm{pDG}^+ (Q_\rho,S)$ where the source term $S$ is unchanged. It leads to,
  \[ \|v_+ \|_{L^\infty (Q_1)} \le \Clmp \left( 9^{\omega_0} \|v_+ \|_{L^2 (Q_\rho)} + \eps_0 \right) .\]
  We thus pick $\eps_0 >0$ such that $\Clmp \eps_0 \le \frac14$ and we estimate the $L^2$-norm of $v_+$ as follows,
  \begin{align*}
    \|v_+ \|_{L^2 (Q_\rho)} & \le |\{ v > 0 \} \cap Q_\rho |^{\frac12} \quad \text{ since } v \le 1 \text{ a.e.  in } Q_\rho \\
                         & = |\{ u < 1 \} \cap Q_\rho |^{\frac12} \\
    & \le \eps_1^{\frac12} |Q_\rho|^{\frac12}. 
  \end{align*}
  Gathering the estimates and using that $\rho \le 2$, we get,
  \[ \|v_+ \|_{L^\infty (Q_1)} \le \Clmp 9^{\omega_0} \eps_1^{\frac12} |Q_2|^{\frac12}   + \frac14.\]
  We thus pick $\eps_1>0$ such that $\Clmp  9^{\omega_0} \eps_1^{\frac12} |Q_2|^{\frac12} = \frac14$ and we get $v \le \frac12$ a.e. in $Q_1$.
  This is equivalent to $u \ge \frac12$ a.e. in $Q_1$, as desired. 
\end{proof}

Arguing like in the elliptic case, we obtain the following refined local maximum principle. 
\begin{cor}[Local maximum principle - again] \label{c:lmp-parab}
  Given a (universal) constant ${p_0} \in (0,2)$, there exists a (universal) constant  $\eClmp >0$, only depending on $d,\lambda,\Lambda$ and ${p_0}$,
  such that for any $u \in \mathrm{pDG}^+ (Q_1,S)$,  
  \[ \| u_+ \|_{L^\infty (Q_{1/2})} \le  \eClmp \left( \|u_+\|_{L^{p_0} (Q_1)} + \|S\|_{L^\infty (Q_1)} \right) \]
  where $\|u_+\|_{L^{p_0} (Q_1)} = \|u_+^{p_0}\|_{L^1 (Q_1)}^{\frac1{p_0}}$. 
\end{cor}

\section{The parabolic De  Giorgi's class pDG\textsuperscript{--}  \&  expansion of positivity}

In this section, we introduce the parabolic De Giorgi's class $\mathrm{pDG}^-$ that contains in particular
all super-solutions of the parabolic equations at stake in this chapter. We then establish that the non-negative elements of this class
enjoy the following property: if they are above $1$ in a set of universal positive measure in the past, they are bounded from
below point-wise in the future. We will follow classical references from the literature by referring to this property as the
\emph{expansion of positivity}.

We first define the parabolic De Giorgi's class that we will work with. 
\begin{defi}[The parabolic De Giorgi's class $\mathrm{pDG}^-$] \label{d:pDG-}
  Let  $I = (a,b]$ with $a,b \in \R$ and $B$ be an open ball of $\R^d$ and $S \in L^2 (I \times B)$.
  A function $u \colon I \times B \to \R$ lies in the \emph{parabolic De Giorgi's class} $\mathrm{pDG}^-(I \times B,S)$
  if $-u \in \mathrm{pDG}^+(I \times B,S)$  and if, for all $X_0 = (t_0,x_0) \in I \times B$ and all $r,R>0$ such that
  $r< R$ and $Q_R (X_0) \subset I \times B$, all $\kappa \in \R$, and all $t_1,t_2 \in (t_0-R^2,t_0]$ with $t_1 <t_2$, we have,
  \begin{multline} \label{e:L2-propa}
    \int_{B_r(x_0)} (u-\kappa)_-^2 (t_2) \dx     \le \int_{B_r(x_0)} (u-\kappa)_-^2 (t_1) \dx \\
    + \frac{\Cpdgm}{(R-r)^2} \int_{t_1}^{t_2} \int_{B_R (x_0)} (u-\kappa)_-^2 \dt \dx + \Cpdgm \int_{t_1}^{t_2} \int_{B_R (x_0)} |S| (u-\kappa)_- \dt \dx
  \end{multline}
  for some universal constant $\Cpdgm \ge 1$.
\end{defi}
Again, we emphasize that functions in De Giorgi's class can be translated and scaled. 
\begin{lemma}[Invariance by scaling and translation of the De Giorgi's class]\label{l:invariance-pdgm}
  If $u \in \mathrm{pDG}^-(I \times B, S)$ and $Q_r (X_0) \subset I \times B$, then the function $v = \lambda u (\frac{t-t_0}{r^2}, \frac{x-x_0}r)$
  lies in $\mathrm{pDG}^-(Q_1, \mathfrak{S})$  with $\mathfrak{S} (t,x)= \frac{\lambda}{r^2} S (\frac{t-t_0}{r^2} , \frac{x-x_0}r )$. 
\end{lemma}
This simple proof is also left to the reader.

We remark that the local energy estimates that we derived in Proposition~\ref{p:lee} implies that weak solutions are in the parabolic De Giorgi's class pDG\textsuperscript{-}. 
  \begin{prop}[Weak solutions are in $\mathrm{pDG}^-$] \label{p:weak-pdgm}
  Let $u \colon I \times B \to \R$ be a weak solution of $\partial_t u = \dive_x ( A \nabla_x u) + S$ in $I \times B$ with $S \in L^2 (I \times B)$.
  Then $u \in \mathrm{pDG}^-(I \times B,S)$.
\end{prop}
\begin{remark}[Super-solutions are in the parabolic De Giorgi's class]
  The parabolic De Giorgi's class $\mathrm{pDG}^-$ contains in particular all super-solutions of parabolic equations mentioned in Remark~\ref{r:sub-super}.
\end{remark}

\begin{figure}[h]
 \centering{\includegraphics[height=5cm]{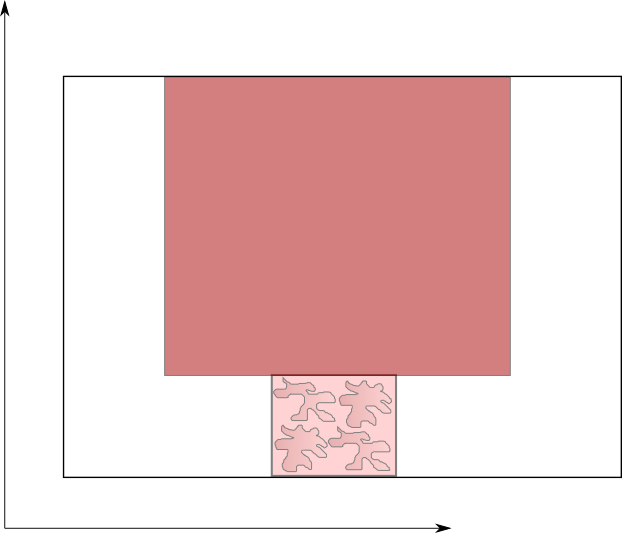}}
 \put(-180,130){\scriptsize $t$}
 \put(-147,111){$\Qexp$}
 \put(-88,90){\huge $Q_1$}
 \put(-88,27){$\Qgu$}
 \put(-65,-10){\scriptsize $(x,v)$}
 \caption{\textit{Geometric setting of the expansion of positivity.}}
   \label{fig:pop-parab}
 \end{figure}

We  now establish the key property of non-negative functions from De Giorgi's class pDG\textsuperscript{-}: they are bounded from below by a universal constant in $Q_1$ 
as soon as  they lie above $1$ in an arbitrary small but universal proportion of a cylinder $\Qgu$ lying in the past, see Figure~\ref{fig:pop-parab}.
\begin{prop}[Expansion of positivity] \label{p:expansion-parab}
  Let $\eta \in (0,1)$ and $\Qext = Q_{\sqrt{1+\eta^2}}$ and $\Qgu = Q_\eta (-1,0,0)$.
  There exist  constants  $\ell_0,\bar \eps_0 \in (0,1)$, depending on $d,\lambda,\Lambda$ and $\eta$, such that for all $u \in \mathrm{pDG}^- (\Qext,S)$ with $S \in L^\infty(\Qext)$
  and $\|S\|_{L^\infty(\Qext)} \le \bar \eps_0$ and $u \ge 0$ a.e. in $\Qext$,
  \[ |\{ u \ge 1 \} \cap  \Qgu | \ge \frac12 |\Qgu| \quad \Rightarrow \quad \{ u \ge \ell_0 \text{ a.e. in } Q_1 \}.\]
\end{prop}

The proof of the previous proposition is split into several lemmas.
We start with proving that the assumption on the measure level-set $\{ u \ge 1\}$ in time and space can be made point-wise in time for some time $t_*$ ``in the past''.
Because we will need to estimate the gradient later on by the local energy inequalities,
we make that that $t_*$ is not too close from the initial time $-1-\eta^2$.
\medskip

Let us consider a radius slightly larger than $1$ in which we aim at getting a lower bound on $u$ ``in a very large proportion of the corresponding ball''.
It is convenient to consider some parameter
\[1< r_\eta < \sqrt{1+\eta^2}\]
for instance $r_\eta =  (1+ \sqrt{1+\eta^2/2})$. 

We first remark that the assumption of Proposition~\ref{p:expansion-parab} implies that
\begin{equation}
    \label{e:inflating}
    |\{ u \ge 1 \} \cap (-1-\eta^2,-1] \times \Breita| \ge \iota_\eta | (-1-\eta^2,-1] \times \Breita|
  \end{equation}
  with $\iota_\eta = \frac{\eta^d}{2 r_\eta^d}$. Indeed, $ \frac12 |\Qgu| = \iota_\eta \eta^2 |\Breita|$. 

  We will first apply the next lemma with $\alpha =\iota_\eta$ but we will then apply it iteratively with a larger $\alpha$.

\begin{remark}[Beyond technicalities]
  The statements of the next three lemmas are a bit technical and we can try to see a bit where we are heading to. 
  
  The first lemma says that if $|\{ u \ge 1\} \cap \Breita|$ fills a  proportion $\iota_\eta$  of a time interval $I$, we pick a time $t_*$ in a sub-interval $I_* \subset I$.
  We do not want $t_*$ to bee too close to beginning of the interval $I$ to make sure that, when applied iteratively, the
  resulting sequence of times increase. But we also want to make sure that $t_*$ is not to close of the end of the interval $I$,
  in order to handle what happens when  the iteration ends. 
\end{remark}

\begin{lemma}[From measure to pointwise in time] \label{l:m-2-p}
  Let $I = (t_I,t_I+\tau_I)$ and  $\tau_{\eta,I} = \frac{2(1-\iota_\eta)}{2-\iota_\eta} \tau_I < \tau_I$ and  $I_*=(t_I + \frac13 \tau_{\eta,I},t_I  + \frac23  \tau_{\eta,I}] \subset I$.

  Suppose that a measurable function $u \colon I \times \Breita \to \R$ satisfies $|\{ u \ge 1 \} \cap I \times \Breita | \ge \iota_\eta |I \times \Breita|$.
Then there exists a time $t_* \in I_*$ with  such that  \( |\{ u (t_*, \cdot)  \ge 1 \} \cap \Breita | \ge \frac{\iota_\eta}2 |\Breita|.\)
\end{lemma}
\begin{remark}
  \label{r:delta-eta}
  When  the time interval is $I = (-1-\eta^2,-1]$, we have $t_* \in ( -1 -\eta^2 + \delta_\eta, -1-\eta^2+ 2 \delta_\eta)$ with $\delta_\eta = \frac{2(1-\iota_\eta)}{3(2-\iota_\eta)} \eta^2$.
  We record  that $ \delta_\eta < \eta^2/2$.
\end{remark}
\begin{proof}
  The assumption is equivalent to $|\{ u <  1 \}\cap I \times \Breita| \le (1-\iota_\eta) |I \times \Breita|$, that we can write,
\[ \int_I |\{ u(t,\cdot) < 1 \} \cap \Breita| \dt \le (1-\iota_\eta) \tau_I |\Breita|.\]
We now a pick a time interval $I_*$ whose length is $2/3$ of $\frac{(1-\iota_\eta)}{1-\iota_\eta/2} \tau_I < \tau_I$,
\[ \fint_{t_I + \frac{2(1-\iota_\eta)}{3(2-\iota_\eta)} \tau_I}^{t_I  + \frac{4(1-\iota_\eta)}{3(2-\iota_\eta)}\tau_I} |\{ u(t,\cdot) < 1 \} \cap \Breita| \dt \le (1-\frac{\iota_\eta}2) |\Breita|.\]
We conclude that there exists  $t_* \in I_*$ such that   \( |\{ u (t_*, \cdot)  < 1 \} \cap \Breita | \le (1-\iota_\eta/2) |\Breita|.\)
\end{proof}
  The next step is to propagate during a time $\tau_\eta$ the point-wise-in-time bound we got from the previous lemma.
  The price to pay for propagation is that we  deteriorate the lower bound from $1$ to $m_\eta$ and the fraction of the ball occupied
  by the level set changes from $\iota_\eta$ to $\iota_\eta/2$.
\begin{remark}[Beyond technicalities, again]
  Following closely the dependence of $\tau_\eta$ and $m_\eta$ with
  respect to all parameters is important to us because we want to make sure that neither the time interval $\tau_\eta$ nor the lower bound $m_\eta$  become too small
  as we iterate the use of these three lemmas. 
\end{remark}

\begin{lemma}[Short time propagation] \label{l:st-propagation}
  Let $\eta \in (0,1)$. There exist $\tau_\eta  \in (0,1)$ (only depending on $d,\lambda,\Lambda,\eta$) and $m_\eta  \in (0,1)$ (only depending on $\iota_\eta$) such that 
  for all  $u \in \mathrm{pDG}^- (\Qext,S)$ with $\|S\|_{L^\infty(\Qext)} \le 1$, if
  \[ \exists t_* \in (-1-\eta^2,0], \qquad |\{ u (t_*, \cdot)  \ge 1 \} \cap \Breita | \ge \frac{\iota_\eta}2 |\Breita|,\]
  then
  \[ \forall t \in I_{\tau_\eta} , \qquad |\{ u (t, \cdot)  \ge m_\eta \} \cap \Breita | \ge \frac{\iota_\eta}4 |\Breita|\]
  where $I_{\tau_\eta} :=(t_*, t_*+ \tau_\eta] \cap (-1-\eta^2,0]$. 
\end{lemma}
\begin{proof}
    On the one hand, we know from the definition of $\mathrm{pDG}^-(\Qext,S)$ -- see \eqref{e:L2-propa} with $\kappa = 1$, $x_0=0$, $\rho \in (0,1)$ -- that for all $t \in (t_*,t_* + \tau] \cap (-1-\eta^2,0]$, we have,
\begin{align*}
  \int_{B_{\rho r_\eta}} (u-1)_-^2 (t) \dx
  \le &\int_{B_{\rho r_\eta}} (u-1)_-^2 (t_*) \dx  + \frac{\Cpdgm}{(1-r)^2} \int_{I_{\tau_\eta}} \int_{\Breita} (u-1)_-^2 \dt \dx \\
  & + \Cpdgm \int_{I_{\tau_\eta}} \int_{\Breita} |S| (u-1)_- \dt \dx \\
  \le &(1-\frac{\iota_\eta}{2}) |\Breita|  + \frac{\Cpdgm}{(1-r)^2} \tau_\eta|\Breita| + \Cpdgm \tau_\eta|\Breita|.
\end{align*}

On the other hand, we have  for all $m \in (0,1)$,
\begin{align*}
  \int_{B_{\rho r_\eta}} (u-1)_-^2 (t) \dx & \ge \int_{\{ u < m\} \cap B_{\rho r_\eta}} (u-1)_-^2 (t) \dx \\
                               &\ge (1-m)^2 |\{ u < m\} \cap B_{\rho r_\eta}| \\
  & \ge (1-m)^2 |\{ u < m\} \cap \Breita| - (1-m)^2 (1-\rho^d)|\Breita|  .
\end{align*}

Combining the resulting inequalities leads for all $\rho, m \in (0,1)$ to,
\[ |\{ u(t,\cdot) < m\} \cap \Breita| \le \left\{ (1-m)^{-2} \left( (1-\frac{\iota_\eta}{2})   + \frac{\Cpdgm}{(1-\rho)^2} \tau_\eta + \Cpdgm \tau_\eta \right) + (1-\rho^d)  \right\}  |\Breita|.\]
We first pick $\rho=\rho_\eta \in (0,1)$ such that $(1-\rho^d) |\Breita| \le \frac{\iota_\eta}{16} |\Breita|$,
\[ |\{ u(t,\cdot) < m\} \cap \Breita| \le  (1-m)^{-2} \left( (1-\frac{\iota_\eta}{2})   + \frac{\Cpdgm}{(1-\rho_\eta)^2} \tau_\eta + \Cpdgm \tau_\eta \right)  |\Breita| + \frac{\iota_\eta}{16} |\Breita|.\]
We then pick $m =m_\eta \in (0,1)$ such that \( (1-m)^{-2} (1-\frac{\iota_\eta}{2}) |\Breita| \le (1 - \frac7{16} \iota_\eta) |\Breita|\)
\[ |\{ u(t,\cdot) < m_\eta \} \cap \Breita| \le  (1-\frac7{16} \iota_\eta) |\Breita| + (1-m_\eta)^{-2} \left( \frac{\Cpdgm}{(1-\rho_\eta)^2} \tau_\eta + \Cpdgm \tau_\eta \right)  |\Breita| + \frac{\iota_\eta}{16} |\Breita|.\]
We finally impose to $\tau_\eta$ to satisfy,
\[ (1-m_\eta)^{-2} \left( \frac{\Cpdgm}{(1-\rho_\eta)^2} \tau_\eta  + \Cpdgm \tau_\eta  \right) |\Breita| \le \frac{\iota_\eta}8 |\Breita|.\]
This leads to
\[ |\{ u(t,\cdot) < m_\eta \} \cap \Breita| \le (1-\frac7{16} \iota_\eta) |\Breita| + \frac{\iota_\eta}8 |\Breita| + \frac{\iota_\eta}{16}|\Breita| = (1-\frac{\iota_\eta}4)|\Breita| \]
as desired. 
\end{proof}
We now want to get a lower bound on a very large proportion of the unit ball in order to be able to apply the upside down local maximum principle  (Corollary~\ref{c:lower-gen}).
While we improve as much as we want the proportion of the truncated cylinder $I_{\tau_\eta} \times \Breita$ where we have a lower bound, we pay it by further deteriorating the lower bound. 

\begin{lemma}[Spreading] \label{l:spreading}
  Let $\eta \in (0,1)$ and $\tau_\eta ,m_\eta \in (0,1)$  be given by Lemma~\ref{l:st-propagation} (short time propagation).

    Suppose that  $u \in \mathrm{pDG}^- (\Qext,S)$ and that there exists $t_* \in (-1-\eta^2+ \delta_\eta,0]$ such that
  \[ \forall t \in I_{\tau_\eta} , \quad |\{ u (t_*, \cdot)  \ge 1 \} \cap \Breita | \ge \frac{\iota_\eta}4 |\Breita|\]
  with $I_{\tau_\eta} :=(t_*,t_*+\tau_\eta] \cap (-1-\eta^2,0]$.
  
  For all $\alpha \in (0,1)$,  there exists $\bar m \in (0,1)$ (only depending on $d,\lambda,\Lambda$ and $\eta,\alpha$ and on any lower bound on $|I_{\tau_\eta}|$) such that
  if $\|S\|_{L^\infty(Q_2)} \le \bar m$, then
  \[ |\{ u \ge \bar m \} \cap I_{\tau_\eta}  \times \Breita | \ge \alpha |I_{\tau_\eta} \times \Breita|.\]
\end{lemma}
\begin{proof}
  Let $N \in \mathbb{N}$ with $N \ge 2$, to be chosen.

\paragraph{Sequence of truncated functions.}
  For all $k \in \{0,\dots, N \}$, we consider
  \[ \kappa_k = m_\eta^k \quad \text{ and } \quad u_k = (\kappa_k- \max (u(t,\cdot), \kappa_{k+1}))_+=
    \begin{cases}
      \kappa_k - \kappa_{k+1} & \text{ if } u (t,\cdot) \le \kappa_{k+1} , \\
      \kappa_k - u & \text{ if } \kappa_{k+1} < u(t,\cdot) < \kappa_k , \\
      0 & \text{ if } u(t,\cdot) \ge \kappa_k .
    \end{cases}
  \]
  By assumption and the fact that $\kappa_k \le 1$, we know that for all $t \in I_{\tau_\eta} = (t_*,t_*+ \tau_\eta] \cap (-1-\eta^2,0]$ and all $k \ge 1$, 
  \begin{equation} \label{e:lwb}
    |\{u (t, \cdot) \ge \kappa_k \} \cap \Breita | \ge \frac{\iota_\eta}4 |\Breita|.
  \end{equation}

\paragraph{Estimate for the truncated gradient.}
We use next the fact that $-u \in \mathrm{pDG}^+ (\Qext,S)$  with $\|S\|_{L^\infty(\Qext)} \le m_\eta^N \le \kappa_k$
and an intermediate cylinder ù$\Qmid$ such that
  \[ I_{\tau_\eta} \times \Breita \subset  (-1-\eta^2 + \delta_\eta , 0] \times \Breita \subset \Qmid \subset \Qext = (-1-\eta^2,0] \times B_{\sqrt{1+\eta^2}}\]
to estimate $\|\nabla_x (u-\kappa_k)_- \|_{L^2 (I_{\tau_\eta} \times \Breita)}$ as follows,
  \begin{align} 
    \nonumber   \|\nabla_x (u-\kappa_k)_- \|_{L^2 (I_{\tau_\eta} \times \Breita)}^2 & \le \Cpdgp \left( \left( \frac1{\delta_\eta} + \frac{1}{(r_\eta -1)^2} \right)
                                                                            \int_{\Qext} (u-\kappa_k)_-^2 + \int_{\Qext} |S|(u-\kappa_k)_- \right) \\
    & \le C_\eta \kappa_k^2 
  \label{e:eek}
  \end{align}
  for some constant $C_\eta$ depending on $d,\lambda,\Lambda$ and $\eta$.

  \paragraph{Intermediate values.}
  We note that $u_k = T_k(u)$ where $T_k$ is a Lipschitz function. In particular, we can apply Poincaré-Wirtinger's inequality from Proposition~\ref{p:poincare} to $u_k$ and
  argue like in the proof of the elliptic intermediate value lemma~\ref{l:ivl-elliptic}. Recalling that $\bar u_k = |\Breita|^{-1} \int_{\Breita} u_k \dx$, we have
  \begin{align*}
    \int_{\Breita} |u_k -\bar u_k | \dx & \le \Cpw \int_{\Breita} |\nabla_x u_k | \dx \\
                                    & = \Cpw \int_{\Breita} |\nabla_x (u-\kappa_k)_- | \un_{\{\kappa_{k+1} < u < \kappa_k\}} \dx.
  \end{align*}
  Integrating with respect to $t \in I_{\tau_\eta} = (t_*,t_*+ \tau_\eta ] \cap (-1,0]$, we get
  \begin{align*}
    \int_{I_{\tau_\eta} \times \Breita} |u_k -\bar u_k | \dt \dx & \le \Cpw \int_{I_{\tau_\eta} \times \Breita} |\nabla_x (u-\kappa_k)_- | \un_{\{\kappa_{k+1} < u < \kappa_k\}} \dt \dx \\
    & \le \Cpw \|\nabla_x (u-\kappa_k)_- \|_{L^2 (I_{\tau_\eta} \times \Breita)} \left| \{ \kappa_{k+1} \le u < \kappa_k\} \cap I_{\tau_\eta} \times \Breita \right|^{\frac12} \\
                                    & \le \Cpw C_\eta ^{\frac12} \kappa_k \left| \{ \kappa_{k+1} \le u < \kappa_k\} \cap I_{\tau_\eta} \times \Breita\right|^{\frac12}
  \end{align*}
where we used \eqref{e:eek} to get the last line. 
  
  We now get a lower bound on the left hand side of the first inequality. Keeping in mind that $u_k \ge 0$ (and consequently $\bar u_k \ge 0$), 
  \begin{align*}
    \int_{\Breita} |u_k -\bar u_k| \dx &\ge \int_{ \{ u_k = 0 \} \cap \Breita} (\bar u_k) \dx \\
    & = (\bar u_k) |\{u_k =0\} \cap \Breita | \\
                                       &= \frac1{|\Breita|} \left( \int_{\Breita} u_k \dx \right)|\{u_k =0\} \cap \Breita | \\
                                   & \ge  \frac{\kappa_k -\kappa_{k+1}}{|\Breita|}  |\{u(t,\cdot) < \kappa_{k+1}\} \cap \Breita| \cdot |\{u (t,\cdot) \ge  \kappa_k \} \cap \Breita | \\
    & \ge (1-m_\eta) \kappa_k \frac{\iota_\eta}4 |\{u(t,\cdot) < \kappa_{k+1}\} \cap \Breita|
  \end{align*}
  where we used \eqref{e:lwb} to get the last line. Integrating with respect to $t \in I_{\tau_\eta}$, we get
 \[ \int_{I_{\tau_\eta} \times \Breita} |u_k -\bar u_k| \dt \dx \ge (1-m_\eta) \kappa_k \frac{\iota_\eta}4 |\{u < \kappa_{k+1}\} \cap I_{\tau_\eta} \times \Breita|.\]
  
  We thus deduce that we have
  \[
    (1-m_\eta) \kappa_k \frac{\iota_\eta}4 |\{u(t,\cdot) < \kappa_{k+1}\} \cap I_{\tau_\eta} \times \Breita| 
    \le \Cpw C_\eta^{\frac12} \kappa_k \left| \{ \kappa_{k+1} \le u < \kappa_k\} \cap I_{\tau_\eta} \times \Breita\right|^{\frac12}
  \]
  and we write this under the following form,
  \[ |\{u < \kappa_{k+1}\} \cap I_{\tau_\eta} \times \Breita| \le \Cexp  \left| \{ \kappa_{k+1} \le u < \kappa_k\} \cap I_{\tau_\eta} \times \Breita\right|^{\frac12}\]
  for some   constant $\Cexp$ depending on $\eta,\eta$ and $d,\lambda,\Lambda$ (dimension and ellipticity constants).

  \paragraph{Iteration.}
  By considering $A_k = |\{u < \kappa_{k+1}\} \cap I_{\tau_\eta} \times \Breita|$, we can write the previous inequality,
  \[ A_{k+1}^2 \le \Cexp^2  [A_k - A_{k+1}]. \]
  We remark that $A_k \ge A_{N}$ and $A_1 \le \tau_\eta |\Breita|$.   We then can sum the previous inequalities over $k \in \{1,\dots, N -1\}$ and get,
  \[ (N-1) A_{N}^2 \le \sum_{k=1}^{N-1} A_{k+1}^2  \le \Cexp^2 A_1 \le \Cexp^2 \tau_\eta |\Breita|.\]
  In particular, we can pick $N$ (only depending on $d,\lambda,\Lambda$ and $\eta,\alpha$ and any lower bound on $|I_{\tau_\eta}|$) such that
  \[ |\{ u < m_\eta^{N} \} \cap I_{\tau_\eta}  \times \Breita | = A_{N} \le (1-\alpha) |I_{\tau_\eta} \times \Breita|. \]
  The precise condition is $N \ge 1+ \frac{\Cexp^2 \tau_\eta}{(1-\alpha)^2 |I_{\tau_\eta}|^2}$. Then we pick $\bar m = m_\eta^N$. 
\end{proof}

The proof of expansion of positivity proceeds basically in two steps.
We first extend easily the initial information we have from the small ball $B_\eta$ to the ``large'' ball $\Breita$.
We then use the three previous lemmas:
\begin{enumerate}[label=(\roman*)]
  \item  to get a time for which we have a lower bound  on the  measure on a super-level set of $f$;
  \item to propagate this lower bound for some short time $\tau_\eta$;
    \item to spread the information and restore the proportion of the truncated cylinder $I \times \Breita$ on which  $f$ is bounded from below.
\end{enumerate}
\begin{proof}[Proof of Proposition~\ref{p:expansion-parab} (expansion of positivity)]
  Let $\eps_1$ be given by the upside down maximum principle from Corollary~\ref{c:lower-gen}. 
  \medskip
  
  \noindent \textsc{Initialization of the iteration.}
  We observed earlier that the assumption of the proposition implies that \eqref{e:inflating} holds true.
  This corresponds to the assumption of Lemma~\ref{l:m-2-p}. The lemma implies that
  there exists $t_* \in (-1-\eta^2 + \delta_\eta, -1]$ (see Remark~\ref{r:delta-eta}) such that
  \[ |\{ u (t_*, \cdot)  \ge 1 \} \cap \Breita | \ge \frac{\iota_\eta}2 |\Breita|.\]

  Now Lemma~\ref{l:st-propagation} implies that 
  \[ \forall t \in I_{\tau_\eta} , \quad |\{ u (t, \cdot)  \ge m_{\iota_\eta} \} \cap \Breita | \ge \frac{\iota_\eta}4 |\Breita|\]
  with $I_{\tau_\eta} = (t_*, t_*+ \tau_{\iota_\eta}) \cap (-1-\eta^2,0]$. The parameter $\tau_{\iota_\eta}$ only depends on $d,\lambda,\Lambda, \eta$.
  
  We next apply Lemma~\ref{l:spreading} with $\alpha = \iota_\eta$,   we get $\bar m_{\iota_\eta}$ such that
  \[ |\{ u \ge \bar m_{\iota_\eta} \} \cap I_{\tau_\eta}  \times \Breita | \ge \iota_\eta |I_{\tau_\eta} \times \Breita|.\]
  The parameter $\bar m_{\iota_\eta}$ only depends on $d,\lambda,\Lambda,\eta$ and a lower bound of
  the length interval $|I_{\tau_\eta}|$.
  \medskip
  
  \noindent \textsc{Time iteration.} In the case where $t_* + \tau_{\iota_\eta} < 0$, we can apply again Lemmas~\ref{l:m-2-p}, \ref{l:st-propagation} and \ref{l:spreading}:
  we construct of sequence of times $t_k$ such that $t_{k+1} \in (t_k + \frac13 \tau, t_k + \frac23 \tau )$ where
  \[ \tau = \frac{2(1-\iota_\eta)}{(2-\iota_\eta)} \tau_\eta \]
  only depends on $d,\lambda,\Lambda,\eta$   (we recall that $\tau_\eta$ comes from Lemma~\ref{l:st-propagation}).
  We  can proceed till the rank $K \ge 1$ for which  $(t_K , t_K + \tau)$ leaks out $(-1-\eta^2,0]$. This translates into,
  \begin{align*}
\text{ if } K =1, \quad &    \begin{cases}
    t_K \in (-1-\eta^2/2, 0 ), \\
    t_K + \tau >0, \\ 
    \end{cases}  \\[1ex]
\text{ if } K \ge 2, \quad &     \begin{cases}
    t_K \in (t_{K-1} + \frac13 \tau, t_{K-1} + \frac23 \tau), \\
    t_{K-1} + \tau \le 0,\\
      t_K + \tau > 0. &
    \end{cases}
  \end{align*}
  Such a maximal rank $K$ exists because $t_{k+1} \ge t_k + \frac13 \tau$ for all $k \ge 2$. 
  At the last rank $K \ge 2$, we have $t_K \in (-\tau,-\frac13 \tau)$  with $\tau$ only depending on $d,\lambda,\Lambda,\eta$. At this final rank, either $K=1$ or $K \ge 2$,
  the time interval $I_{\tau_\eta}$ is given by $(t_K, 0)$, whose length is bounded
  from below by a  constant  only depending on $d,\lambda,\Lambda,\eta$.
\medskip
  
  \noindent \textsc{Collecting the information.}
  By applying iteratively Lemma~\ref{l:st-propagation}, we proved that there exists a  lower bound $\ell_1 >0$  only depending on $d,\lambda,\Lambda,\eta$ and time intervals $I_k=(t_k,t_k+ \tau_k)$ such that
  \begin{equation}
    \label{e:collection}
    \forall k \in \{1,\dots,K\}, \forall t \in I_k, \quad |\{ u \ge \ell_1 \} \cap I_k \times \Breita | \ge \frac{\iota_\eta}2 |\Breita|.
  \end{equation}
  Moreover, $I_k \cap I_{k+1} \neq \emptyset$ and $\cap_{k \ge 1} I_k \cap (-1-\eta^2,0] = (t_*,0] \supset (-1 - \eta^2/2,0]$ (see Remark~\ref{r:delta-eta}). 
  We conclude that,
  \[ \forall t \in (-1 - \eta^2/2,0], \quad |\{ u (t) \ge \ell_1 \} \cap \Breita | \ge \frac{\iota_\eta}2 |\Breita|.\]
  We now see that it is convenient to pick $r_\eta = \sqrt{1+\eta^2/2}$.
  Applying one last time Lemma~\ref{l:spreading} leads to,
  \[ |\{ u \ge 2\ell_0 \} \cap Q_{r_\eta} | \ge (1-\eps_1) |Q_{r_\eta}|.\]
  The previous reasoning works by choosing $\bar \eps_0 = 2 \ell_0$ for the upper bound for source terms.

  Now the upside down maximum principle from Corollary~\ref{c:lower-gen} yields the result. 
 \end{proof}

\section{Improvement of oscillation}

We defined two De Giorgi's classes for parabolic equations: $\mathrm{pDG}^+$ and $\mathrm{pDG}^-$.
The intersection of these two classes contains all weak solutions of the parabolic equations we deal with in this chapter.
\begin{defi}[The parabolic De Giorgi's class $\mathrm{pDG}$] \label{defi:pDG}
  Let  $I = (a,b]$ with $a,b \in \R$ and $B$ be an open ball of $\R^d$ and $S \in L^2 (I \times B)$.
  A function $u \colon I \times B \to \R$ lies in the \emph{parabolic De Giorgi's class} $\mathrm{pDG}(I \times B,S)$
  if it lies both in $\mathrm{pDG}^+(I \times B,S)$ and in $\mathrm{pDG}^-(I \times B,S)$.
\end{defi}
We can now establish that the oscillation of elements of this parabolic class improves when zooming in. 
\begin{prop}[Improvement of oscillation] \label{p:improve-osc-parab}
  Let $\bar \eps_0$  be given by Proposition~\ref{p:expansion-parab} about expansion of positivity. 
  There exists a universal constant  $\mu \in (4^{-1} , 1)$ such that for all $u \in \mathrm{pDG}(Q_2,S)$ with $S \in L^\infty(Q_2)$
  with $\|S\|_{L^\infty(Q_2)} \le \bar \eps_0$ and $u \in L^\infty(Q_2)$,
 \[\osc_{Q_2} u \le 2 \quad \Rightarrow \quad \osc_{Q_1} u \le 2 \mu .\]
\end{prop}
\begin{proof}
 Without loss of generality, we can assume  that $0  \le u \le 2$ a.e. in $Q_2$ by considering $\tilde{u} = u -  \inf{Q_2} u$.
 We still have an essentially bounded source term $S$ and the upper essential bound is not modified. We consider $\Qgu = Q_1 (-1,0,0)$ and distinguish two cases.
  \begin{itemize}
  \item We assume first that $|\{ u \ge 1 \} \cap \Qgu | \ge \frac12 |\Qgu|$. Then expansion of positivity from Proposition~\ref{p:expansion-parab} implies
    that $u \ge \ell_0$ a.e. in $Q_1$. In particular, $\osc_{Q_1} u \le 2 -\ell_0$. 
  \item In the other case, $|\{ u \ge 1 \} \cap \Qgu | < \frac12|\Qgu|$, so that the function $v = 2 -u$ satisfies $|\{ v \le 1 \} \cap \Qgu | < \frac12 |\Qgu|$.
    Equivalently, we have $|\{v > 1 \} \cap \Qgu | > \frac12|\Qgu|$. In particular, $|\{ v \ge 1 \} \cap \Qgu| \ge \frac12 |\Qgu|$. The source term $|S|$ in the
    definition of $\mathrm{pDG}^\pm (\Qgu,S)$ is not changed and we can apply Proposition~\ref{p:expansion-parab} as in the first case and deduce that $v \ge \ell_0$ a.e. in $Q_1$.
    This translates into $u \le 2 -\ell_0$ a.e. in $Q_1$, yielding in this case too that $\osc_{Q_1} u \le 2 -\ell_0$. 
  \end{itemize}
  In both cases, we proved that $\osc_{Q_1} u \le 2 -\ell_0$, reaching the desired conclusion with $\mu = \max (4^{-1}, (1-\ell_0)/2)$.   
\end{proof}
\begin{proof}[Proof of  (De Giorgi \& Nash's) theorem~\ref{t:dg-parab}]
The theorem is a consequence of the local maximum principle (Proposition~\ref{p:lmp-parab}) and of the improvement of oscillation  (Proposition~\ref{p:improve-osc-parab}). 
The local maximum principle applied to $u$ and $-u$ with $X_0=0$ and $r=\frac34$ and $R=1$ implies that
  \begin{equation} \label{e:linfty} \sup_{Q_{\frac34}} |u| \le \Clmp \left( \|u\|_{L^2 (Q_1)} + \|S\|_{L^\infty (Q_1)} \right).
  \end{equation}
  As far as the H\"older semi-norm is concerned, we shall prove that there exists $\alpha \in (0,1]$ and $C_0 \ge 1$ (both universal) such that for all $X_0 \in Q_{\frac12}$, and all $r>0$,
  \[ \osc_{Q_r (X_0) \cap Q_{\frac12}} u \le C_0 \left( \|u\|_{L^\infty (Q_{\frac34})} + \|S\|_{L^\infty (Q_1)} \right) r^\alpha.\]
  This implies that $[u]_{\Cpar^\alpha (Q_{\frac12})} \le C_0 \left( \|u\|_{L^\infty (Q_{\frac34})} + \|S\|_{L^\infty (Q_1)} \right)$ (see Proposition~\ref{p:holder-para}).

  We then consider such a point $X_0 \in Q_{\frac12}$. We infer from \eqref{e:linfty} that $u \in L^\infty (Q_{\frac14} (X_0))$. In order to apply the
  improvement of oscillation (Proposition~\ref{p:improve-osc-parab}), we consider
  \[ \tilde u (X) = \frac{u (X_0+\frac{X}8)}{ \|u\|_{L^\infty (Q_{\frac34})} + \bar \eps_0^{-1} \|S\|_{L^\infty (Q_1)}}.\]
  Then $\| \tilde u \|_{L^\infty (Q_2)} \le 1$  and the source term $\tilde S (X)= \bar \eps_0 8^{-2} \frac{S(X_0 +8^{-1} X)}{\|S\|_{L^\infty (Q_1)}} $ satisfies $\|\tilde S \|_{L^\infty (Q_2)} \le \bar \eps_0$ (see Lemmas~\ref{l:invariance-pdgp} and \ref{l:invariance-pdgm}).
  We thus get from Proposition~\ref{p:improve-osc-parab} that $\osc_{Q_1} \tilde u \le 2 \mu$. We now scale recursively the function $\tilde u$ and consider,
  \[ \forall X \in Q_2, \quad  \tilde{u}_k (X) = \mu^{-k} \tilde{u} (( 1/2)^k X) \]
  whose source term $\tilde S_k (X) = (1 / (4\mu))^k \tilde S ((1/2)^k X)$. We remark that $\|\tilde S_k\|_{L^\infty (Q_1)} \le \| \tilde S\|_{L^\infty (Q_1)} \le \bar \eps_0$ since $\mu \ge 1/4$.
  We conclude that $\osc_{Q_2} \tilde{u}_k \le 2$ for all $k \ge 1$. This translates into,
  \[ \osc_{Q_{r_k}} \tilde{u} \le 2 \mu^k = 2^{1-\alpha} r_k^\alpha \quad \text{ with } \quad r_k = 2 \frac{1}{2^k} \quad \text{ and } \quad (1/2)^\alpha = \mu.\]
  Now for $r \in (0,2]$, there exists $k \ge 0$ such that $r_{k+1} \le r \le r_k$. This implies that
  \[\osc_{Q_r} \tilde{u} \le \osc_{Q_{r_k}} \tilde{u} \le 2^{1-\alpha} r_k^\alpha = 2 r_{k+1}^\alpha \le 2 r^\alpha. \]
  In terms of the function $u$, this implies that for all $r \in (0,2]$,
  \[\osc_{Q_{\frac{r}8} (X_0)} u  \le 8^\alpha 2   \left( \|u\|_{L^\infty (Q_{\frac34})} + \bar \eps_0^{-1} \|S\|_{L^\infty (Q_1)}\right) (r/8)^\alpha. \]
  Since for $s \ge \frac14$, we have
  \[\osc_{Q_s (X_0) \cap Q_{\frac12}} u  \le 2 \|u\|_{L^\infty (Q_{\frac34})} (4s)^\alpha, \]
  we conclude that
  \[[u]_{\Cpar^\alpha (Q_{\frac12})} \le 2^{4\alpha+1} \left( \|u\|_{L^\infty (Q_{\frac34})} + \bar \eps_0^{-1} \|S\|_{L^\infty (Q_1)} \right). \qedhere\]
\end{proof}

\section{(Weak) Harnack's inequality}

In this section, we show that elements of the parabolic De Giorgi's class $\mathrm{kDG}^-$ satisfies a weak Harnack's inequality.
We state it at unit scale.  
\begin{thm}[Weak Harnack's inequality] \label{t:whi-parab}
There exists a (small) universal constant  $\omega \in (0,1)$ and two positive universal constants $\Cwhi$ and ${p_0}$
  such that for $\Qpast = Q_\omega (-1+\omega^2,0,0)$ and $\Qfuture =Q_\omega$, 
  and $u \in \mathrm{pDG}^- (Q_1,S)$ with $S \in L^\infty (Q_1)$ and $f \ge 0$ a.e. in $Q_1$, we have 
  \[ \| u \|_{L^{p_0} (\Qpast)} \le \Cwhi \left ( \inf_{\Qfuture} u + \| S\|_{L^\infty (Q_1)} \right). \]
\end{thm}
\begin{remark}[About the $L^{p_0}$-``norm'' for $p \in (0,1)$]
  We let $\| u \|_{L^{p_0} (\Qpast)}$ denote $\left( \int_{\Qpast} u^{p_0} \right)^{1/p}$ even if $p$ could be smaller than $1$. 
\end{remark}
It is then possible to combine the weak Harnack's inequality with the local maximum principle in order to get Harnack's inequality
for solutions of parabolic equations, and more generally for all elements in the parabolic De Giorgi's class of the domain with essentially bounded source terms.
\begin{thm}[Harnack's inequality] \label{t:harnack-parab}
  There exist two universal constants $R_0>1$ and $\omega \in (0,1)$ and a positive constant $\Charnack$ 
  such that for $\Qharnack = (-1,0] \times B_{R_0} \times B_{R_0}$ and $\Qpast^* = Q_{\omega/2} (-1+\omega^2,0,0)$ and $\Qfuture =Q_\omega$, 
  and $u \in \mathrm{pDG} (\Qharnack,S)$ with $S \in L^\infty (\Qharnack)$ and $u \ge 0$ a.e. in $\Qharnack$, we have 
  \[ \sup_{\Qpast^*} u \le \Charnack \left ( \inf_{\Qfuture} u + \| S\|_{L^\infty (\Qharnack)} \right). \]
\end{thm}
\begin{proof}
  Apply first  Corollary~\ref{c:lmp-parab} between $\Qpast^*$ and $\Qpast$ (from the statement of the weak Harnack's inequality).
  Thanks to Lemma~\ref{l:invariance-pdgp}, it can be scaled and translated in time by considering $\tilde f (t,x) = f(-1+\omega^2 + (t/\omega^2), x/\omega)$.
We can then combine the resulting estimate with the one given by Theorem~\ref{t:whi-parab}.
\end{proof}
\begin{remark}
  The fact that $\Qpast$ is replaced with $\Qpast^*$ in the statement of Harnack's inequality
  is irrelevant since $\inf_{\Qfuture} f \le \inf_{\Qfuture^*} f$ with $\Qfuture^* = Q_{\omega/2}$. 
\end{remark}

\subsection{Generating and propagating a lower bound}

The proof combines the expansion of positivity from Proposition~\ref{p:expansion-parab} with a covering argument.
We aim at  estimating  $\| u \|_{L^{p_0} (\Qpast)}$ by the infimum of $f$ in $\Qfuture$. By linearity, we can reduce to $\inf_{\Qfuture} u \le 1$.
Establishing the estimate amounts to proving that there exists $\eps >0$ (universal) such that for all $t \ge 1$,
\[ |\{ u > t \} \cap \Qpast | \le C t^{-\eps} .\]
A further reduction is to prove that there exists $M>1$ and $\mu \in (0,1)$ such that for all integers $k \ge 1$, we have
\[ |\{ u > M^k \} \cap \Qpast| \le C (1 - \mu)^k .\]
In order to prove this, we consider $U_{k+1} = \{ u > M^{k+1} \} \cap \Qpast$ and we want to prove that $|U_{k+1}| \le (1-\mu) |U_k|$.
In order to prove this inequality, we cover the set $U_{k+1}$ with small cylinders $Q$ where we have a lower bound on $f$ in measure.
By using the expansion of positivity once, we generate a lower bound on a cylinder in the future, with a larger radius. 
By applying iteratively this expansion of positivity, we propagate this lower bound in the future, till the final time. Since we know
that $f$ takes values smaller than $1$ in $\Qfuture$, this gives us some information on the radius of the initial cylinder $Q$.

\subsection{Covering sets with ink spots}

In order to state the covering result that we need in order to establish the weak Harnack's inequality,
we need to introduce the notion of \emph{stacked cylinders}. Given an integer $m \ge 1$ and a cylinder $Q = Q_r (X_0)$,
the stacked cylinder $\bar{Q}^m$ equals $\{ (t,x) \colon 0 < t-t_0 < mr^2,  |v-v_0| < r \}$. \label{p:stack}
\begin{thm}[Leaking ink spots in the wind] \label{t:is-parab}
  Let $E \subset F$ be measurable and bounded sets of $\R \times \R^d$ and $E \subset F \cap Q_1$.
  We assume that there exist two constants $r_0 \in (0,1)$ and an integer $m \ge 1$ such that
  for any cylinder $Q = Q_r (X_0) \subset Q_1$ such that $|Q \cap E | > \frac12|Q|$, we have $\bar{Q}^m \subset F$ and $r < r_0$.
  Then $|E| \le \frac{m+1}m (1-c) \left( |F \cap Q_1| + Cm r_0^2 \right)$. The constants $c \in (0,1)$ and $C>1$ only depend on dimension $d$. 
\end{thm}
\begin{remark}[About the mental picture]
  The two sets $E$ and $F$ are seen as ink spots. They stand in the wind because of the time variable (time delay). And the ink spot $E$ can leak out of the
  reference cylinder $Q_1$. 
\end{remark}
\begin{remark}[About the factor $1/2$]
  The factor $\frac12$ in both assumptions can be replaced with an arbitrary parameter $\mu \in (0,1)$. In this case, the conclusion is
  $|E | \le \frac{m+1}m (1-c\mu)(|F \cap Q_1| + C m r_0^2)$ for some $c\in (0,1)$ and $C>1$ only depending on dimension. 
\end{remark}

We postpone the proof until the next section. 

\subsection{Expansion of positivity for stacked cylinders and  for large times}

In this subsection, we derive the two results that will allow us to use the covering argument from Theorem~\ref{t:is-parab}. There are two assumptions
on cylinders intersecting $E$ in a good proportion: after stacking them, they should lie in $F$, and their radius should be under control. 
\begin{itemize}
\item On the one hand, we check that if we choose  $\eta$ depending on the integer $m$ (from the statement of Theorem~\ref{t:is-parab}) then Proposition~\ref{p:expansion-parab} 
yields  a lower bound in the stacked cylinder $\bar{Q}_1^m$ from an information in measure in $Q_1$.
\item
On the other hand, we apply iteratively Proposition~\ref{p:expansion-parab} with $\eta = 1/2$ in order to estimate
how the lower bound that is generated for small times deteriorates for large ones. 
\end{itemize}

\paragraph{Expansion of positivity for a staked cylinder.} 
We first scale and translate in time the result from Proposition~\ref{p:expansion-parab}  in order to get a statement with $\Qgu$ replaced with $Q_1$. We notice that the stacked cylinder $\bar{Q}_1^m$ equals $(0,m] \times B_{m+2} \times B_1$.
\begin{prop}[Expansion of positivity for a stacked cylinder] \label{p:expansion-parab-s}
  Let $m \ge 1$ be an integer and let   $R_m>1$  be given by Proposition~\ref{p:expansion-parab} for $\eta = 1/\sqrt{m}$.
  Let  $\Qstack = (-1,m] \times B_{\sqrt{2}}$.

  There exists a constant  $M >1$, depending on $d,\lambda,\Lambda$ and $m$, 
  such that, if $u \in \mathrm{pDG}^- (\Qstack,0)$  and $u \ge 0$ a.e. $\Qstack$, then,
  \[ |\{ u \ge M \} \cap Q_1 | \ge \frac12 |Q_1| \quad \Rightarrow \quad \{ u \ge 1 \text{ a.e. in } \bar{Q}_1^m \}.\]
\end{prop}

\begin{figure}[h]
\centering{\includegraphics[height=7cm]{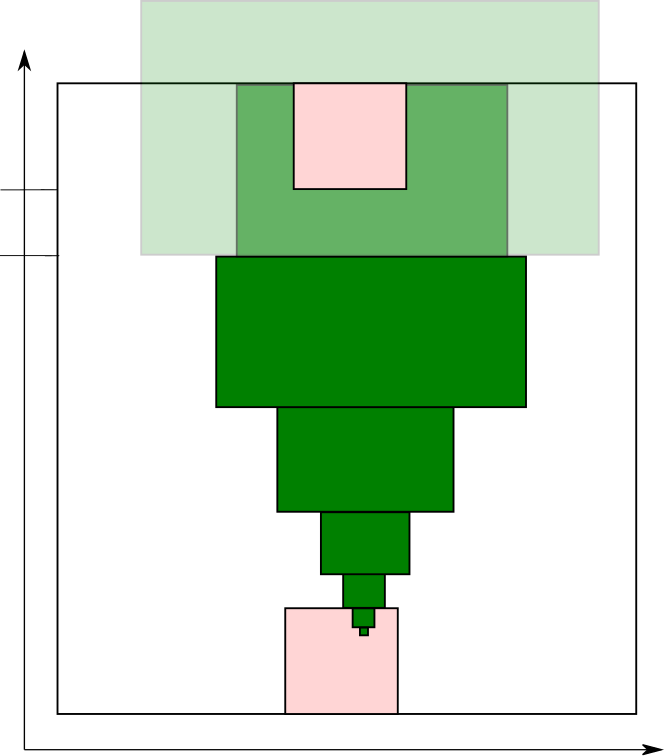}}
\put(-195,147){$-\omega^2$}
\put(-213,130){$t_0+T_N$}
\put(10,0){$(x,v)$}
\put(-180,180){$t$}
\put(-97,20){$\Qpast$}
\put(-96,160){\scriptsize $\Qfuture$}
\put(-90,110){$Q[N]$}
\put(-100,137){$Q[N+1]$}
\caption{\textit{Stacking iteratively cylinders above an initial one  contained in
    $Q_-$.} We see that the stacked cylinder obtained after $N+1$
  iterations by doubling the radius leaks out of the domain. This is
  the reason why $Q[N+1]$ is chosen in a way that it is contained in
  the domain and its ``predecessor'' is contained in $Q[N]$. Notice
  that the cylinders $Q[k]$ are in fact slanted since they are not
  centered at the origin. We also mention that $Q[N+1]$ is chosen
  centered if the time $t_0+T_N$ is too close to the final time~$0$.}
  \label{fig:stackedcylinders-parab}
\end{figure}
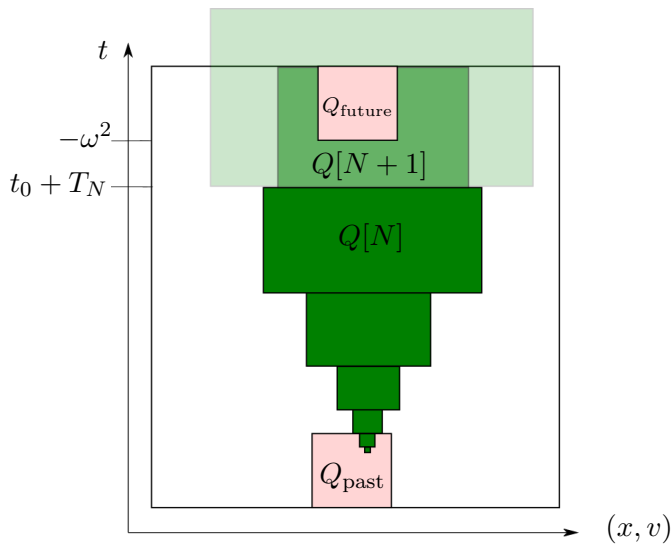

\paragraph{Iteratively stacked cylinders.} We are going to apply iteratively Proposition~\ref{p:expansion-parab} to control the lower bound generated after applying it once.
 We need to make sure that the iterated cylinders do not exhibit the spatial domain  and that their union  captures the cylinder $\Qfuture$, see Figure~\ref{fig:stackedcylinders-parab}.
Recall that $\Qpast = Q_\omega (-1+\omega^2,0)$ and $\Qfuture= Q_\omega$.
\begin{lemma}[Iteratively stacked cylinders]\label{l:stack-parab}
  Let $\omega \in (0,5^{-1/2})$.
  Given $Q = Q_r (X_0) \subset \Qpast$, we define for all $k \ge 1$,  \( T_k = \sum_{j=1}^k (2^j r)^2 \) and pick $N \ge 1$
  the largest integer such that $t_0 +T_N \le 0$. In particular $2^N r \le 1$.

  If $R$ denotes $|t_0+T_N|^{1/2}$, we consider $R_{N+1} =  \max (R,\rho)$ with $\rho = 2 \omega$ and
  \[ \forall k \in \{1,\dots,N\}, \quad X_k = X_0 + (T_k,0) \quad \text{ and } \quad X_{N+1} = \begin{cases} X_N + (R,0) & \text{ if } R \ge \rho, \\
                                                                                                         0 & \text{ if } R < \rho. \end{cases} \]
  We finally define $R_k = 2^k r$ for $k \in \{1,\dots,N\}$ and $Q[k] = Q_{R_k} (X_k)$ for $k \in \{1,\dots, N+1\}$.

  These cylinders $Q[k]$ are such that
  \[ Q[k] \subset (-1,0] \times B_2 \quad  \text{ and } \quad Q[N+1] \supset \Qfuture \quad \text{ and }\quad Q[N] \supset \tilde{Q}[N] \]
  where $\tilde{Q}[N] = Q_{\frac{R_{N+1}}2} (X_{N+1} + (-R_{N+1}^2,0))$. 
\end{lemma}
\begin{proof}
We first check that the sequence of cylinders is well defined for $\omega <5^{-1/2}$. Since
$r \le \omega$, we have $t_0+T_1 \le -1 +\omega^2 + 4r^2 < -1+5 \omega^2<0$. Let $N \ge 1$ be the largest integer such that $t_0 + T_N < 0$.

We check next that $\Qfuture \subset Q[N+1]$.

If $R < \rho$, then $Q[N+1] = Q_\rho$ and we simply remark that $\omega \le \rho$ and recall that $\Qfuture = Q_\omega$ to conclude.

In the other case, that is to say when $R \ge \rho$, we have $Q[N+1] = Q_R (X_{N+1})$ with $X_{N+1} = X_N + (R,0) = (t_0 + T_N + R^2, x_0) = (0,x_0)$.
We have to check that $Q_\omega \subset Q_R (0,x_0)$. We have $\omega \le \rho \le R$ and $\omega +|x_0|\le \omega + r \le 2 \omega \le  \rho \le R$. 
\medskip

Let us now check that for all $k \in \{1,\dots, N+1\},$ we have $Q[k] \subset Q_1$.

By definition of $N$, we have $- (2^{N+1} r)^2 < t_0+T_N <0$. In particular, $R = |t_0+T_N|^{\frac12}\le 2^{N+1} r$ and $R \le 1$ (because $t_0 + T_N \in (-1,0]$) 
$(2^N r)^2 \le T_N \le -t_0 \le 1$. This implies that for all $k \in \{1,\dots, N\}$, $r_k \le 1$. Recalling that $X_k = (t_0+T_k,x_0)$ and $x_0 \in B_r \subset B_\omega$, we thus get
$Q[k] \subset (-1,0] \times B_2$. 

As far as $Q[N+1]$ is concerned, we have $R_{N+1} = \max (R,\rho) \le 1$. Moreover, $x_{N+1} = x_0 \in B_\omega \subset B_1$· So we do have $Q[N+1] = Q_{R_{N+1}} (0,x_{N+1}) \subset (-1,0]\times B_2$. 
\medskip

We are left with proving that $\tilde{Q}[N]  \subset Q[N]$.

If $R \ge \rho$, then $\tilde{Q}_N = Q_{R/2} (0,x_0)$ and $X_{N} = (0,x_0)$ and $R/2 \le 2^N r = R_N$. 

If now $R < \rho$, we have to check that $Q_{\rho/2} (-\rho^2,0) \subset Q_{2^N r}(0,x_0)$. This inclusion holds true if
\[ (2^N r)^2 \ge \rho^2/4 + \rho^2 = 5 \rho^2 /4 \quad \text{ and } \quad \rho/2 + |x_0| \le 2^N r.\]
Because $x_0 \in B_r \subset B_\omega \subset B_\rho$, the second condition holds if $3 \rho/2 \le 2^N r$. And in this case, the first condition is also satisfied. 
In order to prove this inclusion, we first estimate $2^N r$ from below. Since $t_0+T_{N+1} > 0$ and $-t_0 \ge 1 -\omega^2$, we have
$T_{N+1} = (4/3)(4^{N+1}-1)r^2 \ge 1-\omega^2$ and in particular $4^N r^2 \ge  (3/16)(1-\omega^2) \ge 1/8$ (because $\omega^2  \le 2/3$). We conclude
that
\[ 2^N r \ge 1/(2\sqrt2).\]
We finally remark that we do have $3/2 \rho \le 1/(2 \sqrt{2})$ because $\rho \le 1$. 
\end{proof}

\subsection{Iterated expansion of positivity}

\begin{prop}[Iterated expansion of positivity] \label{p:iep-parab}
  There exists a universal constant $\gamma_0 >0$ such that for all $u \in \mathrm{pDG}^- (Q_1,0)$,  all $A>0$ and all cylinder $Q_r (X_0) \subset \Qpast$, 
  \[   |\{ u > A \} \cap Q_r (X_0)| > \frac12 |Q_r(X_0)| \qquad  \Rightarrow \qquad \bigg\{ u \ge A (r/2)^{\gamma_0} \text{ a.e. in } \Qfuture \bigg\}.\]
\end{prop}
\begin{proof}
  We first apply Proposition~\ref{p:expansion-parab} to the function $v = u /A$ after scaling it. This implies that $v \ge A \ell_0$ in $Q[1]$.
  We then apply it iteratively and get $v \ge A \ell_0^k$ in $Q[k]$ for all $k \in \{1,\dots,N\}$. In particular, $v \ge A \ell_0^N$ in $\tilde{Q}[N]$.
  This cylinder is the ``predecessor'' of $Q[N+1]$ and we thus finally get $v \ge A \ell_0^{N+1}$ in $Q[N+1]$. Because $Q[N+1]$ contains $\Qfuture$, we finally
  get $v \ge A \ell_0^{N+1}$ in $\Qfuture$. Now we remember that $2^N r \le 1$ (see Lemma~\ref{l:stack-parab}). We pick $\gamma_0$ such that $\ell_0 = 2^{-\gamma_0}$
  and we write $\ell_0^{N+1} = \left(2^{-(N+1)} \right)^{\gamma_0} \ge (r/2)^{\gamma_0}$.
\end{proof}

\begin{proof}[Proof of Theorem~\ref{t:whi-parab} (weak Harnack's inequality)]
  The proof proceeds in several steps. \medskip

\noindent \textsc{Reduction.}
We first reduce to the case $S=0$ by considering $\tilde u = u + \|S\|_{L^\infty(Q_1)} (t+1)$. 

Second,  we  reduce to the case $\inf_{\Qfuture} u \le 1$ by considering $\tilde u = u / \max(1,\inf_{\Qfuture} u) $. Indeed, $u \le \tilde u \le u + \| S\|_{L^\infty(Q_1)} $ and it is in $\mathrm{pDG}^-(Q_1,0)$.
\medskip

\noindent \textsc{Parameters.} We now aim at proving that there exist two universal constants $p>0$ and $C>0$ such that \( \| u\|_{L^{p_0} (\Qpast)} \le \Cwhi.\)
This is equivalent to prove that there exists three universal constants $M >1$ and $\tilde \mu \in (0,1)$ and $\tilde{C}>1$ such that
\[  \forall k \ge 1, \qquad |\{ u > M^k \} \cap \Qpast | \le \tilde C (1-\tilde \mu)^k .\]
For $k=1$, we simply pick $\tilde \mu \le 1/2$ and $\tilde C \ge 2 |\Qpast|$. We then argue by induction.  
We are going to apply Theorem~\ref{t:is-parab} (about covering with ink spots) for some integer $m \ge 1$ large enough so that $\frac{m+1}m (1-c) <1-c/2$. The parameter $m$ only depends on $c = c(d)$, it is therefore universal.

We are going to use Proposition~\ref{p:expansion-parab-s} (propagation of positivity for stacked cylinders) with $m$ universal as above. Then we obtain another universal parameter $R_m$  from Proposition~\ref{p:expansion-parab}, see the statement of Proposition~\ref{p:expansion-parab-s}. We will also use Proposition~\ref{p:iep-parab} (iterated expansion of
positivity) from which we get yet another universal parameter $R_{1/2}$. Now we choose $R_0 = \max (R_{1/2},23 m^3 R_m)$.
\medskip

\noindent \textsc{The covering argument.}
We are going to apply the ink spots theorem to the sets $E_0 = \{ u > M^{k+1} \} \cap \Qpast$ and $F_0 = \{ u > M^k \} \cap Q_1$ after transforming $\Qpast$ into $Q_1$.
We thus consider a cylinder $Q \subset \Qpast$ such that $|E_0 \cap Q | > \frac12 |Q|$.
We have to check that the stacked cylinder $\bar{Q}^m$ is a subset of $F_0$ and that the radius of $Q$ is controlled by some constant $r_k$. 

We start with checking that $\bar{Q}^m \subset F_0$ for $Q$ such that $|E_0 \cap Q | > \frac12 |Q|$, that is to say
\[ |\{ u > M^{k+1} \} \cap Q| >\frac12 |Q|.\]
 If $Q = Q_r(X_0)$, we consider for $z \in Q_1$ the scaled function $v (X) = M^{-k} u (r (X_0 + X))$, so that
\( |\{ v > M \} \cap Q_1| >\frac12 |Q_1|.\) We have $v \in \mathrm{pDG}^-(\Qstack,0)$.
We deduce from Proposition~\ref{p:expansion-parab-s} that $v \ge 1$ a.e. in $\bar{Q}_1^m$. This means that $u \ge M^k$ a.e. in $\bar{Q}^m$.
We thus proved that $\bar{Q}^m \subset F_0$. 

We now estimate $r$ from above for $Q=Q_r (X_0) \subset \Qpast$ such that \( |\{ u > M^{k+1} \} \cap Q| >\frac12 |Q|\).  Proposition~\ref{p:iep-parab} implies that
$u \ge M^{k+1} (r/2)^{\gamma_0}$ in $\Qfuture$. This implies that $M^{k+1} (r/2)^{\gamma_0} \le 1$ that is to say $r \le 2 M^{-\frac{k+1}{\gamma_0}}=:r_k$.
\medskip

\noindent \textsc{Conclusion.}
Now Theorem~\ref{t:is-parab} implies that
\[ |\{ u > M^{k+1} \} \cap \Qpast| \le  (1-c/2) \left( |\{ u > M^{k} \} \cap \Qpast| + Cm 4 M^{-2\frac{k+1}{\gamma_0}}\right).\]
We use the induction assumption and get
\[ |\{ u > M^{k+1} \} \cap \Qpast| \le  (1-c/2) \left( \tilde C (1-\tilde \mu)^k + Cm 4 M^{-2\frac{k+1}{\gamma_0}}\right).\]
Recall that $M>1$ and $\gamma_0 >0$ are universal. We now choose $\tilde \mu$ so that $(1-\tilde \mu) \ge M^{-2/\gamma_0}$ and
$(1-\tilde \mu)^2 \ge 1- c/2$. 
We get
\begin{align*}
  |\{ u > M^{k+1} \} \cap \Qpast| & \le  (1-\tilde \mu)^2 \left( \tilde C (1-\tilde \mu)^k + Cm 4 (1-\tilde \mu)^{k+1} \right) \\
  & \le \bigg((1-\tilde \mu) \tilde C + 4 Cm\bigg) (1-\tilde \mu)^{k+1} .
\end{align*}
We thus pick $\tilde C$ such that $(1-\tilde \mu) \tilde C + 4 Cm \le \tilde C$ that is to say $\tilde C \ge 4Cm \tilde \mu^{-1}$. 
\end{proof}

\section{Proof of the ink spots theorem}

The assumption of the ink spots theorem asserts that the set $E$ can be covered by cylinders and if more than half the cylinder
lies in $E$, then the corresponding stacked cylinder $\bar Q^m$ is contained in $F$. The conclusion asserts that the volume
of $E$ is bounded from above (up to some multiplicative constant) by the volume of $F$. In order to relate these two volumes,
it is necessary to extract from the original covering another one made of disjoint cylinders, and to make sure that we do not
lose too much by doing so. This is made possible thanks to a parabolic variation of Vitali's lemma with Euclidean balls. 

\subsection{A parabolic Vitali's covering lemma}

As explained in the previous paragraph, Vitali's lemma asserts that a countable disjoint family of cylinders can be extracted from
any covering of a set. We make sure that we do not lose too much by doing so is by imposing that the whole set is recovered
if the radii of cylinders of the sub-covering are multiplied by $5$.

For an arbitrary cylinder $Q \subset \R^{1+d}$, if $Q = Q_r(X_0)$ with $X_0 = (t_0,x_0)$, then $5Q$ denotes $Q_{5r}(t_0 +12r^2,x_0)$.
It is necessary to update the top of the cylinder in order to extract a disjoint sub-cover, see in particular Lemma~\ref{l:overlap-parab}.
\begin{lemma}[Vitali]\label{l:vitali-parab}
  Let $\{Q_j\}_{j \in J}$ be a family of parabolic cylinders whose radii $r_j$ satisfy $\sup_{j \in J} r_j < +\infty$.
  There exists a countable sub-family $\{ Q_{j_i} \}_{i \in \N}$ of disjoint cylinders
  such that
  \[ \cup_{j \in J} Q_j \subset \cup_{i \in \N} 5 Q_{j_i}.\]
\end{lemma}
In order to prove this lemma, we first deal with two overlapping cylinders.
\begin{lemma}[Overlaping parabolic cylinders]\label{l:overlap-parab}
  Let $Q_i = Q_{r_i} (X_i)$ for $i=1,2$ such that $Q_1 \cap Q_2 \neq \emptyset$ and $r_2 \le 2r_1$.
  Then $Q_2 \subset 5 Q_1$. 
\end{lemma}
\begin{proof}
We first reduce to the case $X_1 =0$ by translating both cylinders. By assumption, there exists $X_{1,2} \in Q_1 \cap Q_2$.
This means that there exists $t_{1,2} \in (-r_1^2,0]$ and $x_{1,2} \in B_{r_1}$ such that
\[
  t_2-r_2^2 \le t_{1,2} \le t_2 \quad \text{ and } \quad |x_{1,2} - x_2 | < r_2.
\]
The fact that $Q_2 \subset 5 Q_1$ is equivalent to the following condition
\[
  -13 r_1^2 < t_2 -r_2^2 \le t_2 \le 12 r_1^2 \quad \text{ and } \quad |x_2 | + r_2 < 5r_1.
\]
We check these inequalities one after the other. First, $t_2 \ge t_{1,2} > -r_1^2 \ge r_2^2 - 13 r_1^2 $.
Second, $t_2 \le t_{1,2} + r_2^2 \le r_2^2 \le 4 r_1^2$. Third, $|x_2| \le |x_{1,2}| + |x_{1,2} -x_2| < r_1 + r_2 \le 3 r_1$.
\end{proof}
\begin{proof}[Proof of Lemma~\ref{l:vitali-parab} (Vitali)]
  Let $R = \sup_{j \in J} r_j$ where $r_j$ denotes the radius of the parabolic cylinder $Q_j$.
  Let $\mathcal{F}$ denote the family of cylinders $\{Q_j\}_{j \in J}$ and consider for all $n \ge 1$ the sub-family,
  \[ \mathcal{F}_n = \left\{ Q_j : \; j \in J, \; \frac{R}{2^n}< r_j \le \frac{R}{2^{n-1}} \right\}. \]
  We now construct families $\mathcal{G}_n$ by induction as follows: let $\mathcal{G}_1$ be any maximal disjoint sub-family of $\mathcal{F}_1$.
  Such a sub-family exists because of Zorn's lemma from set theory. If now $n \ge 1$ and $\mathcal{G}_1, \dots, \mathcal{G}_n$ are already constructed,
  then $\mathcal{G}_{n+1}$ is a maximal sub-family of
  \[ \left\{ Q_j \in \mathcal{F}_{n+1} \; : \; Q_j \cap Q_l = \emptyset \text{ for all } Q_l \in \mathcal{G}_1 \cup \dots \cup \mathcal{G}_n\right\}.\] 
  Roughly speaking, we add cylinders with smaller and smaller radii by making sure that they do not intersect the ones we already collected.
  We finally consider \[ \mathcal{G} = \cup_{n=1}^\infty \mathcal{G}_n.\]
  We now verify that this sub-family satisfies the conclusion of the lemma. We consider the sequence of cylinders $Q_{j_i}$ for $i=1,\dots,n$ such that
  $\mathcal{G} = \{ Q_{j_i}\}_{i \in \N}$. Then for $Q_j \in \mathcal{F}$, there exists $n \ge 1$ such that $Q_j \in \mathcal{F}_n$.
  Assume first that $n \ge 1$. By maximality of $\mathcal{F}_1$, there exists $Q_l \in \mathcal{F}_1$ such that $Q_j \cap Q_l \neq \emptyset$. 
  Assume now that $n \ge 2$. Because $\mathcal{G}_n$ is
  maximal, there exists $Q_l \in \mathcal{G}_m$ with $m \in \{1,\dots, n-1\}$ such that $Q_j \cap Q_l \neq \emptyset$. By definition of $\mathcal{F}_n$ and $\mathcal{G}_m$, we have
  $r_j \le \frac{R}{2^{n-1}}$ and $r_l \ge \frac{R}{2^m}$ with either $m=n=1$ or $1 \le m \le n-1$. In both cases, $ r_j \le 2 r_l$. Lemma~\ref{l:overlap-parab} then implies
  that $Q_j \subset 5 Q_l$.   
\end{proof}

\subsection{Lebesgue's differentiation theorem with parabolic cylinders}

\begin{thm}[Lebesgue's differentiation]\label{t:lebesgue-parab}
  Let $u \in L^1 (\R^{1+d})$.  Then for a.e. $X \in \R^{1+d}$,
\[
\lim_{r \to 0+}  \fint_{ Q_r (X)} |u-u(X)| =0
\]
where $\fint_{Q} v= \frac{1}{|Q|}\int_Q v$ for any cylinder $Q \subset \R^{d+1}$ and $v \in L^1 (Q)$.
\end{thm}
The proof of this theorem relies on a functional inequality involving the maximal function.
For $v \in L^1 (\R^{1+d})$, it is defined by,
\[ M v (X) = \sup_{Q \ni X} \fint_Q |v|.\]
\begin{lemma}[The maximal inequality]\label{l:maximal-parab}
  For all $\kappa >0$,
  \[ | \{ M v > \kappa \} \cap \R^{1+d}| \le \frac{C}\kappa \|v \|_{L^1} \]
  for some constant $C$ only depending on dimension. 
\end{lemma}
\begin{proof}
  For every $X \in \R^{1+d}$ such that $Mv (X) > \kappa$, there exists a cylinder $Q$ containing $X$ such that
  \[ \int_{Q} |v| \ge \frac\kappa2 |Q|. \]
  This means that the set $\{ M v > \kappa \}$ is covered with cylinders $\{Q_j\}$ satisfying the previous inequality.
  We know from Vitali's lemma~\ref{l:vitali-parab} that there exists a finite sub-family $\{Q_{j_i}\}_{i \in \N}$ such that
  \[ \{ M v > \kappa \} \subset \cup_{i \in \N} Q_{j_i}.\]
  With such a covering in hand, we can estimate the $L^1$-norm of $v$ as follows:
  \begin{align*}
    \int_{\R^{1+d}} |v| &\ge \sum_{i \in \N} \int_{Q_{j_i}} |v| \\
                               & \ge \frac\kappa2 \sum_{i \in \N}  | Q_{j_i} | \\
                              & = \frac\kappa{ 5^{1+d}2 } \sum_{i \in \N}  | 5 Q_{j_i} | \\
                        & \ge \frac{\kappa}{ 5^{1+d}2 } |\{ M v > \kappa\}|.
  \end{align*}
We thus get the maximal inequality with $C = 5^{1+d} 2$. 
\end{proof}
\begin{proof}[Proof of Theorem~\ref{t:lebesgue-parab} (Lebesgue's differentiation)]
  Let $u_n$ be continuous on $\R^{1+d}$ and such that
  \[ \|u_n - u\|_{L^1} \le \frac1{2^n} .\]
  We can also assume that $u_n \to u$ almost everywhere in $\R^{1+d}$ \cite[Theorem~4.9]{MR2759829}. Let $\mathcal{N}_0$ denote
  the negligible set outside which  point-wise convergence holds. 
  The maximal inequality from Lemma~\ref{l:maximal-parab} tells us that,
  \[ |\{ M (u_n -u) > \kappa \} | \le \frac{C}\kappa 2^{-n}. \]
  This is implies that the non-negative function $\sum_{n\in N} \un_{\{ M (u_n-u) > \kappa\}}$ is integrable on $\R^{1+d}$.
  It is thus finite outside of a negligible set $\mathcal{N}_1 \subset \R^{1+d}$. 
  This implies that there exists $n_\kappa \in \N$ such that for all $n \ge n_\kappa$,
  \[ M (u_n -u) \le  \kappa \quad \text{ outside } \mathcal{N}_1.\]
  For all $i \in N$, we now we pick $\kappa = 1/i$ and construct an increasing sequence $n_i$ such that
  \[ M (u_{n_i} -u) \le \frac1i \quad \text{ outside } \mathcal{N}_1.\]

  With such a sequence of functions in hand, we can write for $X \in \R^{1+d} \setminus (\mathcal{N}_0 \cup \mathcal{N}_1)$ and $i \in \N$,
  \[
    \fint_{Q_r (X)} |u - u(X)|  \le \fint_{Q_r (X)} |u - u_{n_i}| +\fint_{Q_r (X)} |u_{n_i} - u_{n_i}(X)| + |u_{n_i}(X) - u(X)| .
  \] 
  In the right hand side, the first term in bounded from above by $1/i$ because $X \notin \mathcal{N}_1$ and the third term goes to $0$ as $i \to \infty$ because $X \notin \mathcal{N}_0$.
  As far as the second term is concerned, the continuity of $u_{n_i}$ implies that it converges to $0$ too as $i \to \infty$. We thus proved that the left hand side
  tends to $0$ as $i \to \infty$. 
\end{proof}

\subsection{Proof of the ink spots theorem}

The first step of the proof of Theorem~\ref{t:is-parab} is to address the case where the two sets $E$ and $F$ are contained in the cylinder $Q_1$
and in which there is no time delay (no stacked cylinder). 
\begin{lemma}[Crawling ink spots] \label{l:is-parab-nowind}
  Let $E \subset F \subset Q_1$ be measurable  sets of $\R^{1+d}$.
  We assume that $|E| \le \frac12 |Q_1|$ and that
  for any cylinder $Q = Q_r (z_0) \subset Q_1$ such that $|Q \cap E | > \frac12|Q|$, we have $Q \subset F$.
  Then $|E| \le  (1-c)  |F |$. The constant $c \in (0,1)$  only depends on dimension $d$. 
\end{lemma}
\begin{remark}[The factor $1/2$]
  The factor $\frac12$ in both assumptions can be replaced with an arbitrary parameter $\mu \in (0,1)$. In this case, the conclusion is
  $|E | \le (1-c\mu)|F|$ for some $c\in (0,1)$ only depending on dimension. 
\end{remark}
\begin{proof}
    By applying Lebesgue's differentiation theorem~\ref{t:lebesgue-parab} to the indicator function $\un_E$,
  we know that for a.e. $x \in E$, there exists a cylinder $Q^x$ such that $|E \cap Q^x | \ge (1-\iota) |Q^x|$.
  Let us now choose a maximal cylinder $\Qmax^x \subset Q_1$ satisfying $|E \cap Q^x | \ge (1-\iota) |Q^x|$. It is of the form $ \Qmax^x = Q_{\bar r} (\bar t, \bar x)$. 
  By assumption, we know that $\Qmax^x \neq Q_1$ and $\Qmax^x \subset F$.

  We now claim that $|E \cap \Qmax^x | = \frac12 |\Qmax^x|$. 
  If the claim does not hold, then $\Qmax^x \neq Q_1$ and there would be a cylinder $Q^x$ and a $\delta>0$
  such that $\Qmax^x \subset Q^x \subset (1+\delta) \Qmax^x$ with  $Q^x \subset Q_1$ and $|E \cap Q^x| > \frac12 |Q^x|$,
  contradicting the maximality of $\Qmax^x$.
  
  The set $E$ is covered by cylinders $\Qmax^x$. By Vitali's lemma~\ref{l:vitali-parab}, there exists a countable sub-collection of non-overlapping cylinders
  $Q^j = Q_{r_j} (z_j)$, $j \ge 1$, such that $E \subset \cup_{j=1}^\infty 5  Q^j$. Since $Q^j \subset F$ and $|Q^j \cap E| = \frac12 |Q^j|$,
  this implies that $|Q^j \cap (F \setminus E) | = \frac12 |Q^j|$.
  \[
    |F \setminus E|  \ge \sum_{j=1}^\infty |Q^j \cap (F \setminus E)| 
    = \frac12 \sum_{j=1}^\infty  |Q^j| 
    = \frac12 5^{-1-d} \sum_{j=1}^\infty \iota |5 Q^j| 
    \ge \frac12 5^{-1-d} |E|.
  \]
  We conclude that $|F| \ge (1+5^{-1-d}2^{-1}  ) |E|$, from which we get $|E| \le (1-c) |F|$ with $c=5^{-1-d}2^{-2}$.
\end{proof}
We need two preparatory lemmas before proving the ink spot theorem with time delay (wind) and/or leakage. The first one concerns the measure
of a union of time intervals $(a_k -h_k,a_k]$ compared to the measure of the union of their stacked versions $(a_k,a_k + m h_k)$. 
\begin{lemma}[Sequence of time intervals] \label{l:intervals}
  For all $k \ge 1$, let $a_k \in \R$ and $h_k >0$. Then,
  \[ \left| \bigcup_k (a_k,a_k+ m h_k) \right| \ge \frac{m}{m+1} \left| \bigcup_k (a_k-h_k,a_k] \right|.\]
\end{lemma}
\begin{proof}
We consider the open set $\bigcup_{k=1}^N (a_k,a_k+ m h_k)$ of $\R$. Its connected components are open intervals and we can write
\[ \bigcup_{k=1}^N (a_k,a_k+ m h_k) = \cup_l I_l\]
for some disjoint open intervals $I_l$. Each interval $I_l$ is a union of a finite number of intervals $(a_k,a_k+m h_k)$.
Let $a_-- h_-$ be the minimum of the corresponding $a_k- h_k$'s and $a_+ + m h_+$ be the maximum of the corresponding $a_k+ m h_k$'s.

We have in particular $I_l \supset (a_-,a_-+mh_-) \cup (a_+,a_++mh_+)$ with $a_+ + mh_+ \ge a_- + mh_-$.
On the one hand this implies in particular that,
\[ |I_l | \ge a_+ + m h_+ - \min (a_-,a_+)  \ge a_+ + mh_+ - a_-. \]
On the other hand, the fact that $a_-+ m h_- \le a_+ +m h_+$ implies that
\[ a_+ + m h_+ -a_- \ge \frac{m}{m+1} (a_+ + mh_+ - (a_--h_-)).\]
We thus have
\[ |I_l| \ge \frac{m}{m+1} (a_+ + mh_+ - (a_--h_-)).\] 
We remark next that for all $k$ such that $(a_k,a_k + mh_k) \subset I_l$, we have 
\[ a_+ + mh_+ - (a_- -h_-) \ge a_k +m h_k - (a_--h_-) \ge a_k -(a_- - h_-) .\]
Taking the supremum over $k$ yields,
\[ a_+ + mh_+ - (a_- -h_-) \ge \left| \bigcup_{k: (a_k,a_k + mh_k) \subset I_l} (a_k-h_k,a_k] \right|.\]
We thus reach the intermediate conclusion that,
\[ \left| \bigcup_{k=1}^N (a_k,a_k+ m h_k) \right| = \left| \bigcup_l I_l \right| \ge \sum_l |I_l| \ge \frac{m}{m+1} \sum_l \left| \bigcup_{k: (a_k,a_k + mh_k) \subset I_l} (a_k-h_k,a_k] \right|\]
Because the intervals $I_l$ are disjoint, this implies that
\[ \left| \bigcup_{k=1}^N (a_k,a_k+ m h_k) \right| \ge \frac{m}{m+1} \left| \bigcup_{k=1}^N (a_k-h_k,a_k] \right|.\]
Letting $N$ go to $\infty$ allows us to conclude. 
\end{proof}
We can now use this lemma about sequences of intervals to deal with sequence of stacked cylinders.
\begin{lemma}[Overlaping stacked parabolic cylinders] \label{l:stack-overlap-parab}
  Let $\{Q_j\}$ be a family of parabolic cylinders and let $\bar Q_j^m$ be the corresponding stacked cylinders as defined on page~\pageref{p:stack}.
  We have,
  \[ \left| \bigcup_j \bar Q_j^m \right| \ge \frac{m}{m+1} \left| \bigcup_j Q_j\right|.\]
\end{lemma}
\begin{proof}
  We use successively Fubini's theorem, the definition of $\bar Q_j^m$ and  Lemma~\ref{l:intervals} in order to get,
\[  \begin{aligned}[b]
    \left| \bigcup_j \bar Q_j^m \right| &  = \int_{\R^d} \left| \left\{t \in \R: (t,x) \in \bigcup_j \bar Q_j^m \right\}\right| \dx \\
                                        & = \int_{\R^d} \left|  \bigcup_{j : x \in B_{r_j} (x_j)}  (t_j,t_j + m r_j^2) \right| \dx \\
                                        & \ge \frac{m}{m+1}\int_{\R^d} \left|  \bigcup_{j : x \in B_{r_j} (x_j)}  (t_j-r_j,t_j] \right| \dx \\
    & = \frac{m}{m+1} \left| \bigcup_j Q_j\right|.  
    \end{aligned}
    \qedhere
    \]
\end{proof}
Our next task is to get the ink spot theorem in the case where $E$ and $F$ are contained in the cylinder $Q_1$.
In other words, we postpone the treatment of cylinders leaking out of $Q_1$. 
\begin{thm}[Ink spots in the wind] \label{t:is-parab-noleak}
  Let $E \subset F \subset Q_1$ be measurable sets of $\R \times \R^d$.
  We assume that $|E| \le \frac12|Q_1|$ and that
  there exists an integer $m \ge 1$ such that  for any cylinder $Q = Q_r (z_0) \subset Q_1$ such that $|Q \cap E | > \frac12|Q|$, we have $\bar Q^m \subset F$.
  Then $|E| \le  (1-c) \frac{m+1}m |F |$. The constant $c \in (0,1)$  only depends on dimension $d$. 
\end{thm}
\begin{proof}
  We consider the family $\mathcal{Q}$ of parabolic cylinders $Q$ contained in $Q_1$ such that $|Q \cap E| > \frac12|Q|$. We let $G$ denote their union: $G = \bigcup_{Q \in\mathcal{Q}} Q$.
  We know from Lemma~\ref{l:is-parab-nowind} (crawling ink spots) that $E \le (1-c)|G|$. Moreover, the assumption of the theorem implies that $F$ contains the union of
  the corresponding stacked cylinders: $F \supset \bigcup_{Q \in \mathcal{Q}} \bar Q^m$. Using Lemma~\ref{l:stack-overlap-parab} about overlapping stacked cylinders, we obtain the
  following chain of inequalities,
  \[ |F| \ge \left| \bigcup_{Q \in \mathcal{Q}} \bar Q^m \right| \ge \frac{m}{m+1} \left| \bigcup_{Q \in \mathcal{Q}} Q \right| = \frac{m}{m+1}|G| \ge
    \frac{m}{(m+1)(1-c)} |E|. \qedhere\]
\end{proof}
We finally prove the covering result that was used in the derivation of the weak Harnack's inequality.
\begin{proof}[Proof of Theorem~\ref{t:is-parab}]
The assumption of the theorem implies that $|E| \le \frac12 |Q_1|$. Indeed, if this does not true, then $1 \le r_0$, contradicting the fact that $r_0 \in (0,1)$. 

  We consider again the family $\mathcal{Q}$ of parabolic cylinders $Q$ contained in $Q_1$ such that $|Q \cap E| > \frac12|Q|$.
  We let $\bar F$ denote the union of the corresponding stacked cylinders: $\bar F = \bigcup_{Q \in\mathcal{Q}} \bar Q^m$. Theorem~\ref{t:is-parab-noleak}
  implies that
  \[ |E| \le \frac{m+1}m (1-c) |\bar F| = \frac{m+1}m (1-c) \bigg[ |\bar F\cap Q_1| + |\bar F \setminus Q_1| \bigg].\]
  Moreover, the assumptions of Theorem~\ref{t:is-parab} imply that $\bar F \subset F$.   We are thus left we estimating $|\bar F \setminus Q_1|$.
  We claim that for all $Q \in \mathcal{Q}$, we have $\bar Q^m \subset (-1,m r_0^2] \times B_1$. Indeed, $Q = Q_r (z_0)$ for some $z_0 \in Q_1$ and $r< r_0$ and
  $\bar Q^m = (t_0,t_0+ mr^2) \times B_r(x_0)$. In particular, $\bar Q^m \setminus Q_1 \subset (0,mr_0^2) \times B_1$ and thus $\bar F \setminus Q_1 \subset (0,mr_0^2) \times B_1$.
  This implies that $|\bar F \setminus Q_1| \le |Q_1|m r_0^2$. We thus proved,
  \[ |E| \le  \frac{m+1}m (1-c) \bigg[ | F\cap Q_1| + |Q_1| m r_0^2 \bigg]\]
  as desired. 
\end{proof}

\section{Bibliographical notes}
\label{s:biblio-parab}

\paragraph{Parabolic De Giorgi's classes.}
We refer the reader to Section~\ref{s:biblio-elliptic} of the previous chapter (concerned with elliptic equations)  for first definitions of parabolic De Giorgi's classes and first proofs
of H\"older continuity and weak Harnack's inequality for weak solutions and for elements de pDG classes. The difference between $\mathrm{pDG}^+$ and $\mathrm{pDG}^-$ only lies in the propagation in time of $L^2$-norms. We will have to strengthen this assumption in the kinetic setting by imposing a local Poincaré-Wirtinger's inequality. This is the reason why we present in this chapter the original proof by G.~L.~Wang \cite{MR1032780} (see below). 

\paragraph{Kruzhkov's method.}
We already mentioned that S.~N.~Kruzhkov \cite{MR171086} gave an alternative proof of Moser's Harnack inequality for elliptic and parabolic equations.
In particular, is used a different Poincaré inequality, due to S.~L.~Sobolev and V.~P.~Il$'$in, see \cite[Theorem~1.1.]{MR171086} for references. 
He follows J.~Moser by considering the logarithm of the solution. He replaces the logarithm by a smooth approximation of it. His proof is described 
in \cite{zbMATH07750909} with the notation and techniques from this book.
We draw the attention of the reader towards the fact he applies his Poincaré's inequality in the $x$ variable for fixed times $t$.
He thus have to prove a (time propagation) result in the spirit of Lemma~\ref{l:st-propagation}. 

\paragraph{Degenerate equations.}
Let us come back to the techniques built on De Giorgi's original ideas (iterative truncation, gain of integrability, improvement of oscillation). 
They were  used in many subsequent works and we will give just a few references in this paragraph. Degenerate elliptic equations like the porous medium one can be handled. 
This was first observed by L.~Caffarelli and A.~Friedman \cite{MR534112}: they prove the continuity of global solutions. A local version of this result
can be found in \cite{MR684758} by E.~DiBenedetto. Then  E.~DiBenedetto and A.~Friedman addressed in \cite{MR783531} the case of degenerate parabolic \emph{systems}.
All these articles build on De Giorgi's techniques. E.~DiBenedetto contributed to this field by numerous articles, dealing in particular with $p$-Laplace operators.
His book \cite{MR1230384} from 1993 was influential. He wrote another book with U.~Gianazza and V.~Vespri \cite{MR2865434} about equations of $p$-Laplace and porous media type. 

\paragraph{Expansion of positivity.}
The wording ``Expansion of positivity'' first appears in a paper written by E.~DiBenedetto \cite{MR1230384}. He gave a useful and meaningful name to a phenomena
exhibited  in most of the works dealing with De Giorgi's methods that are mentioned in the book you have in hand (or on your screen). In other contexts, it is called
the doubling property. It is related to growth lemmas by E.~Landis, N.~V.~Krylov, M.~V.~Safonov, among other authors. 
In this book, the starting point of the proof of expansion of positivity is the one by G.~L.~Wang contained in \cite{MR1032780} and also presented in \cite{MR1465184}.
But the geometric setting is different, and we argue by following the path that the kinetic proof will traverse.

\paragraph{Ink spots in the wind.}
We quickly come back to the measure result used in the covering argument to derive the weak Harnack's inequality, see Theorem~\ref{t:is-parab}.
It differs from the crawling ink spots lemma for elliptic equations (see Lemma~\ref{l:is} and Section~\ref{s:biblio-elliptic}) in the sense that its statement involves so called
stacked cylinders. Such cylinders are obtained by stacking a finite number of copies of an original cylinder above it, in the future. They have to be considered because of the time variable. Indeed, the expansion of positivity takes place in the future. Moreover, it is necessary to consider cylinders that are going to leak out of the domain.
 E.~Landis \cite{MR1487894} refers to this type of result as crawling ink spots ``in the wind'', see for example \cite[Lemma~2.3]{MR563790}. 
The proof that is presented in this chapter is the parabolic counterpart of the one contained in \cite{MR4049224} about the kinetic case.

\chapter{Kinetic Fokker-Planck equations}
\label{c:kin}

In this chapter, we prove that weak solutions of a class of kinetic equations with rough coefficients
are locally H\"older continuous, in the spirit of De Giorgi's theorem for elliptic equations and
Nash's theorem for parabolic equations. In order to do so, we will proceed as in the two previous chapters about elliptic and parabolic equations:
 we shall first show that these solutions lie in an appropriate kinetic De Giorgi's class, and then derive a local H\"older estimate for elements of this class.

\section{Kinetic Fokker-Planck equations}

Let us first define the kinetic equations that we are going to work with throughout this chapter. 
Let $I$ be a bounded interval of $\R$ of the form $(a,b]$ with $a,b \in \R$, let $\Ox$ and $\Ov$ be two open sets of $\R^d$, and let
\[ \domain = I \times \Ox \times \Ov.\]
Let $\lambda, \Lambda$ be two positive constants with $\lambda \le \Lambda$.
We consider
\[ \mathcal{E} (\lambda, \Lambda) = \{ A \in L^\infty ( \domain, \mathbb{S}_d (\R)),
\text{ a.e. in } \domain, \forall \xi \in \R^d, \lambda |\xi|^2 \le A \xi \cdot \xi \le \Lambda |\xi|^2 \}.\]
For $A \in \mathcal{E} (\lambda, \Lambda)$ and $B \in L^\infty (\domain)$ such that
\[ \|B\|_{L^\infty(\domain)} \le \Lambda,\]
and $S \in L^1 (\domain)$, we consider the following equation,
\begin{equation}
  \label{e:kfp}
  \partial_t f + v \cdot \nabla_x f = \dive_v ( A \nabla_v f) + B \cdot \nabla_v f+ S \quad \text{ in } \quad \domain .
\end{equation}
In the left hand side of the equation, the differential operator is known as the \emph{free transport operator}.
The right hand side contains the diffusion operator that was considered for (local) elliptic and  parabolic equations, together with drift and source terms.
The diffusion operator acts on the velocity variable $v$ only.

\subsection*{Kolmogorov's equation}

When $A$ is constant and equal to the identity matrix and $B$ is constant and equal to zero, \eqref{e:kfp} is called the \emph{Kolmogorov equation},
\begin{equation}
  \label{e:kolm}
  \partial_t f + v \cdot \nabla_x f = \Delta_v f + S.
\end{equation}

\subsection*{Degenerate ellipticity \& hypoellipticity}

The key difficulty  to be addressed in order to get a regularity result \`a la De Giorgi for this class of equations
is the lack of uniformly ellipticity in $(x,v)$. They are only elliptic in the velocity variable $v$. The equation is thus
degenerate elliptic (parabolic), in the sense that some eigenvalues (the ones corresponding to the $x$ diffusion) equal zero. 

But these degenerate parabolic equations enjoy a specific structure: the free transport operator knows how to talk to the diffusion operator in velocity.
More precisely,  the free transport operator $(\partial_t + v \cdot \nabla_x)$ does not commute with the vector field $\nabla_v$ and their commutator
is precisely $\nabla_x$: for a smooth function $f$, we have
\[ [(\partial_t + v \cdot \nabla_x), \nabla_v] f:=(\partial_t + v \cdot \nabla_x) \nabla_v f - \nabla_v (\partial_t + v \cdot \nabla_x) f  =  \nabla_xf.\]
The vector field $\nabla_v$ represents the diffusion since one can write $\dive_v (A \nabla_v f) = (\sqrt{A} \nabla_v)^\ast (\sqrt{A} \nabla_v)$ where $(\sqrt{A} \nabla_v)^\ast$ denotes the adjoint (in $L^2$)
of the differential operator $(\sqrt{A} \nabla_v)$. 

In H\"ormander's hypoelliptic theory \cite{hormander}, the map $z \mapsto A (z)$ is assumed to be smooth.
We do not want to make such an assumption on coefficients because we aim at getting regularity estimates for non-linear problems. 

When proving the kinetic De Giorgi \& Nash's result (Theorem~\ref{t:dg-kinetic}), the regularity in the $x$ variable
will not be recovered by a commutator argument. It will take the form of a transfer of regularity result, saying that
any regularity in $v$ can be transferred in some regularity in $x$. Such a phenomenon was first discovered by H\"ormander \cite{hormander},
it was then further explored by establishing \emph{averaging lemmas}: given a smooth function $\varphi$, the transport operator ensures that any mean
$\int_{\R^d} f(t,x,v) \varphi (v) \dv$ in $v$ is more regular than the solution $f$ itself. It was then possible
to use these lemmas to prove transfer of regularity results, resulting in the translation of $v$ regularity into
some regularity in $x$ of the function itself (and not only its mean). See the bibliographical section~\ref{s:biblio-kin} for references. 

In order to prove such a transfer of regularity result, we will use a trick due to A.~Pascucci and S.~Polidoro~\cite{pp}
that will be described in due time (see Section~\ref{s:x-reg}). 

\subsection*{Kinetic scaling}

\label{p:scaling}
Let $R>0$. For $z = (t,x,v) \in \R\times \R^d \times \R^d$, we define the scaling operator $\sigma_R$ by
\[ \sigma_R (z) =  (R^2t,R^3 x,Rv).  \]
If $f$ is a solution of the kinetic Fokker-Planck equation \eqref{e:kfp} for some $A \in \mathcal{E} (\lambda,\Lambda)$,
then the function $f_R (z) = f (\sigma_R (z))$ satisfies \eqref{e:kfp} with $A$ is replaced with $A_R (z) = A (\sigma_R(z))$.
Notice that $A_R \in \mathcal{E}(\lambda,\Lambda)$. 

We will sometimes abuse notation by simply writing $Rz$ for $\sigma_R (z)$. 

\subsection*{Galilean invariance}

In contrast with elliptic and parabolic equations, the class of kinetic Fokker-Planck equations is not invariant
by translation: given $z_0 =(t_0,x_0,v_0)$ and a solution $f$ of \eqref{e:kfp}, the function $g(z) = f(z_0+z)$ is
not a solution of \eqref{e:kfp} for another $A_0 \in \mathcal{E}(\lambda,\Lambda)$.\label{elambda2} Since kinetic equations
encode a law from statistical physics, it is expected to be invariant under Galilean transformations. More precisely,
we define
\[ z_0 \circ z = (t_0+t, x_0+x + t v_0,v_0+v).\]
Times and velocities are just translated while the position variable is corrected by $t v_0$, encoding the fact that
we look at the equation in a frame moving at constant speed $v_0$ relatively to the reference frame.
It is useful to compute for $z =(t,x,v)$ and $z_i = (t_i,x_i,v_i)$, $i=1,2$,
\[ z^{-1} = (-t,-x+tv,-v) \quad \text{ and } \quad z_1^{-1} \circ z_2 = (t_2 -t_1,x_2 -x_1 - (t_2-t_1) v_1 , v_2 -v_1). \]

\section{Kinetic geometry: cylinders, distance and H\"older regularity}

\subsection*{Kinetic cylinders}
\label{kcyl}
In order to define a kinetic distance, we first introduce cylinders respecting the kinetic scaling and Galilean invariance:
for $R>0$ and $z_0  \in \R^{1+2d}$,
\begin{align*}
  Q_1 &= (-1,0] \times B_1 \times B_1, \\
  Q_R & = (-R^2,0] \times B_{R^3} \times B_R, \\
  Q_R (z_0) & = \{ z_0 \circ z : z \in Q_R \}.
\end{align*}
We will use next a more explicit definition of $Q_R (z_0)$,
\[ Q_R (z_0 ) = \{ (t,x,v) \in \R^{1+2d} : -R^2 < t-t_0 \le 0, |x-x_0 - (t-t_0)v_0 | < R^3, |v-v_0|< R\}.\]
The following lemma is useful.
\begin{lemma} \label{l:add-rad}
If $z_i \in Q_{r_i}$ for $i=1,2$, then $z_1 \circ z_2 \in Q_{r_1+r_2}$.  
\end{lemma}
\begin{proof}
  We recall that $z_1 \circ z_2 = (t_1+t_2, x_1+x_2 + t_2 v_1, v_1+v_2)$.
  \begin{align*}
    |t_1+t_2| &\le r_1^2 + r_2^2 \le (r_1+r_2)^2 \\
    |v_1+v_2| & \le r_1 + r_2 
  \end{align*}
\[    |x_1+x_2+ t_2 v_1 |  \le r_1^3 + r_2^3 + r_2^2 r_1 \le (r_1+r_2)^3.   \qedhere \]
\end{proof}
Let $\bar B_1$ denote the closed unit ball.
We also consider the closed cylinders $\bar Q_R (z_0)$, associated with $\bar Q_1 = [-1,0]\times \bar B_1 \times \bar B_1$. 

\subsection*{Kinetic distance}

We now define a kinetic distance on the space $\R^{1+2d}$.
\begin{defi}[kinetic distance] \label{defi:kin-distance}
For $z_1,z_2 \in \R^{1+2d}$,
the kinetic distance $\dkin (z_1,z_2)$ is  the infimum of the set of real numbers $r>0$ such that
there exists $z \in \R^{1+2d}$ such that $z_1,z_2 \in Q_r (z)$. 
\end{defi}

It is not clear from this definition that $\dkin$ satisfies the triangle
inequality. Before justifying it, we start with collecting elementary properties
immediately obtained from the previous definition. 
\begin{lemma}[Elementary properties of $\dkin$]\label{l:elem}
  The kinetic distance satisfies:
\begin{enumerate}[label=(\roman*)]
\item \label{i:sym}
  \textsc{(Symmetry)}  for all $z_1,z_2 \in \R^{1+2d}$, $\dkin (z_1,z_2) = \dkin (z_2,z_1)$;
\item  \label{i:left-inv}
\textsc{(Left invariance)}  for all $z_1,z_2,z \in \R^{1+2d}$, $\dkin (z_1,z_2) = \dkin (z^{-1} \circ z_1, z^{-1} \circ z_2)$;
\item \textsc{(Scaling)} \label{i:d-scaling}
  for all $z_1, z_2 \in \R^{1+2d}$, $R>0$, we have
  $\dkin (\sigma_R (z_1), \sigma_R (z_2)) = R \dkin (z_1,z_2)$. 
\end{enumerate}
\end{lemma}
The triangle inequality will follow from the following lemma, that  is of independent interest.
\begin{lemma}[Formula for the kinetic distance] \label{l:expression}
  Let $z_1,z_2 \in \R^{1+2d}$. Then
  \[ \dkin (z_1,z_2) = \min_{w \in \R^d} \max \left(|t_1-t_2|^{\frac12},   |v_1-w|,|v_2-w|, 2^{-1/3} |(x_1-x_2) - (t_1-t_2) w |^{1/3} \right). \]
  In particular, there exists $z$ such that $z_1,z_2 \in \bar Q_\rho (z)$ for $\rho = \dkin (z_1,z_2)$.
  More precisely, $z = (s,y,w)$ with $s = \max(t_1,t_2)$, $y = \frac{x_1+x_2}2 - w$ and $w$ realizes the minimum in the previous formula. 
\end{lemma}
\begin{proof}
  We first remark that there exists $z =(t,x,v) \in \R^{1+2d}$ such that $z_1,z_2 \in \bar Q_\rho (z)$ with
  $\rho = \dkin (z_1,z_2)$. Indeed, if one considers a sequence $r_n >0$ realizing the infimum from the definition of
  the kinetic distance, and if $z_n$ denotes the corresponding centers for cylinders, then $r_n$ and $z_n$ are bounded
  sequences, in particular, they converge up to a sub-sequence.

  The fact that $z_1,z_2 \in \bar Q_\rho (z)$ with $z=(s,y,w)$ is equivalent to,
  \[ \text{ for } i =1,2, \quad \begin{cases} -\rho^2 \le t_i - s \le 0, \\ |v_i - w| \le \rho, \\ |x_i - y - (t_i-s)w| \le \rho^3. \end{cases}\]
  This is equivalent to $s \ge \max(t_1,t_2)$ and,
  \[ \rho \ge \max \left( (s-t_1)^{\frac12}, (s-t_2)^{\frac12}, |v_1-w|,|v_2-w|, |x_1 - y - (t_1-s)w|^{\frac13}, |x_2 - y - (t_2-s)w|^{\frac13} \right).\]
  Since $\rho$ is as small as possible, this inequality is in fact an equality. Moreover, this equality holds for any $z$ such that $z_1,z_2 \in \bar Q_\rho (z)$.
  This implies that
  \[ \rho = \min_{i,s,y,w} \left\{ \max \left( (s-t_i)^{\frac12},  |v_i-w|,|x_i - y - (t_i-s)w|^{\frac13} \right) : i=1,2, s \ge \max(t_1,t_2)\right\}.\]
  We next remark that the minimum in $y$ of $\max_{i=1,2} |x_i - y - (t_i-s)w|^{\frac13}$ is reached for $y = \frac12 \sum_{i=1,2} (x_i - (t_i-s)w) =  sw + \frac12 \sum_{i=1,2} (x_i - t_iw) $.
  It is then equal to
  \[ 2^{-1/3}\left|  (x_1-x_2) - (t_1-t_2) w \right|^{\frac13}.\]
  We deduce from this observation that 
  \[ \rho = \min_{i,s,w} \left\{ \max \left( (s-t_i)^{\frac12}, |v_i-w|,2^{-\frac13}|(x_1-x_2) -  (t_1-t_2)w|^{\frac13} \right) : i=1,2, s \ge \max(t_1,t_2)\right\}.\]
  We conclude the proof of the lemma by finally remarking that the infimum in $s$ is reached for $s = \max(t_1,t_2)$. 
\end{proof}

We are now ready to prove that $\dkin$ satisfies the triangle inequality.
\begin{lemma}[Distance property]
Let $z_1,z_2,z_3 \in \R\times \R^d \times \R^d$. We have: $\dkin (z_1,z_2) \le \dkin (z_1,z_3) + \dkin (z_3,z_2)$. 
\end{lemma}
\begin{proof}
  By the left invariance property (Lemma~\ref{l:elem}-\ref{i:left-inv}), we can assume that $z_3 = 0$. For $i=1,2$, let $r_i$ denote $\dkin (z_i,0)$ (recall the symmetry property).
  By Lemma~\ref{l:expression}, we know that there exists $w_i \in \R^d$ such that $r_i = \max (|t_i|^{\frac12}, |v_i-w_i|,|w_i|, 2^{-\frac13}|x_i- t_i w_i|^{\frac13})$ for $i=1,2$.
  We now remark that,
  \begin{align*}
    |t_1 -t_2|^{\frac12} \le (|t_1|+|t_2|)^{\frac12} \le |t_1|^{\frac12} + |t_2|^{\frac12} &\le r_1 + r_2, \\
    |v_i-(w_1+w_2)| \le |v_i-w_i| + |w_j| &\le r_1 + r_2, \qquad (j\neq i), 
  \end{align*}
  and
  \begin{align*}
    |(x_1-x_2) - (t_1-t_2)( w_1+w_2) |^{\frac13} &\le \left| |x_1-t_1 w_1| + |x_2 -t_2 w_2| + |t_1||w_2|+|t_2||w_1| \right|^{\frac13} \\
                                                 & \le \left| 2 r_1^3 + 2 r_2^3 + r_1^2 r_2 + r_2^2 r_1 \right|^{\frac13} \\
                                                 & \le 2^{\frac13} (r_1 + r_2).
  \end{align*}
  Then the result follows from the characterization of $\dkin$ given by Lemma~\ref{l:expression}. 
\end{proof}

We finally give upper and lower bounds of the kinetic distance in terms
of a supremum norm of kinetic type. For $z \in \R^{1+2d}$,
we define
\[ |z|_\infty = \max \left(|t|^{\frac12},|x|^{\frac13}, |v|\right).\]
\begin{remark} \label{rem:stupid}
We notice that for $z_1,z_2 \in \R^{1+2d}$, we have $z_1 \in Q_{r_{1,2}} (z_2)$ with $r_{1,2} = |z_2^{-1} \circ z_1|_\infty$. 
\end{remark}
\begin{lemma}[Upper and lower bounds] \label{l:up-low}
  For all $z_1,z_2 \in \R^{1+2d}$, we have: \( \frac12 |z_2^{-1} \circ z_1 |_\infty \le \dkin (z_1,z_2) \le |z_2^{-1} \circ z_1|_\infty.\)
  The constants $1$ and $\frac12$ are optimal. 
\end{lemma}
\begin{proof}
  By left invariance, see Lemma~\ref{l:elem}-\ref{i:left-inv}, we reduce to the case $z_2=0$.  Let $d_1$ denote $\dkin (z_1,0)$ and $D_1 = |z_1|_\infty$. Let also $w_1$ denote
  the $w$ realizing the minimum in the explicit formula of the kinetic distance,see  Lemma~\ref{l:expression}
  From this lemma, we can pick $w=0$ and get that $d_1 \le D_1$.
  Let us prove the lower bound. If $D_1 = |t_1|^{\frac12}$, then $d_1 \ge D_1$
  thanks to the explicit form of $\dkin$, see Lemma~\ref{l:expression}.
  
  If $D_1 = |v_1|$, then we distinguish two sub-cases.
If $|w_1|\le |v_1|/2$,
  then $d_1 \ge |v_1 -w_1| \ge |v_1|/2 = D_1/2$. If now $|w_1| \ge |v_1|/2$,
  then $d_1 \ge |w_1| \ge |v_1|/2 = D_1/2$.
  
  We are left with assuming that $D_1 = |x_1|^{\frac13}$. We also distinguish two
  sub-cases. If $|w_1| \le |x_1|^{\frac13}/2$, then $|t_1w_1| \le D_1^2 |x_1|^{\frac13} \le |x_1|/2$. In turn, $|x_1-t_1w_1| \ge |x_1|/2$ and $d_1 \ge \frac1{4^{\frac13}} |x_1|^{\frac13} \ge \frac12 D_1$. If $|w_1| \ge |x_1|^{\frac13}/2$, then
  $d_1 \ge |w_1| \ge |x_1|^{\frac13}/2 = D_1/2$. 

  The optimality of the constant $\frac12$ in the lower bound is obtained
  by choosing $t_1=0$ and $x_1=0$. In this case, $d_1 = |v_1|/2 = D_1/2$.
  The optimality of the constant $1$ in the upper bound is easily obtained
  for instance with $x_1=x_2=w_1=w_2$ and $t_1 \neq t_2$. 
\end{proof}

\subsection*{Kinetic H\"older regularity}

We now characterize the H\"older regularity of a function by its oscillation.
It is convenient to use $\osc (f \;|\; E)$ for  $\osc_E f$ when the writing of $E$ is too long to be written as a subscript.  
\begin{prop}[H\"older regularity and oscillation] \label{p:holder-kin}
  Let $f \in L^\infty (Q_R (z_0))$ be such that for all $z \in Q_R (z_0)$ and all $r>0$, we have
  \( \osc (f \;|\; Q_r (z) \cap Q_R (z_0))  \le C_0 r^\alpha.\)
  Then for all $z_1,z_2 \in \bar Q_R (z_0)$, we have
  \( | f(z_1) - f(z_2) | \le C_0 \dkin (z_1,z_2)^\alpha.\)
\end{prop}
\begin{remark}
If $z \notin Q_R(z_0)$, then for $r>0$ small enough, $Q_r (z) \cap Q_R (z_0)$ is empty. 
\end{remark}
\begin{proof}
  We first prove the result by further assuming that $f$ is continuous.
  Let $z_1, z_2 \in \bar Q_R (z_0)$ and $\rho = \dkin (z_1,z_2)$. The definition of the kinetic distance implies that there exists $z \in \R^{1+2d}$
  such that $z_1,z_2 \in \bar Q_R (z)$. 
  The assumption implies that
  \begin{align*}
    |f(z_2) - f(z_1) | &\le \osc (f \; | \; Q_{\rho} (z_1) \cap Q_R (z_0) ) \\
                       & \le C_0 \rho^\alpha .
  \end{align*}
  
  We now address the case where $f$ is only essentially bounded in $Q_R (z_0)$.
  We  regularize it by a convolution procedure respecting the Lie group structure: for $\eps >0$,
  we consider a mollifier $\rho^\eps (z) = \eps^{-2 - 4d} \rho (S_{1/\eps} (z^{-1}))$ with $\rho$ smooth, non-negative, supported in $Q_1$ and $\int_{\R^{1+2d}} \rho (z) \dd z =1$,
  and we define for $z \in Q_{R-\eps} (z_0)$,
  \[ f^\eps (z) = \int_{\R^{1+2d}}  f(\zeta^{-1} \circ z) \rho_\eps (\zeta) \dd \zeta .\]
  This function $f^\eps$ is continuous in $Q_{R-\eps_0} (z_0)$ because $f$ is essentially bounded.
  We now verify  that $\bar f^\eps$ inherits the property satisfied by the oscillation of $f$.
  Let $\eps_0>0$ and $\eps \in (0,\eps_0)$ and $z \in \R^{1+2d}$. Then for $z_1 \in Q_{R-\eps_0} (z_0)$ and 
   $\zeta^{-1} \in Q_\eps$, we have $\zeta^{-1} \circ z_1 \in Q_R (z_0)$ (see Lemma~\ref{l:add-rad}). In particular,
  \begin{align*}
    \osc ( f^\eps \; | \; {Q_r (z) \cap Q_{R-\eps_0} (z_0)} )  & \le \int_{\R^{1+2d}}  \osc ( f \; | \; Q_r (\zeta^{-1} \circ z) \cap Q_R (z_0) )   \rho_\eps (\zeta) \dd \zeta \\
    & \le C_0 r^\alpha.
  \end{align*}
  From the continuous case, we deduce that for all $z_1,z_2 \in Q_{R-\eps_0} (z_0)$, we have
  \[ |f^\eps(z_2)-f^\eps(z_1)| \le C_0^\alpha \dkin (z_1,z_2)^\alpha. \]
  By dominated convergence, we see that $f^\eps \to f$ a.e. in $Q_{R-\eps_0} (z_0)$.  We thus conclude that for a.e. $z_1,z_2 \in Q_{R-\eps_0} (z_0)$,
  \[ |f(z_2)-f(z_1)| \le C_0^\alpha \dkin (z_1,z_2)^\alpha.\]
  In particular, $f$ can be extended into a (H\"older) continuous function in $Q_{R-\eps_0} (z_0)$.
  Since $\eps_0>0$ is arbitrarily small, the proof is complete. 
\end{proof}

\subsection*{Kinetic Young's inequality}

We now introduce the \emph{kinetic convolution} of two functions $f$ and $g$  defined in $\R^{1+2d}$,
 \[ f \astkin g (z) = \int_{\R^{1+2d}} f (\zeta^{-1} \circ z) g (\zeta) \dd \zeta = \int_{\R^{1+2d}} f (\zeta) g (z \circ \zeta^{-1}) \dd \zeta .\]
 It can be useful to write the formula with $z=(t,x,v)$ and $\zeta = (s,y,w)$,
 \begin{equation} \label{e:kin-convol-txv}
   f \astkin g (t,x,v) = \int_{\R^{1+2d}} f (t-s,x-y - (t-s) w , v-w) g(s,y,w) \ds \dy \dw  .
 \end{equation}
 
 Notice that we used it in the previous proof to regularize $f$ by considering a mollifier $\rho$. 
The following elementary lemma will be useful in the sequel.
\begin{lemma}[kinetic convolution and $L^2$ duality] \label{l:convol-L2}
  For $f,g,h \colon \R^{1+2d} \to \R$, we have
  \[ \int_{\R^{1+2d}} ( f\astkin g) h = \int_{\R^{1+2d}} g (\check{f} \astkin h) = \int_{\R^{1+2d}} f (h \astkin \check{g}) \]
  where $\check{f} (z) = f (z^{-1})$ and $\check{g} (z) = g (z^{-1})$. 
\end{lemma}

The group $(\R^{1+2d}, \circ)$  equipped with the Lebesgue measure is uni-modular and the classical Young's inequality holds  \cite{MR2768550}. 

In order to state it, we let $M_1^+( \R^{1+2d})$ denote the set of non-negative measures on $\R^{1+2d}$. For $\mathfrak{m} \in M_1^+ (\R^{1+2d})$,
the total mass $\mathfrak{m} (\R^{1+2d})$ is denoted by $\| \mathfrak{m}\|_{M_1^+ (\R^{1+2d})}$.

We also consider weak Lebesgue spaces $L^{p,\weak}(\R^N)$. It is the set of functions $f \colon \R^N \to \R$ such that their weak $L^p$-norm is finite,
\[ \|f\|_{L^{p,\weak}(\R^N)} := \sup_{\alpha >0}  \hspace{2ex}  \alpha |\{ \xi \in \R^N : f(\xi) > \alpha \}|^{\frac1p} . \]
\begin{remark}
  The reader that is unfamiliar with this notion can check that Bienaymé-Chebyshev's inequality implies that
  classical $L^p$ functions are in $L^{p, \weak}$ and that their weak $L^p$ norm is bounded by their classical $L^p$ norm.
  See also \cite[p.~25]{boyer2012mathematical}.
\end{remark}
\begin{lemma}[Young's inequality for the kinetic convolution] \label{l:young-kin}
\begin{enumerate}[label=(\roman*)] 
\item  For $f \in L^p (\R^{1+2d})$ and $g \in L^q (\R^{1+2d})$, with $p,q \in [1,\infty]$, we have $f \astkin g \in L^r(\R^{1+2d})$
  with $1+\frac1{r} = \frac{1}p + \frac{1}q$, and the following estimate holds,
  \[
    \| f \astkin g \|_{L^r (\R^{1+2d})} \le \| f \|_{L^p (\R^{1+2d})} \|g \|_{L^q (\R^{1+2d})} .
  \]
\item For $\mathfrak{m} \in M_1^+ (\R^{1+2d})$ and $f \in L^{p} (\R^{1+2d})$ with $p \in (1,\infty)$ and $g \in L^{p,\weak} (\R^{1+2d})$, 
  \begin{align*}
    \| f \astkin \mathfrak{m} \|_{L^{p} (\R^{1+2d})} & \le \| f \|_{L^{p} (\R^{1+2d})} \|\mathfrak{m} \|_{M_1^+ (\R^{1+2d})}, \\
    \| g \astkin \mathfrak{m} \|_{L^{p,\weak} (\R^{1+2d})} & \le \| g \|_{L^{p,\weak} (\R^{1+2d})} \|\mathfrak{m} \|_{M_1^+ (\R^{1+2d})}.
  \end{align*}
\end{enumerate}
\end{lemma}

\section{The Kolmogorov equation}

In order to study the regularity of solutions $f$ and their truncation $(f-\kappa)_\pm$,
we use a trick due to A.~Pascucci and S.~Polidoro \cite{pp} and write the kinetic Fokker-Planck equation
$(\partial_t + v \cdot \nabla_x ) f = \dive_v (A \nabla_v f) + B \cdot \nabla_v f + S$ as
\[(\partial_t + v \cdot \nabla_x) f - \Delta_v f = \dive_v ((A-I) \nabla_v f) + B \cdot \nabla_v f + S \]
where $I$ denotes the $d \times d$ identity matrix. We thus turn the kinetic Fokker-Planck equation
into the Kolmogorov equation with a singular source term. We remark that local energy estimates ensure
that $\nabla_v f$ is square integrable. In particular, the source term in the Kolmogorov equation can
be written as $\dive_v \mathfrak{S} + \bar S$ for some $\mathfrak{S}, \bar S \in L^2$. 

It is convenient to write the Kolmogorov operator $\mathcal{K}$,
\[ \mathcal{K} f = (\partial_t + v \cdot \nabla_x) f - \Delta_v f.\]

\subsection{Fundamental solution}

A.~Kolmogorov computed in \cite{MR1503147} solutions of the linear Fokker-Planck equation
in the case where the diffusion matrix $A$ is constant and equal to the identity.
They can be represented thanks to the counterpart of the heat kernel,
\[
  \Gamma (t,x,v) =
  \begin{cases}
\left(\frac{3}{4 \pi^2}\right)^{\frac{d}2} \frac{1}{t^{2d}}\exp \left[ - \frac{3 \left|x+\frac{t}2 v \right|^2}{t^3} - \frac{|v|^2}{4t} \right]  & \text{ if } t > 0, \\     
    0 & \text{ if } t \le 0.
  \end{cases}
\]
Such a function is referred to as the \emph{fundamental solution} of the Kolmogorov equation. 

 Proofs will use the following functions associated with $\Gamma$: for $t>0$ and $x,v \in \R^d$,
 \begin{align*}
   \Gamma_1 (x,v) &= \Gamma (1,x,v) \\
   \Gammax (t,x,v) &= t \nabla_x \Gamma (t,x,v) = \frac1{t^{2d+\frac12}} \nabla_x \Gamma_1 \left( \frac{x}{t^{\frac32}}, \frac{v}{t^{\frac12}}\right) , \\
      \Gammav (t,x,v) &=  \nabla_v \Gamma (t,x,v) = \frac1{t^{2d+\frac12}} \nabla_v \Gamma_1 \left( \frac{x}{t^{\frac32}}, \frac{v}{t^{\frac12}}\right) .
 \end{align*}
The second formulas for $\Gammax$ and $\Gammav$ are consequences of the scale invariance stated in the next proposition. 
 \begin{prop}[Properties of the fundamental solution] \label{p:fundamental}
Let $\Gamma$ be the fundamental solution of the Kolmogorov equation.
   \begin{enumerate}[label=(\roman*), nosep]
   \item \label{i:scale} \textsc{(Scale invariance)} If $\Gamma_1 (x,v)$ denotes $\Gamma(1,x,v)$, then $\Gamma(t,x,v) = t^{-2d} \Gamma_1 (t^{-\frac32} x, t^{-\frac12} v)$
     and $\int_{\R^{2d}} \Gamma_1 (x,v) \dx \dv =1$. 
     \item \label{i:equation} \textsc{(Solution of the equation)} For any $\tau >0$, the function $\Gamma$ is smooth in all variables in $[\tau, + \infty) \times \R^d \times \R^d$
       and its partial derivatives satisfy the Kolmogorov equation.
     \item \label{i:integrability}
       \textsc{(Integrability)} For all $T>0$, the function $\Gamma$ is in $L^p((0,T) \times \R^{2d})$ for any $p \in [1,1 + \frac{1}{2d})$ and it is in $L^{1+\frac1{2d},\weak}((0,T) \times \R^{2d})$.
       The functions $\Gammax$ and $\Gammav$ are in $L^q((0,T) \times \R^{2d})$ for any $q \in [1, 1 + \frac1{4d+1})$ and they are in $L^{1+\frac1{4d+1},\weak}((0,T) \times \R^{2d})$.
     \end{enumerate}
 \end{prop}
 In order to prove this proposition, we need the following lemma.
 \begin{lemma}[Integrability of a time-scaled function]\label{l:integ-weak}
   Let $F (t,x,v) = t^{-2d - {\beta_0}} G (t^{-3/2}x, t^{-1/2} v)$ with $G \in L^1 \cap L^\infty (\R^{2d})$.
   Then for all $T>0$, we have $F \in L^p ((0,T) \times \R^{2d})$ for $p \in [1, \frac{1+2d}{2d+{\beta_0}})$ and $F \in L^{\frac{1+2d}{2d+{\beta_0}},\weak} ((0,T) \times \R^{2d})$. 
 \end{lemma}
\begin{proof}
\begin{align*} \int_0^T \int_{\R^{2d}}F^p (t,x,v) \dt \dx \dv
  & = \int_0^T \frac1{t^{(2d+{\beta_0})p-2d}} \left\{
    \int_{\R^{2d}}G^p \left( \frac{x}{t^{3/2}},\frac{v}{t^{1/2}} \right)  t^{-2d} \dx \dv \right\} \dt \\
  & = \| G \|_{L^p (\R^{2d})}^p \int_0^T \frac1{t^{(2d+{\beta_0}) p-2d}} \dt.
\end{align*}
This integral is finite if and only if $(2d+{\beta_0})p-2d < 1$.

For the end point case, we write
\begin{align*} |\{ F (t,x,v) > \alpha \} \cap (0,T) \times \R^{2d}\}|
  &= |\{ G (t^{-3/2} x,t^{-1/2}v) > \alpha t^{2d+{\beta_0}} \} \cap (0,T) \times \R^{2d}\}|  \\
  &= \int_0^T |\{ G (t^{-3/2} x,t^{-1/2}v) > \alpha t^{2d+{\beta_0}} \}\}| \dt \\
  & = \int_0^T |\{ G (y,w) > \alpha t^{2d+{\beta_0}} \}  | t^{2d} \dt \\
  & = \int_0^\alpha |\{ G (y,w) > s \} | \left(\frac{s}\alpha\right)^{\frac{2d}{2d+{\beta_0}}} \frac1{2d+{\beta_0}} \left(\frac{s}\alpha \right)^{\frac1{2d+{\beta_0}}-1} \frac{\ds}\alpha \\
  & = \alpha^{-\frac{1+2d}{2d+{\beta_0}}} (2d+{\beta_0})^{-1} \int_0^\alpha |\{ G (y,w) > s \} |   s^{\frac1{2d+{\beta_0}}} \ds. 
\end{align*}
We conclude that $\sup_{\alpha >0} \alpha^{\frac{1+2d}{2d+{\beta_0}}} |\{ F (t,x,v) > \alpha \} \cap (0,T) \times \R^{2d}\}| < +\infty$. 
\end{proof}
\begin{proof}[Proof of Proposition~\ref{p:fundamental}]
The scale invariance stated in \ref{i:scale} is straightforward.
The constant appearing in the formula of $\Gamma$ (or equivalently, the fact that $\Gamma_1$ has mass $1$ in $(x,v)$) will be justified together with \ref{i:equation}. 
\medskip

Proof of \ref{i:equation}. We compute the Fourier transform in $(x,v)$ of the Schwartz function $\Gamma(t,\cdot,\cdot)$,
\begin{align*}
  g (t,\varphi,\xi) &= \int_{\R^{2d}} \Gamma (t,x,v) e^{-i 2 \pi \varphi \cdot x - i 2 \pi \xi \cdot v} \dx \dv \\
                    & = t^{-2d} \int_{\R^{2d}} \Gamma_1 (t^{-3/2} x,t^{-1/2} v) e^{-i 2\pi t^{3/2} \varphi \cdot t^{-3/2} x -  2\pi i t^{1/2} \xi \cdot t^{-1/2} v} \dx \dv \\
    & = \int_{\R^{2d}} \Gamma_1 (\bar x,\bar v) e^{-i (t^{3/2} \varphi) \cdot \bar x - i (t^{1/2} \xi) \cdot \bar v} \dd \bar x \dd \bar v.
\end{align*}
We can write
\[ g(t,\varphi,\xi) = \hat \Gamma_1 (t^{3/2} \varphi, t^{1/2} \xi).\]

We remark  that $\Gamma_1$ can be written as,
\[ \Gamma_1 (x,v) = C_\Gamma  e^{ - \frac12 C  (x,v) \cdot (x,v) }\]
with $C_\Gamma = (\frac{3}{4 \pi^2})^{d/2}$ and where the $2d \times 2d$ real symmetric matrix $C$, its determinant and its inverse are  given by
\[ C  = \begin{pmatrix} 6 I &  3  I \\  3  I & 2  I \end{pmatrix} \quad \text{ and } \quad  |\det C | = 3^d
  \quad \text{ and } \quad C^{-1}  = \begin{pmatrix} \frac{2}3  I &   -I  \\   -I & 2  I\end{pmatrix}.\]
In particular,
\begin{align*}
  \hat \Gamma_1 (\varphi,\xi)
  &= C_\Gamma \int_{\R^{2d}} e^{ - \frac12 C  (x,v) \cdot (x,v) } e^{- 2\pi i (x,v)  \cdot (\varphi,\xi)} \dx \dv \\
  & = \frac{(2 \pi)^d C_\Gamma}{|\det C|^{1/2}} \int_{\R^{2d}} (2\pi)^{-d} e^{ - \frac12   |(\bar x,\bar v)|^2 } e^{- 2\pi i (\bar x,\bar v)  \cdot \sqrt{C}^{-1}(\varphi,\xi)} \dx \dv .
\end{align*}
We now remark that $\frac{(2 \pi)^d C_\Gamma}{|\det C|^{1/2}}=1$ and $\hat \Gamma_1 (\varphi,\xi)$ coincides with the Fourier transform $\hat \mu$ of the Gaussian
$\mu = (2\pi)^{-d} e^{ - \frac12   |(\bar x,\bar v)|^2 }$ at $\sqrt{C}^{-1}(\varphi,\xi)$.
We know  that  $\hat \mu (\varphi,\xi) = e^{-\frac12 (|\varphi|^2 + |\xi|^2)}$ \cite[Eq.~(1.20)]{MR2768550}. We thus conclude that,
\[
  g(t,\varphi,\xi) = \hat \Gamma_1 (t^{3/2} \varphi, t^{1/2} \xi) = e^{-\frac12 C^{-1} (t^{3/2} \varphi, t^{1/2} \xi)\cdot (t^{3/2} \varphi, t^{1/2} \xi)} =
  e^{- \left( \frac13 t^3 |\varphi|^2 -  t^2 \varphi \cdot \xi +   t |\xi|^2\right) }.
\]
In particular,
\[ g(t,\varphi,\xi) = \exp \left(- \int_0^t ( s^2 |\varphi|^2 - 2 s \varphi \cdot \xi + |\xi|^2) \ds \right) = \exp \left(- \int_0^t  |s \varphi-\xi|^2 \ds \right) .\]  
We now check that $g$ satisfies the differential equation,
\[ \partial_t g + \varphi \cdot \nabla_\xi g  = - |\xi|^2 g .\]
In order to do so, we simply compute the right hand side as follows,
\[ \partial_t g + \varphi \cdot \nabla_\xi g = \left( - |t \varphi - \xi|^2 - \int_0^t 2 [(\xi - s \varphi) \cdot \varphi]  \ds \right) g .\]
Thanks to the inverse Fourier transform,
\[ \Gamma (t,x,v) = \int_{\R^{2d}} g (t,\varphi,\xi) e^{2\pi i (x,v) \cdot (\varphi,\xi)} \dd \varphi \dd \xi, \]
we finally deduce that $(\partial_t + v \cdot \nabla_x - \Delta_v) \Gamma =0$ in $(0,+\infty) \times \R^{2d}$. 
\medskip

Proof of \ref{i:integrability}. We only prove the integrability of $\Gamma$ since the cases of $\Gammax$ and $\Gammav$ are treated similarly.  Now Lemma~\ref{l:integ-weak} leads to the result. 
\end{proof}

The following lemma will be used in order to represent solutions by the fundamental solution.
\begin{lemma}[Kinetic convolution of $\check\Gamma$ with smooth functions] \label{l:gamma-varphi}
  Let $\varphi \in C^\infty_c (\R^{1+2d})$. Then
  \[
    \check \Gamma \astkin (\partial_t + v \cdot \nabla_x + \Delta_v) \varphi = -\varphi
  \]
  where $\check \Gamma (z) = \Gamma (z^{-1})$. 
\end{lemma}
\begin{proof}
Let us consider the adjoint of the Kolmogorov operator $\mathcal{K}^* \varphi = (\partial_t + v \cdot \nabla_x + \Delta_v) \varphi$. We then write for $z = (t,x,v)$,
\begin{align*}
  \left( \check \Gamma \astkin \mathcal{K}^* \varphi \right) (z)
  &= \int \check \Gamma (\zeta^{-1} \circ z) \mathcal{K}^* \varphi (\zeta) \dd \zeta \\
  &= \int  \Gamma (z^{-1} \circ \zeta) \mathcal{K}^* \varphi (\zeta) \dd \zeta \\
  &= \int_t^{+\infty} \iint_{\R^{2d}} \Gamma (s-t, y-x-(s-t) v, w-v) \mathcal{K}^* \varphi (s,y,w) \ds \dy \dw.
\end{align*}
The integral being singular at $s=t$, we consider $\eps >0$ and write 
\[   \left( \check \Gamma \astkin \mathcal{K}^* \varphi \right)(z)  = \mathcal{I}_\eps (z) + \mathcal{R}_\eps (z) \]
with
\begin{align*}
  \mathcal{I}_\eps (z) = \int_{t+\eps}^{+\infty} \iint_{\R^{2d}} \Gamma (s-t, y-x-(s-t) v, w-v) \mathcal{K}^* \varphi (s,y,w) \ds \dy \dw, \\
  \mathcal{R}_\eps (z) = \int_t^{t+\eps} \iint_{\R^{2d}} \Gamma (s-t, y-x-(s-t) v, w-v) \mathcal{K}^* \varphi (s,y,w) \ds \dy \dw.
\end{align*}

We first deal with $\mathcal{I}_\eps(z)$. We can integrate by parts in time $\partial_t \varphi$ and integrate by parts in space $v \cdot \nabla_x \varphi$ and in velocity $\Delta_v \varphi$.
If we let $\mathcal{K} \Gamma$ denote $(\partial_t + v \cdot \nabla_x - \Delta_v) \Gamma$, this leads to,
\begin{align*}
\mathcal{I}_\eps (z) &= \int_{t+\eps}^{+\infty} \iint_{\R^{2d}} \mathcal{K} \Gamma (s-t, y-x-(s-t) v, w-v)  \varphi (s,y,w) \ds \dy \dw \\
& - \int_{\R^{2d}} \Gamma (\eps, y-x - \eps v, w-v) \varphi (t+\eps,y,w) \dy \dw.
\end{align*}
We now use that $\mathcal{K} \Gamma =0$ in $(0,+\infty) \times \R^{2d}$ and $\Gamma (t,x,v) = t^{-2d} \Gamma_1 (t^{-3/2}x,t^{-1/2}v)$ in order to get,
\begin{align*}
  \mathcal{I}_\eps (z) &= - \int_{\R^{2d}} \Gamma (\eps, y-x - \eps v, w-v) \varphi (t+\eps,y,w) \dy \dw \\
                       &= - \int_{\R^{2d}} \Gamma (\eps, \bar y , \bar w) \varphi (t+\eps, x + \bar y  + \eps v,v + \bar w) \dy \dw \\
                       &= - \int_{\R^{2d}} \eps^{-2d} \Gamma_1 (\eps^{-3/2} \bar y , \eps^{-1/2} \bar w) \varphi (t+\eps, x+ \bar y  + \eps v,v + \bar w ) \dd \bar y \dd \bar w \\
  & = - \int_{\R^{2d}} \Gamma_1 ( \tilde y , \tilde w) \varphi (t+\eps,  x + \eps^{3/2} \tilde y  + \eps v,   v+\eps^{1/2} \tilde w) \dd \tilde y \dd \tilde w .
\end{align*}
The constant in $\Gamma$ is chosen such that $\iint_{\R^{2d}} \Gamma_1 (x,v) \dx \dv = 1$. This allows us to prove that
\[ \mathcal{I}_\eps (z) \to - \varphi (z) \quad \text{ as } \quad \eps \to 0.\]

We now prove that $\mathcal{R}_\eps (z)$ vanishes at the limit. 
\begin{align*}
\mathcal{R}_\eps (z) &= \int_0^{\eps} \iint_{\R^{2d}} \Gamma (\tau, y-x-\tau v, w-v) (\mathcal{K}^* \varphi) (t+\tau,y,w) \ds \dy \dw \\
                     & = \int_0^{\eps} \iint_{\R^{2d}} \Gamma (\tau, \bar y, \bar  w) (\mathcal{K}^* \varphi) (t+\tau,x+ \bar y + \tau v, v + \bar w) \ds \dd \bar y \dd \bar w \\
  & = \int_0^{\eps} \iint_{\R^{2d}} \Gamma_1 ( \tilde y, \tilde  w) (\mathcal{K}^* \varphi) (t+\tau,x+ \tau^{3/2} \tilde y + \tau v, v + \tau^{1/2} \tilde w) \ds \dd \tilde y \dd \tilde w.
\end{align*}
This implies that
\[ | \mathcal{R}_\eps (z)  | \le \|\mathcal{K}^* \varphi\|_{L^\infty (\R^{1+2d})} \eps .\]
In particular, $\mathcal{R}_\eps (z) \to 0$ as $\eps \to 0$. 
\end{proof}

\subsection{Uniqueness}

In order to obtain the representation of sub-solutions (including truncated solutions),
the following uniqueness result is used for the Kolmogorov equation in Lebesgue spaces. 
\begin{prop}[Uniqueness for Kolmogorov]\label{p:kolm-unique}
  Assume that a function $F \colon \R^{1+2d} \to \R$ is supported in $(0,+\infty) \times \R^{2d}$, that there exists $p \in (1,+\infty)$ such that for all $T>0$,  $F \in L^p ((-T,T) \times \R^{2d})$,  and that it is a distributional solution of $(\partial_t + v \cdot \nabla_x- \Delta_v) F = 0$ in $\R^{1+2d}$. Then $F \equiv 0$ in $\R^{1+2d}$.
\end{prop}
\begin{proof}
  The formal argument is very simple. It is enough to multiply the equation by $p |F|^{p-2}F$ and integrate in $(x,v)$ to get,
  \[ \frac{\dd}{\dt} \int_{\R^{2d}} |F|^p (t,x,v) \dx \dv = - p \int_{\R^{2d}} \nabla_v F \cdot \nabla_v |F|^{p-2} F \dx \dv \le 0.\]
  The non-negative function $t \mapsto \int_{\R^{2d}} |F|^p (t,x,v) \dx \dv$ lies in $L^1 ((-T,T))$, vanishes for $t<0$ and is  decreasing: it has to be $0$ for almost every $t$. 
\medskip
  
\noindent In what follows, given $F\in \R$ and $q >0$, we write $F^q$ for $|F|^{q-1}F$ in order to clarify computation. 

\medskip
\noindent \textsc{Difficulties and sketch of proof.}  To make the previous formal argument rigorous, we shall mollify $F$. We do it by using the kinetic convolution in order to simplify integration by parts.
  Since we may have $p-1 \in (0,1)$, the function $r \mapsto r^{p-1}$ is not Lipschitz and the meaning of  $\nabla_v F^{p-1}$ is to be made precise when $F$ vanishes.
  We circumvent this difficulty by mollifying by a function $\rho$ that is positive everywhere, so that $F^\eps$ is also positive and smooth everywhere. 
  Then we have to truncate the test function to ensure that it is compactly supported. We use two smooth functions $\varphi$ and $\psi$: we will first let $\psi$ converge
  to $1$; second, we will integrate by parts, rearrange terms and discard dissipation. Third,  we will  let $\varphi$  converge to a function only depending on time and
  conclude that the aforementioned $L^1$ function of the time variable is non-increasing. 
  
  \medskip
\noindent \textsc{Mollification and truncation.}
  It is convenient to use the fundamental solution $\Gamma$ at time $\eps$ as a mollifier, even if this is not mandatory. 
  Let $\rho_\eps = \Gamma (\eps,x,v)$   and let $F^\eps = F \astkin \rho_\eps$. This function $F^\eps$ is $C^\infty$, positive everywhere if $F$ is not identically equal to $0$,
  and converges to $F$ in $L^p ((-T,T) \times \R^{2d})$
  as $\eps \to 0$ for all $T>0$.  We then consider two smooth and compactly supported functions $\varphi, \psi \in C_c^\infty (\R^{1+2d})$  and define,
   \[ G_\psi =  \left[ \left[ (F^\eps)^{p-1} \varphi \right] \astkin \check{\rho}_\eps \right] \psi \]
   where $\check{\rho}_\eps (z) = \rho_\eps (z^{-1})$. 
   We are going to use that $\int_{\R^{1+2d}} ( f \astkin g) h = \int_{\R^{1+2d}} f ( h \astkin \check{g})$ (see Lemma~\ref{l:convol-L2}).
   Taking $G_\psi$ as a test function, we have
   \[ \int_{\R^{1+2d}} F (\partial_t + v \cdot \nabla_x ) G_\psi = - \int_{\R^{1+2d}} F \Delta_v G_\psi.\]

  \medskip
\noindent \textsc{Limit as $\psi \to 1$.}
   We now  write $G_\psi = G_1 \psi$ where $G_1$ denotes $G_\psi$ with $\psi \equiv 1$.
   For $\eps>0$ fixed,
   all functions in the weak formulation are smooth and we can pass to the limit by using dominated convergence in order to get,
   \[ \int_{\R^{1+2d}} F (\partial_t + v \cdot \nabla_x ) G_1 = - \int_{\R^{1+2d}}  F \Delta_v G_1.\]

  \medskip
\noindent \textsc{Integration by parts.}
   We now consider the smooth and compactly supported function,
   \[f^\eps = (F^\eps)^{p-1} \varphi \]
   that allows us to write $G_1 =  f^\eps \astkin \check{\rho}_\eps$. In particular,
   \[ (\partial_t + v \cdot \nabla_x ) G_1 = \left[ (\partial_t + v\cdot \nabla_x ) f^\eps \right] \astkin \check{\rho}_\eps
   \quad \text{ and } \quad \Delta_v G_1 = \left[ \Delta_v f^\eps \right] \astkin \check{\rho}_\eps.\]
   We use these formulas and  $\int_{\R^{1+2d}} ( f \astkin g) h = \int_{\R^{1+2d}} f ( h \astkin \check{g})$ from Lemma~\ref{l:convol-L2} and get,
   \[ \int_{\R^{1+2d}} F^\eps (\partial_t + v \cdot \nabla_x ) f^\eps = - \int_{\R^{1+2d}} F^\eps \Delta_v f^\eps = \int_{\R^{1+2d}} \nabla_v F^\eps \cdot \nabla_v f^\eps .\]   
   In view of the definition of $f^\eps$, this leads to
   \begin{multline*}
     \int_{\R^{1+2d}} \varphi F^\eps (\partial_t + v \cdot \nabla_x ) (F^\eps)^{p-1} + \int_{\R^{1+2d}}  (F^\eps)^p (\partial_t + v \cdot \nabla_x ) \varphi \\
     = \int_{\R^{1+2d}} \varphi \nabla_v F^\eps \cdot \nabla_v (F^\eps)^{p-1} + \int_{\R^{1+2d}}  (F^\eps)^{p-1} \nabla_v F^\eps \cdot \nabla_v \varphi .
   \end{multline*}
   Since $F^\eps$ is smooth and positive, we have that
   \begin{align*}
     F^\eps (\partial_t + v \cdot \nabla_x ) (F^\eps)^{p-1} & = p^{-1}  (\partial_t + v \cdot \nabla_x ) |F^\eps|^p, \\
   (F^\eps)^{p-1} \nabla_v F^\eps &= p^{-1} \nabla_v |F^\eps|^p.
   \end{align*}
   We then integrate by  parts and get,
   \begin{multline*}
     - p^{-1} \int_{\R^{1+2d}}  |F^\eps|^p (\partial_t + v \cdot \nabla_x ) \varphi  + \int_{\R^{1+2d}}  |F^\eps|^p (\partial_t + v \cdot \nabla_x ) \varphi \\
     = (p-1) \int_{\R^{1+2d}} \varphi |F^\eps|^{p-2} \nabla_v F^\eps \cdot \nabla_v F^\eps  - p^{-1} \int_{\R^{1+2d}}  |F^\eps|^p \Delta_v \varphi .
   \end{multline*}

   \medskip
\noindent \textsc{Conclusion.}
   We next remark that the first term in the right hand side is non-negative if $\varphi \ge 0$: it corresponds to dissipation. We thus get,
   \[
      (1 - p^{-1}) \int_{\R^{1+2d}}  |F^\eps|^p (\partial_t + v \cdot \nabla_x ) \varphi \ge  -  p^{-1} \int_{\R^{1+2d}}  |F^\eps|^p \Delta_v \varphi .
   \]
   We finally consider $\varphi (t,x,v) = \Phi (x,v) \Theta (t)$. Since $F^\eps \to F$ in $L^p (\R^{1+2d})$, we can pass to the limit $\Phi \to 1$ and get,
   \[
      \int_{\R} \left\{ \int_{\R^{2d}}  |F^\eps|^p \dx \dv \right\} \Theta' (t) \dt \ge  0.
   \]
   This means that the integrable function $t \mapsto \int_{\R^{2d}}  |F^\eps|^p (t,x,v)\dx \dv$ is non-increasing. It is thus equal to $0$ for almost every time $t$. 
This implies that $F^\eps =0$ as a function of $L^p((-T,T) \times \R^{2d})$ for all $T>0$. Letting $\eps \to 0$ allows us to conclude that $F \equiv 0$.  
 \end{proof}

\section{Weak solutions and the kinetic De Giorgi \& Nash's theorem}

In order to state the kinetic counterpart of De Giorgi \& Nash's theorem, the notion of weak solutions is first made precise.
The rest of the section is then devoted to the derivation of the local energy estimates. 

\subsection{Weak solutions of kinetic Fokker-Planck equations}

We recall that $I$ is a bounded interval of $\R$ of the form $(a,b]$ with $a,b \in \R$ and 
that $\Ox$ and $\Ov$ are two open sets of $\R^d$. 

\begin{defi}[Weak solutions]
  Let $S \in L^2 (\domain)$. A function $f \colon \domain \to \R$ is a \emph{weak solution} of $(\partial_t + v \cdot \nabla_x) f = \dive_v (A \nabla_v f) + B \cdot \nabla_v f +S$ in $\domain$ if  $f, \nabla_v f \in L^2 (\domain)$ and  for all $\varphi \in C^\infty_c (\domain)$,
  \[
    \int_{\domain} f  (\partial_t + v \cdot \nabla_x) \varphi  
    = \int_{\domain}  (A \nabla_v f \cdot \nabla_v \varphi) - \int_{\domain} (B \cdot \nabla_v f +S) \varphi.
  \]
\end{defi}
It is also useful to introduce the notion of weak sub-solutions (resp. weak super-solutions). Their a priori regularity coincides with the one of weak solutions, but the equation
is only satisfied (in the sense of distributions) with a non-negative (resp. non-positive)  measure in the right hand side. 
\begin{defi}[Weak sub/super-solutions] \label{defi:sub-super}
  Let $S \in L^2 (\domain)$. A function $f \colon \domain \to \R$ is a \emph{weak sub-solution} (resp.\emph{weak super-solution}) of $(\partial_t + v \cdot \nabla_x) f = \dive_v (A \nabla_v f) + B \cdot \nabla_v f +S$ in $\domain$ if  $f,\nabla_v f \in L^2 (\domain)$ and for all $\varphi \in C_c^\infty(\domain)$ with $\varphi \ge 0$ in $\domain$, 
  \begin{align*}
    \int_\domain f (\partial_t + v \cdot \nabla_x) \varphi &\ge \int_{\domain} A \nabla_v f \cdot \nabla_v \varphi - \int_\domain (B \cdot \nabla_v f +S) \varphi \\
    \bigg(\text{resp.} \quad      \int_\domain f (\partial_t + v \cdot \nabla_x) \varphi &\le \int_{\domain} A \nabla_v f \cdot \nabla_v \varphi - \int_\domain (B \cdot \nabla_v f +S) \varphi  \bigg).
\end{align*}
\end{defi}

\subsection{Kinetic De Giorgi \& Nash's theorem}

The goal of this chapter is to establish the  counterpart of De Giorgi \& Nash's theorem for kinetic Fokker-Planck equations by following the general strategy presented in Chapter~\ref{c:dg-big}. 
\begin{thm}[kinetic De Giorgi \& Nash]\label{t:dg-kinetic}
  Let  $A \in \mathcal{E}(\lambda,\Lambda)$ with $\Omega = Q_1$ and $\lambda, \Lambda >0$.
  There exist two universal constants $\alpha \in (0,1]$ and $\Cdg >0$ such that any
  weak solution $f$ of $(\partial_t + v \cdot \nabla_x )f = \dive_v (A \nabla_v f) + B \cdot \nabla_v f + S$ in $\domain$ with $S \in L^\infty (\domain)$ satisfies, 
  \[ \|f\|_{\Ckin^\alpha (Q_{1/2})} \le \Cdg \left( \|f\|_{L^2 (Q_1)} + \|S\|_{L^\infty (Q_1)} \right).\]
\end{thm}
\begin{remark}[Local maximum principle, intermediate value principle and expansion of positivity]
  Like in the elliptic and parabolic cases, the proof of this theorem consists in establishing first a local maximum principle and then
  proving that the oscillation of solutions improves by a universal factor when zooming in with a universal factor. This improvement of
  oscillation derives from the expansion of positivity of solutions. Such a property is established thanks to an intermediate value principle,
  that can be seen as the kinetic counterpart of the elementary intermediate value lemma~\ref{l:ivl-elliptic}. 
\end{remark}
\begin{remark}[Kinetic De Giorgi's class]
We state the theorem for weak solutions but we will prove it for functions that are in a De Giorgi's class of kinetic type. More precisely,
the local maximum principle  holds  for functions in a class $\mathrm{kDG}^+$ while expansion of positivity  holds  for functions in a
class $\mathrm{kDG}^-$. Finally, the kinetic class $\mathrm{kDG}$ is made of the intersection of these two classes $\mathrm{kDG}^\pm$. 
\end{remark}

\subsection{Local energy estimates}

We recall that $I$ is a bounded interval of $\R$ of the form $(a,b]$ with $a,b \in \R$ and 
that $\Ox$ and $\Ov$ are two open sets of $\R^d$.
\begin{prop}[Local energy estimates]\label{p:EE-kinetic}
  Let $f$ be a weak sub-solution of $(\partial_t + v \cdot \nabla_x) f = \dive_v (A \nabla_v f) + B \cdot \nabla_v f +S$ in $\domain$.
  There exists a set $\mathcal{N} \subset I$  of null measure such that for all $x_0,v_0 \in \R^d$ and  $0 < r_x < R_x$ and $0 < r_v < R_v$,
  for all $t_1,t_2 \in (t_0 -r^2, t_0] \setminus \mathcal{N}$ with $t_1 < t_2$ and all  $\kappa \in \R$,
  \begin{align*}
    \iint_{B_{r_x}(x_0) \times B_{r_v}(v_0)} & (f-\kappa)_+^2 (t_2,x,v) \dx \dv
     + \frac\lambda4 \int_{\Qint} | \nabla_v (f-\kappa)_+ |^2     \\
    & \le     \iint_{B_{r_x}(x_0) \times B_{r_v}(v_0)} (f-\kappa)_+^2 (t_1,x,v) \dx \dv     \\
    & + \left[ \frac{8 \Lambda}{(R_v-r_v)^2} + \frac{4R_v}{R_x-r_x}
      + \frac{\Lambda^2}{\lambda}\right] \int_{\Qext}   (f-\kappa)_+^2 + 2 \int_{\Qext}  S (f-\kappa)_+
  \end{align*}
  with $\Qint = [t_1,t_2] \times B_{r_x}(x_0) \times B_{r_v}(v_0)$ and  $\Qext = [t_1,t_2]\times B_{R_x}(x_0) \times B_{R_v}(v_0)$. 
\end{prop}
In order to derive the local energy estimates, we need to use $(u-\kappa)_+ \varphi$ as a test function in the definition of weak solutions, for any smooth and compactly
supported function $\varphi$. The following lemma allows us to do so.
\begin{lemma}[A non-smooth test function] \label{l:trans-trunc}
Let $f$ be a weak sub-solution of $(\partial_t + v \cdot \nabla_x) f = \dive_v (A \nabla_v f) + B \cdot \nabla_v f +S$ in $\domain$.
  Then for all $\varphi \in C^\infty_c (\domain)$, 
  \begin{align*}
  - \int_{\domain}  (f-\kappa)_+^2 (\partial_t + v \cdot \nabla_x) \varphi^2
  &+ \frac\lambda2 \int_{\domain}  |\nabla_v (f-\kappa)_+|^2 \varphi^2 \\
  &\le   \int_{\domain} \left(4\Lambda | \nabla_v  \varphi |^2 + 2 \frac{\Lambda^2}\lambda \varphi^2 \right) (f-\kappa)_+^2     + 2\int_{\domain}  |S |   (f-\kappa)_+  \varphi^2. 
\end{align*}
\end{lemma}
\begin{proof}
  We consider $\rho \in C^\infty_c (\R^{1+2d})$ with $\rho \ge 0$ and $\int_{\R^{1+2d}} \rho (z) \dz = 1$. Then for any $\eps \in (0,1)$, we define
  $\rho_\eps (z) := \eps^{-4d-2} \rho (\eps^{-1} z)$ where $\eps^{-1} z$ is a short hand notation for $\sigma_{\eps^{-1}} (z)$ (scaling operator, see page~\pageref{p:scaling}).
  We also define  $f^\eps := f \astkin \rho_\eps$.

  We next consider another smooth function $\gamma \in C^\infty_c (\R)$ with $\gamma \ge 0$ and $\int_\R \gamma (r) \dr =1$
  and $\supp \gamma \in [-1,0]$.
  Then for any  $\nu \in (0,1)$, we define the function $P_\nu (r):= (r-\kappa)_+^2 \ast \gamma_\nu$.

  We finally define $\check{\rho}_\eps (z) = \rho_\eps (z^{-1})$ and use $\psi = (P'_\nu(f^\eps) \varphi^2) \astkin \check{\rho}_\eps$ as a test function,
  \begin{align*}
-     \int_{\R^{1+2d}} P_\nu (f^\eps) (\partial_t + v \cdot \nabla_x) \varphi^2 & =  \int_{\R^{1+2d}} \bigg( (\partial_t + v \cdot \nabla_x) f^\eps \bigg) \bigg( P_\nu' (f^\eps) \varphi^2 \bigg) \\
    & = - \int_{\R^{1+2d}}    f (\partial_t + v \cdot \nabla_x) \left[  \bigg( P'_\nu (f^\eps) \varphi^2 \bigg) \astkin \check{\rho}_\eps  \right].
  \end{align*}
  We now use the fact that $f$ is a weak sub-solution in the sense of Definition~\ref{defi:sub-super} in order to get,
  \begin{align*}
-     \int_{\domain}  P_\nu (f^\eps) (\partial_t + v \cdot \nabla_x) \varphi^2 
   \le &   -\int_{\domain} \bigg( A \nabla_v f^\eps \cdot \nabla_v P_\nu'(f^\eps) \bigg)\varphi^2   - \int_{\domain} \bigg( A \nabla_v f^\eps \cdot  \nabla_v  \varphi^2 \bigg) P_\nu'(f^\eps)  \\
                               & + \int_{\domain} \big[ (B \cdot \nabla_v f) \astkin \rho_\eps \big] P'_\nu (f^\eps) \varphi^2 + \int_{\domain} (S \astkin \rho_\eps) P_\nu'(f^\eps) \varphi^2. 
  \end{align*}
Because $P_\nu$ is a convex function, we remark that the first term  in the right hand side is non-positive and we move it to the left,
  \begin{multline*}
-     \int_{\domain}  P_\nu (f^\eps) (\partial_t + v \cdot \nabla_x) \varphi^2 + \int_{\domain} P_\nu''(f^\eps) \bigg( A \nabla_v f^\eps \cdot \nabla_v f^\eps \bigg)\varphi^2
   \le       - 2 \int_{\domain} \bigg( \sqrt{A} \nabla_v f^\eps \cdot  \sqrt{A} \nabla_v  \varphi \bigg) P_\nu'(f^\eps) \varphi \\
         + \int_{\domain} \big[ (B \cdot \nabla_v f) \astkin \rho_\eps \big] P'_\nu (f^\eps) \varphi^2 + \int_{\domain} (S \astkin \rho_\eps) P_\nu'(f^\eps) \varphi^2. 
       \end{multline*}       
   In the two next lines, we use that $\supp \gamma \subset [-1,0]$ so that we have  $P''_\nu (r)\ge 2 \un_{r \ge \kappa}$ and we estimate $P'_\nu (r)$ from above in a rough way,
   \begin{align*}
  \forall r \ge \kappa, \quad   P''_\nu (r) &= 2 ( \un_{r \ge \kappa}) \ast \gamma_\nu = 2 \int_\kappa^{\infty} \gamma_\nu (r-s) \ds  = 2 \int_{-\infty}^{(r-\kappa)} \gamma_\nu (\tau) \dd \tau =2, \\
     \forall r \in \R, \quad   P'_\nu (r) &= 2 (r-\kappa)_+ \ast \gamma_\nu = 2 \int_0^{\infty} (s-\kappa)_+ \gamma_\nu (r-s) \ds \\
     & \le 2 \int_0^{\infty} ((r-\kappa)_+ + \nu) \gamma_\nu (t-s) \ds
                               = 2 ((r-\kappa)_+ +\nu). 
   \end{align*}
We can use these estimates to deduce from the weak formulation the following inequality,
\begin{align*}
    - \int_{\domain}  P_\nu (f^\eps) (\partial_t + v \cdot \nabla_x) \varphi^2 &+ 2\int_{\domain} \un_{\{f^\eps \ge \kappa\}} \bigg( A \nabla_v f^\eps \cdot \nabla_v f^\eps \bigg) \varphi^2 \\
  &\le   4 \int_{\domain} \bigg| \sqrt{A} \nabla_v f^\eps \cdot  \sqrt{A} \nabla_v  \varphi \bigg| \left[ (f^\eps-\kappa)_+ + \nu \right] |\varphi| \\
  &+ 2\int_{\domain} \bigg( \left| (B \cdot \nabla_v f) \astkin \rho_\eps \right| +  \left|S \astkin \rho_\eps\right| \bigg) \left[ (f^\eps-\kappa)_+ + \nu \right] \varphi^2. 
\end{align*}
We now use dominated convergence to pass to the limit in the left hand side as $\nu \to 0$  and we finally get,
\begin{align*}
    - \int_{\domain}  (f^\eps-\kappa)_+^2 (\partial_t + v \cdot \nabla_x) \varphi^2 &+ 2\int_{\domain} \un_{\{f^\eps \ge \kappa\}} \bigg( A \nabla_v f^\eps \cdot \nabla_v f^\eps \bigg) \varphi^2 \\
  &\le 4  \int_{\domain} \bigg| \sqrt{A} \nabla_v f^\eps \cdot  \sqrt{A} \nabla_v  \varphi \bigg| (f^\eps-\kappa)_+  |\varphi|\\
  &+ 2\int_{\domain} \bigg( \left| (B \cdot \nabla_v f) \astkin \rho_\eps \right| +  \left|S \astkin \rho_\eps\right| \bigg)  (f^\eps-\kappa)_+  \varphi^2. 
\end{align*}
Like in the elliptic and parabolic cases, we use next the fact that $\nabla_v (f^\eps-\kappa)_+ = \un_{\{ f^\eps \ge \kappa\}} \nabla_x f^\eps$
and $\un_{\{ f^\eps \ge \kappa\}} = \un_{\{ f^\eps \ge \kappa\}}\un_{\{ f^\eps \ge \kappa\}}$,
\begin{align*}
    - \int_{\domain}  (f^\eps-\kappa)_+^2 (\partial_t + v \cdot \nabla_x) \varphi^2 &+ 2\int_{\domain}  \bigg( A \nabla_v (f^\eps-\kappa)_+ \cdot \nabla_v (f^\eps-\kappa)_+ \bigg) \varphi^2 \\
  &\le 4  \int_{\domain} \bigg| \sqrt{A} \nabla_v (f^\eps-\kappa)_+ \cdot  \sqrt{A} \nabla_v  \varphi \bigg| (f^\eps-\kappa)_+  |\varphi| \\
  &+ 2 \int_{\domain} \bigg( \left| (B \cdot \nabla_v (f-\kappa)_+) \astkin \rho_\eps \right| +  \left|S \astkin \rho_\eps\right| \bigg)  (f^\eps-\kappa)_+  \varphi^2. 
\end{align*}
We continue following what we did in the parabolic case by using Cauchy-Schwarz's inequality to get rid of the first error term in the right hand side,
\begin{align*}
    - \int_{\domain}  (f^\eps-\kappa)_+^2 (\partial_t + v \cdot \nabla_x) \varphi^2 &+  \int_{\domain}  \bigg( A \nabla_v (f^\eps-\kappa)_+ \cdot \nabla_v (f^\eps-\kappa)_+ \bigg) \varphi^2 \\
  &\le   4 \int_{\domain} \bigg| \sqrt{A} \nabla_v  \varphi \bigg|^2 (f^\eps-\kappa)_+^2  \\
  &+ 2\int_{\domain} \bigg( \left| (B \cdot \nabla_v (f-\kappa)_+) \astkin \rho_\eps \right| +  \left|S \astkin \rho_\eps\right| \bigg)  (f^\eps-\kappa)_+  \varphi^2. 
\end{align*}
We finally use ellipticity constants and $|B|\le \Lambda$ and get,
\begin{align*}
  - \int_{\domain}  (f^\eps-\kappa)_+^2 (\partial_t + v \cdot \nabla_x) \varphi^2
  &+ \lambda \int_{\domain}  |\nabla_v (f^\eps-\kappa)_+|^2 \varphi^2 \\
  &\le  4\Lambda \int_{\domain} | \nabla_v  \varphi |^2 (f^\eps-\kappa)_+^2 \\
&    +2 \int_{\domain} \bigg( \Lambda \left| \nabla_v (f-\kappa)_+ \right| \astkin \rho_\eps +  \left|S \astkin \rho_\eps\right| \bigg)  (f^\eps-\kappa)_+  \varphi^2. 
\end{align*}
We can now use the fact that $(f^\eps -\kappa)_+^2$ converges to $(f-\kappa)_+^2$ in $L^2(I \times \Ox , H^1 (\Ov))$ and obtain,
  \begin{align*}
  - \int_{\domain}  (f-\kappa)_+^2 (\partial_t + v \cdot \nabla_x) \varphi^2
  &+ \lambda \int_{\domain}  |\nabla_v (f-\kappa)_+|^2 \varphi^2 \\
  &\le  4\Lambda \int_{\domain} | \nabla_v  \varphi |^2 (f-\kappa)_+^2
    + 2 \int_{\domain} \bigg( \Lambda \left| \nabla_v (f-\kappa)_+ \right|  +  |S | \bigg)  (f-\kappa)_+  \varphi^2. 
\end{align*}
We conclude by using one more time Cauchy-Schwarz's inequality. 
\end{proof}
With this technical lemma in hand, we can establish the local energy estimates.
\begin{proof}[Proof of Proposition~\ref{p:EE-kinetic}]   
The proof proceeds in two steps. We first localize the estimate in $(x,v)$, we then get the result for rational $r,R,\kappa$'s and conclude by a monotonicity argument. 
\medskip

\noindent \textsc{Step~1.}
  Let $\rho$ be the truncation function from Lemma~\ref{l:trunc} that is supported in $B_R $ and equal $1$ in $B_r$.
  Similarly, we consider $\bar \rho$ be the truncation function from Lemma~\ref{l:trunc} that is supported in $B_{R^3} $ and equal $1$ in $B_{r^3} $. In particular, $|\nabla_x \bar \rho|\le 2 (R^3-r^3)^{-1}$.
  Given $ t_1,t_2 \in I$ with $t_1 < t_2$,   we consider a 1D mollifier $\theta \colon \R \to \R$ (smooth and of unit mass) that is supported in $[-1,0]$, and
  the smooth function $\Theta_\eps \colon (0,+\infty) \to \R$ such that $\Theta_\eps (0) = 0$ and for all $t \in \R$, we have $\Theta_\eps'(t) = \theta_\eps (t-t_1) - \theta_\eps (t-t_2)$.
  We choose $\eps>0$ small enough so that $t_1-\eps \in I$. 
  
  We now use Lemma~\ref{l:trans-trunc} with $\varphi (t,x,v) = \Theta_\eps^2 (t) \bar \rho^2 (x) \rho^2 (v) = (\Theta_\eps \otimes \bar \rho \otimes \rho)^2 (t,x,v)$. Since
  \[\partial_t \varphi (t,x,v) = 2 (-\theta_\eps (t-t_2) +\theta_\eps (t-t_1)) \Theta_\eps (t) \bar \rho^2 (x) \rho^2 (v), \]
  we can rearrange terms and get,
  \begin{align*}
   2 \int_{\domain}  &(f-\kappa)_+^2 (\bar \rho \otimes \rho)^2  \Theta_\eps (t) \theta_\eps (t-t_2)
  + \frac\lambda2 \int_{\domain}  |\nabla_v (f-\kappa)_+|^2 (\Theta_\eps \otimes \bar \rho \otimes \rho)^2   \\
                   &\le 2 \int_{\domain} (f-\kappa)_+^2 (\bar \rho \otimes \rho)^2  \Theta_\eps (t) \theta_\eps (t-t_1) \\
                   &  +   2 \int_{\domain}  (f-\kappa)_+^2 (v \cdot \nabla_x \bar \rho (x)) \bar \rho(x)  (\Theta_\eps \otimes \rho)^2 \\
    & +  \int_{\domain} \left(4 \Lambda | \nabla_v  \rho |^2 + 2 \frac{\Lambda^2}\lambda  \rho^2 \right) (\Theta_\eps\otimes \bar \rho)^2 (f-\kappa)_+^2 
       + 2 \int_{\domain}  |S |   (f-\kappa)_+  ( \Theta_\eps \otimes \bar \rho \otimes \rho)^2. 
\end{align*}
  We use next that $\bar \rho \otimes \rho =1$ in $B_{r^3} \times B_r$ and $|\nabla_v \rho| \le 2(R_v-r_v)^{-1}$ and $ |\nabla_x \bar \rho| \le 2(R_x-r_x)^{-1}$ and $\bar \rho \otimes \rho  \le \un_{B_{R^3} \times B_R}$ and get
  \begin{align*}
   \int_{I \times B_{r^3} \times B_r}  &(f-\kappa)_+^2 \Theta_\eps (t) \theta_\eps (t-t_2)
  + \frac\lambda4 \int_{I \times B_{r^3} \times B_r}  |\nabla_v (f-\kappa)_+|^2\Theta_\eps^2  \\
  &\le \int_{\domain} (f-\kappa)_+^2  (\bar \rho \otimes \rho)^2 \Theta_\eps (t) \theta_\eps (t-t_1)   + \int_{I \times B_{R^3} \times B_R}  |S |   (f-\kappa)_+   \Theta_\eps^2\\
    & +  \int_{I \times B_{R^3} \times B_R} \left(\frac{8 \Lambda}{(R_v-r_v)^2} + \frac{4R_v}{R_x-r_x} +  \frac{\Lambda^2}\lambda  \right) (f-\kappa)_+^2 \Theta_\eps^2. 
\end{align*}
We can now pass to the limit $\bar \rho \otimes \rho \to \un_{B_{r^3}} \otimes \un_{B_r}$ by dominated convergence and get,  
  \begin{align*}
   \int_{I \times B_{r^3} \times B_r}  &(f-\kappa)_+^2 \Theta_\eps (t) \theta_\eps (t-t_2)
  + \frac\lambda4 \int_{I \times B_{r^3} \times B_r}  |\nabla_v (f-\kappa)_+|^2 \Theta_\eps^2  \\
  &\le \int_{I \times B_{r^3} \times B_r} (f-\kappa)_+^2   \Theta_\eps (t)  \theta_\eps (t-t_1)        + \int_{I \times B_{R^3} \times B_R}  |S |   (f-\kappa)_+   \Theta_\eps^2  \\
    & +  \int_{I \times B_{R^3} \times B_R} \left(\frac{8 \Lambda}{(R_v-r_v)^2} + \frac{4R_v}{R_x-r_x}+  \frac{\Lambda^2}\lambda  \right) (f-\kappa)_+^2 \Theta_\eps^2. 
\end{align*}
\medskip

\noindent \textsc{Step~2.}
For $r,R,\kappa$ fixed, we consider Lebesgue points of the $L^1$ function
\[t \mapsto \int_{B_{r^3} \times B_r} (f(t,x,v)-\kappa)_+^2 \dx \dv.\]
This provides a set $\mathcal{N}_{r,R,\kappa}\subset I$ of null measure such that for all $t_1,t_2 \in I \setminus \mathcal{N}_{r,R,\kappa}$,
the announced inequality holds (after passing to the limit as $\eps \to 0$).

We then consider the set of null measure $\mathcal{N}$ corresponding to rational $r,R, \kappa$'s.
For other $r,R,\kappa$, we remark that all integrals are non-increasing in $\kappa$ and non-decreasing in $r$ and $R$. 
We thus consider an increasing sequence $\kappa_n$ of rational numbers and decreasing sequences $r_n,R_n$ or rational numbers,
write the corresponding inequality, and pass to the limit thanks to the monotone convergence theorem. 
\end{proof}

\section{Gain of integrability of sub-solutions}
\label{s:x-reg}

This section is devoted to \emph{the} key estimate yielding the local maximum principle for sub-solutions: their gain of integrability.
They are various ways to establish this property. We choose to follow A.~Pascucci and S.~Polidoro \cite{pp} by representing sub-solutions with
the help of the fundamental solution of the Kolmogorov equation. 
 
\subsection{Representation of sub-solutions}

The goal of this subsection is to prove that truncated solutions $(f-\kappa)_+$, and more generally weak sub-solutions,
can be represented thanks to the kernel of the Kolmogorov equation.
The definition of weak sub-solutions is given on page~\pageref{defi:sub-super}
and the reader is reminded that $M_1^+ (\domain)$ denotes the set of non-negative  measures on the set $\domain$. 

\begin{prop}[Representation of sub-solutions] \label{p:truncated}
  Let $S \in L^2 (\domain)$ and $f \colon \domain \to \R$ be a weak sub-solution of $(\partial_t + v \cdot \nabla_x) f = \dive_v (A \nabla_v f) + B \cdot \nabla_v f+ S$
  in $\domain$ and let $\phitrunc$ be $C^\infty$ and compactly supported in  $\domain$. Then for all $\kappa \in \R$, 
  \[
     f \phitrunc  = (\Gammax+  \Gammav) \astkin \Ttrunc + \Gamma \astkin (\Strunc-\mtrunc) 
  \]
  with    $\mtrunc =  \phitrunc \mathfrak{m}$ for some $\mathfrak{m} \in M_1^+ (\R^{1+2d})$ and $\Ttrunc, \Strunc \in L^2 (\R^{1+2d})$ are given
  by,
  \begin{align*}
    \Ttrunc  =& \phitrunc (A-I) \nabla_v  f  , \\
    \Strunc  =& (B \phitrunc  - (A+I) \nabla_v \phitrunc )\cdot \nabla_v f  +  S \phitrunc  
     +   f (\partial_t  + v \cdot \nabla_x - \Delta_v) \phitrunc  .
  \end{align*}
\end{prop}
We write an equation for the product of a local sub-solution with a cut-off function. 
\begin{lemma}[Localization of solutions]\label{l:truncation}
  Let $f$ be a weak sub-solution of $(\partial_t + v \cdot \nabla_x) f = \dive_v (A \nabla_v f) +  S$ in $\domain$ with $S \in L^1 (\domain)$
  and $\mathfrak{m} \in M_1^+ (\domain)$ be given by
  \[ \mathfrak{m} = \dive_v (A \nabla_v f) +  S - (\partial_t + v \cdot \nabla_x) f  .\]
   Let $Q$ be a kinetic cylinder $Q_R (z_0)$
  contained in $\domain$ and let $\phitrunc \colon \R^{1+2d} \to \R$ be $C^\infty$ (in all variables)
  and compactly supported in $Q$. Then the function $\ftrunc = f \phitrunc$ satisfies
  \[ (\partial_t + v \cdot \nabla_x -\Delta_v) \ftrunc = \dive_v \Ttrunc + \Strunc - \mtrunc \quad \text{ in } \R \times \R^d \times \R^d\]
  (in the sense of distributions) where  \( \mtrunc =  \phitrunc \mathfrak{m}\) and
  \begin{align*}
    \Strunc &= S \phitrunc - (A+I) \nabla_v f  \cdot \nabla_v \phitrunc + f (\partial_t + v \cdot \nabla_x -\Delta_v ) \phitrunc, \\
    \Ttrunc &= \phitrunc (A-I) \nabla_v f.
  \end{align*}
\end{lemma}
\begin{proof}
Given $\varphi \in C^\infty_c (\R^{1+2d})$, we use $\varphi \phitrunc$ as a test function for the equation satisfied by $f$ and reach the desired conclusion for $\ftrunc$. 
\end{proof}
\begin{lemma}[Representation formula]\label{l:reg-rep}
  Let $\mathfrak{S}, S \in L^2 (\R^{1+2d})$ and $\mathfrak{m} \in M_1^+ (\R^{1+2d})$. We assume that $\mathfrak{S},S$ and $\mathfrak{m}$ are compactly supported   in $(0, +\infty) \times \R^{2d}$. Then for all $T >0$, the function 
  \[ F = \Gamma \astkin (\dive_v \mathfrak{S} + S-\mathfrak{m})\] 
  is supported in $(0,+\infty) \times \R^{2d}$, it is in $L^p ((-T,T) \times \R^{2d})$ for any $p \in (1,1+\frac1{2d})$ and in $L^{1+\frac1{2d},\weak} ((-T,T) \times R^{2d})$.

  Moreover, \( (\partial_t + v \cdot \nabla_x ) F = \Delta_v F + \dive_v \mathfrak{S} + S - \mathfrak{m}\) in $\R^{1+2d}$ (in the sense of distributions). 
\end{lemma}
\begin{proof}
  The function $F$ can be written as follows,
  \[ F = (\Gammax + \Gammav) \astkin \mathfrak{S} + \Gamma \astkin (S-\mathfrak{m})\]
  Thanks to Young's inequality (Lemma~\ref{l:young-kin}) and the integrability properties of the functions $\Gamma,\Gammax,\Gammav$ established in Proposition~\ref{p:fundamental}-\ref{i:integrability},
  we deduce that: for all $p \in (1, 1 + \frac1{2d})$, all $\sigma \in [1,2 + \frac{4}{2d-1})$ and all $\tau \in [1,2+\frac1d)$,
\[ \Gamma \astkin \mathfrak{m} \in L^p ((-T,T) \times \R^{2d}), \; \Gamma \astkin S \in L^\sigma ((-T,T) \times \R^{2d}), \; (\Gammax + \Gammav) \astkin \mathfrak{S} \in L^\tau ((-T,T) \times \R^{2d})  .\]
In particular, $F \in L^p ((-T,T) \times \R^{2d})$ for all $p \in [1, 1+ \frac1{2d})$ and all $T>0$.
The end point case is treated similarly since the convolution of a finite Radon measure with a function in $L^{p,\weak}$ is in $L^{p,\weak}$. 

Moreover, 
    \[
    (\partial_t + v \cdot \nabla_x - \Delta_v) F = \left( (\partial_t + v \cdot \nabla_x - \Delta_v) \Gamma \right) \astkin (\dive_v \mathfrak{S} + S-\mathfrak{m}) 
         =\dive_v \mathfrak{S} + S-\mathfrak{m} \text{ in } \R^{1+2d}
    \]
    in the sense of distributions. Indeed, consider a smooth and compactly supported test function $\varphi \in C_c^\infty$ and write,
    \begin{align*}
      &\int_{\R^{1+2d}} F  (\partial_t + v \cdot \nabla_x + \Delta_v) \varphi \\
                         & = \int_{\R^{1+2d}} \bigg[ (\Gammax + \Gammav) \astkin \mathfrak{S} + \Gamma \astkin (S - \mathfrak{m}) \bigg] (\partial_t + v \cdot \nabla_x + \Delta_v) \varphi \\
      & = \int_{\R^{1+2d}}   \mathfrak{S} \bigg[ (\check\Gammax + \check\Gammav) \astkin (\partial_t + v \cdot \nabla_x + \Delta_v) \varphi \bigg] + \int_{\R^{1+2d}} (S - \mathfrak{m}) \bigg[ \check\Gamma  \astkin (\partial_t + v \cdot \nabla_x + \Delta_v) \varphi \bigg].
    \end{align*}
We conclude thanks to Lemma~\ref{l:gamma-varphi}.
  \end{proof}

We can now prove the representation formula.
\begin{proof}[Proof of Proposition~\ref{p:truncated} (representation of sub-solutions)]
Without loss of generality, we can assume that $\Omega \subset (0,+\infty) \times \R^{2d}$. 
By Lemmas \ref{l:truncation}, we know that $f \phitrunc$ is a distributional solution of
  \begin{equation} \label{e:pp}
    (\partial_t + v \cdot \nabla_x - \Delta_v) g = \dive_v \Ttrunc + \Strunc - \mtrunc \qquad \text{ in } \quad \R^{1+2d}
  \end{equation}
with $\Ttrunc,\Strunc,\mtrunc$ given in the statement of the proposition. 
Now Lemma~\ref{l:reg-rep} ensures that $F=\Gamma \astkin (\dive_v \Ttrunc + \Strunc - \mtrunc)$ is $L^p ((-T,T) \times \R^{2d})$ for all $p \in (1,1+\frac1{2d})$ and all $T>0$,
and  is a distributional solution of \eqref{e:pp} in $\R^{1+2d}$.
In particular, the function $ h = f \phitrunc - F$ lies in $L^p ((-T,T) \times \R^{2d})$ for all $T>0$,
it is a distributional solution of $(\partial_t + v \cdot \nabla_x - \Delta_v) h =0$ in $\R^{1+2d}$, and it is supported in $(0,+\infty) \times \R^{2d}$.
  Uniqueness for the Kolmomgorov equation (Proposition~\ref{p:kolm-unique}) implies that $h =0$, leading to the representation formula. 
\end{proof}

\subsection{Gain of integrability for sub-solutions}

We are now ready to prove that sub-solutions are locally better than square integrable. 
\begin{prop}[Gain of integrability of sub-solutions]\label{p:subsol-integrability}
  Let $f$ be a weak sub-solution of $(\partial_t + v \cdot \nabla_x ) f = \dive_v (A \nabla_v f) + B \cdot \nabla_v f +S$
  in $\domain$ with $A \in \mathcal{E}(\lambda,\Lambda)$ and $B \in L^\infty (\domain)$ and $S \in L^2 (\domain)$.
  For all $p_c \in (2, 2+ \frac1d)$, there exists $\Cc$ (only depending on $\Lambda,d$ and $p_c$) such that for all $Q_r (z_0) \subset Q_R (z_0) \subset \domain$,
  \begin{multline*}
    \| (f-\kappa)_+ \|_{L^{p_c} (Q_r (z_0))} 
    \le C_\tau \bigg( \frac{ \Lambda +1}{ R-r} \| \nabla_v (f-\kappa)_+\|_{L^2 (Q_R (z_0))} \\
    +   \frac{1}{R^2-r^2}  \| (f-\kappa)_+ \|_{L^2 (Q_R (z_0))} + \|S \un_{\{ f \ge \kappa\}}\|_{L^2 (Q_R (z_0))} \bigg).
  \end{multline*}
\end{prop}

This gain of integrability is in particular true for truncated solutions $(f-\kappa)_+$. 
\begin{lemma}[Truncated sub-solutions are sub-solutions] \label{l:eq-trunc}
  Let $S \in L^2 (\domain)$ and $f \colon \domain \to \R$ be a weak sub-solution of $(\partial_t + v \cdot \nabla_x) f = \dive_v (A \nabla_v f) + B \cdot \nabla_v f + S$
  in $\domain$ and $\kappa \in \R$. Then the truncated function $(f-\kappa)_+$ satisfies 
  \[ (\partial_t + v \cdot \nabla_v) (f-\kappa)_+ = \dive_v (A \nabla_v (f-\kappa)_+)  + S^\kappa_+ -\mathfrak{m}^\kappa_+ \quad \text{ in }  \domain\]
(in the sense of distributions)  with \( S_+^\kappa = B \cdot \nabla_v (f-\kappa)_+ +  S \un_{\{f > \kappa\}} \) and some non-negative measure $\mathfrak{m}^\kappa_+ \in M_1^+ (\domain)$.
\end{lemma}
\begin{proof}
  Let us first construct a $1$-Lipschitz approximation of the function $r_+$. 
  Given $\eps \in (0,1)$, we consider the function $P^\eps_+ (r)$ such that $(P^\eps_+)'' (r) = \theta_\eps (r)$ for some mollifier $\theta \colon \R \to [0,+\infty)$
  and $P^\eps_+ (0) = (P^\eps_+)'(0) = 0$.

  Arguing as in the proof of Lemma~\ref{l:trans-trunc} with $P^\eps_+$ instead of $r_+^2/2$, we get,
  \begin{align*}
 - \int_{\domain}  & P^\eps_+(f) (\partial_t + v \cdot \nabla_x) \varphi \, \dt \dx \dv \\
\le & \iint_{I \times \Ox} \langle  (\partial_t + v \cdot \nabla_x) f, (P^\eps_+)'(f) \varphi \rangle_{H^{-1},H^1_0} \dt \dx  \\
                     = &  -\int_{\domain} (P^\eps_+)''(f) \bigg[ A \nabla_v  f \cdot  \nabla_v f \bigg] \varphi
                       - \int_{\domain} \bigg[ A \nabla_v  f \cdot \nabla_v \varphi \bigg]  (P^\eps_+)'(f) \\
&  + \iiint_\domain B \cdot \nabla_v f (P^\eps_+)'(f) + \int_{\domain} S  (P^\eps_+)'(f)  \varphi \\
    \le &  - \int_{\domain} \bigg[ A \nabla_v  f \cdot \nabla_v \varphi \bigg] (P^\eps_+)'(f) + \iiint_\domain (B \cdot \nabla_v f) (P^\eps_+)'(f)  + \int_{\domain} S  (P^\eps_+)'(f)  \varphi  .
  \end{align*}
Passing to the limit as $\eps \to 0$ by dominated convergence yields
  \begin{multline*}
 - \int_{\domain}   (f-\kappa)_+ (\partial_t + v \cdot \nabla_x) \varphi \, \dt \dx \dv \\
 \le - \int_{\domain} \bigg[ A \nabla_v  f \cdot \nabla_v \varphi \bigg] \un_{\{f \ge \kappa \}} + \iiint_\domain (B \cdot \nabla_v f) \un_{\{ f \ge \kappa\}}
 + \int_{\domain} S  \un_{\{f \ge \kappa \}}  \varphi  .
   \end{multline*}
   We conclude by using that $\nabla_v (f-\kappa)_+ = \un_{\{f \ge \kappa\}} \nabla_v f$ (Proposition~\ref{p:composition}) and by recalling that
   a non-negative distribution is a Radon measure. Such a fact is an easy consequence of Riesz's representation theorem \cite[Theorem 1.25,~p.~39]{MR3409135}. 
 \end{proof}

The proof of the gain of integrability uses the (easy) construction of a cut-off function. 
\begin{lemma}[Cut-off function]\label{l:cutoff}
  For  $r,R>0$ such that $r < R$, there exists a smooth function  $\phitrunc$
  such that $\phitrunc =1$ in $Q_r$ and $\phitrunc=0$ outside $Q_R$. Moreover,
  \[ | (\partial_t +  v \cdot \nabla_x ) \phitrunc |\le 4  (R^2-r^2)^{-1}  \quad \text{ and } \quad |\nabla_v \phitrunc| \le 2 (R-r)^{-1}
    \quad \text{ and } \quad |\Delta_v \phitrunc| \le 2 d (R-r)^{-2}. \]
\end{lemma}
\begin{proof}
 On the one hand,
  we consider $\bar \rho (x)$ supported in $B_{R^3}$ and equal to $1$ in $B_{r^3}$ whose gradient satisfies $|\nabla \bar \rho | \le 2 (R^3 -r^3)^{-1}$. 
  On the other hand, we consider $\rho (v)$ supported in $B_R$ and equal to $1$ in $B_r$, whose gradient satisfies $|\nabla \rho | \le 2 (R-r)^{-1}$.
  As far as the time variable is concerned, we simply take $\Theta$ non-decreasing, equal to $1$ in $(-r^2,0]$ and vanishes in $(-\infty, -R^2]$. Its derivative satisfies
  $0 \le \Theta' \le 2 (R^2 -r^2)^{-1}$.   Then we consider $\phitrunc (z) = (\Theta \otimes \bar \rho \otimes \rho) (z)$.
  \begin{align*}
    | (\partial_t +  v \cdot \nabla_x ) \phitrunc |& \le 2  (R^2-r^2)^{-1} +  (2 R)  (R^3-r^3)^{-1}  \\
    & \le 4  (R^2-r^2)^{-1} 
  \end{align*}
  The bound on $\nabla_v \phitrunc$ corresponds to the bound on $\nabla \rho$. 
\end{proof}  

We can now prove the gain of integrability for sub-solutions.
\begin{proof}[Proof of Proposition~\ref{p:subsol-integrability}]  
  We first deal with the case $z_0=0$.   Let $\phitrunc$ be given by Lemma~\ref{l:cutoff}.  We use the representation formula from Proposition~\ref{p:truncated} 
  and use $\mtrunc \ge 0$ to get,
  \[ 0 \le f \phitrunc \le (\Gammax+  \Gammav) \astkin \Ttrunc + \Gamma \astkin \Strunc.\]
  We now can repeat the computations from the proof of Lemma~\ref{l:reg-rep} and get that for all $\tau \in (1,2+ \frac1d)$, there exists $C_\tau$ (only depending on dimension and $\tau$),
  \[ \| f \phitrunc\|_{L^\tau (\R^{1+2d})} \le C_\tau \left( \|\Ttrunc\|_{L^2 (\R^{1+2d})} + \|\Strunc\|_{L^2 (\R^{1+2d})} \right)\]
  where  $\Ttrunc$ and $\Strunc$ are given by the following formulas, see Proposition~\ref{p:truncated},
  \begin{align*}
    \Ttrunc = &  \phitrunc (A-I) \nabla_v f , \\
    \Strunc = &  ( B \phitrunc - (A+I) \nabla_v \phitrunc) \cdot \nabla_v f   + S \un_{\{ f \ge \kappa\}}\phitrunc
      + f (\partial_t + v \cdot \nabla_x-\Delta_v) \phitrunc  .
  \end{align*}
  \begin{align*}
    \| \Ttrunc\|_{L^2 (\R^{1+2d})} \le &(\Lambda +1 ) \| \nabla_v f\|_{L^2 (Q_R)} \\
    \| \Strunc\|_{L^2 (\R^{1+2d})} \le &( \Lambda + 2 (\Lambda+1) (R-r)^{-1}) \|\nabla_v f\|_{L^2 (Q_R)} + \| S \un_{\{ f \ge \kappa\}} \|_{L^2 (Q_R)} \\
     &   + (4+2d) (R-r)^{-2} \|f \|_{L^2 (Q_R)}.
  \end{align*}
  Combining the estimates  leads to
  \begin{multline*}
    \| f \|_{L^\tau (Q_r (z_0))} 
    \le C_\tau \bigg( \bigg[2 \Lambda +1   +2 (\Lambda+1) (R-r)^{-1} \bigg]\| \nabla_v f\|_{L^2 (Q_R (z_0))} \\
     + (4+2d) (R^2-r^2)^{-1} \| f \|_{L^2 (Q_R (z_0))} + \|S \un_{\{ f \ge \kappa\}}\|_{L^2 (Q_R (z_0))} \bigg).
  \end{multline*}

    Now we reduce to the case $z_0 =0$ by considering $g(z) = f (z_0 \circ z)$ and by applying the previous reasoning to $g$.
\end{proof}

\section{The kinetic De Giorgi's class kDG\textsuperscript{+} and the  maximum principle}

In this section, we introduce the kinetic De Giorgi's class $\mathrm{kDG}^+$ that ensures that the local maximum principle holds.
This class contains in particular all sub-solutions of kinetic Fokker-Planck equations. It is made of functions satisfying
the gain of integrability (from $2$ to $p_c >2$)  in the three variables $(t,x,v)$. 

\subsection{The kinetic De Giorgi's class kDG\textsuperscript{+}}

\begin{defi}[The kinetic De Giorgi's class $\mathrm{kDG}^+$]\label{d:kDG+}
  Let  $I = (a,b]$ with $a,b \in \R$ and $\Ox, \Ov$ be open sets of $\R^d$ and $S \in L^2 (\domain)$.
  
  A function $f \colon \domain \to \R$ lies in the \emph{kinetic De Giorgi's class} $\mathrm{kDG}^+(\domain,S)$ if 
  $f \in  L^\infty (I, L^2 (\Ox \times \Ov))$ and $\nabla_v f \in L^2 (\domain)$ and there exist $p_c >2$ and  $\Ckdgp \ge 1$ such that for all
  for all $z_0  \in \domain$, all $\kappa \in \R$ and all $r,R>0$ such that  $r< R <1$ and $Q_R (z_0) \subset \domain$, 
  \begin{equation} \label{e:pdgk}
    \begin{aligned}
      \| (f-\kappa)_+ \|^2_{L^{p_c} (Q_r (z_0))}  &   \le \Ckdgp (R-r)^{-4}  \int_{Q_R (z_0)} (f-\kappa)_+^2 \\
      & + \Ckdgp (R-r)^{-2} \int_{Q_R (z_0)} |S|^2 \un_{\{u \ge \kappa\}} .
    \end{aligned}
  \end{equation}
\end{defi}
\begin{remark}[Restriction to small radii]
We restrict ourselves to radii $r,R \in (0,1)$ to get cleaner formulas. We do not lose generality by doing so because we can always reduce to this case by scaling. 
  \end{remark}
\begin{remark}[Universal constants]
We recall again that a constant is \emph{universal} if it only depends on dimension and the constants $\Ckdgp, \Ckdgm$ appearing in the definition of
the De Giorgi's classes $\mathrm{kDG}^\pm$. 
\end{remark}
  
  By definition, this class of functions is invariant by scaling and left composition with $z_0 \in \R^{1+2d}$. 
\begin{lemma}[Invariance by scaling and translation of the De Giorgi's class]\label{l:invariance-kdgp}
  If $u \in \mathrm{kDG}^+(\domain, S)$ and $Q_r (z_0) \subset \domain$ and $r \in (0,1)$, then the function $v = \lambda u (\sigma_{r^{-1}} (z_0^{-1} \circ z))$
  lies in $\mathrm{kDG}^+(Q_1, \mathfrak{S})$  with $\mathfrak{S} (z)= \frac{\lambda}{r^2} S (\sigma_{r^{-1}} (z_0^{-1} \circ z) )$. 
\end{lemma}
  The integrability estimate from Proposition~\ref{p:subsol-integrability} and local energy estimates from Proposition~\ref{p:EE-kinetic}
  can be combined to prove that weak solutions of kinetic Fokker-Planck equations in a domain $\domain$ are in the kinetic De Giorgi's class $\mathrm{kDG}^+ (\domain,S)$. 
\begin{prop}[Sub-solutions are in the kinetic De Giorgi's class] \label{p:weak-kdgp}
  Let $f$ be a weak sub-solution of $(\partial_t + v \cdot \nabla_x ) f = \dive_v (A \nabla_v f) + B \cdot \nabla_v f +S$
  in $\domain$ with $A \in \mathcal{E}(\lambda,\Lambda)$ and $B,S \in L^\infty (\domain)$.
  Then $f \in \mathrm{kDG}^+ (\domain,S)$ for some universal $p_c>2$ and a constant $\Ckdgp$ depending on the largest $R$ such that there exists $z_0 \in \domain$ such that
  $Q_R (z_0) \subset \domain$.
\end{prop}
\begin{proof}
  On the one hand, Lemma~\ref{l:eq-trunc} and Proposition~\ref{p:subsol-integrability} can be combined in order to get the following estimate,
  \begin{multline*}
    \| (f-\kappa)_+ \|_{L^{p_c} (Q_r (z_0))} 
    \le C_\tau \bigg( \bigg[\frac{(\Lambda +1 +  \Lambda)R + 2 \Lambda}{ R-r} \bigg]\| \nabla_v (f-\kappa)_+\|_{L^2 (Q_R (z_0))} \\
    + 4 \bigg[ \frac{(\Lambda+1)R +1}{R^2-r^2} \bigg] \| (f-\kappa)_+ \|_{L^2 (Q_R (z_0))} + \|S \un_{\{ f \ge \kappa\}}\|_{L^2 (Q_R (z_0))} \bigg)
  \end{multline*}
 for $p_c \in (2,2+\frac1d)$.
In particular,
  \begin{multline} \label{e:1}
    \| (f-\kappa)_+ \|_{L^{p_c} (Q_r (z_0))}^2     \le C_1 \frac{(R+1)^2}{(R-r)^2} \| \nabla_v (f-\kappa)_+\|_{L^2 (Q_R (z_0))}^2 \\ 
    + C_1 \frac{(R+1)^2}{(R^2-r^2)^2} \| (f-\kappa)_+ \|_{L^2 (Q_R (z_0))}^2 + C_1 \|S \un_{\{ f \ge \kappa\}}\|_{L^2 (Q_R (z_0))}^2 
  \end{multline}
with $C_1$ only depending on $\Lambda$, $\Lambda$ and $p_c$. 
  
 On the other hand, we know from Proposition~\ref{p:EE-kinetic} with $t_2 = t_0$ and $t_1 \in (t_0-R^2, t_0 -r^2]$ that the following local energy estimate hold,
  \begin{align*}
  \frac\lambda2   \| \nabla_v (f-\kappa)_+ \|_{L^2 (Q_r (z_0))}^2     
     \le &      \iint_{B_{r^3}(x_0) \times B_r(v_0)} (f-\kappa)_+^2 (t_1,x,v) \dx \dv  + 2 \int_{Q_R (z_0)}  S (f-\kappa)_+   \\
     & + \left[ \frac{32 \Lambda+4}{(R-r)^2}   + \frac{\Lambda^2}{\lambda}\right] \|(f-\kappa)_+\|_{L^2 (Q_R (z_0))}^2. 
  \end{align*}
  Taking a mean in $t_1$ yields
    \begin{multline*}
      \frac\lambda2   \| \nabla_v (f-\kappa)_+ \|_{L^2 (Q_r (z_0))}^2 \\
      \le  \left[ \frac1{R^2-r^2} +  \frac{32 \Lambda+4}{(R-r)^2}   + \frac{\Lambda^2}{\lambda}\right] \|(f-\kappa)_+\|_{L^2 (Q_R (z_0))}^2 + 2 \int_{Q_R (z_0)}  S (f-\kappa)_+  \\
      \le  \left[   \frac{32 \Lambda+5}{(R-r)^2}   + \frac{\Lambda^2}{\lambda} +1 \right] \|(f-\kappa)_+\|_{L^2 (Q_R (z_0))}^2 +  \int_{Q_R (z_0)}  |S|^2\un_{\{f \ge \kappa\}}.
  \end{multline*}
   We write the last inequality as
  \begin{equation}
    \label{e:2}
    \| \nabla_v (f-\kappa)_+ \|_{L^2 (Q_r (z_0))}^2 \le C_2 ((R-r)^{-2} +1 )\|(f-\kappa)_+\|_{L^2 (Q_R (z_0))}^2 + C_2 \|S \un_{\{ f \ge \kappa\}}\|_{L^2 (Q_R (z_0))}^2
  \end{equation}
  for some constant $C_2$ only depending on $\lambda,\Lambda$ and $p_c$. 
  
  Combining \eqref{e:1} and \eqref{e:2} yields
    \begin{align*} 
      \| (f-\kappa)_+&  \|_{L^{p_c} (Q_r (z_0))}^2     \\
      & \le C_1 \frac{(R+1)^2}{(R-r)^2} \left( C_2 ((R-r)^{-2} +1 )\|(f-\kappa)_+\|_{L^2 (Q_R (z_0))}^2 + C_2 \|S \un_{\{ f \ge \kappa\}}\|_{L^2 (Q_R (z_0))}^2\right) \\
      & + C_1 \frac{(R+1)^2}{(R^2-r^2)^2} \| (f-\kappa)_+ \|_{L^2 (Q_R (z_0))}^2 + C_1 \|S \un_{\{ f \ge \kappa\}}\|_{L^2 (Q_R (z_0))}^2 \\
      & \le  C_3 \left(  \frac{(R+1)^2}{(R-r)^2} \left(\frac{1}{(R-r)^2} +1 \right) + \frac{(R+1)^2}{(R^2-r^2)^2} \right) \|(f-\kappa)_+\|_{L^2 (Q_R (z_0))}^2 \\
      & + C_3 \left(  \frac{(R+1)^2}{(R-r)^2} +1 \right) \|S \un_{\{ f \ge \kappa\}}\|_{L^2 (Q_R (z_0))}^2 .
    \end{align*}
    In particular, for $R<1$, we get $(R-r)<1$ and the previous inequality simplifies into,
    \[
      \| (f-\kappa)_+  \|_{L^{p_c} (Q_r (z_0))}^2 \le  C_4 (R-r)^{-4} \|(f-\kappa)_+\|_{L^2 (Q_R (z_0))}^2 
       + C_4 (R-r)^{-2} \|S \un_{\{ f \ge \kappa\}}\|_{L^2 (Q_R (z_0))}^2 . \qedhere
    \]
\end{proof}

\subsection{The local maximum principle}

We are now ready to state and prove the local maximum principle for kinetic Fokker-Planck equations and the associated kinetic De Giorgi's class. 
\begin{prop}[Local maximum principle] \label{p:lmp-kinetic}
  Let $Q \subset \R^{1+2d}$ be a kinetic cylinder and $p_c>2$ and $\Ckdgp \ge 1$ defining a class $\mathrm{kDG}^+ (Q,S)$. There exists $\Clmp$,
  only depending on $p_c$ and $\Ckdgp$, such that for all  $f \in \mathrm{kDG}^+(Q,S)$, all $Q_R(z_0) \subset Q$ and all $r< R  <1$,
  \[ \| f_+ \|_{L^\infty (Q_r(z_0))} \le \Clmp \left( (R-r)^{-\omega_0} \|f_+ \|_{L^2 (Q_R(z_0))} + \| S\|_{L^\infty (Q_R(z_0))} \right) \]
  with $\omega_0=\frac{p_c (3p_c -2)}{2 (p_c-1)^2}$. 
\end{prop}
\begin{remark}[Universal constants]
We recall again for the reader's convenience that a constant is \emph{universal} if it only depends on dimension $d$ and the constants $\Ckdgp,\Ckdgm$ appearing in the definition of
the De Giorgi's classes $\mathrm{kDG}^\pm$. 
\end{remark}
\begin{proof}
    We first assume that $z_0=0$ and  prove that there exists some universal constant $\delta_0 \in (0,1)$ such that, 
  if $\|S\|_{L^\infty(Q_R)} \le 1$  and if  $\|f_+ \|_{L^2 (Q_R)} \le \delta_0$,  then $f \le 2$ a.e. in $Q_r$. 

\paragraph{Iterative truncation.}
We follow the reasoning from the elliptic and parabolic cases by considering an increasing sequence $\kappa_k$
  and by integrating $(f-\kappa_k)^2$ on shrinking cylinders $Q^k = Q_{r_k}$. Precisely, we consider
  \[ A_k = \int_{Q^k} (f-\kappa_k)_+^2 \dz \]
  with
  \[ \forall k \ge 0, \quad \kappa_k = 2 - 2^{-k}, \quad r_k = r + (R-r)2^{-k} .\]

  In order to obtain an upper bound on $f$ in $Q_{r}$, we have to find two  constants $\beta>1$ and $C>0$ such that, for all $k \ge 0$, we have
  $A_{k+1} \le C^{k+1} A_k^\beta$. Indeed, in this case, Lemma~\ref{l:induc} implies that $A_k \to 0$ as soon as $A_0 < C^{-\frac\beta{\beta-1}}$.
  Since
  \[ A_0 = \int_{Q_R} (f-1)_+^2 \dx \le \|f_+ \|_{L^2 (B_1)}^2 \le \delta_0^2, \]
  we see that we can choose $\delta_0 = (1/2)  C^{-\frac{\beta}{2(\beta-1)}}$. 
  And because the limit of $A_k$ as $k \to +\infty$ is $\|(f-2)_+\|^2_{L^2 (Q_r)}$, the fact that $A_k \to 0$ yields $f \le 2$ almost everywhere in $Q_r$.

  \paragraph{Local gain of integrability.}
  We use the definition of the kinetic De Giorgi's class to write the local energy estimate for $f$ with $z_0=0$, $R= r_k$ and $r = r_{k+1}$.
  In particular, the difference of radii is $r_k - r_{k+1}= (R-r) 2^{-k-1}$, and recalling that $\|S\|_{L^\infty (B_1)} \le 1$, we obtain,
  \begin{align*}
    \| (f-\kappa_{k+1})_+  \|^2_{L^{p_c} (Q^{k+1})}  
    & \le \Ckdgp (R-r)^{-4} 2^{4(k+1)} \int_{Q^k} (f-\kappa_{k+1})_+^2  \\
    & + \Ckdgp (R-r)^{-2} 2^{2(k+1)} | \{ f \ge \kappa_{k+1} \} \cap Q^k| . \\
    \intertext{Use  Bienaymé-Chebyshev's inequality to get $| \{ f \ge \kappa_{k+1} \} \cap Q^k| \le (\kappa_{k+1}-\kappa_k)^{-2} A_k \le 2^{2(k+1)} A_k$,}
    & \le \Ckdgp (R-r)^{-4} 2^{4(k+1)} A_k  + \Ckdgp (R-r)^{-2} 2^{4(k+1)}A_k  \\
    & \le  \Ckdgp (R-r)^{-4} 2^{4k+5} A_k
  \end{align*}
  since $(R-r) \le R \le 1$.

\paragraph{Nonlinear iteration.}
  We now estimate $A_{k+1}$ from above by using H\"older's inequality with $q_c \in (1,2)$ such that $\frac12 = \frac1{p_c} + \frac1{q_c}$,
  \begin{align*}
    A_{k+1} &\le \|(f - \kappa_{k+1})_+\|_{L^{p_c}(Q^{k+1})}^2 \left\| \un_{\{ f \ge \kappa_{k+1}\}}\right\|_{L^{q_c} (Q^{k+1})}^2 \\
            & \le  \|(f - \kappa_{k+1})_+\|_{L^{p_c}(Q^{k+1})}^2 \left|\{ (f-\kappa_k) \ge \kappa_{k+1} - \kappa_k\} \cap B^k\right|^{\frac2{q_c}} \\
            & \le \bigg( \Ckdgp (R-r)^{-4} 2^{4k+5} A_k  \bigg) \bigg( 2^{2(k+1)} A_k\bigg)^{\frac{2}{q_c}} \\
    & =  2^{5 + \frac{4}{q_c}} \Ckdgp (R-r)^{-4} \left(2^{4 +\frac{4}{q_c}} \right)^k  A_k^{1+\frac{2}{q_c}} .
  \end{align*}

This implies in particular that $A_{k+1} \le C^{k+1} (R-r)^{-4 (k+1)} A_k^\beta$ with the universal exponent $\beta = 1 +\frac2{q_c} = 3-2p_c^{-1}>1$
and the universal constant $C \ge 1$ only depending on $q_c$ and $\Ckdgp$. In particular, Lemma~\ref{l:induc}
implies that $A_k$ converges to $0$ as soon as $A_0 < \left[C(R-r)^{-4} \right]^{-\frac{\beta}{(\beta-1)^2}}$. Since $A_0 \le \delta_0^2$ (see the beginning of the proof)
  we  pick $\delta_0 \in (0,1)$ such that
\[ \delta_0^2 = \frac12  \left[C (R-r)^{-4} \right]^{-\frac{\beta}{(\beta-1)^2}} = \frac12  C ^{-\frac{\beta}{(\beta-1)^2}} (R-r)^{\frac{4\beta}{(\beta-1)^2}} .\]

\paragraph{The general case.} We now  remark that if we do not assume anymore that $\|S\|_{L^\infty(Q_R)} \le 1$ and $\|f_+\|_{L^2 (Q_R)} \le \delta_0$,
either $f \le 0$ a.e. in $Q_1$ or $\|f_+ \|_{L^2 (Q_R)} > 0$. In the latter case, we consider
\[ \tilde f  = \frac{f}{ \delta_0^{-1} \|f_+ \|_{L^2 (Q_R)}+\|S\|_{L^\infty(Q_R)} }.\]
This function $\tilde f \in \mathrm{kDG}^+ (Q_R,S)$ with $S$ replaced with
\[ \tilde S = \frac{S}{ \delta_0^{-1} \|f_+ \|_{L^2 (Q_R)}+\|S\|_{L^\infty(Q_R)} } \le 1. \]
Since $\|\tilde f_+ \|_{L^2(Q_R)} \le \delta_0$, we conclude that \( \|\tilde f_+ \|_{L^\infty (Q_r)} \le 2 \), that is to say
 \[ \| f_+ \|_{L^\infty (Q_r)} \le 2 \delta_0^{-1} \|f_+\|_{L^2 (Q_R)} + 2 \|S\|_{L^\infty (Q_R)} . \]
 Since $\delta_0^{-1} = \sqrt{2} C^{\frac{\beta}{2(\beta-1)^2}} (R-r)^{-\frac{2\beta}{(\beta-1)^2}}$, we have $\omega_0 = \frac{2\beta}{(\beta-1)^2} = \frac{q_c(q_c +2)}2 = \frac{p_c (3p_c -2)}{2 (p_c-1)^2}$. 
\end{proof}

\begin{cor}[Upside down local maximum principle] \label{c:lower-gen-kin}
  Let $r,R \in (0,1]$. There exist universal constants $\eps_{0,1}, \eps_1 \in (0,1)$, $\eps_{0,1}$ only depending on $\Ckdgp$ and $p_c$ and $\eps_1$ also depending on $R-r$,
  such that if $-f \in \mathrm{kDG}^+ (Q_R,S)$ with $\|S \|_{L^\infty (Q_R)} \le \eps_{0,1}$ and $f \ge 0$ a.e. in $Q_R$, then
  \[ |\{ f \ge 1 \} \cap Q_R | \ge (1-\eps_1) |Q_R| \quad \Rightarrow \quad \left\{ f \ge \frac12 \text{ a.e. in } Q_r \right\}.\]
\end{cor}
\begin{proof}
  The function $g = 1-f$ also belongs to $\mathrm{kDG}^+(Q_R,S)$ and $g \le 1$ a.e. in $Q_R$. The local maximum principle from Proposition~\ref{p:lmp-kinetic} applied to $g$
  implies that
  \[ \text{for a.e. } z \in Q_r, \qquad 0 \le g (z) \le \Clmp \bigg(  (R-r)^{-\omega_0} \eps_1 |Q_1|   + \eps_{0,1} \bigg) . \]
  We now simply pick $\eps_{0,1}$ and $\eps_1$ such that $(R-r)^{-\omega_0} \eps_1 |Q_1| \le (4 \Clmp)^{-1}$ and $\eps_{0,1} \le (4 \Clmp)^{-1}$
  and deduce that $g \le 1/2$ a.e. in $Q_r$. This means that $f \ge 1/2$ in the small cylinder $Q_r$. 
\end{proof}

\begin{cor}[Local maximum principle - again] \label{p:lmp-again-kinetic}
  Given a (universal) constant $p \in (0,2)$, there exists a (universal) constant  $\pClmp >0$, only depending on $d,\lambda,\Lambda$ and $p$,
  such that for any $f \in \mathrm{kDG}^+ (Q_1,S)$,  then
  \[ \| f_+ \|_{L^\infty (Q_{1/2})} \le  \pClmp \left( \|f_+\|_{L^p (Q_1)} + \|S\|_{L^\infty (Q_1)} \right) \]
  where $\|f_+\|_{L^p (Q_1)} = \|f_+^p\|_{L^1 (Q_1)}^{\frac1p}$. 
\end{cor}
\begin{proof}
  We argue like we did in the elliptic setting. We give details for the reader's convenience. 
  The result   is a consequence of the interpolation of $L^2$ between $L^\eps$ and $L^\infty$. If $\eps < 1$, we interpolate $L^{2/\eps}$ between $L^1$ and $L^\infty$.
Let us make this precise. We start by applying Proposition~\ref{p:lmp-kinetic} for $r,R \in (0,1)$, 
\begin{align*}
  \| f_+ \|_{L^\infty (Q_r)}
  & \le \Clmp \left( (R-r)^{-\omega_0}  \|f_+\|_{L^2 (Q_R)} + \|S\|_{L^\infty (Q_R)} \right) \\
  & \le \Clmp \left( (R-r)^{-\omega_0}  \|f_+\|_{L^\eps (Q_R)}^{\eps/2} \|f_+\|_{L^\infty (Q_R)}^{1-\eps/2} + \|S\|_{L^\infty (B_R)} \right) \\
    & \le \frac12 \|f_+\|_{L^\infty (Q_R)} + K_\eps  (R-r)^{-\omega_\eps}  
\end{align*}
with  \[ K_\eps = 2^{\frac{2-\eps}{\eps}} \Clmp^{\frac2{\eps}}  \|f_+\|_{L^\eps (Q_1)} + \Clmp  \|S\|_{L^\infty (Q_1)}\]
and $\omega_\eps = \frac{2\omega_0}{\eps}$. We now consider $r_0 = \frac12$ and $r_{n+1}= r_n + \delta (n+1)^{-2}$ with
$\delta = \frac12 \left(\sum_{k=1}^\infty k^{-2}\right)^{-1} =  \frac{3}{\pi^2}$. In particular, $\frac12 \le r_n \le 1$ for all $n \ge 0$.
Letting $N_n$ denote $  \| f_+ \|_{L^\infty (Q_{r_n})}$, we thus have,
\[ N_n \le \frac12 N_{n+1} + K_\eps (\delta^{-1} (n+1)^2)^{\omega_\eps} \le \frac12 N_{n+1} + K_\eps \delta^{-\omega_\eps} (n+1)^{2 \omega_\eps}.\]
By induction, we thus get for all $n \ge 1$,
\[ N_0 \le \frac1{2^{n}} N_{n} +  K_\eps \delta^{-\omega_\eps} \left( \sum_{k=1}^n \frac{k^{2\omega_\eps}}{2^{k-1}} \right).\]
Letting $n \to \infty$, we conclude that
\begin{align*}
\| f_+ \|_{L^\infty (Q_{1/2})} & = N_0 \\
   &\le K_\eps \delta^{-\omega_\eps} \left( \sum_{k=1}^\infty \frac{k^{2\omega_\eps}}{2^{k-1}} \right) \\
  & =  \eClmp \left(  \|f_+\|_{L^\eps (B_1)} + \Clmp  \|S\|_{L^\infty (B_1)} \right)
\end{align*}
with $\eClmp = \delta^{-\frac{2\omega_0}{\eps}} \left( \sum_{k=1}^\infty \frac{k^{\frac{4\omega_0}{\eps}}}{2^{k-1}} \right)  \left( 2^{\frac{2-\eps}{\eps}} \Clmp^{\frac2{\eps}} + \Clmp \right)$ and  $\delta =  \frac{3}{\pi^2}$.
Since $\omega_0$ and $\Clmp$ are universal, the constant $\eClmp$ only depends on $d,\lambda,\Lambda$ and $\eps$. 
\end{proof}

\section{The  De Giorgi's class kDG\textsuperscript{--} \& the intermediate value principle}

\subsection{The kinetic De Giorgi's class kDG\textsuperscript{--}}

We introduce next the kinetic De Giorgi's class corresponding to DG\textsuperscript{--} and pDG\textsuperscript{--}.
The elements of this class satisfies the intermediate value principle, that is key to reach H\"older continuity and to prove the weak Harnack's inequality.
\begin{defi}[The kinetic De Giorgi's class $\mathrm{kDG}^-$] \label{d:kDG-}
  Let  $I = (a,b]$ with $a,b \in \R$ and $\Ox, \Ov$ be open sets of $\R^d$ and $S \in L^\infty(\domain)$.
  
  A function $f \colon \domain \to \R$ lies in the \emph{kinetic De Giorgi's class} $\mathrm{kDG}^-(\domain,S)$ if 
  \begin{itemize}
    \item \textsc{(Local gain of integrability)} $-f \in \mathrm{kDG}^+(\domain,S)$;
    \item \textsc{(Local gradient estimate)} For any (not necessarily kinetic) cylinders $\Qint = (T-\tau_-,T] \times B_{r_x} \times B_{r_v}$,
      $\Qext = (T-\tau_+, T] \times B_{R_x} \times B_{R_v}$ and $\Qint \subset \Qext$, and any $\kappa \in \R$, 
      \[ \| \nabla_v (f-\kappa)_-\|_{L^2(\Qint)} \le \Ckdgm \left( e^{-1} \| (f-\kappa)_-\|_{L^2 (\Qext)} + \| S \un_{\{ f \le \kappa\}} \|_{L^2 (\Qext)}\right) \]
      where $e = \min ((\tau_+-\tau_-)^{1/2}, R_v^{-1/2}(R_x-r_x)^{1/2}, R_v -r_v)$. 
    \item \textsc{(local Poincaré-Wirtinger's inequality)}
      For any radius $R>1$ and  $\eta \in (0,1)$ and $\Tmid \in (-1-2\eta^2,-1-\eta^2)$,
      let $\Qmid=(\Tmid,0] \times B_{8R} \times B_{2R}$ and  $Q_- = Q_\eta (-1-\eta^2,0,0)$ and $Q_+ = Q_1$. 
      For any $z_0, \rho$ such that $z_0 \circ \sigma_\rho (\Qmid) \subset \domain$ and any $\kappa \in \R$,
      the function $g = (f-\kappa)_- \left(\rho^{-1} (z_0^{-1} \circ z) \right)$ satisfies,
      \[ \left\| \left( g - \langle g \rangle_{Q_-} + \omega (R) \right)_- \right\|_{L^1 (Q_+)} \le \Ckdgm \left( \| \nabla_v g \|_{L^1 (\Qmid)} + \|S\|_{L^\infty (\Qmid)} \right) \]
      where $\langle g \rangle_{Q_-} = \Gamma \astkin \big( g \un_{Q_- \cap \Qmid} \big)$ and $\omega (R) \to 0$ as $R \to \infty$.
    \end{itemize}

\end{defi}
\begin{remark}[Universal constants]
We recall again that a constant is called \emph{universal} if it only depends on the constant appearing in the definition of
the De Giorgi's classes $\mathrm{kDG}^\pm$. 
\end{remark}

By definition, this class is invariant under scaling and left composition with any $z_0 \in \R^{1+2d}$. 
\begin{lemma}[Invariance by scaling and translation of the De Giorgi's class]\label{l:invariance-kdgm}
  If $u \in \mathrm{kDG}^-(\domain, S)$ and $Q_r (z_0) \subset \domain$ and $r \in (0,1)$, then the function $v = \lambda u (\sigma_{r^{-1}} (z_0^{-1} \circ z))$
  lies in $\mathrm{kDG}^-(Q_1, \mathfrak{S})$  with $\mathfrak{S} (z)= \frac{\lambda}{r^2} S (\sigma_{r^{-1}} (z_0^{-1} \circ z) )$. 
\end{lemma}

We will prove in the next subsection that weak super-solutions of kinetic Fokker-Planck equations are in the kinetic De Giorgi's class kDG\textsuperscript{--}.
\begin{prop}[Weak super-solutions are in kDG\textsuperscript{--}] \label{p:weak-kDG-}
  Let $f$ be a weak super-solution of $(\partial_t + v \cdot \nabla_x) f = \dive_v (A \nabla_vf ) + B \cdot \nabla_v f + S$ in $\domain$
  with $S \in L^\infty(\domain)$ and $A \in \mathcal{E}(\lambda,\Lambda)$ and $\|B\|_{L^\infty(\domain)} \le \Lambda$. 
  Then $f \in \mathrm{kDG}^- (\domain,S)$. 
\end{prop}

\subsection{Weak super-solutions are in the kinetic De Giorgi's class kDG\textsuperscript{--}}

Because we already derived local energy estimates for super-solutions (see Proposition~\ref{p:EE-kinetic}) that provides the local gradient estimate
(through a mean in time, see the proof of Proposition~\ref{p:weak-pdg} from Chapter~\ref{c:parabolic}), 
and that we can deduce the gain of integrability (Proposition~\ref{p:weak-kdgp}) from them,
the proof of Proposition~\ref{p:weak-kDG-} reduces to establishing the weak Poincaré-Wirtinger's inequality for weak super-solutions. Let us state it. 
\begin{prop}[Weak Poincaré-Wirtinger's inequality for weak super-solutions] \label{p:weak-PW}
  Let $\Tmid \in (-1-2_\eta^2, -1-\eta^2)$ and $\Qmid=(\Tmid,0] \times B_{8R} \times B_{2R}$ and  $Q_- = Q_\eta (-1-\eta^2,0,0)$ and $Q_+ = Q_1$. 
Let $f$ be a weak super-solution of $(\partial_t + v\cdot \nabla_x)f \ge \dive_v (A \nabla_v f) + B \cdot \nabla_v f + S$ in $\Qmid$
with $S \in L^\infty (\Qmid)$. Then
      \[  \left\| \left( f - \langle f \rangle_{Q_-} + \omega (R) \right)_- \right\|_{L^1 (Q_+)} \le \Ckdgm \left( \| \nabla_v f \|_{L^1 (\Qmid)} + \|S\|_{L^\infty (\Qmid)} \right) \]
      where $\langle f \rangle_{Q_-} = \Gamma \astkin \big( f \un_{Q_- \cap \Qmid} \big)$ and $\omega (R) = CR^{-2}$ for some constant $C$ only depending on dimension $d$. 
\end{prop}
The proof of the weak Poincaré-Wirtinger's inequality starts with establishing
a local estimate by using the representation of weak sub-solutions. 
\begin{lemma}[A local estimate] \label{l:lee-kin}
  Let $\Psi \colon \R^{1+2d} \to [0,1]$ be $C^\infty$, supported in $\Qmid$ and identically equal to $1$ in $Q_+$. 
  Let $f \in L^2(\Qmid)$ such that $\nabla_v f \in L^2 (\Qmid)$ and 
   $(\partial_t + v \cdot \nabla_x) f \ge \dive_v (A \nabla_v f) + B \cdot \nabla_v f + S$ in $\mathcal{D}'(\Qmid)$. Then
  \[ \| (f- \bar f)_- \|_{L^1 (Q_+)} \le C  \bigg( \|\nabla_v f \|_{L^1 (\Qmid)}  +  \| S\|_{L^1 (\Qmid)} \bigg) \]
  where $C = C(d,\Lambda)$ and   \[ \bar f = \Gamma \astkin (f \mathcal{K} \Psi).\]
\end{lemma}
\begin{proof}
  We first use the representation of sub-solutions (Proposition~\ref{p:truncated} applied to $-\ftrunc$) in order to get that 
  \[ \ftrunc = f \Psi = (\Gammax + \Gammav) \astkin \Ttrunc + \Gamma \astkin (\Strunc + \mtrunc) \]
  with
  \[
    \begin{cases}
    \Ttrunc & = \Psi (A-I) \nabla_v f, \\
      \Strunc & = (B \Psi - (A+I) \nabla_v \Psi) \cdot \nabla_v f + S \Psi + f \mathcal{K} \Psi.
    \end{cases}
  \]
  We then consider $g = \ftrunc - \bar f$  with $\bar f = \Gamma \astkin (f \mathcal{K} \Psi)$ and define
  \[ \bStrunc = (B \Psi - (A+I) \nabla_v \Psi) \cdot \nabla_v f + S \Psi.\]
  Then
  \[ -g \le (\Gammax + \Gammav) \astkin (-\Ttrunc) + \Gamma \astkin (-\Strunc). \]
  In particular,
  \[ g_- = \max(0,-g) \le |(\Gammax + \Gammav) \astkin \Ttrunc| + |\Gamma \astkin \Strunc|. \]
  We now repeat the reasoning from the proof of Lemma~\ref{l:reg-rep}.
  Thanks to Young's inequality (Lemma~\ref{l:young-kin}) and the integrability properties of the functions $\Gamma,\Gammax,\Gammav$ established
  in Proposition~\ref{p:fundamental}-\ref{i:integrability},
  we deduce that for all $T>0$, 
  \[  \Gamma \astkin \bStrunc \text{ and }  (\Gammax + \Gammav) \astkin \Ttrunc \text{ are in }  L^1 ((-T,T) \times \R^{2d})  .\]
  Moreover, there exists $C_d$ only depending on dimension such that,
  \[    \|g_-\|_{L^1 (Q_+)}  \le C_d \left( \|\Ttrunc\|_{L^1 (\R^{1+2d})} + \|\bStrunc\|_{L^1 (\R^{1+2d})}  \right).\] 
  The constant $C_d$ only depends on dimension because it is related to Lebesgue
  norms of $\Gamma$ on a time interval that is related $Q_+$ and $\Qmid$, both contained in $(-1,0]$. 
  We conclude by using the formulas for $\Ttrunc$ and $\bStrunc$. 
\end{proof}
The next step in the derivation of the weak Poincaré-Wirtinger's inequality from Proposition~\ref{p:weak-PW}
is to construct a cut-off function whose free transport part is non-negative, and bounded from below by $1$ in the past. Here is a precise statement.
Recall that $\Tmid = -1- 2{\etad}^2+ e$ with $e = \delta_1 \etad^2/2$.
\begin{lemma}[Cut-off function] \label{l:cut-off}
  Given $\delta_1,\etad \in (0,1)$, there exists a $C^\infty$ function $\Psi_1 \colon [\Tmid,0] \times \R^{2d}$
  that is supported in $[\Tmid,0] \times B_8 \times B_2$, identically equal to $1$ in $Q_+$, and such that
  $(\partial_t + v \cdot \nabla_x) \Psi_1  \ge 0$ and   $(\partial_t + v \cdot \nabla_x) \Psi_1  \ge 1$ in $[\Tmid, -1- \etad^2] \times B_1 \times B_1$. 
\end{lemma}
\begin{proof}
  Consider $\Psi_1 (t,x,v) = \varphi_1 (t) \varphi_2 (x-tv) \varphi_3 (v)$ for $C^\infty$ functions $\varphi_i$ valued in $[0,1]$ and such that
  \begin{itemize}
  \item $\varphi_1 (\Tmid)=0$, $\varphi_1' \ge 0$ in $[\Tmid,0]$, $\varphi_1' = 1$ in $[\Tmid,-1- \etad^2]$ and $\varphi_1 \equiv 1$ in $[-1,0]$;
  \item $\varphi_2$ is supported in $B_4$ and equal to $1$ in $B_3$;
    \item $\varphi_3$ is supported in $B_2$ and equal to $1$ in $B_1$. 
    \end{itemize}
    Let us check that the conditions are satisfied. As far as the free transport is concerned, we remark that $(\partial_t + v \cdot \nabla_x)(\varphi_2 (x-tv))=0$ and get,
    \[ (\partial_t + v \cdot \nabla_x) \Psi_1 (t,x,v) = \varphi_1' (t) \varphi_2 (x-tv) \varphi_3 (v). \]
    In particular, for $(t,x,v) \in [\Tmid,-1-\etad^2] \times B_1 \times B_1$, we have $|x-tv| \le 2$ and
    \[ (\partial_t + v \cdot \nabla_x) \Psi_1 (t,x,v) = \varphi_1' (t) \ge 1. \]
    
    Moreover, for $v \in B_2$ and $t \in [\Tmid,0]$ and $|x| \ge 8$, we have $|x-tv| \ge 4$ so that $\varphi_2 (x-tv)=0$. This implies that $\Psi_1$ is supported
     in $[\Tmid, 0] \times B_8 \times B_2$. To finish with, when $(t,x,v) \in Q_+$, we have $|x-tv| \le 3$ and we conclude that $\Psi_1 (t,x,v)=1$. 
\end{proof}
We can now study the local mean $\bar f$ associated with this test function $\Psi_1$, after re-scaling it.
\begin{lemma}[Control of the local mean] \label{l:local-mean}
  Let $f \colon \Qext \to [0,1]$ and $\Psi_1$ be the cut-off function from Lemma~\ref{l:cut-off}. 
  There exist $R>1$ and $\theta \in (0,1)$, only depending on $d,\etad,\delta_1,\delta_2$, such that the local mean
  \[ \bar f = \Gamma \astkin (f \mathcal{K} \Psi) \quad \text{  with } \quad \Psi (t,x,v) = \Psi_1 (t, x/R, v/R) \]
  satisfies,
  \[ \bar f \ge \Gamma \astkin \left( f \un_{Q_- \cap \Qmid} \right) - C R^{-2}\]
  for some constant $C$ only depending on the dimension. 
\end{lemma}
\begin{proof}
   We remark that the scaling implies that for all $(t,x,v) \in (\Tmid, -1 -\etad^2) \times B_{8R} \times B_{2R}$,
  \[ (\partial_t + v \cdot \nabla_x ) \Psi (t,x,v) \ge 1 \quad \text{ and } \quad \Delta_v \Psi (t,x,v) = R^{-2} \Delta_v \Psi_1 (t,x/R,v/R).\]
  
  We split $\bar f$ into two pieces: $\bar f = \bar f_{(1)} + \bar f_{(2)}$ with
  \[
    \bar f_{(1)} = \Gamma \astkin (f (\partial_t + v \cdot \nabla_x) \Psi) \quad \text{ and } \quad 
    \bar f_{(2)}  = - \Gamma \astkin (f \Delta_v \Psi).
  \]
  
  We first estimate $\bar f_{(1)}(z)$ from below for $z=(t,x,v) \in Q_+ = Q_1$.
  In order to do so, we first remark that for $\zeta =(s,y,w) \in Q_- \cap \Qmid \subset (-2,-1] \times B_1 \times B_1$, we have
  \[ |t-s| < 3, \quad |x-y - (t-s)w| \le 5, \quad |v-w| \le 2.\]
  This means that $\zeta^{-1} \circ z \in (-2,2) \times B_5 \times B_2$. 
  
  Keeping this remark in mind and using \eqref{e:super-level-lap}, we write,
  \begin{align*}
    \bar f_{(1)} (z) & = \int  \Gamma (\zeta^{-1} \circ z) f (\zeta) (\partial_t + v \cdot \nabla_x)\Psi (\zeta)  \dd \zeta. \\
    \intertext{Remark that the three terms in the product forming the integrand are non-negative,} 
    \bar f_{(1)} (z)   & \ge \int_{Q_- \cap \Qmid} \Gamma (\zeta^{-1} \circ z) f(\zeta) \dd \zeta.
  \end{align*}
  
We now turn our attention to $\bar f_{(2)}$. Since we have $\Delta_v \Psi (t,x,v) = R^{-2} \Delta_v \Psi_1 (t,x/R,v/R)$ and $f \in [0,1]$ in $\Qext$ and $\Psi$ is supported in $\Qext$, we can write
\[ |\bar f_{(2)} | \le \frac{\| \Delta_v \Psi_1\|_\infty}{R^2}  \left( \Gamma \astkin \un_{\Qext} \right) \le \frac{\| \Delta_v \Psi_1\|_\infty}{R^2}  \left( \Gamma \astkin \un_{\{ -1-2\etad^2\le t \le 0 \}} \right) \le \frac{3 \| \Delta_v \Psi_1\|_\infty}{R^2} .\]
We used that $1+2 \eta_2^2 \le 3$ and the fact that for $t$ fixed, $\Gamma(t,\cdot,\cdot)$ has mass $1$ to get the last inequality. 
\end{proof}

Then the weak Poincaré-Wirtinger's inequality follows from the combination of the three previous lemmas. 
\begin{proof}[Proof of Proposition~\ref{p:weak-kDG-}]
Combine Lemmas~\ref{l:lee-kin}, \ref{l:cut-off} and \ref{l:local-mean}. 
\end{proof}

\subsection{The intermediate value principle}

 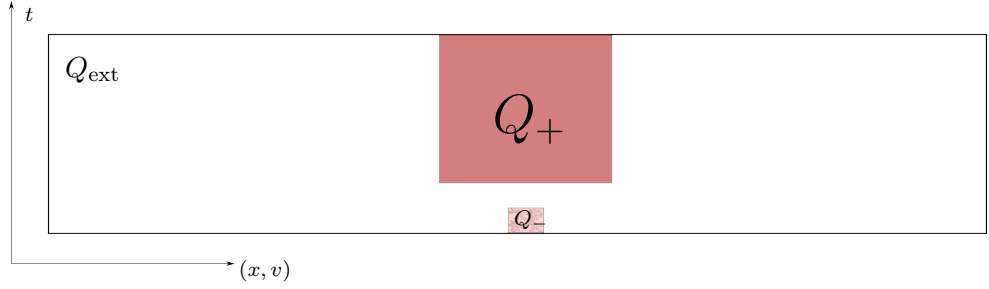
\begin{figure}[h]
 \centering{\includegraphics[height=3.5cm]{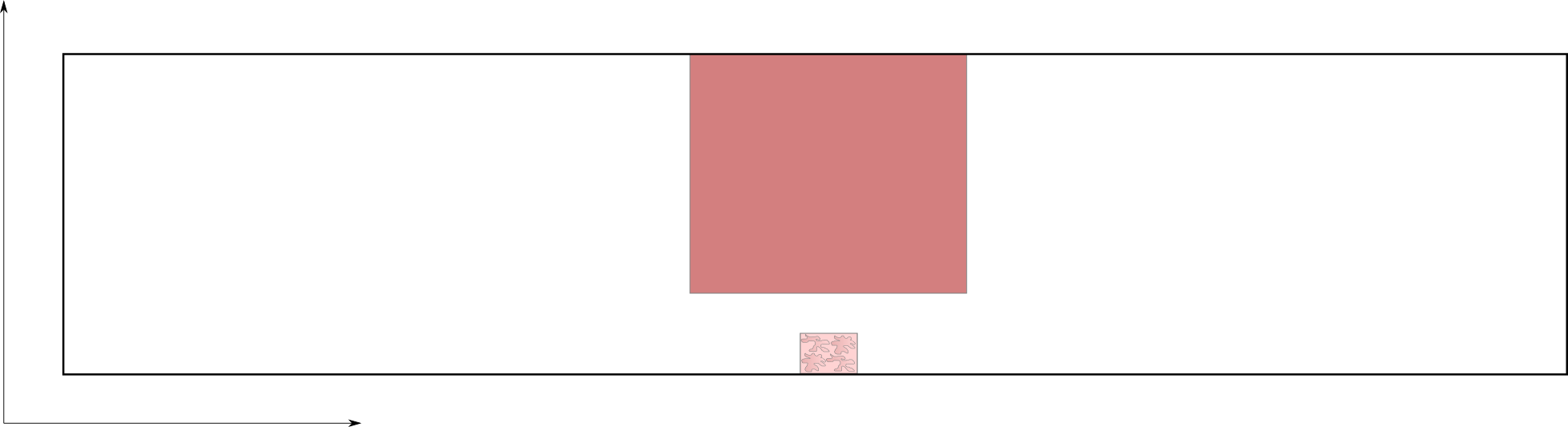}}
 \put(-360,92){\scriptsize $t$}
 \put(-345,71){$\Qext$}
 \put(-185,50){\huge $Q_+$}
 \put(-177,15){\scriptsize $Q_-$}
 \put(-280,-4){\scriptsize $(x,v)$}
 \caption{\textit{Geometric setting of the intermediate value principle.}}
   \label{fig:iv-principle}
 \end{figure}
\begin{prop}[Intermediate value principle] \label{p:iv-principle}
Given $\delta_1,\delta_2,\etad \in (0,1)$, there exist $R>2$ and $\theta,\delta_{1,2},\eps_{0,2} \in (0,1)$, only depending on $d,\lambda,\Lambda,\delta_1,\delta_2,\eta_2$, such that
for $\Qext = (-1-2\eta_2^2, 0] \times B_{\cyril{8}R} \times B_{\cyril{2}R}$ and $Q_- = Q_\etad (-1-\etad^2,0,0)$ and $Q_+ = Q_1$, for any $f \in \mathrm{kDG}^- (\Qext,S)$
with $\|S\|_{L^\infty(\Qext)} \le \eps_{0,2}$,  if
  \[ | \{ f \ge 1 \} \cap Q_-| \ge \delta_1 |Q_-| \quad \text{ and } \quad |\{ f \le \theta \} \cap Q_+| \ge \delta_2 |Q_+| ,\]
  then $|\{ \theta < f < 1 \} \cap \Qext | \ge \delta_{1,2} |\Qext|$.
\end{prop}
\begin{proof}
We choose $\Tmid = -1- 2{\etad}^2+ e$ with $e=\delta_1 \etad^2 /2$.
  The proof proceeds in three short steps. 

\noindent \textsc{Cropping the cylinder $Q_-$.}
We start with proving that the cropped cylinder \(Q_- \cap \Qmid\)  satisfies
  \begin{equation}
    \label{e:super-level-lap}
    |\{ f \ge 1 \} \cap Q_- \cap \Qmid | \ge \frac{\delta_1}2 |Q_-|.
  \end{equation}
  Recall that $e = \delta_1 \etad^2 /2$ and $\Tmid = -1- 2{\etad}^2+ e$ and that we have by assumption that
  \[ \int_{-1-2\etad^2}^{-1-\etad^2} |\{ f(t) \ge 1 \} \cap B_{\etad^3} \times B_\etad | \dt \ge \delta_1 |Q_-|.\]
  This implies that
  \[ | \{ f \ge 1 \} \cap Q_- \cap \Qmid| = \int_{-1-2\etad^2+e}^{-1-\etad^2} |\{ f(t) \ge 1 \} \cap B_{\etad^3} \times B_\etad | \dt \ge \delta_1 |Q_-| -  e \etad^{-2} |Q_-| = \frac{\delta_1}2 |Q_-|. \]
\medskip

\noindent \textsc{Lower bound on the local mean if $f \in [0,1]$.}
We next assume that $f$ takes values in $[0,1]$ (recall that it is only assumed to be non-negative).
We claim that, if  \( |\{ f = 1 \} \cap Q_-| \ge \delta_1 |Q_-|,\)  then
there exists $\theta >0$ and $R>1$ such that the local mean appearing in the definition~\ref{d:kDG-} of the kinetic De Giorgi's class $\mathrm{kDG}^-(\Qext,S)$ satisfies,
\begin{equation}
  \label{e:lower-mean}
  \langle f \rangle_{Q_-} - \omega (R) \ge  \sqrt{\theta} \text{ in } Q_+. 
\end{equation}
Indeed, we can use the definition of $\mathrm{kDG}^-(\Qext,S)$ and estimate the local mean $\langle f \rangle_{Q_-}$ as follows,
\begin{align*}
\langle f \rangle_{Q_-} & =   \int_{Q_- \cap \Qmid} \Gamma (\zeta^{-1} \circ z) f(\zeta) \dd \zeta\\
  & \ge \left(\inf_{(0,2) \times B_5 \times B_2} \Gamma  \right) | \{ f = 1 \} \cap Q_- \cap \Qmid  | \\
                     & \ge \bar C_d \frac{\delta_1}2  |Q_-|\\
                      & \ge \bar C_d \frac{\delta_1}2 \etad^{4d+2}|Q_1|
  \end{align*}
with $\bar C_d = \inf_{(0,2) \times B_5 \times B_2} \Gamma$. We now choose $R>1$ such that
\[ \omega (R) \le \bar C_d \frac{\delta_1}4 \etad^{4d+2}|Q_1| \]
and we obtain the claim by choosing $\sqrt{\theta} =   \bar C_d \frac{\delta_1}4  \etad^{4d+2}|Q_1| $.
\medskip

\noindent \textsc{Reaching the conclusion when $f \in [0,1]$.}
We combine the local estimate from Lemma~\ref{l:lee-kin} with the lower bound on the local mean from \eqref{e:lower-mean} in order to write,
\begin{align*}
  (\sqrt{\theta}-\theta) \delta_2 |Q_+| & \le (\sqrt{\theta} - \theta) |\{ f \le \theta \} \cap Q_+ | \\
                        & \le \| (\sqrt{\theta} -f)_+\|_{L^1 (Q_+)} \\
                        & \le \| (\bar f - f)_+ \|_{L^1 (Q_+)} \\
                        & \le C \bigg( \|\nabla_v f \|_{L^1 (\Qmid)} + \| S \|_{L^1 (\Qmid)} \bigg).
\end{align*}
We now make appear the intermediate value set. We use that $f \in [0,1]$ and $\nabla_v f = 0 $ in $\{f =1\}$ (see Proposition~\ref{p:composition}) in order to write
\[ \nabla_v f = \un_{\{ \theta < f < 1 \}} \nabla_v (f-1)_- +  \nabla_v (f-\theta)_-.\]
This decomposition leads to the following estimate,
\begin{align*}
  \|\nabla_v f \|_{L^1 (\Qmid)} \le & |\{ \theta < f < 1 \}\cap \Qmid|^{\frac12} \|\nabla_v (f-1)_-\|_{L^2 (\Qmid)}+ |\Qmid|^{\frac12} \| \nabla_v (f-\theta)_-\|_{L^2 (\Qmid)} \\
\intertext{we use now the local energy estimate for $f$ from $\Qmid$ to $\Qext$,}
                                 \le & \Ckdgm |\{ \theta < f < 1 \}\cap \Qext|^{\frac12} \left( e^{-1} \|(f-1)_- \|_{L^2 (\Qext)}  + \|S \|_{L^2 (\Qext)}\right)\\
                                    & + \Ckdgm  \left( e^{-1} \|(f-\theta)_- \|_{L^2 (\Qext)}  + \|S  \|_{L^2 (\Qext)}\right)\\
  \le &  (1+e^{-1}) \Ckdgm |\Qext|^{\frac12}  |\{ \theta < f < 1 \}\cap \Qext|^{\frac12} + \Ckdgm  ( |\Qext|^{\frac12} \theta + \eps_{0,2})
\end{align*}
where we recall that we chose $e =\delta_1 \etad^2/2$. 

We now combine the consequence of the local estimate with the estimate of the gradient and we get,
\[ \left(\sqrt{\theta} - \theta \right) \delta_2 |Q_+| \le (1+e^{-1}) \Ckdgm |\Qext|^{\frac12}  |\{ \theta < f < 1 \}\cap \Qext|^{\frac12} + \Ckdgm  ( |\Qext|^{\frac12} \theta + \eps_{0,2}). \]
We now pick $\theta$ and $\eps_{0,2}$ such that 
\[ \Ckdgm  ( |\Qext|^{\frac12} \theta + \eps_{0,2}) \le \frac12 \left(\sqrt{\theta}-\theta \right) \delta_2 |Q_+|\]
and we get the result with $\delta_{1,2} = \frac{(\sqrt{\theta}-\theta) \delta_2}{2(1+e^{-1})\Ckdgm|\Qext|^{\frac12}}$. 
\medskip

\noindent \textsc{Removing the upper bound condition on $f$.}
If now $f$ takes values larger than $1$, we  replace it with $\tilde f = \min (f,1) = 1 -(1-f)_+$. Thanks to Lemma~\ref{l:eq-trunc}, we know that it is  a super-solution of $(\partial_t + v \cdot \nabla_x) f \ge \dive_v (A \nabla_v f) + B \cdot \nabla_v f + \tilde S$ in $\mathcal{D}'(\Qext)$ with $\tilde S = S \un_{\{ \tilde f =1\}}$. Since their intermediate value and sub-level sets coincide and $\{ \tilde f = 1 \} = \{ f \ge 1\}$, the previous conclusion for $\tilde f$ implies the conclusion from the statement for $f$. 
\end{proof}

\section{Expansion of positivity and improvement of oscillation}

\subsection{Expansion of positivity}

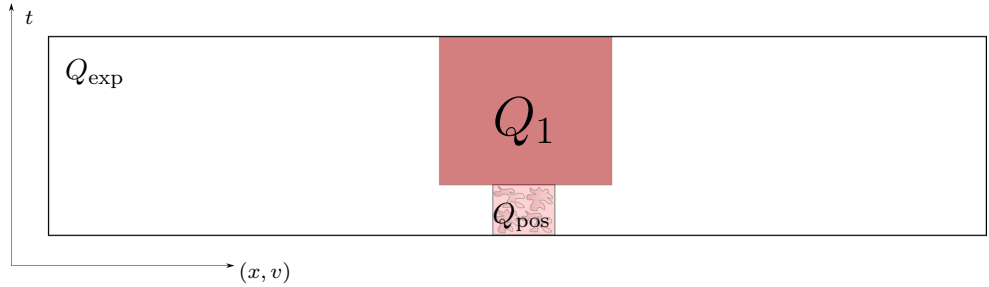
\begin{figure}[h]
 \centering{\includegraphics[height=3.5cm]{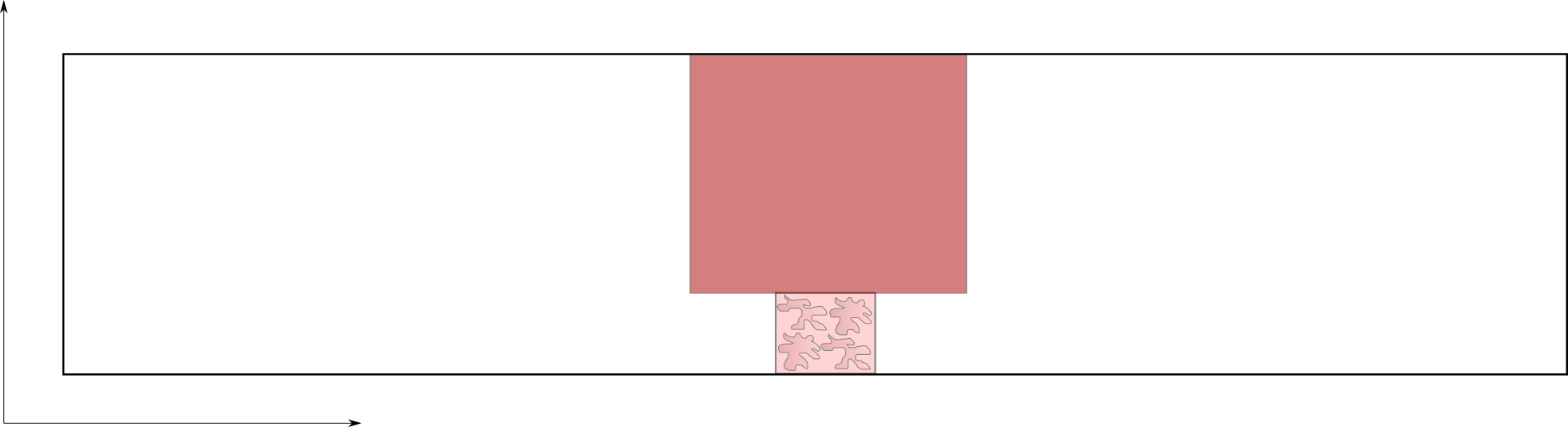}}
 \put(-360,92){\scriptsize $t$}
 \put(-345,71){$\Qexp$}
 \put(-185,50){\huge $Q_1$}
 \put(-185,17){$\Qgu$}
 \put(-280,-4){\scriptsize $(x,v)$}
 \caption{\textit{Geometric setting of the expansion of positivity.}}
   \label{fig:pop}
 \end{figure}

In this subsection, we prove that if $\{ f \ge 1\}$ is of positive measure in the cylinder $\Qgu$ lying in the past, a point-wise lower bound is generated in the future (in $Q_1$), see Figure~\ref{fig:pop}. 
 The parameter $\eta_0 \in (0,1)$ is useful when proving weak Harnack's inequality. We will use $\eta_0 =1/2$ and $\eta_0 = 1/\sqrt{m}$ where
 $m$ is an integer related to the covering argument, see Theorem~\ref{t:ink-kin} and Proposition~\ref{p:expansion-kin-s}.
 In order to get De Giorgi \& Nash's theorem, $\eta_0 = 1/2$ is enough. 
\begin{prop}[Expansion of positivity] \label{p:expansion-kin}
  Let   $\eta_0 \in (0,1)$. There exist  constants  $\ell_0,\eps_0 \in (0,1)$ and $R>1$, depending on $d,\lambda,\Lambda$ and $\eta_0$,
  such that for all $f \in \mathrm{kDG}^- (\Qexp,S)$  with $S \in L^\infty(\Qexp)$  such that $\|S\|_{L^\infty(\Qexp)} \le  \eps_0$ and $f \ge 0$ a.e. $\Qexp$, 
  \[ |\{ f \ge 1 \} \cap \Qgu | \ge \frac12 |\Qgu| \quad \Rightarrow \quad \{ f \ge \ell_0 \text{ a.e. in } Q_1 \}\]
  where $\Qgu = Q_{\eta_0} (-1,0,0)$ and  $\Qexp = (-1-\eta_0^2,0] \times B_{23R} \times B_{3R}$.
\end{prop}
\begin{proof}
The proof consists in applying Corollary~\ref{c:lower-gen-kin} to the function $f$ after scaling it. 
The intermediate value principle from Proposition~\ref{p:iv-principle} ensures that one of the scaling functions $f_k = \theta^{-k} f$
necessarily satisfies the assumption of the corollary. 
\medskip

\noindent \textsc{Geometric setting.}
Because we want to get the lower bound on $Q_1$, we need to get an information in measure from Proposition~\ref{p:iv-principle} in
a  cylinder $Q_{\eta_1}$ slightly bigger than $Q_1$. In order to do so, we first produce a time gap between $Q_1$ and $\Qgu$ (changing the latter in $\bQgu$),
to the cost of a factor $2$ on the lower bound on the measure of  the super-level set of $f$. Then we adjust the scaling factor $\eta_1$ in the intermediate value
principle so that, after scaling, we recover  the geometric setting of the expansion of positivity: $Q_-$ coincides with $\bQgu$. 
\begin{figure}[h]
 \centering{\includegraphics[height=2.5cm]{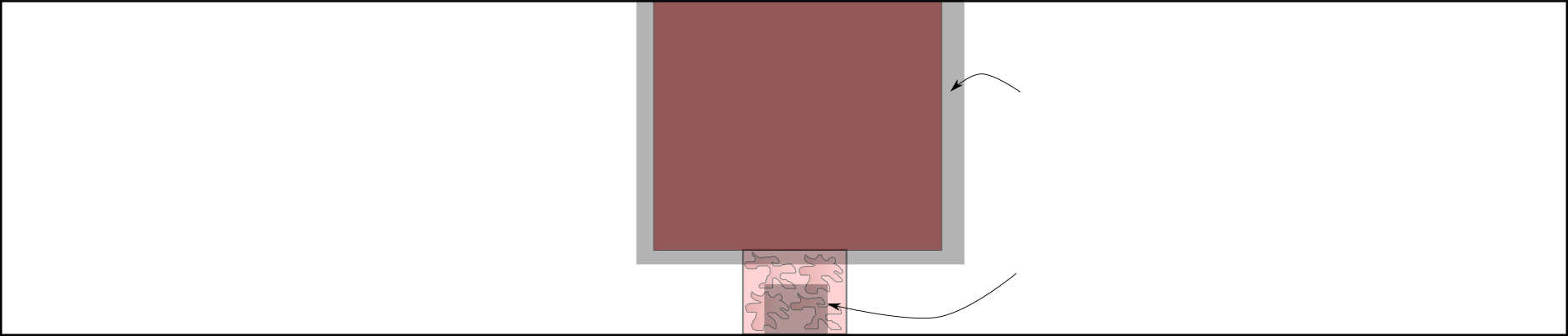}}
 \put(-325,56){$\Qexp$}
 \put(-115,45){\huge $Q_{_{\eta_1}}$ \scriptsize ($Q_+$ after scaling)}
 \put(-115,10){$\bQgu$ \text{\scriptsize ($Q-$ after scaling)}} 
 \caption{\textit{Geometric setting: from expansion of positivity to intermediate values.}}
   \label{fig:pop-proof}
 \end{figure}

We now make this precise. We first pick $\eta_1 \in (0,\eta_0)$, only depending on $\eta_0$, such that $|Q_{\eta_0}| - |Q_{\eta_1}| = \frac14 |Q_{\eta_0}|$. 
Such a choice  ensures that
\begin{equation}
  \label{e:cond-uls}
  |\{ f \ge 1 \} \cap  \bQgu | \ge \frac14 | \bQgu| 
\end{equation}
with $\bQgu = Q_{\eta_1}(-1-\eta_0^2+\eta_1^2,0,0)$.

Then we consider $\etad$ such that the scaling of $Q_-$ from the intermediate value principle (Proposition~\ref{p:iv-principle}) coincides with $\bQgu$.
\[ 1 + \eta_0^2 = 2 \eta_1^2 + \frac{\eta_1^2}{\etad^2} \Leftrightarrow \etad = \frac{\eta_1}{\sqrt{1+ \eta_0^2 - \eta_1^2}}.\]
The scaling factor $\eta_1/\etad$ is less than $\sqrt{2}$, then $\Qexp$ contains $\Qext$ after scaling. In particular, $2 \sqrt{2} \le 3$ and $8 \sqrt{2}^3 \le 23$. 

We next consider $\tilde{f}_k (z) = f_k ((\eta_1/\eta_2) z) = \theta^{-k} f ((\eta_1/\eta_2) z)$ where $r z = \sigma_r (z)$ is the kinetic scaling, see page~\pageref{p:scaling}.
We still have to fix the scaling parameter $\theta \in (0,1)$. Now $\tilde{f}_k$ can be studied in the geometric setting of the intermediate value principle (see Proposition~\ref{p:iv-principle}). 
\medskip

\noindent \textsc{Parameters from the local maximum and intermediate value principles.}
 We have to fix a parameter $\eps_0$ measuring the size of source terms in such a way  that all scaled functions $f_k$ are in a kinetic De Giorgi's class with a source term $S_k = \theta^{-k} S$ satisfying Conditions from Corollary~\ref{c:lower-gen-kin} and Proposition~\ref{p:iv-principle}.
In order to do so, we first get $\eps_1$ from Corollary~\ref{c:lower-gen-kin} with $R=1$, $r =\eta_2/\eta_1<1$.
Then we get $\delta_{1,2}$ from Proposition~\ref{p:iv-principle} for $\delta_1=\frac14$ and $\delta_2 = \eps_1$.
We  consider next the largest integer $N \ge 1$ such that $N \delta_{1,2} \le 1$.
We finally take $\eps_0 = \theta^{N+1} \min (\eps_{0,1},\eps_{0,2})$.

\medskip

\noindent \textsc{Finite iteration.}
We remark that for all $k \in \{1, \dots, N+1\}$, we have
\[ |\{ \tilde{f}_k \ge 1\} \cap Q_-| \ge |\{ \tilde{f} >1\} \cap Q_-| \ge \frac14|Q_-|\]
because $\tilde{f}_k \ge \tilde{f}_1$ and \eqref{e:cond-uls} corresponds to the previous estimate with $k=1$. 

Moreover, we consider the set $E \subset \{1,\dots, N+1\}$ of integers $k$ such that $|\{ \tilde{f}_k \le \theta \} \cap Q_1| \ge \eps_1 |Q_1|$.
For those $k$'s, the intermediate value principle implies that $|\{ \theta < \tilde{f}_k < 1 \} \cap \Qext| \ge \delta_{1,2}|\Qext|$.
This means that $|\{ \theta^{k+1} < \tilde{f} < \theta^k \} \cap \Qext| \ge \delta_{1,2} |\Qext|$. The sets $\{ \theta^{k+1} < \tilde{f} < \theta^k \} \cap \Qext$ are
disjoint in $\Qext$ and we conclude that
\[ |\Qext| \ge \sum_{k \in E} |\{ \theta^{k+1} < \tilde{f} < \theta^k \} \cap \Qext| \ge (\#E) \delta_{1,2} |\Qext|.\]
In particular, $(\# E)\delta_{1,2} \le 1$ and this implies that $\#E \le N$. 
\medskip

\noindent \textsc{Conclusion. }
We thus proved that there exists $k_0 \in \{1,\dots, N+1\} \setminus E$. This means that
\(|\{ \tilde{f}_{k_0} \le \theta \} \cap Q_1| < \eps_1 |Q_1|\) or equivalently,
\[ |\{ \tilde{f}_{k_0+1} > 1 \} \cap Q_1| >  (1-\eps_1) |Q_1|.\]
The upside down maximum principle from Corollary~\ref{c:lower-gen-kin} with $R=1$ and $r= \eta_2/\eta_1$ then implies that
\( \tilde{f}_{k_0+1} \ge \frac12  \) almost everywhere in $Q_r$, that is to say $f \ge \frac12 \theta^{k_0+1}$ in $Q_1$.
We thus get the desired estimate with $\ell_0 = \frac12 \theta^{N+1}$. 
\end{proof}

\subsection{Improvement of oscillation}

An easy consequence from the propagation of positivity is the improvement of oscillation and, in turn,  De Giorgi \& Nash's theorem. 
\begin{prop}[Improvement of oscillation] \label{p:improve-osc-kin}
  Let $\eps_0$  be given by Proposition~\ref{p:expansion-kin} about expansion of positivity.
  There exists a universal constant  $\mu \in (0 , 1)$ such that for all $f \in \mathrm{kDG}^+(Q_2,S) \cap \mathrm{kDG}^-(Q_2,S)$
  with $S \in L^\infty(Q_2)$ with $\|S\|_{L^\infty(Q_2)} \le \bar \eps_0$ and $u \in L^\infty(Q_2)$,
  \[\osc_{Q_2} f \le 2 \quad \Rightarrow \quad \osc_{Q_\omega} f \le 2 \mu .\]
\end{prop}
\begin{proof}
  We first embed $\Qexp$ from Proposition~\ref{p:expansion-kin} into a large kinetic cylinder. Since $\eta_0 \in (0,1)$ and $R >2$,
  we see that $\Qexp \subset Q_{6R}$. Then we consider  for $z \in Q_{6R}$,
  \[\tilde f (z) = f ((3R)^{-1} z)- \essinf_{Q_2} f.\]
  The function $\tilde f$ takes values in $[0,2]$ and the corresponding source term $\tilde S$ is such that
  \[\|\tilde S\|_{L^\infty(Q_{6R})} \le \eps_0.\]
  We distinguish two cases.
  \begin{itemize}
  \item If $|\{ \tilde f \ge 1 \} \cap \Qgu| \ge \frac12 |\Qgu|$, then  expansion of positivity from Proposition~\ref{p:expansion-kin}
    yields that $\tilde f \ge \ell_0$ a.e. in $Q_1$.
  \item If $|\{ \tilde f \ge 1 \} \cap \Qgu| < \frac12 |\Qgu|$, then we consider $g = 2 -\tilde f$ and we have $|\{ g \le 1 \} \cap \Qgu| < \frac12 |\Qgu|$
    or equivalently, $|\{ g > 1 \} \cap \Qgu| > \frac12 |\Qgu|$. In particular, $|\{ g \ge 1 \} \cap \Qgu| \ge \frac12 |\Qgu|$.
    Then  expansion of positivity from Proposition~\ref{p:expansion-kin} yields that $g \ge \ell_0$ a.e. in $Q_1$, or equivalently $f \le 2-\ell_0$ a.e. in $Q_1$.
  \end{itemize}
  In both cases, we obtain that $\osc_{Q_1} \tilde f  \le 2 -\ell_0$, that is to say $\osc_{Q_\omega} f \le 2 -\ell_0$ with $\omega = (6R)^{-1}$. 
\end{proof}

\subsection{Proof of the kinetic De Giorgi \& Nash's theorem}

\begin{proof}[Proof of Theorem~\ref{t:dg-kinetic} (De Giorgi \& Nash's theorem)]
The theorem is a consequence of the local maximum principle (Proposition~\ref{p:lmp-kinetic}) and of the improvement of oscillation (Proposition~\ref{p:improve-osc-kin}). 
The local maximum principle applied to $f$ and $-f$ with $z_0=0$ and $r=\frac34$ and $R=1$ implies that
  \begin{equation} \label{e:linfty-kin} \sup_{Q_{\frac34}} |f| \le \Clmp \left( 4^{\omega_0}\|f\|_{L^2 (Q_1)} + \|f\|_{L^\infty (Q_1)} \right).
  \end{equation}  
  As far as the H\"older semi-norm is concerned, we prove that there exists $\alpha \in (0,1]$ and $C_0 \ge 1$ (both universal) such that for all $z_0 \in Q_{\frac12}$, and all $r>0$,
  \[ \osc_{Q_r (z_0) \cap Q_{\frac12}} f \le C_0 \left( \|f\|_{L^\infty (Q_{\frac34})} + \|f\|_{L^\infty (Q_1)} \right) r^\alpha.\]
  This implies that $[f]_{\Cpar^\alpha (Q_{\frac12})} \le C_0 \left( \|f\|_{L^\infty (Q_{\frac34})} + \|S\|_{L^\infty (Q_1)} \right)$ (see Proposition~\ref{p:holder-kin}).

  We thus consider an arbitrary point $z_0 \in Q_{\frac12}$. We infer from \eqref{e:linfty-kin} that $f \in L^\infty (Q_{\frac14} (z_0))$. In order to invoke the
 improvement of the oscillation of $f$ (Proposition~\ref{p:improve-osc-kin}), we introduce
  \[ \tilde f (z) = \frac{f (z_0 \circ 8^{-1} z)}{ \|f\|_{L^\infty (Q_{\frac34})} + \bar \eps_0^{-1} \|S\|_{L^\infty (Q_1)}}.\]
  Then $\| \tilde f \|_{L^\infty (Q_2)} \le 1$  and the source term $\tilde S (z)= \bar \eps_0 8^{-2} \frac{S(z_0 \circ 8^{-1} z)}{\|S\|_{L^\infty (Q_1)}} $ satisfies
  $\|\tilde S \|_{L^\infty (Q_2)} \le  \eps_0$.
  We thus get from Proposition~\ref{p:improve-osc-kin} that $\osc_{Q_\omega} \tilde f \le 2 \mu$. We now scale recursively the function $\tilde f$ and consider,
  \[ \forall z \in Q_2, \quad  \tilde{f}_k (z) = \mu^{-k} \tilde{f} (( \omega/2)^k  z) \]
  whose source term $\tilde S_k (z) = (\omega^2 / (4\mu))^k \tilde S ((\omega/2)^k z)$. 
  We remark that $\|\tilde S_k\|_{L^\infty (Q_1)} \le \| \tilde S\|_{L^\infty (Q_1)} \le \bar \eps_0$ if we assume (without loss of generality) that $\mu \ge \omega^2/4$.
  We conclude that $\osc_{Q_2} \tilde{f}_k \le 2$ for all $k \ge 1$. This translates into,
  \[ \osc_{Q_{r_k}} \tilde{f} \le 2 \mu^k = 2^{1-\alpha} r_k^\alpha \quad \text{ with } \quad r_k = 2 \frac{\omega^k}{2^k} \quad \text{ and } \quad (\omega/2)^\alpha = \mu.\]
  Now for $r \in (0,2]$, there exists $k \ge 0$ such that $r_{k+1} \le r \le r_k$. This implies that
  \[\osc_{Q_r} \tilde{f} \le \osc_{Q_{r_k}} \tilde{f} \le 2^{1-\alpha} r_k^\alpha = 2 r_{k+1}^\alpha \le 2 r^\alpha. \]
  In terms of the function $f$, this implies that for all $r \in (0,2]$,
  \[\osc_{Q_{\frac{r}8} (z_0)} f  \le 8^\alpha 2   \left( \|f\|_{L^\infty (Q_{\frac34})} + \bar \eps_0^{-1} \|S\|_{L^\infty (Q_1)}\right) (r/8)^\alpha. \]
  Since for $s \ge \frac14$, we have
  \[\osc_{Q_s (z_0) \cap Q_{\frac12}} f  \le 2 \|f\|_{L^\infty (Q_{\frac34})} (4s)^\alpha, \]
  we conclude that
  \[[f]_{\Cpar^\alpha (Q_{\frac12})} \le 2^{4\alpha+1} \left( \|f\|_{L^\infty (Q_{\frac34})} + \bar \eps_0^{-1} \|S\|_{L^\infty (Q_1)} \right). \qedhere\]
\end{proof}

\section{(Weak) Harnack's inequality}

In this section, we show that elements of the kinetic De Giorgi's class $\mathrm{kDG}^-$ satisfies a weak Harnack's inequality.
We state it at unit scale. We already saw when deriving the expansion of positivity that it is necessary to have some room
in $(x,v)$ around the unit cylinder if the time interval is constrained to be $(-1,0]$. This impacts the geometric setting of the weak Harnack's inequality as well. 
\begin{thm}[Weak Harnack's inequality] \label{t:whi-kin}
  There exists two universal constants $R_0>1$ and $\omega \in (0,1)$ and two positive universal constants $\Cwhi$ and $p$
  such that for $\Qwhi = (-1,0] \times B_{R_0} \times B_{R_0}$ and $\Qpast = Q_\omega (-1+\omega^2,0,0)$ and $\Qfuture =Q_\omega$, 
  and $f \in \mathrm{kDG}^- (\Qwhi,S)$ with $S \in L^\infty (\Qwhi)$ and $f \ge 0$ a.e. in $\Qwhi$, we have 
  \[ \| f \|_{L^p (\Qpast)} \le \Cwhi \left ( \inf_{\Qfuture} f + \| S\|_{L^\infty (\Qwhi)} \right). \]
\end{thm}
\begin{remark}
  We let $\| f \|_{L^p (\Qpast)}$ denotes $\left( \int_{\Qpast} f^p \right)^{1/p}$ even if $p$ could be smaller than $1$. 
\end{remark}
It is then possible to combine the weak Harnack's inequality with the local maximum principle in order to get Harnack's inequality for solutions of kinetic Fokker-Planck equations.
\begin{thm}[Harnack's inequality] \label{t:harnack-kin}
  There exists two universal constants $R_0>1$ and $\omega \in (0,1)$ and a positive constant $\Charnack$ 
  such that for $\Qharnack = (-1,0] \times B_{R_0} \times B_{R_0}$ and $\Qpast^* = Q_{\omega/2} (-1+\omega^2,0,0)$ and $\Qfuture =Q_\omega$, 
  and $f \in \mathrm{kDG}^- (\Qharnack,S) \cap \mathrm{kDG}^+ (\Qharnack,S)$ with $S \in L^\infty (\Qharnack)$ and $f \ge 0$ a.e. in $\Qharnack$, we have 
  \[ \sup_{\Qpast^*} f \le \Charnack \left ( \inf_{\Qfuture} f + \| S\|_{L^\infty (\Qharnack)} \right). \]
\end{thm}
\begin{proof}
Apply first  Proposition~\ref{p:lmp-again-kinetic} between $\Qpast^*$ and $\Qpast$ (from the statement of the weak Harnack's inequality). 
Then combine the estimate with the one given by Theorem~\ref{t:whi}.
\end{proof}
\begin{remark}
  The fact that $\Qpast$ is replaced with $\Qpast^*$ in the statement of  Harnack's inequality
  is irrelevant since $\inf_{\Qfuture} f \le \inf_{\Qfuture^*} f$ with $\Qfuture^* = Q_{\omega/2}$. 
\end{remark}

\subsection{Generating and propagating a lower bound}

The proof combines the expansion of positivity from Proposition~\ref{p:expansion-kin} with a covering argument.
We aim at  estimating  $\| f \|_{L^p (\Qpast)}$ by the infimum of $f$ in $\Qfuture$. By linearity, we can reduce to $\inf_{\Qfuture} f \le 1$.
Establishing the estimate amounts to proving that there exists $\eps >0$ (universal) such that for all $t \ge 1$,
\[ |\{ f > t \} \cap \Qpast | \le C t^{-\eps} .\]
A further reduction is to prove that there exists $M>1$ and $\mu \in (0,1)$ such that for all integers $k \ge 1$, we have
\[ |\{ f > M^k \} \cap \Qpast| \le C (1 - \mu)^k .\]
In order to prove this, we consider $U_{k+1} = \{ f > M^{k+1} \} \cap \Qpast$ and we want to prove that $|U_{k+1}| \le (1-\mu) |U_k|$.
In order to prove this inequality, we cover the set $U_{k+1}$ with small cylinders $Q$ where we have a lower bound on $f$ in measure.
By using the expansion of positivity once, we generate a lower bound on a cylinder in the future, with a larger radius. Then
by applying iteratively this expansion of positivity, we propagate this lower bound in the future, till the final time. Since we know
that $f$ takes values smaller than $1$ in $\Qfuture$, this gives us some information on the radius of the initial cylinder $Q$.

\subsection{The Ink spots theorem in the kinetic geometric setting}

We first review the covering result that we will be using in the derivation of the weak Harnack's inequality. 

In order to state the covering result that we need in order to establish the weak Harnack's inequality,
we need to introduce the notion of \emph{stacked cylinders}. Given an integer $m \ge 1$ and a cylinder $Q = Q_r (z_0)$,
the stacked cylinder $\bar{Q}^m$ equals $\{ (t,x,v) \colon 0 < t-t_0 < mr^2, |x-x_0-(t-t_0)v_0|< (m+2)r^3, |v-v_0| < r \}$. 
\begin{thm}[Leaking ink spots in the wind]\label{t:ink-kin}
  Let $E$ and $F$ be two bounded measurable sets of $\R^{1+2d}$ such that $E \subset F \cap Q_1$.
  We assume that there exist two constants $r_0 \in (0,1)$ and an integer $m \ge 1$ such that
  for any cylinder $Q = Q_r (z_0) \subset Q_1$ such that $|Q \cap E | > \frac12|Q|$, we have $\bar{Q}^m \subset F$ and $r < r_0$.
  Then $|E| \le \frac{m+1}m (1-c) \left( |F \cap Q_1| + C m r_0^2 \right)$. The constant $c \in (0,1)$  and $C>1$ only depend on dimension $d$.
\end{thm}
The proof of this theorem is an easy adaptation of the parabolic one. We postpone it until Section~\ref{s:proof-ink}.

\subsection{Expansion of positivity for stacked cylinders and  for large times}

In this subsection, we derive the two results that will allow us to use the covering argument from Theorem~\ref{t:ink-kin}. There are two assumptions
on cylinders intersecting $E$ in a good proportion: after stacking them, they should lie in $F$, and their radius should be under control. 
\begin{itemize}
\item On the one hand, we check that if we choose  $\eta$ depending on the integer $m$ (from the statement of Theorem~\ref{t:ink-kin}) then Proposition~\ref{p:expansion-kin} 
yields  a lower bound in the stacked cylinder $\bar{Q}_1^m$ from an information in measure in $Q_1$.
\item
On the other hand, we apply iteratively Proposition~\ref{p:expansion-kin} with $\eta = 1/2$ in order to estimate
how the lower bound that is generated for small times deteriorates for large ones. 
\end{itemize}

\paragraph{Expansion of positivity for a staked cylinder.} 
We first scale and translate in time the result from Proposition~\ref{p:expansion-kin}  in order to get a statement with $\Qgu$ replaced with $Q_1$. We notice that the stacked cylinder $\bar{Q}_1^m$ equals $(0,m] \times B_{m+2} \times B_1$. \label{p:stack-kin}
\begin{prop}[Expansion of positivity for a stacked cylinder] \label{p:expansion-kin-s}
  Let $m \ge 1$ be an integer and let   $R_m>1$  be given by Proposition~\ref{p:expansion-kin} for $\eta = 1/\sqrt{m}$.
  Let  $\Qstack = (-1,m] \times B_{23m^3 R_m} \times B_{3m R_m}$.

  There exists a constant  $M >1$, depending on $d,\lambda,\Lambda$ and $m$, 
  such that, if $f \in \mathrm{kDG}^- (\Qstack,0)$  and $f \ge 0$ a.e. $\Qstack$, then,
  \[ |\{ f \ge M \} \cap Q_1 | \ge \frac12 |Q_1| \quad \Rightarrow \quad \{ f \ge 1 \text{ a.e. in } \bar{Q}_1^m \}.\]
\end{prop}

\begin{figure}[h]
\centering{\includegraphics[height=7cm]{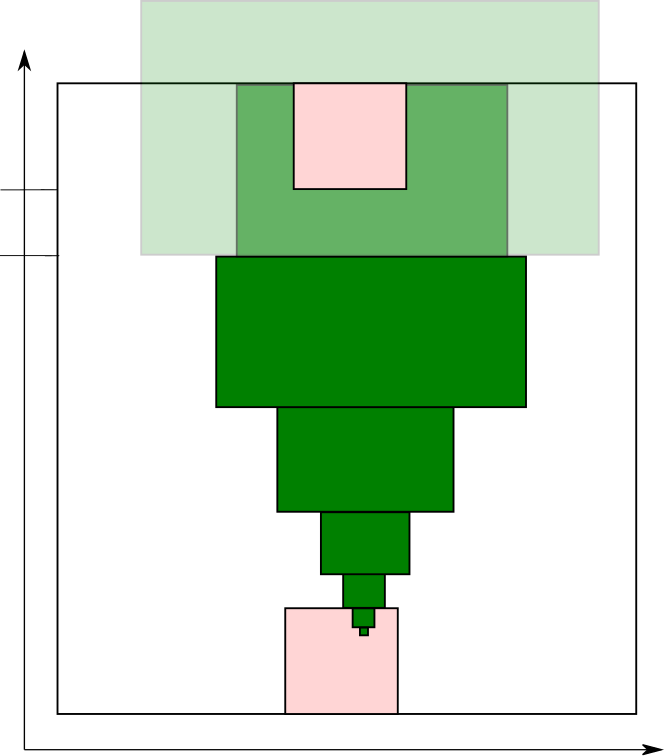}}
\put(-195,147){$-\omega^2$}
\put(-213,130){$t_0+T_N$}
\put(10,0){$(x,v)$}
\put(-180,180){$t$}
\put(-97,20){$\Qpast$}
\put(-96,160){\scriptsize $\Qfuture$}
\put(-90,110){$Q[N]$}
\put(-100,137){$Q[N+1]$}
\caption{\textit{Stacking iteratively cylinders above an initial one  contained in
    $Q_-$.} We see that the stacked cylinder obtained after $N+1$
  iterations by doubling the radius leaks out of the domain. This is
  the reason why $Q[N+1]$ is chosen in a way that it is contained in
  the domain and its ``predecessor'' is contained in $Q[N]$. Notice
  that the cylinders $Q[k]$ are in fact slanted since they are not
  centered at the origin. We also mention that $Q[N+1]$ is chosen
  centered if the time $t_0+T_N$ is too close to the final time~$0$.}
  \label{fig:stackedcylinders}
\end{figure}

\paragraph{Iteratively stacked cylinders.} We are going to apply iteratively Proposition~\ref{p:expansion-kin} to control the lower bound generated after applying it once.
 We need to make sure that the iterated cylinders do not exhibit the domain in $(x,v)$ and that their union  captures the cylinder $\Qfuture$, see Figure~\ref{fig:stackedcylinders}.
Recall that $\Qpast = Q_\omega (-1+\omega^2,0,0)$ and $\Qfuture= Q_\omega$.
\begin{lemma}[Iteratively stacked cylinders]\label{l:stack}
  Let $\omega \in (0,10^{-2})$. Given $Q = Q_r (z_0) \subset \Qpast$, we define for all $k \ge 1$,  \( T_k = \sum_{j=1}^k (2^j r)^2 \) and pick $N \ge 1$
  the largest integer such that $t_0 +T_N \le 0$. In particular $2^N r \le 1$.

  If $R$ denotes $|t_0+T_N|^{1/2}$, we consider $R_{N+1} =  \max (R,\rho)$ with $\rho =(4\omega)^{1/3}$ and
  \[ \forall k \in \{1,\dots,N\}, \quad z_k = z_0 \circ (T_k,0,0) \quad \text{ and } \quad z_{N+1} = \begin{cases} z_N \circ (R,0,0) & \text{ if } R \ge \rho, \\
                                                                                                         0 & \text{ if } R < \rho. \end{cases} \]
  We finally define $R_k = 2^k r$ for $k \in \{1,\dots,N\}$ and $Q[k] = Q_{R_k} (z_k)$ for $k \in \{1,\dots, N+1\}$.

  These cylinders $Q[k]$ are such that
  \[ Q[k] \subset (-1,0] \times B_2 \times B_2 \quad  \text{ and } \quad Q[N+1] \supset \Qfuture \quad \text{ and }\quad Q[N] \supset \tilde{Q}[N] \]
  where $\tilde{Q}[N] = Q_{\frac{R_{N+1}}2} (z_{N+1} \circ (-R_{N+1}^2,0,0))$. 
\end{lemma}
\begin{proof}
We first check that the sequence of cylinders is well defined for $\omega <10^{-2}$. Since
$r \le \omega$, we have $t_0+T_1 \le -1 +\omega^2 + 4r^2 < 0$. Let $N \ge 1$ be the largest integer such that $t_0 + T_N < 0$.

We check next that $\Qfuture \subset Q[N+1]$.

If $R < \rho$, then $Q[N+1] = Q_\rho$ and we simply remark that $\omega \le \rho = (4 \omega)^{1/3}$ to conclude.

In the other case, that is to say when $R \ge \rho$, we have $Q[N+1] = Q_R (z_{N+1})$ with $z_{N+1} = z_N \circ (R,0,0) = (t_0 + T_N + R^2, x_0,v_0) = (0,x_0,v_0)$.
In particular, $z_{N+1}^{-1} = (0,-x_0,-v_0)$ with $z_0 = (t_0,x_0,v_0) \in \Qpast = Q_\omega (-1+\omega^2,0,0)$.  We have to check that $Q_\omega (z_{N+1}^{-1}) \subset Q_R$. In this case,  for $z=(t,x,v) \in Q_\omega$,
\[ z_{N+1}^{-1} \circ z  = (t,-x_0+x-t v_0 , v-v_0) \in Q_R\]
if $\omega^2 \le R^2$ and  if $\omega^3 + \omega^3 + \omega^3 \le R^3$ and if $2 \omega \le R$.
This is true as soon as $4\omega \le R^3$, that holds true because we are dealing with the case $\rho \le R$. 
\medskip

Let us now check that for all $k \in \{1,\dots, N+1\},$
$Q[k] \subset (-1,0] \times B_{2} \times B_2$.

As far as $Q[N+1]$ is concerned, we use the fact that $R = |t_0+T_N|^{\frac12} \le 1$ and
$\rho = (4\omega)^{\frac13} \le 1$ to get $R_{N+1} \le 1$. Moreover
$z_{N+1} \in Q_1$ (because it is obtained by translating in time $z_0 \in \Qpast$) and thus $Q[N+1] \subset (-1,0]\times B_2 \times B_2$.

We remark
$(2^N r)^2 \le T_N \le -t_0 \le 1$.  If $\bar z_k =(t_k,x_k,v_k) \in Q[k]$ for
$k \le N$ then there exists $(t,x,v) \in Q_1$ such that
$\bar z_k = z_0 \circ (T_k,0,0) \circ ((2^kr)^2t,(2^k r)^3x,2^k r v)$. This
implies that $x_k = x_0 +T_k v_0+ (2^kr)^2 t v_0 + (2^kr)^3 x$ and
$v_k = v_0 + 2^k r v$ and since $z_0 \in \Qpast$,
\[ |x_k| \le |x_0| + |v_0| + |v_0| + |x|\le \omega^3 +  2\omega + 1 \le 2 \quad \text{ and } \quad |v_k| \le \omega + 1 \le 2.\]
In particular $Q[k] \subset (-1,0] \times B_2 \times B_2$. 
\medskip

We are left with proving that $\tilde{Q}[N]  \subset Q[N]$.

If $R \ge \rho$, then $\tilde{Q}_N = Q_{R/2} (z_{N+1})$ and $z_{N+1} = (0,x_0 - R^{1/2} v_0,v_0)$. We remark that $R/2 \le 2^N r \le 1$ (since $T_{N+1} >0$). 
In particular $|x_{N+1}| \le \omega^3 + R^{1/2} \omega \le \omega^3 + \sqrt{2} \omega \le 1$ and $|v_{N+1}| \le 1$ and $\tilde{Q}_N \subset Q_1 (z_{N+1}) \subset (-1,0] \times B_2 \times B_2$. 

Let us deal with the case $R \le \rho$. In view of the definitions of these cylinders, this is equivalent to
\[Q_{\rho/2} (\bar z) \subset Q_{2^N r} \text{ with } \bar z = (-T_N,0,0) \circ z_0^{-1} \circ (-\rho^2,0,0).\]

In order to prove this inclusion, we first estimate $2^N r$ from below. Since $t_0+T_{N+1} > 0$ and $-t_0 \ge 1 -\omega^2$, we have
$T_{N+1} = (4/3)(4^{N+1}-1)r^2 \ge 1-\omega^2$ and in particular $4^N r^2 \ge  (3/16)(1-\omega^2) \ge 1/8$ (since $\omega^2 \le 10^{-4} \le 2/3$). We conclude
that
\begin{equation}\label{e:lower}
  2^N r \ge 1/(2\sqrt2).
\end{equation}

With such a lower bound in hand, we now compute $\bar z = (R^2-\rho^2, -x_0 + (t_0+\rho^2)v_0, -v_0)$ and get
for $z \in Q_{\rho/2}$,
\[ \bar z \circ z = (R^2-\rho^2 + t,  -x_0 + (t_0+\rho^2)v_0 + x -tv_0,v-v_0) \in Q_{2\rho}. \]
Indeed, recalling that $\rho^3 = 4 \omega = 0,04$, we have $-2\rho^2 < R^2 -\rho^2 + t \le 0$ and $|-x_0 + (t_0+\rho^2-t)v_0 + x | \le \omega^3 +3 \omega + (\rho/2)^3 \le (2 \rho)^3$ and $|v-v_0|\le (\rho/2) + \omega \le 2 \rho$.
\end{proof}

\subsection{Iterated expansion of positivity}

\begin{prop}[Iterated expansion of positivity] \label{p:iep-kin}
  Let $R_{1/2}$ be the universal constant given by Proposition~\ref{p:expansion-kin} with $\eta_0 =1/2$
  and let $R_0 \ge R_{1/2}$ and $\Qwhi = (-1,0] \times B_{R_0} \times B_{R_0}$.
  There exists a universal constant $\gamma_0 >0$ such that for all $f \in \mathrm{kDG}^- (\Qwhi,0)$,  all $A>0$ and all cylinder $Q_r (z_0) \subset \Qpast$, 
  \[   |\{ f > A \} \cap Q_r (z_0)| > \frac12 |Q_r(z_0)| \qquad  \Rightarrow \qquad \bigg\{ f \ge A (r/2)^{\gamma_0} \text{ a.e. in } \Qfuture \bigg\}.\]
\end{prop}
\begin{proof}
  We first apply Proposition~\ref{p:expansion-kin} to the function $g = f /A$ after scaling it. This implies that $g \ge A \ell_0$ in $Q[1]$.
  We then apply it iteratively and get $g \ge A \ell_0^k$ in $Q[k]$ for all $k \in \{1,\dots,N\}$. In particular, $g \ge A \ell_0^N$ in $\tilde{Q}[N]$.
  This cylinder is the ``predecessor'' of $Q[N+1]$ and we thus finally get $g \ge A \ell_0^{N+1}$ in $Q[N+1]$. Because $Q[N+1]$ contains $\Qfuture$, we finally
  get $g \ge A \ell_0^{N+1}$ in $\Qfuture$. Now we remember that $2^N r \le 1$ (see Lemma~\ref{l:stack}). We pick $\gamma_0$ such that $\ell_0 = 2^{-\gamma_0}$
  and we write $\ell_0^{N+1} = \left(2^{-(N+1)} \right)^{\gamma_0} \ge (r/2)^{\gamma_0}$.
\end{proof}

\begin{proof}[Proof of Theorem~\ref{t:whi-kin} (weak Harnack's inequality)]
  The proof proceeds in several steps. \medskip

\noindent \textsc{Reduction.}
We first reduce to the case $S=0$ by considering $\tilde f = f + \|S\|_{L^\infty(\Qwhi)} (t+1)$. 

Second,  we  reduce to the case $\inf_{\Qfuture} f \le 1$ by considering $\tilde f = f / \max(1,\inf_{\Qfuture} f) $. Indeed, $f \le \tilde f \le f + \| S\|_{L^\infty(\Qwhi)} $ and it is in $\mathrm{kDG}^-(\Qwhi,0)$.
\medskip

\noindent \textsc{Parameters.} We now aim at proving that there exist two universal constants $p>0$ and $C>0$ such that \( \| f\|_{L^p (\Qpast)} \le \Cwhi.\)
This is equivalent to prove that there exists three universal constants $M >1$ and $\tilde \mu \in (0,1)$ and $\tilde{C}>1$ such that
\[  \forall k \ge 1, \qquad |\{ f > M^k \} \cap \Qpast | \le \tilde C (1-\tilde \mu)^k .\]
For $k=1$, we simply pick $\tilde \mu \le 1/2$ and $\tilde C \ge 2 |\Qpast|$. We then argue by induction.  
We are going to apply Theorem~\ref{t:ink-kin} (about covering with ink spots) for some integer $m \ge 1$ large enough so that $\frac{m+1}m (1-c) <1-c/2$. The parameter $m$ only depends on $c = c(d)$, it is therefore universal.

We are going to use Proposition~\ref{p:expansion-kin-s} (propagation of positivity for stacked cylinders) with $m$ universal as above. Then we obtain another universal parameter $R_m$  from Proposition~\ref{p:expansion-kin}, see the statement of Proposition~\ref{p:expansion-kin-s}. We will also use Proposition~\ref{p:iep-kin} (iterated expansion of
positivity) from which we get yet another universal parameter $R_{1/2}$. Now we choose $R_0 = \max (R_{1/2},23 m^3 R_m)$.
\medskip

\noindent \textsc{The covering argument.}
We are going to apply the ink spots theorem to the sets $E_0 = \{ f > M^{k+1} \} \cap \Qpast$ and $F_0 = \{ f > M^k \} \cap \Qwhi$ after transforming $\Qpast$ into $Q_1$.
We thus consider a cylinder $Q \subset \Qpast$ such that $|E_0 \cap Q | > \frac12 |Q|$.
We have to check that the stacked cylinder $\bar{Q}^m$ is a subset of $F_0$ and that the radius of $Q$ is controlled by some constant $r_k$. 

We start with checking that $\bar{Q}^m \subset F_0$ for $Q$ such that $|E_0 \cap Q | > \frac12 |Q|$, that is to say 
\[ |\{ f > M^{k+1} \} \cap Q| >\frac12 |Q|.\]
We recall that $\sigma_r$ denotes the scaling operator. If $Q = Q_r(z_0)$, we consider for $z \in Q_1$ the scaled function $g (z) = M^{-k} f (\sigma_r (z_0 \circ z))$, so that
\( |\{ g > M \} \cap Q_1| >\frac12 |Q_1|.\) We have $g \in \mathrm{kDG}^-(\Qstack,0)$.
We deduce from Proposition~\ref{p:expansion-kin-s} that $g \ge 1$ a.e. in $\bar{Q}_1^m$. This means that $f \ge M^k$ a.e. in $\bar{Q}^m$.
We thus proved that $\bar{Q}^m \subset F_0$. 

We now estimate $r$ from above for $Q=Q_r (z_0) \subset \Qpast$ such that \( |\{ f > M^{k+1} \} \cap Q| >\frac12 |Q|\).  Proposition~\ref{p:iep-kin} implies that
$f \ge M^{k+1} (r/2)^{\gamma_0}$ in $\Qfuture$. This implies that $M^{k+1} (r/2)^{\gamma_0} \le 1$ that is to say $r \le 2 M^{-\frac{k+1}{\gamma_0}}=:r_k$.
\medskip

\noindent \textsc{Conclusion.}
Now Theorem~\ref{t:ink-kin} implies that
\[ |\{ f > M^{k+1} \} \cap \Qpast| \le  (1-c/2) \left( |\{ f > M^{k} \} \cap \Qpast| + Cm 4 M^{-2\frac{k+1}{\gamma_0}}\right).\]
We use the induction assumption and get
\[ |\{ f > M^{k+1} \} \cap \Qpast| \le  (1-c/2) \left( \tilde C (1-\tilde \mu)^k + Cm 4 M^{-2\frac{k+1}{\gamma_0}}\right).\]
Recall that $M>1$ and $\gamma_0 >0$ are universal. We now choose $\tilde \mu$ so that $(1-\tilde \mu) \ge M^{-2/\gamma_0}$ and
$(1-\tilde \mu)^2 \ge 1- c/2$. 
We get
\begin{align*}
  |\{ f > M^{k+1} \} \cap \Qpast| & \le  (1-\tilde \mu)^2 \left( \tilde C (1-\tilde \mu)^k + Cm 4 (1-\tilde \mu)^{k+1} \right) \\
  & \le \bigg((1-\tilde \mu) \tilde C + 4 Cm\bigg) (1-\tilde \mu)^{k+1} .
\end{align*}
We thus pick $\tilde C$ such that $(1-\tilde \mu) \tilde C + 4 Cm \le \tilde C$ that is to say $\tilde C \ge 4Cm \tilde \mu^{-1}$. 
\end{proof}

\section{Proof of the ink spots theorem}
\label{s:proof-ink}

This section is devoted to the proof of the ink spots theorem for kinetic cylinders. The proof follows very closely
the one presented for parabolic cylinders. Still, we repeat every reasoning for the reader's convenience, and to allow them
to read both chapters independently. \medskip

The assumption of the ink spots theorem asserts that the set $E$ can be covered by cylinders and if more than half the cylinder
lies in $E$, then the corresponding stacked cylinder $\bar Q^m$ is contained in $F$. The conclusion asserts that the volume
of $E$ is bounded from above (up to some multiplicative constant) by the volume of $F$. In order to relate these two volumes,
it is necessary to extract from the original covering another one made of disjoint cylinders, and to make sure that we do not
lose too much by doing so. This is made possible thanks to a kinetic variation of Vitali's lemma with Euclidean balls. 

\subsection{A kinetic Vitali's covering lemma}

As explained in the previous paragraph, Vitali's lemma asserts that a countable disjoint family of cylinders can be extracted from
any covering of a set. We make sure that we do not lose too much by doing so is by imposing that the whole set is recovered
if the radii of cylinders of the sub-covering are multiplied by $5$.

For an arbitrary cylinder $Q \subset \R^{1+2d}$, if $Q = Q_r(z_0)$ with $z_0 = (t_0,x_0,v_0)$, then $5Q$ denotes $Q_{5r}(t_0 +12r^2,x_0,v_0)$.
It is necessary to update the top of the cylinder in order to extract a disjoint sub-cover, see in particular Lemma~\ref{l:overlap-kin}.
\begin{lemma}[Vitali]\label{l:vitali-kin}
  Let $\{Q_j\}_{j \in J}$ be a family of kinetic cylinders whose radii $r_j$ satisfy $\sup_{j \in J} r_j < +\infty$.
  There exists a countable sub-family $\{ Q_{j_i} \}_{i \in \N}$ of disjoint cylinders  such that
  \[ \cup_{j \in J} Q_j \subset \cup_{i \in \N} 5 Q_{j_i}.\]
\end{lemma}
In order to prove this lemma, we first deal with two overlapping cylinders.
\begin{lemma}[Overlaping kinetic cylinders]\label{l:overlap-kin}
  Let $Q_i = Q_{r_i} (z_i)$ for $i=1,2$ such that $Q_1 \cap Q_2 \neq \emptyset$ and $r_2 \le 2r_1$.
  Then $Q_2 \subset 5 Q_1$. 
\end{lemma}
\begin{proof}
  We first reduce to the case $z_1 =0$ by translating both cylinders. By assumption, there exists $z_{1,2} \in Q_1 \cap Q_2$.
This means that there exist $t_{1,2} \in (-r_1^2,0]$ and $x_{1,2} \in B_{r_1^3}$ and $v_{1,2} \in B_{r_1}$ such that
\[
  t_2-r_2^2 \le t_{1,2} \le t_2 \quad \text{ and } \quad |x_{1,2} - x_2 - (t_{1,2}-t_2) v_2| < r_2^3 \quad \text{ and } \quad |v_{1,2} - v_2 | < r_2.
\]
The fact that $Q_2 \subset 5 Q_1$ is equivalent to the following condition
\[
  -13 r_1^2 < t_2 -r_2^2 \le t_2 \le 12 r_1^2 \quad \text{ and } \quad |x_2| + r_2^2 |v_2| + r_2^3 \le (5r_1)^3\quad \text{ and } \quad |v_2 |+r_2 < 5r_1.
\]
We check these inequalities one after the other. First, $t_2 \ge t_{1,2} > -r_1^2 \ge r_2^2 - 13 r_1^2 $.
Second, $t_2 \le t_{1,2} + r_2^2 \le r_2^2 \le 4 r_1^2$.
Third,  $|v_2| \le |v_{1,2}| + |v_{1,2} -v_2| < r_1 + r_2 \le 3 r_1$. This implies directly the last inequality and
it also allows us to justify the third one because we also have:
\[ |x_2| < |x_{1,2}|+ |t_{1,2}-t_2||v_2| + r_2^3 \le r_2^3 + r_2^2 r_2 + r_2^3 = 3r_2^3. \]
These estimates for $x_2$ and $v_2$ imply that $|x_2| + r_2^2 |v_2| + r_2^3 \le (5r_1)^3\le 5r_2^3 \le 25 r_2^3$.
\end{proof}
The proof of Vitali's lemma is copied/pasted from the previous chapter but we include it here for the reader's convenience. 
\begin{proof}[Proof of Lemma~\ref{l:vitali-kin} (Vitali)] 
  Let $R = \sup_{j \in J} r_j$ where $r_j$ denotes the radius of the kinetic cylinder $Q_j$.
  Let $\mathcal{F}$ denote the family of cylinders $\{Q_j\}_{j \in J}$ and consider for all $n \ge 1$ the sub-family,
  \[ \mathcal{F}_n = \left\{ Q_j : \; j \in J, \; \frac{R}{2^n}< r_j \le \frac{R}{2^{n-1}} \right\}. \]
  We now construct families $\mathcal{G}_n$ by induction as follows: let $\mathcal{G}_1$ be any maximal disjoint sub-family of $\mathcal{F}_1$.
  Such a sub-family exists because of Zorn's lemma from set theory. If now $n \ge 1$ and $\mathcal{G}_1, \dots, \mathcal{G}_n$ are already constructed,
  then $\mathcal{G}_{n+1}$ is a maximal sub-family of
  \[ \left\{ Q_j \in \mathcal{F}_{n+1} \; : \; Q_j \cap Q_l = \emptyset \text{ for all } Q_l \in \mathcal{G}_1 \cup \dots \cup \mathcal{G}_n\right\}.\] 
  Roughly speaking, we add cylinders with smaller and smaller radii by making sure that they do not intersect the ones we already collected.
  We finally consider \[ \mathcal{G} = \cup_{n=1}^\infty \mathcal{G}_n.\]
  We now verify that this sub-family satisfies the conclusion of the lemma. We consider the sequence of cylinders $Q_{j_i}$ for $i=1,\dots,n$ such that
  $\mathcal{G} = \{ Q_{j_i}\}_{i \in \N}$. Then for $Q_j \in \mathcal{F}$, there exists $n \ge 1$ such that $Q_j \in \mathcal{F}_n$.
  Assume first that $n \ge 1$. By maximality of $\mathcal{F}_1$, there exists $Q_l \in \mathcal{F}_1$ such that $Q_j \cap Q_l \neq \emptyset$. 
  Assume now that $n \ge 2$. Because $\mathcal{G}_n$ is
  maximal, there exists $Q_l \in \mathcal{G}_m$ with $m \in \{1,\dots, n-1\}$ such that $Q_j \cap Q_l \neq \emptyset$. By definition of $\mathcal{F}_n$ and $\mathcal{G}_m$, we have
  $r_j \le \frac{R}{2^{n-1}}$ and $r_l \ge \frac{R}{2^m}$ with either $m=n=1$ or $1 \le m \le n-1$. In both cases, $ r_j \le 2 r_l$. Lemma~\ref{l:overlap-kin} then implies
  that $Q_j \subset 5 Q_l$.   
\end{proof}

\subsection{Lebesgue's differentiation theorem with kinetic cylinders}

We now turn to Lebesgue's differentiation theorem in the kinetic geometry. Because we have Vitali's lemma, proofs are exactly the same as in the parabolic chapter. Only
notation slightly changes.
\begin{thm}[Lebesgue's differentiation]\label{t:lebesgue-kin}
  Let $f \in L^1 (\R^{1+2d})$.  Then for a.e. $z \in \R^{1+2d}$,
\[
\lim_{r \to 0+}  \fint_{ Q_r (z)} |f-f(z)| =0
\]
where $\fint_{Q} g= \frac{1}{|Q|}\int_Q g$ for any cylinder $Q \subset \R^{1+2d}$ and $g \in L^1 (Q)$.
\end{thm}
The proof of this theorem relies on a functional inequality involving the maximal function.
For $g \in L^1 (\R^{1+2d})$, it is defined by,
\[ M g (t,x) = \sup_{Q \ni (t,x)} \fint_Q |g|.\]
\begin{lemma}[The maximal inequality]\label{l:maximal-kin}
  For all $\kappa >0$,
  \[ | \{ M g > \kappa \} \cap \R^{1+2d}| \le \frac{C}\kappa \|g \|_{L^1} \]
  for some constant $C$ only depending on dimension. 
\end{lemma}
\begin{proof}
  For every $z \in \R^{1+2d}$ such that $Mg (t,x) > \kappa$, there exists a cylinder $Q$ containing $z$ such that
  \[ \int_{Q} |g| \ge \frac\kappa2 |Q|. \]
  This means that the set $\{ M g > \kappa \} \cap \R^{1+2d}$ is covered with cylinders $\{Q_j\}$ satisfying the previous inequality.
  We know from Vitali's lemma~\ref{l:vitali-kin} that there exists a finite sub-family $\{Q_{j_i}\}_{i \in \N}$ such that
  \[ \{ M g > \kappa \} \subset \cup_{i \in \N} Q_{j_i}.\]
  With such a covering in hand, we can estimate the $L^1$-norm of $g$ as follows:
  \begin{align*}
    \int_{\R^{1+2d}} |g| &\ge \sum_{i \in \N} \int_{Q_{j_i}} |g| \\
                               & \ge \frac\kappa2 \sum_{i \in \N}  | Q_{j_i} | \\
                              & = \frac\kappa{ 5^{1+2d}2 } \sum_{i \in \N}  | 5 Q_{j_i} | \\
                        & \ge \frac{\kappa}{ 5^{1+2d}2 } |\{ M g > \kappa\}|.
  \end{align*}
We thus get the maximal inequality with $C = 5^{1+2d} 2$. 
\end{proof}
\begin{proof}[Proof of Theorem~\ref{t:lebesgue-kin} (Lebesgue's differentiation)]
  Let $f_n$ be continuous on $\R^{1+2d}$ and such that
  \[ \|f_n - f\|_{L^1} \le \frac1{2^n} .\]
  We can also assume that $f_n \to u$ almost everywhere in $\R^{1+2d}$ \cite[Theorem~4.9]{MR2759829}. Let $\mathcal{N}_0$ denote
  the negligible set outside which  point-wise convergence holds. 
  The maximal inequality from Lemma~\ref{l:maximal-kin} tells us that,
  \[ |\{ M (f_n -f) > \kappa \} | \le \frac{C}\kappa 2^{-n}. \]
  This is implies that the non-negative function $\sum_{n\in N} \un_{\{ M (f_n-f) > \kappa\}}$ is integrable over $\R^{1+2d}$ (Borel-Cantelli).
  It is thus finite outside of a negligible set $\mathcal{N}_1 \subset \R^{1+2d}$. 
  This implies that there exists $n_\kappa \in \N$ such that for all $n \ge n_\kappa$,
  \[ M (f_n -f) \le  \kappa \quad \text{ outside } \mathcal{N}_1.\]
  For all $i \in N$, we now we pick $\kappa = 1/i$ and construct an increasing sequence $n_i$ such that
  \[ M (f_{n_i} -f) \le \frac1i \quad \text{ outside } \mathcal{N}_1.\]

  With such a sequence of functions in hand, we can write for $z \in \R^{1+2d} \setminus (\mathcal{N}_0 \cup \mathcal{N}_1)$ and $i \in \N$,
  \[
    \fint_{Q_r (z)} |f - f(z)|  \le \fint_{Q_r (z)} |f - f_{n_i}| +\fint_{Q_r (z)} |f_{n_i} - f_{n_i}(z)| + |f_{n_i}(z) - f(z)| .
  \] 
  In the right hand side, the first term in bounded from above by $1/i$ because $z \notin \mathcal{N}_1$ and the third term goes to $0$ as $i \to \infty$ because $z \notin \mathcal{N}_0$.
  As far as the second term is concerned, the continuity of $f_{n_i}$ implies that it converges to $0$ too as $i \to \infty$. We thus proved that the left hand side
  tends to $0$ as $i \to \infty$. 
\end{proof}

\subsection{Proof of the ink spots theorem}

We continue following the reasoning from the parabolic chapter. 
The first step of the proof of Theorem~\ref{t:ink-kin} is to address the case where the two sets $E$ and $F$ are contained in the cylinder $Q_1$
and in which there is no time delay (no stacked cylinder). Again, the proof of this lemma is copied/pasted from the parabolic chapter but parabolic cylinders
are replaced with kinetic cylinders. 
\begin{lemma}[Crawling ink spots] \label{l:is-kin-nowind}
  Let $E \subset F \subset Q_1$ be measurable  sets of $\R^{1+2d}$.
  We assume that $|E| \le \frac12 |Q_1|$ and that
  for any cylinder $Q = Q_r (z_0) \subset Q_1$ such that $|Q \cap E | > \frac12|Q|$, we have $Q \subset F$.
  Then $|E| \le  (1-c)  |F |$. The constant $c \in (0,1)$  only depends on dimension $d$. 
\end{lemma}
\begin{remark}[The factor $1/2$]
  The factor $\frac12$ in both assumptions can be replaced with an arbitrary parameter $\mu \in (0,1)$. In this case, the conclusion is
  $|E | \le (1-c\mu)|F|$ for some $c\in (0,1)$ only depending on dimension. 
\end{remark}
\begin{proof}
    By applying Lebesgue's differentiation theorem~\ref{t:lebesgue-kin} to the indicator function $\un_E$,
  we know that for a.e. $x \in E$, there exists a cylinder $Q^x$ such that $|E \cap Q^x | \ge (1-\iota) |Q^x|$.
  Let us now choose a maximal cylinder $\Qmax^x \subset Q_1$ satisfying $|E \cap Q^x | \ge (1-\iota) |Q^x|$. It is of the form $ \Qmax^x = Q_{\bar r} (\bar t, \bar x)$. 
  By assumption, we know that $\Qmax^x \neq Q_1$ and $\Qmax^x \subset F$.

  We now claim that $|E \cap \Qmax^x | = \frac12 |\Qmax^x|$. 
  If the claim does not hold, then $\Qmax^x \neq Q_1$ and there would be a cylinder $Q^x$ and a $\delta>0$
  such that $\Qmax^x \subset Q^x \subset (1+\delta) \Qmax^x$ with  $Q^x \subset Q_1$ and $|E \cap Q^x| > \frac12 |Q^x|$,
  contradicting the maximality of $\Qmax^x$.
  
  The set $E$ is covered by cylinders $\Qmax^x$. By Vitali's lemma~\ref{l:vitali-kin}, there exists a countable sub-collection of non-overlapping cylinders
  $Q^j = Q_{r_j} (z_j)$, $j \ge 1$, such that $E \subset \cup_{j=1}^\infty 5  Q^j$. Since $Q^j \subset F$ and $|Q^j \cap E| = \frac12 |Q^j|$,
  this implies that $|Q^j \cap (F \setminus E) | = \frac12 |Q^j|$.
  \[
    |F \setminus E|  \ge \sum_{j=1}^\infty |Q^j \cap (F \setminus E)| 
    = \frac12 \sum_{j=1}^\infty  |Q^j| 
    = \frac12 5^{-1-d} \sum_{j=1}^\infty \iota |5 Q^j| 
    \ge \frac12 5^{-1-d} |E|.
  \]
  We conclude that $|F| \ge (1+5^{-1-d}2^{-1}  ) |E|$, from which we get $|E| \le (1-c) |F|$ with $c=5^{-1-d}2^{-2}$.
\end{proof}
We need two preparatory lemmas before proving the ink spot theorem with time delay (wind) and/or leakage. The first one was proved in the parabolic chapter
and it concerns the measure of a union of time intervals $(a_k -h_k,a_k]$ compared to the measure of the union of their stacked versions $(a_k,a_k + m h_k)$. 
\begin{lemma}[Sequence of time intervals] \label{l:intervals-kin}
  For all $k \ge 1$, let $a_k \in \R$ and $h_k >0$. Then,
  \[ \left| \bigcup_k (a_k,a_k+ m h_k) \right| \ge \frac{m}{m+1} \left| \bigcup_k (a_k-h_k,a_k] \right|.\]
\end{lemma}
We can now use this lemma about sequences of intervals to deal with sequence of stacked kinetic cylinders. The proof of this lemma
has to be adapted to take into account the $x$ variable. Here, there is a substantial difference  with the parabolic proof. 
\begin{lemma}[Overlaping stacked kinetic cylinders] \label{l:stack-overlap-kin}
  Let $\{Q_j\}$ be a family of kinetic cylinders and let $\bar Q_j^m$ be the corresponding stacked cylinders as defined on page~\pageref{p:stack-kin}.
  We have,
  \[ \left| \bigcup_j \bar Q_j^m \right| \ge \frac{m}{m+1} \left| \bigcup_j Q_j\right|.\]
\end{lemma}
\begin{proof}
  We use  Fubini's theorem in order to write,
\[  
  \left| \bigcup_j \bar Q_j^m \right|   = \int_{\R^d} \left| \left\{(t,x) \in \R^{1+d}: (t,x,v) \in \bigcup_j \bar Q_j^m \right\}\right| \dv .\]
We are thus left with proving that for $v$ fixed,
\[ \left| \left\{(t,x) \in \R^{1+d}: (t,x,v) \in \bigcup_j \bar Q_j^m \right\}\right|\ge \frac{m}{m+1}
  \left| \left\{(t,x) \in \R^{1+d}: (t,x,v) \in \bigcup_j Q_j \right\}\right|.\]
Let $r_j>0$ and $z_j=(t_j,x_j,v_j) \in \R^{1+2d}$ be such that $Q_j = Q_{r_j} (z_j)$.
If $|v-v_j| < r_j$, then
\begin{multline*}
  \left\{(t,x) \in \R^{1+d}: (t,x,v) \in \bar Q_j^m \right\} \\
  = \left\{ (t,x) \in \R^{1+d}: 0<t-t_j<mr_j^2, |x-x_j - (t-t_j) v_j | < (m+2)r_j^3 \right\}.
\end{multline*}
If $|v-v_j|\ge r_j$, then the set in the left hand side is empty. In the other case, $|(t-t_j) (v_j-v)| < m r_j^3$.
In particular, the right hand side contains the set
\[\left\{ (t,x) \in \R^{1+d}: 0<t-t_j<mr_j^2, |x-x_j - (t-t_j) v | < 2 r_j^3 \right\}.\]
Because $v$ is fixed, we can make the change of variables $(t,x) \mapsto (t,x+tv)$. It has Jacobian $1$.
Let $z_j = x_j -t_j v$. We thus have,
\[
   \left| \left\{(t,x) \in \R^{1+d}: (t,x,v) \in \bigcup_j \bar Q_j^m \right\}\right| 
   \ge \left| \bigcup_{j : |v_j -v|< r_j} (t_j,t_j+mr_j^2) \times B_{2r_j^3}\right|.
\]
We use  Fubini's theorem again,
\[
   \left| \left\{(t,x) \in \R^{1+d}: (t,x,v) \in \bigcup_j \bar Q_j^m \right\}\right| 
   \ge \int_{\R^d} \left| \bigcup_{j : |v_j -v|< r_j, |z-z_j|<r_j^3} (t_j,t_j+mr_j^2) \right| \dz.
\]
We now use Lemma~\ref{l:intervals-kin} and the change of variables $(t,z) \mapsto (t,z-tv)$ (of Jacobian $1$)  in order to get, 
\begin{align*}
   \left| \left\{(t,x) \in \R^{1+d}: (t,x,v) \in \bigcup_j \bar Q_j^m \right\}\right| 
   &\ge \frac{m}{m+1}\int_{\R^d} \left| \bigcup_{j : |v_j -v|< r_j, |z-z_j|<2 r_j^3} (t_j-r_j^2,t_j) \right| \dz \\
   & \ge \frac{m}{m+1}\int_{\R^d} \left| \bigcup_{j : |v_j -v|< r_j, |x - x_i -(t-t_j)v|<2r_j^3} (t_j-r_j^2,t_j) \right| \dz \\
  & \ge \frac{m}{m+1}\int_{\R^d} \left| \bigcup_{j : |v_j -v|< r_j, |x - x_i -(t-t_j)v_j|<r_j^3} (t_j-r_j^2,t_j] \right| \dz 
\end{align*}
where we used in the last line again that $|(t-t_j)(v-v_j)|< r_j^3$.  We now recognize in the right hand side
the desired lower bound. 
\end{proof}
Our next task is to get the ink spot theorem in the case where $E$ and $F$ are contained in the cylinder $Q_1$.
In other words, we postpone the treatment of cylinders leaking out of $Q_1$. Once again, no change for this proof compared to the parabolic setting. 
\begin{thm}[Ink spots in the wind] \label{t:is-kin-noleak}
  Let $E \subset F \subset Q_1$ be measurable sets of $\R^{1+2d}$.
  We assume that $|E| \le \frac12|Q_1|$ and that
  there exists an integer $m \ge 1$ such that  for any cylinder $Q = Q_r (z_0) \subset Q_1$ such that $|Q \cap E | > \frac12|Q|$, we have $\bar Q^m \subset F$.
  Then $|E| \le  (1-c) \frac{m+1}m |F |$. The constant $c \in (0,1)$  only depends on dimension $d$. 
\end{thm}
\begin{proof}
  We consider the family $\mathcal{Q}$ of kinetic cylinders $Q$ contained in $Q_1$ such that $|Q \cap E| > \frac12|Q|$. We let $G$ denote their union: $G = \bigcup_{Q \in\mathcal{Q}} Q$.
  We know from Lemma~\ref{l:is-kin-nowind} (crawling ink spots) that $E \le (1-c)|G|$. Moreover, the assumption of the theorem implies that $F$ contains the union of
  the corresponding stacked cylinders: $F \supset \bigcup_{Q \in \mathcal{Q}} \bar Q^m$. Using Lemma~\ref{l:stack-overlap-kin} about overlapping stacked cylinders, we obtain the
  following chain of inequalities,
  \[ |F| \ge \left| \bigcup_{Q \in \mathcal{Q}} \bar Q^m \right| \ge \frac{m}{m+1} \left| \bigcup_{Q \in \mathcal{Q}} Q \right| = \frac{m}{m+1}|G| \ge
    \frac{m}{(m+1)(1-c)} |E|. \qedhere\]
\end{proof}
We finally prove the covering result that was used in the derivation of the weak Harnack's inequality. The parabolic proof is slightly changed because kinetic cylinders
can leak out not only on the top (times after $t=0$) but also in the spatial variable ($x$ outside $B_1$).
\begin{proof}[Proof of Theorem~\ref{t:ink-kin}]
The assumption of the theorem implies that $|E| \le \frac12 |Q_1|$. Indeed, if this does not true, then $1 \le r_0$, contradicting the fact that $r_0 \in (0,1)$. 

  We consider again the family $\mathcal{Q}$ of kinetic cylinders $Q$ contained in $Q_1$ such that $|Q \cap E| > \frac12|Q|$.
  We let $\bar F$ denote the union of the corresponding stacked cylinders: $\bar F = \bigcup_{Q \in\mathcal{Q}} \bar Q^m$. Theorem~\ref{t:is-kin-noleak}
  implies that
  \[ |E| \le \frac{m+1}m (1-c) |\bar F| = \frac{m+1}m (1-c) \bigg[ |\bar F\cap Q_1| + |\bar F \setminus Q_1| \bigg].\]
  Moreover, the assumptions of Theorem~\ref{t:ink-kin} imply that $\bar F \subset F$.

  We are thus left we estimating $|\bar F \setminus Q_1|$. We first remark that because $|E| < \frac12 |Q_1|$, the
  conclusion of the theorem is trivial if $mr_0^2 \ge 1$ (take $C=\frac12$). We now assume that $mr_0^2 \le 1$. 

  We claim that for all $Q \in \mathcal{Q}$, we have $\bar Q^m \subset (-1,m r_0^2] \times B_{1+mr_0^2}\times B_1$. Indeed, $Q = Q_r (z_0)$ for some $z_0 \in Q_1$ and $r< r_0$ and
  $\bar Q^m = \{ (t,x,v) : t_0< t < t_0+ mr^2, |x-x_0 - (t-t_0)v_0|< (m+2) r^3 , |v-v_0|<r\}$.
  For $(t,x,v) \in \bar Q^m$, we have $t < mr_0^2$ because $t_0 \le 0$. We also have $|x| \le |x_0| + |t-t_0||v_0| \le 1+mr_0^2$.
  
  We thus proved that for all $Q \in \mathcal{Q}$,
  \[ \bar Q^m \setminus Q_1 \subset (0,mr_0^2) \times B_1 \times B_1 \cup (-1,mr_0^2) \times (B_{1+mr_0^2} \setminus B_1) \times B_1.\]
  This implies that
  \[ \bar F \setminus Q_1 \subset (0,mr_0^2) \times B_1 \times B_1 \cup (-1,mr_0^2) \times (B_{1+mr_0^2} \setminus B_1) \times B_1.\]
  We deduce from this inclusion that
  \begin{align*}
    |\bar F \setminus Q_1| & \le |B_1|^2 m r_0^2 + 2 \left[(1+m r_0^2)^d - 1 \right] |B_1|^2 \\
    & \le (1 +  2 d 2^{d-1}   )|B_1|^2 m r_0^2.
  \end{align*}
  In turn, this leads to,
  \[ |E| \le  \frac{m+1}m (1-c) \bigg[ | F\cap Q_1| + C m r_0^2 \bigg]\]
with $C= (1 + d 2^d     )|B_1|^2 m$. 
\end{proof}

\section{Transfer of regularity \& regularity of sub-solutions}

The main goal of this  section is  to study the regularity with respect to the spatial variable of sub-solutions of kinetic Fokker-Planck equations. This spatial regularity
is obtained from the natural energy estimates for sub-solutions. Indeed, these estimates ensure that sub-solutions are $H^1$ in the  variable $v$. The free  transport operator transfers
this regularity in $v$ into regularity in $x$.

\subsection{Regularity of the fundamental solution of the Kolmogorov equation}

In order to state a regularity property of the fundamental solution of the Kolmogorov equation, we first give the
definition of the fractional Laplacian. This singular integral operator is helpful when measuring the regularity of fractional order of a function.
\medskip

The fractional Laplace operator can be easily defined via Fourier transform, but it also has a singular
integral representation. We will use the latter definition. \medskip

   Let $\alpha \in (0,1)$.   For $f \in C^2 (\R^d)$ and $f$ bounded in $\R^d$, the fractional Laplacian of $f$ is defined by the following singular integral,
   \begin{align*}
     (-\Delta)^{\frac\alpha2} f (x) & = \PV c_{d,\alpha} \int_{\R^d} (f(x) -f (y)) \frac{\dy}{|y-x|^{d+\alpha}} ,
                                      \intertext{it is understood in the principal value sense \cite{MR1232192},}
                                      & = \lim_{\eps \to 0} c_{d,\alpha} \int_{|y-x|> \eps} (f(x) -f (y)) \frac{\dy}{|y-x|^{d+\alpha}}.
   \end{align*}
   The normalization constant $c_{d,\alpha}$ only depends on dimension and $\alpha$.

   Here are some classical properties of the fractional Laplacian that we will use when studying the integrability
   and the regularity of the fundamental solution of the Kolmogorov equation. 
   \begin{prop}[Fractional Laplacian]\label{p:fractional}
     Let $\alpha \in (0,1)$ and $p \in [1,2)$. 
   \begin{itemize}
 \item
   For $R >0$, if $f_R (x) = f (Rx)$, then  \[(-\Delta)^{\frac\alpha2} (f_R)(x) = R^\alpha (-\Delta)^{\frac\alpha2} (f) (Rx).\]
 \item
   There exist a constant $\Calp >1$, only depending on $s,p$ and dimension $d$, such that the norm and the semi-norm,
   \begin{align*}
     \|f\|_{W^{\alpha,p}(\R^d)} &:= \|f\|_{L^p (\R^d)} + \| (-\Delta)^{\frac\alpha2} f \|_{L^p (\R^d)} \\
     \|f\|_{\dot{W}^{\alpha,p} (\R^d)}  & := \left( \iint_{\R^d \times \R^d} \frac{|f(x) - f(y)|^p}{|x-y|^{d+ \alpha p}} \dx \dy \right)^{\frac1p}
   \end{align*}
satisfy  $\|f\|_{\dot{W}^{\alpha,p}(\R^d)} \le \Calp \|f\|_{W^{\alpha,p}(\R^d)}$.
 \end{itemize}
\end{prop}
We are now ready to state the regularity properties of $\Gamma$. 
\begin{prop}[Additional property of the fundamental solution]
   Let $\Gamma$ be the fundamental solution of the Kolmogorov equation.

   For all $T>0$ and $\eps \in (0,\frac13)$, we have  $(-\Delta_x)^{\frac\eps{2}} \Gamma \in L^{p_\eps} ((0,T) \times \R^{2d})$ for $p_\eps \in [1, 1 + \frac{2-3 \eps}{4d + 3 \eps })$
       and in $L^{1+\frac{2-3\eps}{4d+3 \eps}, \infty} ((0,T) \times \R^{2d})$. 
       The functions  $(-\Delta_x)^{\frac\eps{2}} \Gammax,(-\Delta_x)^{\frac\eps{2}} \Gammav$ are in $L^{q_\eps} ((0,T) \times \R^{2d})$ for $q_\eps \in [1,1 + \frac{1-3\eps}{4d+ 1 + 3\eps})$
       and in $L^{1+\frac{1-3\eps}{4d+1+3 \eps},\weak}((0,T) \times \R^{2d})$.
     \end{prop}
\begin{proof}
Thanks to the properties of the fractional Laplacian recalled in Proposition~\ref{p:fractional}, 
we know that
\begin{align*}
  (-\Delta)_x^{\frac\eps2} \Gamma (t,x,v) = \frac1{t^{2d+ \frac{3 \eps}2}} \left[ (-\Delta)_x^{\frac\eps2} \Gamma_1 \right] \left(\frac{x}{t^{\frac32}},\frac{v}{t^{\frac12}} \right), \\
  (-\Delta)_x^{\frac\eps2} \Gammax (t,x,v) = \frac1{t^{2d+ \frac{1+3 \eps}2}} \left[ (-\Delta)_x^{\frac\eps2} \Gamma_1 \right] \left(\frac{x}{t^{\frac32}},\frac{v}{t^{\frac12}} \right), \\
(-\Delta)_x^{\frac\eps2} \Gammav (t,x,v) = \frac1{t^{2d+ \frac{1+3 \eps}2}} \left[ (-\Delta)_x^{\frac\eps2} \Gamma_1 \right] \left(\frac{x}{t^{\frac32}},\frac{v}{t^{\frac12}} \right).
\end{align*}
We now compute,
\[ \frac{1+2d}{2d+ \frac{3 \eps}2} = 1 + \frac{2- 3\eps}{4d + 3\eps}\qquad \text{ and } \qquad \frac{1+2d}{2d+ \frac{1+3 \eps}2} = 1 + \frac{1-3\eps}{4d+ 1+3 \eps}. \]
The conclusion follows from Lemma~\ref{l:integ-weak}. 
\end{proof}

\subsection{Regularity in the space variable of sub-solutions}

In the next proposition, we obtain the regularity in the $x$ variable of $(f-\kappa)_+$, and more generally of sub-solutions.
Given a cylinder $Q = Q_R (z_0)$ with $z_0 = (t_0,x_0,v_0)$, we define,
\begin{multline*}     \|g \|_{L^p_{t,v}\dot{W}_x^{\alpha,p} (Q)} \\
  := \left( \iint_{(t_0-R^2,t_0] \times B_R (v_0)} \left\{ \iint_{Q^{(t,v)} \times Q^{(t,v)}}  \frac{|f(t,x,v) - f(t,y,v)|^p}{|x-y|^{d+\alpha p}} \dx \dy \right\}
    \dt \dv \right)^{\frac1p}
\end{multline*}
with $Q^{(t,v)} = \{ x \in \R^d : (t,x,v) \in Q\}= B_{R^3} (x_0 + (t-t_0)v_0)$.
\begin{prop}[$x$-regularity of sub-solutions]\label{p:x-reg}
  Let $f$ be a weak sub-solution of $(\partial_t + v \cdot \nabla_x ) f = \dive_v (A \nabla_v f) + B \cdot \nabla_v f +S$
  in $\domain$ with $A \in \mathcal{E}(\lambda,\Lambda)$ and $B \in L^\infty (\domain)$ and $S \in L^2 (\domain)$.
  For all $\eps \in (0,\frac13)$ and $p_\eps \in (1, 1 + \frac{2-3 \eps}{4d + 3 \eps })$,  all $Q_R (z_0) \subset \domain$ with $R<1$ and $r \in (0,R)$
  and $\mu_+^\kappa \in \R$, 
  \begin{align*}
    \| f \|_{L^{p_\eps}_{t,v}\dot{W}_x^{\eps,p_\eps} (Q_r (z_0))}      \le &C_c  (R-r)^{-1} \| \nabla_v f\|_{L^2 (Q_R (z_0))} \\
    + & C_c (R-r)^{-2} \| f \|_{L^2 (Q_R (z_0))} + C_c\|S \|_{L^2 (Q_R (z_0))}
  \end{align*}
  for some $C_c$ that only depends on $\Lambda$ and $d$. 
\end{prop}
\begin{remark}
  For $\eps \simeq \frac13$, we conclude that $(-\Delta)_x^{(1/3-0)/2}f \in L_{t,x,v}^{p}$ with $p = 1 + (4d+1)^{-1}-0$.
  This estimate is a slight improvement of the one  obtained    by J.~Guerand and C.~Mouhot in \cite{MR4453413}. 
\end{remark}
\begin{remark}
    One could use the Sobolev's inequality to gain some integrability in $x$ from this regularity estimate. But they do not provide
    integrability above $L^2$. More precisely, for all $\eps$ and $p_\eps$ as in the statement, $W^{\eps,p_\eps}(\R^d)$ does not embed into $L^2(\R^d)$. 
\end{remark}

\begin{proof}[Proof of Proposition~\ref{p:x-reg}]
  Let us assume that $z_0=0$. 
  Let $\phitrunc$ be given by Lemma~\ref{l:cutoff}. We use the representation formula from Proposition~\ref{p:truncated} in order to get
    \[
     f \phitrunc  = (\Gammax+  \Gammav) \astkin \Ttrunc + \Gamma \astkin (\Strunc-\mtrunc) 
  \]
  with $\Ttrunc,\Strunc \in L^2 (\R^{1+2d})$ given by
  \[  \begin{cases}
    \Ttrunc = & \phitrunc (A-I) \nabla_v f  ,\\
        \Strunc = &  ( B \phitrunc - A \nabla_v \phitrunc ) \cdot \nabla_v f  + S \phitrunc  + f (\partial_t + v \cdot \nabla_x - \Delta_v) \phitrunc    
  \end{cases}\]
  and some measure $\mtrunc \in M_1^+ (\R^{1+2d})$.
  In particular, if $|D|_x^\eps$ denotes $(-\Delta_x)^{\frac\eps2}$, then we have from Proposition~\ref{p:fundamental} that for all $\eps \in (0,1/3)$,
    \[
      |D|_x^\eps f \phitrunc  = \underset{L^{q_\eps}}{\underbrace{(|D|_x^\eps\Gammax+  |D|_x^\eps\Gammav)}} \astkin \underset{L^2}{\underbrace{\Ttrunc}} + \underset{L^{p_\eps}}{\underbrace{(|D|_x^\eps\Gamma)}} \astkin (\Strunc-\mtrunc) 
    \]
    with $p_\eps \in (1, 1 + \frac{2-3 \eps}{4d + 3 \eps })$ and $q_\eps \in [1,1 + \frac{1-3\eps}{4d+ 1 + 3\eps})$.

    This implies that the first term is in $L^{r_\eps} (\R^{1+2d})$ with $r_\eps$ such that $1+ \frac1{r_\eps} = \frac12 + \frac1{q_\eps}$.
    If $q_\eps = 1 + \frac{1-3\eps}{4d+ 1 + 3\eps}$, then $r_\eps = \frac{4d+2}{2d + 3 \eps} > 1 + \frac{2-3 \eps}{4d+3\eps}$. We conclude that
    $|D|_x^\eps f \phitrunc \in L^{p_\eps} (\R^{1+2d})$ and
    \[ \| (-\Delta_x)^{\frac{p_\eps}2} f \phitrunc \|_{L^{p_\eps} (\R^{1+2d})} \le \bar C_\eps \left( \|\Ttrunc\|_{L^2 (\R^{1+2d})} + \|\Strunc\|_{L^2 (\R^{1+2d})}
        + \|\mtrunc\|_{M^1_+ (\R^{1+2d})} \right)  \]
    for some constant $\bar C_\eps\ge 1$ only depending on dimension $d$ and $\eps$. 
    Moreover, from Lemma~\ref{l:reg-rep}, we have
    \[ \|  f \phitrunc \|_{L^{p_\eps} (\R^{1+2d})} \le \tilde C_\eps \left( \|\Ttrunc\|_{L^2 (\R^{1+2d})} + \|\Strunc\|_{L^2 (\R^{1+2d})}
        + \|\mtrunc\|_{M^1_+ (\R^{1+2d})} \right)  \]
    for some constant $\tilde C_\eps \ge 1$ only depending on dimension $d$ and $\eps$.
    Thanks to Proposition~\ref{p:fractional}, we conclude that,
        \[ \|  f \phitrunc \|_{L^{p_\eps}_{t,v} \dot{W}^{\eps,p_\eps}_x (\R^{1+2d})} \le  C_\eps \left( \|\Ttrunc\|_{L^2 (\R^{1+2d})} + \|\Strunc\|_{L^2 (\R^{1+2d})}
        + \|\mtrunc\|_{M^1_+ (\R^{1+2d})} \right)  \]
for some $C_\eps \ge 1$ only depending on $\eps$ and $d$. 
    
We now estimate $\|\mtrunc\|_{M^1_+ (\R^{1+2d})}$ by simply integrating \eqref{e:pp} against $1$,
\[ \|\mtrunc\|_{M^1_+ (\R^{1+2d})} \le \|\Strunc\|_{L^1 (\R^{1+2d})} \le |Q_R|^{\frac12} \|\Strunc\|_{L^2 (\R^{1+2d})}.\]
We finally obtain
\[
  \|  f \phitrunc \|_{L^{p_\eps}_{t,v} \dot{W}^{\eps,p_\eps}_x (\R^{1+2d})} 
  \le  C_\eps (1+ |Q_R|^{\frac12}) \left( \|\Ttrunc\|_{L^2 (\R^{1+2d})} +  \|\Strunc\|_{L^2 (\R^{1+2d})}  \right)  .
\]

Arguing as in the proof of Proposition~\ref{p:subsol-integrability}, we get estimates for $\Ttrunc$ and $\Strunc$ that leads to,
  \begin{align*}
    \|f \|_{L^{p_\eps}_{t,v}\dot{W}_x^{\eps,p_\eps} (Q_r)}  
    \le &C_c  (R-r)^{-1} \| \nabla_v f\|_{L^2 (Q_R )}  \\
     + &C_c (R-r)^{-2} \|  f  \|_{L^2 (Q_R )} + \|S \un_{\{ f \ge \kappa\}}\|_{L^2 (Q_R )}.
   \end{align*}

In the case where $z_0 \neq 0$, we applying the previous reasoning to $g (z) = f(z_0 \circ z)$. 
\end{proof}

\section{Bibliographical notes}
\label{s:biblio-kin}

\paragraph{Kinetic geometry.}
We explained in the introductory chapter (see page~\pageref{c:history}) that the regularity theory  for
kinetic Fokker-Planck equations took a new direction with the work by E.~Lanconelli and S.~Polidoro \cite{MR1289901}
about ultraparabolic equations. The geometric setting is the one obtained from H\"ormander's hypoellipticity theory where non-smooth coefficients are considered.
Roughly speaking, the free transport operator equals a finite sum of squares of vector fields 
In particular, in the papers mentioned in this chapter, the authors used systematically the  geometry associated with these ultraparabolic equations.
Among other things, they introduced functional spaces respecting this geometry, by studying for instance the H\"older regularity of solutions with respect to the norm
$\| z_2^{-1} \circ z_1 \|_\infty$ (see Lemma~\ref{l:up-low}). The kinetic distance from Definition~\ref{defi:kin-distance} was introduced more recently by L.~Silvestre
and the author in \cite{zbMATH07365668}. Kinetic H\"older spaces where introduced and used in the context of Schauder estimates. See in particular the definition of
kinetic H\"older spaces in \cite{zbMATH07365668,zbMATH07480708}. 

\paragraph{Kinetic De Giorgi \& Nash's theorem and Harnack inequalities.}
Theorem~\ref{t:dg-kinetic} was first proved by W.~Wang and L.~Zhang \cite{wz09}. The strong Harnack inequality was proved in \cite{zbMATH07050183}
by a compactness argument (for the intermediate value principle). Then constructive proofs where given in \cite{zbMATH07750909} and \cite{MR4453413},
by Kruzhkov's method \cite{MR171086} and through a trajectory argument, respectively. The result was then extended in various directions, for instance
by dealing with ultraparabolic equations \cite{anceschi2025poincareinequalityquantitativegiorgi}. The optimal dependency of the constant appearing in
Harnack's inequality with respect to $\lambda,\Lambda$ (ellipticity) is obtained in \cite{dietert2025criticaltrajectorieskineticgeometry}. 

\paragraph{Weak solutions.} Various notions of weak solutions were used in the works mentioned in this section.
In \cite{pp,wz09}, the authors impose that $(\partial_t + v \cdot \nabla_x)f$ is square integrable, which is too strong.
In \cite{zbMATH07050183}, the definition imposes $f \in L^\infty_t L^2_{x,v}$ because of the natural energy estimates.
In this book, we follow  a point of view that aligns with the classical parabolic one \cite{L1,MR1465184} by replacing the condition $f \in L^\infty_t L^2_{x,v}$ with $f$ square integrable
(in all variables).
More recently, P.~Auscher, L.~Niebel and the author introduced in \cite{auscher2025weaksolutionskolmogorovfokkerplanckequations} an even weaker notion of solution and
showed that time continuity with values in $L^2_{x,v}$ can be obtained,
in the spirit of Lions's embedding theorem. It is worth pointing out that the study of weak sub- and super-solutions is
out of the scope of this study. This being said, the section dedicated to the representation of weak sub-solutions borrows ideas coming from \cite{auscher2025weaksolutionskolmogorovfokkerplanckequations}. This is true in particular for the uniqueness proof, even if both statements and proofs differ.

\paragraph{Functional analysis framework and kinetic Poincaré inequalities.}
Kinetic functional spaces can be found in the first contributions by L.~H\"ormander \cite[p.~152]{hormander}  and they were used for instance by P.~Degond in \cite{MR875086}. 
More recently,   D.~Albritton, S.~Armstrong, J.-C.~Mourrat and M.~Novack contributed with \cite{MR4776290} to clarify the functional analysis framework for the study of kinetic
(including Fokker-Planck) equations. They  also established some functional inequalities, including some of Poincaré and H\"ormander types, and introduced new techniques to
establish them.

\paragraph{Pascucci-Polidoro's trick.} A.~Pascucci and S.~Polidoro \cite{pp} first obtained the local maximum principle for ultraparabolic equations by a Moser iteration procedure.
The gain of integrability was obtained in three steps: by first deriving local $v$-gradient estimates à la Caccioppoli; second, by artificially adding and subtracting a Laplacian in $v$;
third, by using the fundamental solution of the Kolmogorov equation, the diffusion operator with rough coefficients being treated as a source term. More generally, many of the papers
of the Italian school and later of W.~Wang and L.~Zhang \cite{wz09} use this trick. In contrast, averaging lemmas were used in the article by F.~Golse, C.~Mouhot, A.~F.~Vasseur and the
author \cite{zbMATH07050183} in order to gain  integrability and establish the intermediate value principle. This latter result is obtained by a compactness argument. 
Pascucci-Polidoro's trick was used by L.~Silvestre and the author \cite{MR4049224} when deriving a local H\"older estimate for a class of kinetic equations with integral diffusion.

\paragraph{Weak Poincaré's inequalities.} 
Weak Poincaré's inequalities are Poincaré's inequalities where the $L^2$-norm of the function (minus its local mean) is replaced with
the $L^2$-norm of its positive part. They are typically satisfied by sub-solutions. They first appeared in \cite{wz09} where W.~Wang and L.~Zhang
established H\"older regularity of solutions of so-called ultraparabolic equations. This class of equations contains in particular kinetic Fokker-Planck equations
that are studied in this chapter.  J.~Guerand and the author \cite{zbMATH07750909} established
an inequality directly inspired from \cite{wz09}, but relating the cylinder $\bar Q_1$ from the ``past'' to the cylinder $Q_1$ in the ``future''.
At the same time, J.~Guerand and C.~Mouhot \cite{MR4453413} also established some weak Poincaré's inequalities by using the $x$-regularity of sub-solutions and
by relating points in the past to point in the future by trajectories. Their trajectories are piece-wise smooth and follow alternatively the vector field
associated with free transport $(\partial_t + v \cdot \nabla_x)$ and the one related to $\nabla_v$. It relies on the idea coming from \cite{zbMATH07480708} that hypoellipticity
with rough coefficients can be recovered by commuting trajectories rather than vector fields (like in H\"ormander's far reaching hypoelliptic theory). 
A key observation where then made by L.~Niebel and R.~Zacher \cite{MR4875497}: it is possible to construct directly kinetic trajectories, without trying to commute
trajectories along vector fields. 
 This approach reached maturity with the paper \cite{anceschi2025poincareinequalityquantitativegiorgi} by further simplifying the construction from \cite{MR4875497}
and by dealing with ultraparabolic equations and integral diffusion. We finally refer the reader that is interested in this trend of research to
the very recent contribution \cite{dietert2025criticaltrajectorieskineticgeometry}.

\paragraph{Gain of integrability.}
The proof of the gain of integrability (see Proposition~\ref{p:subsol-integrability}) using the fundamental solution of the Kolmogorov equation
originates from the work A.~Pascucci and S.~Polidoro \cite{pp}
and it is also inspired by the reasoning by J.~Guerand and C.~Mouhot \cite{MR4453413}. In the former work, the authors use the fundamental solution as a test function in the
weak formulation while in the latter one, the representation of sub-solutions contained in Proposition~\ref{p:truncated} is used (without proof). 
It is explained in  \cite{dietert2025criticaltrajectorieskineticgeometry} that the gain of integrability can be established by using kinetic trajectories. 

\paragraph{Intermediate value principle.}
The proof of the intermediate principle from Proposition~\ref{p:iv-principle} relies on ideas from various works. The idea of retrieving to the function $f$
the function $\bar f = \Gamma \astkin (f \mathcal{K} \psi)$ (where $\mathcal{K}$ denotes the Kolmogorov operator) comes from \cite{wz09}, while the control
of the mean slightly departs from this work and follows the reasoning from \cite{zbMATH07750909}. In particular, information from the past is used in \cite{zbMATH07750909}
while \cite{wz09} followed the parabolic proof by propagating for short times some bounds on the super-level set of $f$ (see Lemma~\ref{l:st-propagation} from Chapter~\ref{c:parabolic}).

\paragraph{The kinetic ink spots theorem.} For general comments about ink-spots covering results, the reader is referred to the bibliographical section~\ref{s:biblio-parab}
of Chapter~\ref{c:parabolic}. 
Theorem~\ref{t:ink-kin} and its proof are extracted from \cite[Corollary~10.2]{MR4049224}. They rely on ideas introduced in the parabolic setting
in \cite{zbMATH06233951}. This being said, the parabolic proof contained in this book makes use of a Caldéron-Zygmund decomposition instead of a covering lemma à la Vitali.
Such an approach seems impossible to adapt to the kinetic setting because kinetic cylinders are slanted and overlap in the spatial variable. 

\paragraph{Regularity of sub-solutions and transfer of regularity.}
The last section of this chapter contains regularity results with respect to the space variable for sub-solutions of kinetic Fokker-Planck equations. A similar result
where first proved in \cite{zbMATH07050183} by using classical transfer of regularity properties due to F.~Bouchut \cite{MR1949176}. The use of the fundamental
solution was used in \cite{MR4453413}. The reasoning contained in this book is new. In particular, the estimates that are obtained are sharp in view of the sharpness
of Young's inequality on uni-modular groups. Moreover, this reasoning can be used to prove some sharp results from \cite{MR1949176} by changing the diffusion in the
Kolmogorov operator. For instance, \cite[Proposition~1.1]{MR1949176} can be proved by considering $\mathcal{K} = (\partial_t + v \cdot \nabla_x) + (-\Delta_v)^{\frac12}$. 

\paragraph{Kinetic De Giorgi's classes.} The notion of kinetic De Giorgi's class (see Definitions~\ref{d:kDG+} and \ref{d:kDG-}) is new.
Their definition can be modified in different ways. For instance the perturbed weak Poincaré-Wirtinger's inequality from Definition~\ref{d:kDG-}
can be proven without an error term $\omega_R$ and with $\langle f \rangle_{Q_-}$ replaced with $\fint_{Q_-} f$. 

\paragraph{Nash's approach via Fundamental solutions.} This book focuses on De Giorgi's approach to local regularity. In the parabolic case, Nash's proof \cite{nash} relies on the
existence of fundamental solutions of parabolic equations with variable coefficients and on an estimate on their logarithm, known as Nash's $G$ bound.
P.~Auscher, L.~Niebel and the author  constructed in \cite{MR4898687}
fundamental solutions for kinetic Fokker-Planck equations with variable coefficients and H.~Dietert together with L.~Niebel established in \cite{dietert2025nashsgboundkolmogorov}
the kinetic counterpart of Nash's $G$ bound. 

\paragraph{Conditional regularity program for Landau.} An important motivation for the study of such a class of equations is the one derived by Lev Landau -- see \cite{landau1936kinetische}. 
The Landau equation is nonlinear and describes the interaction between charged particles in a plasma, see for instance \cite[\S~4]{landau1980statistical}.
It can be written as follows,
\[ \partial_t f + v \cdot \nabla_x f = \nabla_v (A_f \nabla_v f -   b_f f)  \]
where $A_f (v) = \int_{\R^d} f(v-w) a (v-w) \dd w$ with $a_{ij} (z) = |z|^\gamma(\delta_{ij} |z|^2 - z_iz_j)$ for some $\gamma \in [-d,1]$
and $b_f = \dive_v A_f$. Under some conditions on the following density functions,
\[ \rho (t,x) = \int_{\R^d} f(t,x,v) \dv, \quad E (t,x) = \int_{\R^d} f(t,x,v) |v|^2 \dv, \quad H(t,x) = \int_{\R^d} f \ln f (t,x,v) \dv, \]
it has been known for a long time that the matrix $A_f$ is uniformly elliptic, see for instance the article by L.~Desvillettes and C.~Villani \cite{MR1737548} and more recently
the contribution by L.~Silvestre \cite{MR3582250}. It was conjectured in \cite{zbMATH07050183} that, as long as the three previous density functions are ``under control'', solutions
of the Landau equation (at least for $\gamma \ge -2$ in dimension $3$) remain smooth. The local H\"older estimate from \cite{zbMATH07050183} was the first step in this program.
It was completed in \cite{MR3778645} and \cite{zbMATH07177445}.

\bibliographystyle{plain}
\bibliography{dg.bib}

\end{document}